\documentclass{book}
\usepackage[a4paper, top=3cm, bottom=3cm]{geometry}
\usepackage[latin1]{inputenc}
\usepackage{setspace}
\usepackage{fancyhdr}
\usepackage{tocloft}
\usepackage{geometry}                
\usepackage[parfill]{parskip}    % Activate to begin paragraphs with an empty line rather than an indent
\usepackage{graphicx}
\usepackage{amssymb}
\usepackage{epstopdf}
\usepackage[russian, english]{babel}
\usepackage[all, 2cell]{xy}
\usepackage{epigraph}
\begin{document}

\pagestyle{empty}
%\pagenumbering{}
% Set book title
\title{\textbf{Axiomatic Method and Category Theory}}
% Include Author name and Copyright holder name
\author{Andrei Rodin}

% 1st page for the Title
%-------------------------------------------------------------------------------
\maketitle

% 2nd page, thanks message
%-------------------------------------------------------------------------------
\thispagestyle{empty}

\newpage

\pagestyle{plain}

% General definitions for all Chapters
%-------------------------------------------------------------------------------

% Define Page style for all chapters

% Delete the current section for header and footer

% Set arabic (1,2,3...) page numbering
\pagenumbering{arabic}

% Set double spacing for the text
\doublespacing

% Not enumerated chapter
%-------------------------------------------------------------------------------

% If the chapter ends in an odd page, you may want to skip having the page
%  number in the empty page
\frontmatter
\newpage
\chapter*{Preface}
I first learned about category theory about 20 years ago from Yuri I. Manin's course on algebraic geometry \cite{Manin:1970} when I was preparing my dissertation on Euclid's \emph{Elements} and was focused on studying Greek mathematics and classical Greek philosophy. Then I convinced myself that the mathematical category theory is philosophically relevant not only because of its name but also because of its content and because of its special role in the contemporary mathematics, which I privately compared to the role of the notion of \emph{figure} in Euclid's geometry. Today I have more to say about these matters. The broad historical and philosophical context, in which I studied category theory, is made explicit throughout the present book. My interest to the Axiomatic Method stems from my work on Euclid and extends through Hilbert and axiomatic set theories to Lawvere's axiomatic topos theory to the  Univalent Foundations of mathematics recently proposed by Vladimir Voevodsky. This explains what the two subjects appearing in the title of this book share in common.  

The next crucial biographical episode took place in 1999 when I was a young scholar visiting Columbia University on the Fulbright grant working on ontology of events under the supervision of Achille Varzi. As a part of my Fulbright program I had to make a presentation in a different American university, and I decided to use this opportunity for talking about the philosophical significance of category theory (I cannot now remember how exactly I married then this subject with the event ontology). Achille Varzi kindly arranged for me the invitation from Barry Smith to give a talk at his seminar on formal ontology in the SUNY in Buffalo. When I sent to Barry Smith my abstract he replied that nobody except probably Bill Lawvere will be able to understand my paper, and suggested to make the paper more accessible to the general audience. By that time I had already read some of Lawvere's papers but was wholly unaware about the fact that Lawvere worked in the same university and could attend my planned talk. So I took Smith's words for a joke. When I realized that this was not a joke I was very excited and, as it turned out, not without a reason because my meeting with Lawvere during this visit indeed determined the direction of my research for many years to come. This book is a summary of what I have achieved so far working in this direction. 

\textbf{Acknowledgement}. My main intellectual debt is to Bill Lawvere. I am also very grateful to all those friends and colleagues with whom I discussed the content of this book at various occasions and who gave me valuable advices and opportunities for its presentation. Leaving too many people out I mention Samson Abramsky, Vladimir Arshinov, Sergei Artemov, Mark van Atten, Steven Awodey, Andrej Bauer, Jean B\'enabou, Jean-Yves B\'eziau, Olivia Caramello, Pierre Cartier,  Tatiana Chernigovskaya, Anatoly Chussov, Bob Coecke, Maxim Djomin, Andreas D\"oring, Haim Gaifman, Ren\'e Guitart, Brice Halimi, Geoffrey Hellman, Jaakko Hintikka, Christian Houzel, Daniel Isaacson, Valery Khakhanjan, Anatole Khelif, Anders Kock, Roman Kossak, Anatoly Kritchevets, Marc Lachieze-Rey, Michiel van Lambalgen, Vladislav Lektorsky, Elena Mamchur, Yuri Manin, Per Martin-L\"of, Jean-Pierre Marquis, Fred Muller, John Mayberry, Colin McLarty, Arkady Nedel, Marco Panza, Vasily Perminov, Alberto Peruzzi, Richard Pettigrew, Alain Prout\'e, Oleg Prozorov, David Rabouin, Mehrnoosh Sadrzadeh, Gabriel Sandu, Dirk Schlimm, Valdislav Shaposhnikov, Ivahn Smadja, Sergei Soloviev, John Stachel, Jean-Jacques Szczeciniarz, Achille Varzi, Vladimir Vasyukov, Vladimir Voevodsky, Michael Wright and Noson Yanofsky.  My special personal debt is to Marina Brudastova without whose moral support this work could not be completed. 

\newpage
\thispagestyle{empty}
\tableofcontents

\chapter{Introduction}
\epigraph{Logical and mathematical concepts must no longer produce instruments for building a metaphysical ``world of thought'': their proper function and their proper application is only within the empirical science.}{\emph{Ernest Cassirer}}
\epigraph{Mathematics is a part of physics. It is a part of physics where experiments are cheap. [..] In the middle of the 20th century there were attempts to separate mathematics from physics. The results turned to be catastrophic.}{\emph{Vladimir Arnold}}

The main motivation of writing this book is to develop the view on mathematics described in the above epigraphs. Some 200 years ago this view used to be by far more common and easier to justify than today. It is sufficient to say that it made part of Kant's view on mathematics, and that Kant's view on mathematics remained extremely influential until the very end of the 19th century. When Cassirer defended this Kantian view in the beginning of the 20th century it was already  polemical. When Arnold defended it in the end of the 20th century and in the beginning of this current century it already sounded as an intellectual provocation, and so his words sound today. Kant, Cassirer and Arnold do not speak about the same mathematics: each speaks about mathematics of his own time. So the growing polemical attitude to their shared view reflects not only a change of the common opinion about the subject but a change of this subject itself. It is a common place that the modern mathematics is more abstract and more detached from physical experience than it used to be in Euclid's times and in Kant's times. When I say that I nevertheless want to defend the view on mathematics as a part of physics this means that I also want to contribute to changing the character of current mathematics, but not only to changing the common views about it. 

The above is a motivation behind this book but not its purpose. The purpose is much more limited. In order to justify the view on mathematics as part of physics I would need to write at least as much about physics as about mathematics. But this book is mainly about mathematics and about logic; physics is mentioned in it only occasionally. Yet more specifically I shall focus on the Axiomatic Method and Category Theory (including the categorical logic, which is a part of modern logic using category-theoretic methods). Let me explain why.

When Arnold talks about recent attempts to separate mathematics from physics he has in mind \emph{Elements of Mathematics} by Nicolas Bourbaki  \cite{Bourbaki:1939-1988} that aims at developing the whole of mathematics systematically from the first principles, i.e., on an axiomatic basis. Bourbaki's \emph{Elements} continue the long tradition of presenting renewed foundations of mathematics in the form of \emph{Elements}: this tradition begins with Euclid's \emph{Elements} (and earlier versions of Greek \emph{Elements} that have been lost) and continues through the whole history of mathematics until today. (I say a bit more about this tradition in the introductory part of Part \textbf{I}). Arnold sees the key to the problem in Bourbaki's Axiomatic Method, and takes a notoriously hostile attitude towards the Axiomatic Method in general. I observe on my part that the problem of separating mathematics from physics concerns the specific form of the Axiomatic Method used by Bourbaki  rather the the Axiomatic Method in general. It is clear, in particular, that Euclid's method does not produce the same effect. And I further observe that Bourbaki's Axiomatic Method is a version of Hilbert's Axiomatic Method presented in Hilbert's \emph{Foundations of Geometry} of 1899 \cite{Hilbert:1899}, which is another example of renewed  mathematical \emph{Elements} playing a more special but perhaps even more important role in the 20th century mathematics than Bourbaki's \emph{Elements}. So I conclude that the origin of Arnold's problem should be traced back at least to the beginning rather than only to the middle of the 20th century. This explains my focus on Axiomatic Method and its history.

Why Category Theory? The mathematical notion of category (which has no immediate relation to the philosophical notion widely known under this name) was invented in 1945 by Eilenberg and Mac Lane \cite{Eilenberg&MacLane:1945} for general purposes, some of which I explain in Chapter \textbf{8}, see also \cite{Kromer:2007} for details. In his thesis defended in 1963 \cite{Lawvere:1963} and a series of papers based on this thesis \cite{Lawvere:1964}, \cite{Lawvere:1966a}, \cite{Lawvere:1966b}, \cite{Lawvere:1967} Lawvere put forward a program of categorical (i.e., category-theoretic) foundations of mathematics and opened a new research field known today under the name of categorical logic, see \cite{Marquis&Reyes:2012} for the most recent historical account. Although Lawvere and other people who pursued the program of categorical foundations have never explicitly challenged Hilbert's Axiomatic Method (albeit they did and do challenge some special applications of this method, most importantly its applications in the standard axiomatic set theories) I shall try to show in this book that some recent works in categorical logic and new foundations of mathematics effectively modify Hilbert's Axiomatic Method and develop it in a wholly new direction. As it always happens in the intellectual history this new development continue some earlier developments, which I shall also take into account. In the last Chapter of this book I generalize upon these tendencies and describe a hypothetical New Axiomatic Method, which admittedly does not yet exist in the form of precise logical and mathematical procedure. I hope that my proposed general philosophical vision of this new method will contribute to its future technical development and also help to use it outside the pure mathematics and its philosophy.  

As the reader shall see the New Axiomatic Method establishes closer relationships between mathematics and physics and so suggests a solution of Arnold's problem. Although I cannot fully justify this claim in this book (because I am not going to discuss physics systematically) I do prepare a philosophical background for such a justification. The issue of relationships between mathematics and physics is a hardcore philosophical issue, and I believe that Arnold's problem cannot be solved without taking this philosophical issue seriously. Another hardcore philosophical issue that comes into the play as soon as one discusses the use of Axiomatic Method in mathematics is the relationships between mathematics and logic. This latter philosophical issue unlike the former is in the focus of this book. The main philosophical dilemma that I consider is, roughly, this: either (i) logic is fundamental in the sense that it gives us an independent access to an ideal space of logical possibilities where the actual world exists side-by-side with plenty of other possible worlds, which can be explored only mathematically, or as Cassirer insists in the above epigraph, (ii) logic and mathematics must stick to the actual world as we know it through empirical sciences, and by all means must avoid producing possible ``metaphysical worlds of thought'' even if these appear more logically coherent and more mathematical beautiful than our actual world. With many important reservations that this rough formulation requires I shall defend the latter view. The former view (which also obviously needs a more precise formulation) I call \emph{logicism}, and when it is applied to mathematics I call it \emph{mathematical logicism}. Beware that this meaning of ``mathematical logicism'' is broader than Russell's radical version of mathematical logicism according to which mathematics \emph{is} logic \cite{Russell:1903}. So a central purpose of this book is to refute mathematical logicism and defend an alternative way of thinking about logic and mathematics. 

Talking about these philosophical issues I would like to stress that I study primarily their implementation in mathematics. When in the beginning of the 20th century Cassirer, Russell and other people discussed hot philosophical issues concerning mathematics and logic they not only made general philosophical arguments but also referred to the actual state of affairs in their contemporary science and to the history of these subjects. They also often contributed themselves to the ongoing research in mathematics and logic. In this book I follow the same pattern of philosophical discussion paying a lot of attention to some recent mathematical works and to the history of the subject but without trying to make any mathematical contribution.        

Before I summarize the content of this book chapter by chapter let me say a few more words about its style and its methodology. I stick to the traditional idea  according to which philosophy and its history naturally combine together. When this view is applied to the philosophy of science and mathematics the result is sometimes called  the \emph{historical epistemology} \cite{Rheinberger:2010}. So what I am doing in this book can be described as a historical epistemology of logic and mathematics. However one important reservation is here in order. In my understanding the past history, the present state of affairs and the anticipated future of a given discipline are parts of the same whole. This whole can be described as the current state of affairs in a broader sense of the word, which includes both the historical reflection upon the past and the projection towards the future of the given discipline. When I talk in this book about mathematics and its philosophy I think about these subjects in this way. When such a view is called historical this should mean the attention to development of the given discipline but not the exclusive attention to its past.  
 
 Although I write about logic and mathematics I don't use myself any formal logical or other mathematical means for expressing and justifying my arguments. A century ago this point would be hardly worth mentioning but since using formal methods in philosophy in general and in philosophy of mathematics in particular is nowadays popular (particularly in the philosophical school that calls itself \emph{Analytic philosophy}) this point requires some explanations. Without going into a long discussion on this sensitive issue let me boldly express my believe that the natural language and the philosophical prose remain so far the best instruments for historical and philosophical work, or at least for the kind of such work that I want to do. The clarity and the exactness that formal methods bring to philosophy come with a price, which for my purposes is unacceptable. This price amounts to certain philosophical assumptions, without which these formal methods cannot work. I am not prepared to pay this price until I can see clearly these assumptions and thus know the price exactly. A philosophical and historical analysis of the notion of logical formalization is a part of my present project (see particularly Chapters \textbf{2} and \textbf{9}). Even if a formal theory of formalization is possible I cannot see that it can be useful for this purpose.  I shall not return to the question of using formal methods in philosophy in what follows but the reader will see that my analysis of the idea of logical formalization hardly supports the idea of using it as an universal instrument for philosophizing.  
 
Although I am not going to use formal methods for philosophical purposes the reader will find below a lot of rudimentary mathematics. Since this book is about mathematics, and a part of this book is about very recent mathematics, which still remains a work in progress (see \textbf{6.9 - 6.10}), this is not surprising. So let me explain my strategy of presenting the relevant mathematical content and mention some mathematical prerequisites for reading this book. My intention is to make this book readable both for a working mathematician interested in philosophy and history of this discipline and for a philosopher like myself, who studies (or wants to study) mathematics and its history, and finds a broad philosophical inspiration in this discipline. To present a fragment of modern mathematics to a wider audience is a very challenging task, which normally should not be combined with any philosophical agenda. I certainly do have a philosophical agenda, which I have already outlined earlier in this Introduction. This is why writing this book I have tried to reduce the burden of explaining mathematics to minimum. At the same time I tried to avoid any \emph{metaphoric} talk about mathematical concepts - even if some people would argue that any talk about mathematics outside the pure mathematics is doomed to be metaphoric. So I could not avoid the burden of explaining some mathematics completely but tried to use the most elementary examples and also tried to use some existing introductory expositions when such were available.  In each particular case I refer to the existing mathematical literature and chose this literature accordingly to my specific purpose.  

For the first superficial reading the given book is self-sustained and, as I hope, it gives a right idea of what I am after. A more attentive critical reading is by far more demanding. The ideal judge of this book is a working mathematician who is also a working philosopher and working historian of mathematics having some broader philosophical and scientific interests, which include some interest in physics, its history and its philosophy. I know several people who at some  degree of approximation fit this description but I rather imagine an average reader of this book as a person like myself who during these recent years has learnt  some philosophy, some mathematics and some history of both subjects, and who tries to make these ends meet. I shall say more about the mathematical prerequisites and give some suggestions for reading (in addition to references found in the main text) in the following summary of the Chapters. 

Part \textbf{I} of this book treats the history of Axiomatic Method. As I have already explained this history is not only about the past. Only Chapter  \textbf{1} on Euclid concerns what is indeed in the past (albeit in \textbf{1.5} I show that even in this case the past continues to live in the present); Chapter  \textbf{2} on Hilbert treats (in the original historical context) what remains today the standard notion of Axiomatic Method; Chapter  \textbf{4} on Lawvere treats what I suggest as a conceptual basis of the New Axiomatic Method. So these three Chapters of this book present, roughy, the past, the present and the anticipated future of the Axiomatic Method. Chapter  \textbf{3} is reserved for studying the fate of Hilbert's Axiomatic Method in the 20th century mathematics. 

Instead of trying to reconstruct a general history of Axiomatic Method, I decided to chose these three key figures and look at the relevant parts of their work more attentively. Although a historical discussion on Euclid found in Chapter  \textbf{1} may appear out of place in a book about today's mathematics it is important for me for several reasons. According to a common view (supported by Hilbert himself at some occasions), Hilbert's Axiomatic Method improves upon Euclid's method in terms of logical rigor and logical clarity. Of course, in such a general formulation this view can hardly be challenged. However in order to see how exactly this improvement on rigor and clarity has been achieved in the 20th century we need first to study Euclid's method on its own rights. This requires some special hermeneutical techniques, which are well-known to historian of mathematics but are less familiar to logicians, mathematicians and philosophers who also write about this subject. We shall see that in some respects Euclid's and Hilbert's method are different in principle, so that the difference between these methods does not reduce to differences in degrees of continuous magnitudes like rigor and clarity. In addition to my attempt to reconstruct Euclid's mathematical reasoning in its proper terms (and in some terms borrowed from Greek philosophy) I explain in this Chapter the relevance of Euclid's geometry to Kant's philosophy of mathematics. In the end of this Chapter I point to some Euclidean patterns of reasoning in the recent mathematics. The main textual reference in this Chapter is obviously Euclid's \emph{Elements}, which is now available in a new English translation  \cite{Euclid:2011}. An interested reader who would like to study the history of Greek mathematics more broadly and would like to better understand Euclid's special place in this history (this is an important subject that I wholly skip in this book) is advised to begin with \cite{Heath:1981}, \cite{Heath:2003} and then study more recent secondary literature.  

Chapter  \textbf{2} on Hilbert is also written in a historical style and contains extended quotes from Hilbert's writings. Although I leave outside the scope of my discussion most of the contemporary context of Hilbert's work  I follow the development of Hilbert's own ideas rather closely and distinguish in it several stages. In its narrow historical aspect my treatment of Hilbert's work contains nothing original. However I also make an attempt to reconstruct the history of some relevant notions (or at least to keep track of their changing meaning) including the notion of being formal. This historical discussion is combined with an explanation of Hilbert's Formal Axiomatic Method, which can be used by a non-mathematical reader for the first acquaintance  with this basic method of modern mathematical reasoning. Someone well acquainted with this method will find here an analysis of certain assumptions required by this method, which remain tacit when this method becomes an intellectual habit and is used automatically. I shall pay a lot of attention to philosophical remarks made by Hilbert in his presentations of Axiomatic Method trying to reconstruct Hilbert's thinking and its philosophical motivation. I also discuss in this Chapter some related subjects including the notion of logicality, diagrammatic and symbolic thinking and some others. This Chapter presents (in its historical original form) the core notion of modern Formal Axiomatic Method, which I contrast in what follows to more traditional Euclid's method, on the one hand, and to some later versions of Axiomatic Method including the anticipated New Axiomatic Method, on the other hand. 

The main suggested reading for Chapter  \textbf{2} is Hilbert's \emph{Foundations of Geometry}, which exist in multiple editions including the English edition \cite{Hilbert:1950} and some later English editions. I highly recommend this reading also to a non-mathematical reader of this book because the real subject-matter of this short masterpiece is the Axiomatic Method itself rather than geometry, and so this short book can be used as a shortcut to the modern style of mathematical thinking. For a later more developed systematic presentation of Formal Axiomatic Method and its underlying philosophy I refer the reader to Tarski's textbook \cite{Tarski:1941}. This  textbook presents in a very clear form a philosophical view on logic and mathematics that I discuss in my present book.  

In Chapter  \textbf{3} I talk about applications of Hilbert's Axiomatic Method in the 20th century mathematics and stress the fact that it has  hardly ever been used in its original form and for its originally intended purpose. I discuss from this point view some formal studies of axiomatic set theories, Bourbaki's \emph{Elements of Mathematics} \emph{Bourbaki:1939-1988} and more specifically an unpublished Bourbaki's draft \cite{Bourbaki:19??}. My main observation amounts to saying that both the modern set theory and Bourbaki's structural mathematics can be described in Hilbert's terms as a \emph{metatheory} or in Tarski's terms as a \emph{model theory} of certain Hilbert-style axiomatic theory or, more typically, of a number of such theories. Since this metatheory or model theory itself is developed by some other means (i.e., \emph{not} axiomatically in Hilbert's sense) one can say that the mainstream mathematics widely applies Hilbert's Formal Axiomatic Method only with a pinch of salt. In the mainstream structural mathematics of the 20th century this method serves as a method of definition and constructing new concepts rather than method of building deductive theories. On the basis of this observations I criticize Hilbert's Axiomatic Method arguing that it is not apt to support mathematical theories useful in the modern physics. Finally I consider in this Chapter Tarski's topological model of intuitionistic propositional logic \cite{Tarski:1956} and stress its unusual character: although, technically speaking, there is no big difference between modeling a given formal theory and modeling a given logical calculus, philosophically it makes a huge difference and requires a rethinking of the whole idea of Axiomatic Method. Although Tarski himself does not draw from this work such far-reaching conclusions I use this example in the following Chapter as a historical prototype of the New Axiomatic Method. 

In addition to the literature referred to in Chapter  \textbf{3} I suggest reading the classical introduction \cite{Fraenkel&Bar-Hillel&Levy:1973} to the modern axiomatic set theory including its last philosophical chapter, and Galileo's \emph{Two New Sciences} \cite{Galilei:1974} where the author stresses the constructive experimental character of the New Science against the background of the earlier Scholastic patterns of doing science.  

Chapter \textbf{4} plays a central role in this book because here I first introduce the notion of category and discuss a new notion of Axiomatic Method, which emerges in category theory and, more specifically, in categorical logic. Although  categorical logic is already a well established subject (see \cite{Marquis&Reyes:2012} for a historical introduction) I decided to follow here the pattern of the first two Chapters and focus my attention on the work of one particular person, namely Lawvere, who founded this discipline in 1960-ies; as before I combine here a historical and a systematic orders of presentation and pay a minute attention to Lawvere's philosophical comments found throughout his writings. After presenting Lawvere's categorical axiomatization of (the category of) sets \cite{Lawvere:1964} and of the category of categories \cite{Lawvere:1966a}, which gives the first idea of using the category theory for axiomatization, I turn to Lawvere's critique of the standard Formal Axiomatic Method as ``subjective'' and explain his idea of \emph{objective} conceptual logic realized by category-theoretic means. I begin this latter discussion by considering two Lawvere's papers \cite{Lawvere:1966b}, \cite{Lawvere:1967} that mark the birth of the categorical logic, and in the same context explain Lawvere's notion of quantifiers as adjoint functors to the substitution functor. Then I make a digression on Curry's \emph{combinatorial logic}, type theory and the so-called \emph{Curry-Howard correspondence}, and show how these conceptual developments combine in Lawvere's notion of Cartesian closed category. Then after a brief discussion on Lawvere's notions of hyperdoctrine (that conceptually connects to the discussion on homotopy type theory found in  \textbf{6.9}) and functorial semantics (further discussed in  \textbf{9.2}) I turn to philosophical issues and discuss the role of Hegel's dialectical logic in Lawvere's thinking, which Lawvere stresses himself at many instances. Here I provide a philosophical reconstruction of Hegel's distinction between the \emph{objective} and the \emph{subjective} logic and then describe how this philosophical distinction is realized by Lawvere with the technical means of categorical logic. This discussion helps me then for interpreting the groundbreaking paper \cite{Lawvere:1970a}  where Lawvere suggests his axiomatization of topos theory and demonstrates the strength of his notion of internal logic of a given category. In the last Chapter  \textbf{9} I use Lawvere's axiomatization of topos theory as a basic example of the new axiomatic approach, which I try to describe in general terms under the title of New Axiomatic Method.

For a better understanding of Chapter  \textbf{4} it would be useful if the reader get some knowledge of basic category theory beforehand (albeit this is not an absolutely necessary requirement and the reader can also follow references during the reading). For a non-mathematical reader or a reader with a modest mathematical background I recommend \cite{Lawvere&Schanuel:1997} and   \cite{Lawvere&Rosebrugh:2003} co-authored by Lawvere as a very accessible introduction into the subject. For a mathematical reader not familiar with categorical logic I recommend \cite{MacLane&Moerdijk:1992} that covers most of the mathematical material that I discuss in this Chapter (but unfortunately skips hyperdoctrines). There is a huge gap in terms of required mathematical skills between these two suggested readings and by the present day this gap has not been yet filled in spite of many very valuable attempts such as \cite{Reyes&Reyes&Zolfaghari:2004}. I believe that there is a principle and not only technical and pedagogical difficulty involved with the project of writing a fairly elementary introduction to category, topos theory and categorical logic. The problem is that the elementary introductions like \cite{Lawvere&Schanuel:1997},   \cite{Lawvere&Rosebrugh:2003} and \cite{Reyes&Reyes&Zolfaghari:2004} begin with considering the category of finite sets, which are first introduced naively as bags of dots and then are treated in terms of their maps. Although such an introduction is geometrical in its character the basic geometry reduces here to the geometry of bags of dots, which is a geometry of a very special sort. A genuine continuous geometry appears then only at the much more advanced level and in a much more abstract form of Grothendieck topology and Grothendieck topos, which are systematically treated in \cite{MacLane&Moerdijk:1992} and other books of the same advanced level. So it still remains, in my view, a challenging task to follow Hilbert's example and rewrite Euclidean or other simple intuitive geometry in new categorical terms. Voevodsky Univalent Foundations discussed in \textbf{6.10} appear to be a step in this direction. 

Talking about elementary introductions to category theory and topos theory I would like also to mention \cite{McLarty:1992} by McLarty. The expression ``elementary' theory'' in the title does not stand for being easy to grasp by a beginner but is used in the technical sense of being a first-order theory in the sense of modern logic and the standard Formal Axiomatic Method. This book is a systematic presentation of category and topos theory which fully complies with the requirement of Formal Axiomatic Method and at the same time treats the internal logic of a given topos and the idea of internal description of a given topos with its internal language. So for a logically-minded philosopher habituated to formal methods this book may also serve as an introduction into the subject. I would like to stress however that since in the present book I discuss specific features of Lawvere's axiomatic thinking, which fall apart from the standard Formal Axiomatic Method, studying McLarty's book does not replace studying Lawvere's original works even if, formally speaking, McLarty's book fully covers the same subject.

Part \textbf{II} is devoted to the notion of identity (in mathematics). This may appear as a side subject with respect to the general theme of this book but it is actually not. A mathematical logicist argues like this: in order to build a mathematical theory in an axiomatic form one needs first to fix some basic logical notions like that of being the \emph{same} (or being equal). Unless this is done beforehand and quite independently from the content of any particular mathematical theory, so the argument goes, no axiomatic construction of mathematical theories is possible. A similar point can be made, of course, about other logical notions including logical connectives ``and'', ``or'', the notion of logical inference, of truth-value, etc. This standard logicist argument does not go through in the case of categorical logical, or at least it does not go through immediately, because the categorical logic \emph{internalizes} the logical notions, i.e., reconstructs them in terms of a given mathematical theory (see   \textbf{4.9} and  \textbf{9.3}). This applies to logical connectives, the relation of inference, quantifiers, truth-values and to some other logical notions. It also applies to the logical identity relation but this case turns to be both more difficult and more mathematically and philosophically interesting than other cases. So I treat it systematically in the two consequent Chapters making the Part  \textbf{II}. 

In Chapter  \textbf{5} I consider the question of identity/equality in mathematics in general beginning with some naive observations and historical examples. In particular, I briefly consider Plato's view according to which the mathematical equality is a weak form of strict identity: while the latter applies only the ideal world of Forms the former applies in the world of mathematics, which takes an intermediate position between the world of immutable Forms and the world of changing material beings. Plato's theory is an echo of the modern mathematical structuralism discussed later in Chapter  \textbf{8}. In Chapter  \textbf{5} I also show the significance of discussions about identity in mathematics in Frege's and Russell's works for establishing the logicist view on mathematics in the end of the 19th and the beginning of the 20th century. Then I turn to more theoretical subjects including a discussion on classes and individuals, and a discussion of the distinction between logical extension and logical intension. This Chapter resumes with a discussion on Martin-L\"of's intuitionistic type theory \cite{Martin-Lof:1984} that provides a theory of identity types, which is very non-trivial in the intensional case. I compare Martin-L\"of's approach to identity with Frege's approach and reconsider Frege's famous \emph{Venus} example through the optics of Martin-L\"of's type theory.

Chapter \textbf{6} continues to treat the issue of identity but this time with new approaches coming from category theory and some related fields. In the beginning of this Chapter I stress the conceptual similarity and the conceptual difference between the logical notion of relation and geometrical notion of transformation aka mapping or simply map. On this basis I re-introduce the notion of category with a naive geometrical example, stress the geometrical origin of categorical thinking and the relationships between category theory and Klein's \emph{Erlangen Program}. (I come back to this topic in  \textbf{8.6}). Then I turn to more advanced geometrically motivated categories and show how they realize the idea of identity as a map (rather than a relation). In particular, I consider B\'enabou's \emph{fibered categories} \cite{Benabou:1985} and higher categories (aka $n$-categories) - first in an abstract form and then in the geometrical form of homotopy categories. So I approach the hot subject of \emph{homotopy type theory}, which brings together identity types of Martin-L\"of's type theory and the geometrical approaches to identity and the homotopical higher category theory. When I began to study these two subjects about ten years ago the precise mathematical connection between them was not yet established and the mathematical discipline of homotopy type theory did not yet exist. So it was for me a great relief to learn that these ideas combine not only at the level of speculative philosophy but also in precise mathematical terms. I conclude this Chapter with a presentation of Voevodsky's new foundations of mathematics that he calls Univalent Foundations \cite{Voevodsky:2010},\cite{Voevodsky:2011}. In Chapter \textbf{9} I refer to the Univalent Foundations as an example of a new form of axiomatic presentation along with the example of Lawvere's axiomatic topos theory. 

As a general mathematical reading for Part  \textbf{II} I recommend Leinster's book \cite{Leinster:2004} on higher category-theory, which has great pedagogical advantages, Granstrom's book \cite{Granstrom:2011} on type theory, which also provides a philosophical perspective on this theory, Jacob's book \cite{Jacobs:1999} that stresses the link between categorical logic and type theory. The homotopy type theory has been not yet exposed in textbooks but there are very clear expository papers \cite{Awodey&Warren:2009} and \cite{Awodey:2010}.  

Last  Part \textbf{III} of the book treats two different subjects, which fall under the scope of Hegel-Lawvere's  distinction between objective and subjective features of logic and mathematics. In Chapter  \textbf{7} I discuss the issue of mathematical intuition from a historical perspective and argue using some historical examples that mathematical intuitions change through the historical time at least as rapidly as do mathematical concepts. The main purpose of this Chapter is to refute the popular view according to which mathematics always develops by increasing its degree of abstractness and according to which the highly abstract character of modern mathematical concepts does not allow for a faithful intuitive representation in principle. I suggest an alternative picture of the historical development of mathematics where concepts and intuitions develop side-by-side but sometimes the conceptual development takes over the intuitive development and sometimes, on the contrary, the intuitive development takes over the conceptual one. 

I expect that a phenomenologically-minded philosophical reader may object that what I discuss is not the strict philosophical notion of intuition but rather a commonsensical meaning of the word ``intuition'' as a bunch of helpful analogies borrowed from the everyday life or elsewhere. I argue in this Chapter that the changing mathematical intuition that I describe qualifies at least as intuition in Kant's sense of the term. The lack of discussion of Husserl's views is indeed a significant lacuna of this Chapter that I cannot easily fix.  So I leave it for a future work. 

Although I wholly share Lawvere's Hegelian view concerning the objective character of scientific logic (which perfectly squares with Cassirer's view on the place and the role of mathematics and logic expressed in the above epigraph) I also stress the role of the subjective intuition because it provides the necessary link that connects the pure mathematics to the individual sensual experience to the scientific empirical methods to the whole body of empirical science. Without such a link Hegel's objective dialectical logic too easily turns into a new form of speculative dogmatic metaphysics wholly detached from reality. One may suggest that since the dogmatic dialectics is an obvious oxymoron it cannot refer to anything real. But the dialectical logic quite rightly protects one from such naive conclusions made on abstract logical grounds: as a matter of painful historical fact the examples of dogmatic misuse of philosophical dialectics are abound
\footnote{
Unlike the older forms of dogmatism the more recent dialectical dogmatism does not use any fixed system of beliefs but enforces a permanent organized change of one's beliefs on changing pragmatic grounds (political, economical, etc.).}. 

In Chapter  \textbf{8} I discuss structuralism including its mathematical variety. Considering structuralism as a suggestive idea rather than a system of stable philosophical views I argue against the received view according to which category theory brings about a new variety of structuralism and provides a new framework for developing structural mathematics. I recognize the role of structural thinking in the development of category theory and describe this role in this Chapter.  In particular, I elaborate on Eilenberg and Mac Lane's idea of category theory    as a continuation of Klein's \emph{Erlangen Program} \cite{Eilenberg&MacLane:1945}. This very analogy allows me to specify the crucial difference between Klein's structural thinking and new categorical thinking: when groups are generalized up to categories the notion of invariant structure is replaced by the notion of covariant or contravariant functor. I argue that the structuralist thinking about functoriality in terms of preservation of invariant structures is, generally, inappropriate; then I suggest a different philosophical view (or rather another suggestive idea) where the notion of functoriality (i.e., of co- and contravariance) becomes central. Although this conceptual development begins with a mere generalization of the structuralist  \emph{Erlangen Program} it brings about a new view, which is very unlike the structuralist view. In the end of this Chapter  (\textbf{8.8}) I suggest a purely geometrical way of thinking about categories alternative to the more convenient way of thinking about categories as categories of structures. The basic idea here is thinking of geometrical objects as maps from types (of geometrical objects) to spaces. I demonstrate this approach with some elementary examples from the 19th century geometry. Thus in my suggested post-structuralist picture the notion of object (this time understood as a map) becomes once again central. 

The conceptual change described in Chapter  \textbf{8} affects not only the choice of structures explored with the Formal Axiomatic Method but also this method itself. So in the concluding Chapter  \textbf{9} I make the long-promised attempt to describe the New Axiomatic Method more systematically. I first describe the two basic functions of Axiomatic Method, which Lawvere calls the \emph{unification} and the \emph{concentration}. Here I contrast the unificatory strategy of the New Method to the more traditional unificatory strategy of Formal Axiomatic Method, which has a structuralist and a logicist underpinning. Then I describe the \emph{concentration} part, which turns to be more traditional and in a new form reproduces some features of Euclid's Axiomatic Method. The most original part of the New Axiomatic Method is, of course, its logical part, which involves the notion of internal logic. Generalizing on works of Lawvere and Voevodsky I describe here in general terms a way of using the internal logic of some given category (which is construed in intuitive geometrical terms at the first step of the axiomatic construction) for improving upon the construction of this very category and providing it with some deductive structure. This way of using logic for building mathematical theories suggests a new way of thinking about the role of logic in mathematical theories, which is very unlike Hilbert's and Tarski's. 

In my suggested approach logic is designed along with the rest of conceptual construction rather than used as a ready-made foundation for making further mathematical constructions. One may think that the freedom of making up logical calculi added to the freedom of making up new axiomatic mathematical theories (assured already by Hilbert) only reinforce the inflation of the ``metaphysical world of thought''. In fact the New Axiomatic Method prevents this inflation in two different ways. First, by taking into account the objective meaning of the category of interest (which can be, for example, a spatiotemporal category used in physics) and, second, by requiring the relevant logic to be the internal logic of this given category. While the former feature is at some degree also compatible with the standard Formal Axiomatic Method the latter feature is a genuinely original contribution of the New Method.The New Method no longer reduces the function of  logical formalization to a logical censorship; instead logic is used here as a flexible tool for the internal conceptual reconstruction. 

 An important part of my argument consists of pointing to Lawvere's and Voevodsky's works as applications of  
this New Method, and stressing the fact that in both cases it allows for a remarkable conceptual simplification and clarification of otherwise difficult and conceptually problematic theories. Since in both cases the relevant logic is internal with respect to its base category this logic inherits the objective meaning of this base category. This allows me to suggest that the New Axiomatic Method may help to bridge the gap between mathematics and physics created and justified by the standard Formal Axiomatic Method and by the logicist view on mathematics that underpins this standard method. Notwithstanding my critique of Hilbert's version of Axiomatic Method developed throughout in this book, I believe (contra Arnold) that Hilbert was perfectly right when he described this method as ``the basic instrument of all research'' (\cite{Hilbert:1927}, p. 467) and when he said that ``[t]o proceed axiomatically means [..] nothing else than to think with consciousness'' (\cite{Hilbert:1922} p. 1120)

\mainmatter
%Finally, include the ToC
% First enumerated chapter
%-------------------------------------------------------------------------------
\part{A Brief History of the Axiomatic Method}

In his famous address ``Axiomatic Thought'' delivered before the Swiss Mathematical Society in Zurich in 1917 Hilbert says:

\begin{quote}
If we consider a particular theory more closely, we always
see that a few distinguished propositions of the field of
knowledge underlie the construction of the framework of
concepts, and these propositions then suffice by
themselves for the construction, in accordance with
logical principles, of the entire framework. [..]  These fundamental propositions can be regarded [..]  as
the \underline{axioms} of the individual fields of knowledge : the
progressive development of the individual field of
knowledge then lies solely in the further logical
construction of the already mentioned framework of
concepts. This standpoint is especially predominant in
pure mathematics. [.. A]nything at all that can be the object of scientific
thought becomes dependent on the axiomatic method,
an thereby indirectly on mathematics. (\cite{Hilbert:1918}, p. 1108-1115)
\end{quote}

In a different paper  the author makes a further epistemological claim:

\begin{quote}
The axiomatic method is and remains the indispensable
tool, appropriate to our minds, for all exact research in
any field whatsoever: it is logically incontestable and at
the same time fruitful. [..] To proceed axiomatically means
in this sense nothing else than to think with
consciousness. (\cite{Hilbert:1922} p. 1120)
\end{quote}

Although Hilbert's enthusiasm about the Axiomatic Method and his high esteem of the role of this method in science may be not  universally accepted today, the modern notion of axiomatic theory remains shaped by Hilbert's works; his \emph{Grundlagen der Geometrie} (Foundations of Geometry) first published in 1899 \cite{Hilbert:1899} still serves as a paradigm of axiomatic mathematical theory. As soon as this method is understood in the above general terms one may think that it has been practiced by mathematicians since the early days of their discipline. Indeed in the \emph{Introduction} to his \emph{Foundations of Geometry} of 1899  \cite{Hilbert:1899} Hilbert states the following:

\begin{quote}
Geometry, like arithmetic, requires for its logical development only a small number of
simple, fundamental principles. These fundamental principles are called the axioms of
geometry. The choice of the axioms and the investigation of their relations to one another
is a problem which, since the time of Euclid, has been discussed in numerous excellent
memoirs to be found in the mathematical literature. This problem is tantamount to the
logical analysis of our intuition of space. (Hereafter \cite{Hilbert:1899} is quoted in English translation \cite{Hilbert:1950}) 
\end{quote}

Notice Euclid's name is the above quote. Evidently Hilbert had in mind Euclid's \emph{Elements} when he prepared his \emph{Foundations of Geometry} for publication.  Hilbert aims at developing Euclidean geometry on a wholly new conceptual basis. \footnote {I agree with David Rowe when he says that ``The reform of geometry that [Hilbert] envisaged in \emph{Grundlagen der Geometrie} was primarily conceived as a renewal of the fundamental structures of classical Euclidean geometry.'' (\cite{Rowe:2000}, p.71)} 
In this sense Hilbert's \emph{Foundations} of 1899 qualifies as a fairly revolutionary work. However one should not forget  that rewriting geometrical chapters of Euclid's \emph{Elements} in new terms is itself an old and well-establish tradition in the history of mathematical thought. Hilbert's \emph{Foundations of Geometry} (as well as Bourbaki's open-ended \emph{Elements of Mathematics} \cite{Bourbaki:1939-1988} produced later in the 20th century) make part of this long-term tradition and can be compared with such groundbreaking works of earlier generations as, for example,  \emph{Restored Euclid}  by Borelli (1658) \cite{Borelli:1658}, \emph{New Elements of Geometry} by Arnauld (1667)\cite{Arnauld:1683} and  \emph{Euclid Freed from All Flaws} by Saccheri (1733)\cite{Saccheri:1733}. Thus the Hilbertian revolution that still strongly influences today's  mathematical practice is certainly not the first revolution of this sort and hopefully not the last one. 

Hilbert thinks of his new version of Axiomatic Method as a development of and improvement over Euclid's method of theory-building. Surely Hilbert is aware about the fact that his method is not the same as Euclid's; we shall see that Hilbert in fact quite precisely points to the key difference (see \textbf{2.6}). The purpose of Chapter \textbf{1} is to describe this difference more precisely and more systematically.  In Chapter \textbf{2}, I focus on Hilbert's work and compare Hilbert's approach to Euclid's. In Chapter \textbf{3}, I consider applications of Hilbert's Axiomatic Method in the 20th century mathematics and, in particular, in Bourbaki. In Chapter \textbf{4}, I discuss Lawvere's work and show how some basic features of Euclid's approach deliberately ignored by Hilbert get a new life in today's categorical logic.

\chapter{Euclid: Doing and Showing}
Reading older mathematical texts always involves a hermeneutical dilemma: in order to make sense of the mathematical content of a given old text one wants to interpret it in modern terms; in order to see the difference between the modern mathematical thinking and older ways of mathematical thinking one wants to avoid anachronisms and understand the old text in its own terms \cite{Unguru:1975}. Any scholar studying older mathematics needs to find a way between the Scylla of ``antiquarianism'' that seeks the scholar's conversion into a person living during a different historical epoch, and the Charybdis of radical ``presentism'' that finds in older texts nothing but a minor part of today's standard mathematical curricula and wholly ignores the historical change of basic patterns of mathematical thinking
\footnote{Being between Scylla and Charybdis is an idiom deriving from Greek mythology. Scylla and Charybdis were mythical sea monsters noted by Homer. Scylla was rationalized as a rock shoal (described as a six-headed sea monster) on the Italian side of the strait and Charybdis was a whirlpool off the coast of Sicily. They were regarded as a sea hazard located close enough to each other that they posed an inescapable threat to passing sailors; avoiding Charybdis meant passing too close to Scylla and vice versa. (after WikipediA)}
. My way through the channel is the following. I read Euclid's text verbatim (relying on Heiberg's edition of the original Greek \cite{Euclides:1883-1886} and using Fitzpatrick's new English translation \cite{Euclid:2011}), consider its most important modern interpretations (including overtly anachronistic ones), criticize some of these interpretations on the basis of textual evidences, and finally suggest some alternative interpretations. In order to prevent the risk of losing the main argument behind the following historical details I formulate now my general conclusion. Contrary to popular opinion Euclid's geometry is not a system of propositions some of which have a special status of axioms while some other are derived from the axioms according to certain rules of logical inference. It can be rather described after Friedman as ``a form of rational argument'' (\cite{Friedman:1992}, p. 94)
\footnote{See the full quote from Friedman in the end of Section \textbf{2.5}.}, where some non-propositional content (including non-propositional first principles) is indispensable. Precipitating what follows (see particularly \textbf{2.6}) let me mention that certain non-propositional principles also make part of modern formal theories in the form of \emph{syntactic rules}. As we shall now see in Euclid the non-propositional aspect of mathematical reasoning plays a more prominent role.

\section{Demonstration and ``Monstration''} 

All Propositions of Euclid's \emph{Elements} (with few easily understandable exceptions) fit into the scheme described by Proclus in his \emph{Commentary} \cite{Proclus:1970} as follows: 

\begin{quote}
Every Problem and every Theorem that is furnished with all its parts should contain the following elements: an \emph{enunciation}, an \emph{exposition}, a \emph{specification}, a \emph{construction}, a \emph{proof}, and a \emph{conclusion}. Of these \emph{enunciation} states what is given and what is being sought from it, a perfect \emph{enunciation} consists of both these parts. The \emph{exposition} takes separately what is given and prepares it in advance for use in the investigation. The \emph{specification} takes separately the thing that is sought and makes clear precisely what it is. The \emph{construction} adds what is lacking in the given for finding what is sought. The \emph{proof} draws the proposed inference by reasoning scientifically from the propositions that have been admitted. The \emph{conclusion} reverts to the \emph{enunciation}, confirming what has been proved. (\cite{Proclus:1970}, p.203, italic is mine)
\end{quote}

It is appropriate to notice here that the term ``proposition'', which is traditionally used in translations as a common name of Euclid's problems and theorems, is not found in the original text of the \emph{Elements}: Euclid numerates these things throughout each Book without naming them by any common name. (The reader will shortly see why this detail is important.) The difference between problems and theorems is explained in \textbf{1.4}  below. Let me now show how Proclus' scheme applies to Proposition 5 of the First Book (Theorem 1.5), which is a well-known theorem about angles of the isosceles triangle. References in square brackets are added by the translator; some of them will be discussed later on. Words in round brackets are added by the translator for stylistic reason. Words in angle brackets are borrowed from the above Proclus' quote. Throughout this Chapter I write these words in italic when I use them in Proclus' specific sense. 

\paragraph{[\emph{enunciation:}]}
\begin{quote}
For isosceles triangles, the angles at the base are equal to one another, and if the equal straight lines are produced then the angles under the base will be equal to one another.
\end{quote}

\paragraph{[\emph{exposition}]:}

\begin{quote}
Let $ABC$ be an isosceles triangle having the side $AB$ equal to the side $AC$; and let the straight lines $BD$ and $CE$ have been produced further in a straight line with $AB$ and $AC$ (respectively). [Post. 2].
\end{quote}

\begin{center}
\includegraphics [scale=0.5]{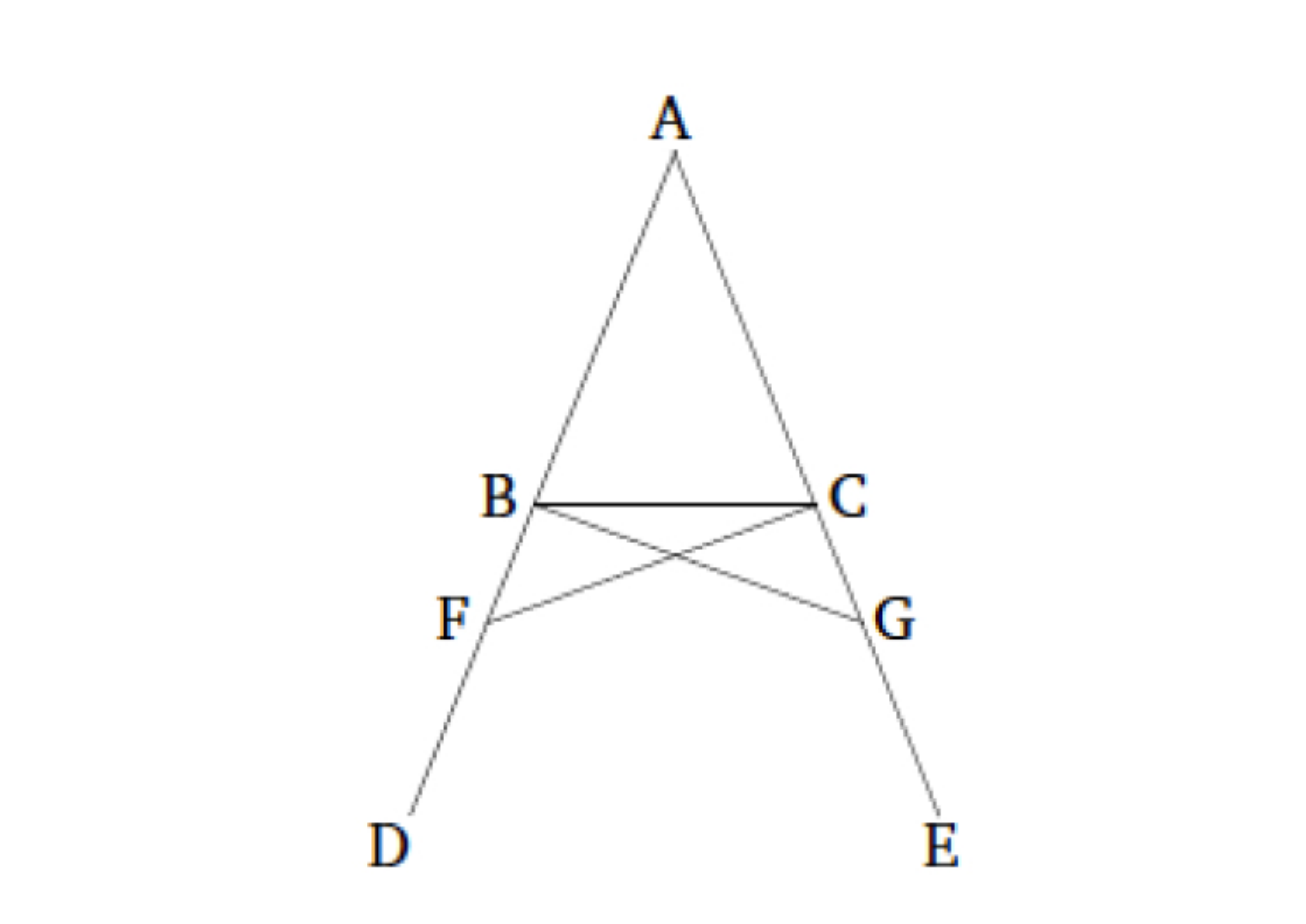}
\end{center}

\begin{center}
Fig. 1.1
\end{center}

\paragraph{[\emph{specification}:]}

\begin{quote}
I say that the angle $ABC$ is equal to $ACB$, and (angle) $CBD$ to $BCE$.
\end{quote}

\paragraph{[\emph{construction}:]}

\begin{quote}
For let a point $F$ be taken somewhere on $BD$, and let $AG$  have been cut off from the greater  $AE$,  equal to the lesser $AF$ [Prop. 1.3]. Also, let the straight lines $FC$, $GB$ have been joined. [Post. 1] 
\end{quote}

\paragraph{[\emph{proof}:]}

\begin{quote}
In fact, since $AF$ is equal to $AG$, and $AB$ to $AC$,
the two (straight lines) $FA$, $AC$ are equal to the two (straight lines) $GA$, $AB$, respectively. They also encompass a common angle $FAG$.
Thus, the base $FC$ is equal to the base $GB$, and the triangle $AFC$ will be equal to the triangle $AGB$,
and the remaining angles subtended by the equal sides will be equal  to the corresponding remaining angles [Prop. 1.4]. (That is) $ACF$ to $ABG$, and $AFC$ to $AGB$.
And since the whole of $AF$ is equal to the whole of $AG$,
within which $AB$ is equal to $AC$, the remainder $BF$ is thus equal to the remainder $CG$ [Ax.3].
But $FC$ was also shown (to be) equal to $GB$.
So the two (straight lines) $BF$, $FC$ are equal to the two (straight lines) $CG$, $GB$ respectively, and the angle $BFC$ (is) equal to the angle $CGB$,
while the base $BC$ is common to them. Thus
the triangle $BFC$ will be equal to the triangle $CGB$,
and the remaining angles subtended by the equal sides will be equal to the corresponding remaining angles [Prop. 1.4]. Thus $FBC$ is equal to $GCB$, and $BCF$ to $CBG$. Therefore, since the whole angle $ABG$ was shown (to be) equal to the whole angle $ACF$, within which 
$CBG$ is equal to $BCF$,
the remainder $ABC$ is thus equal to the remainder $ACB$ [Ax. 3].
And they are at the base of triangle $ABC$.
And $FBC$ was also shown (to be) equal to $GCB$.
And they are under the base.\end{quote}

\paragraph{[\emph{conclusion}:]}

\begin{quote}
Thus, for isosceles triangles, the angles at the base are
equal to one another, and if the equal sides are produced
then the angles under the base will be equal to one another.
(Which is) the very thing it was required to show.
\end{quote}

An obvious difference between Proclus' analysis of the above theorem and its usual modern analysis is the following. For a modern reader the proof of this theorem begins with Proclus' \emph{exposition} and includes Proclus' \emph{specification}, \emph{construction} and \emph{proof}. Thus for Proclus the \emph{proof} is only a part of what we call today the proof of this theorem. Also notice that Euclid's theorems conclude with the words ``which ... was required to \emph{show}'' (as correctly translates Fitzpatrick) but not with the words ``what it was required to \emph{prove}'' (as inaccurately translates Heath \cite{Heath:1926}). The standard Latin translation of this Euclid's formula as \emph{quod erat demonstrandum} is also inaccurate. These inaccurate translations conflate two different Greek verbs: ``apodeiknumi'' (English ``to prove'', Latin ``demonstrare'') and ``deiknumi'' (English ``to show'', Latin ``monstrare''). The difference between the two verbs can be clearly seen in the two Aristotle's \emph{Analytics}: Aristotle uses the verb ``apodeiknumi'' and the derived noun ``apodeixis'' (proof) as technical terms in his syllogistic logic, and he uses the verb ``deiknumi'' in a broader and more informal sense when he discusses epistemological issues (mostly in the \emph{Second Analytics}). Without trying to trace here the history of Greek logical and mathematical terminology and speculate about possible influences of some Greek writers on some other writers, I would like to stress the remarkable fact that Aristotle's use of verbs ``deiknumi'' and ``apodeiknumi'' agrees with Euclid's and Proclus'. In my view this fact alone is  sufficient for taking  seriously the difference between the two verbs and distinguishing between \emph{proof} and ``showing'' (or otherwise between \emph{demonstration} and \emph{monstration})\footnote{As far as mutual influences are concerned two things are certain: (i) Proclus read Aristotle and (ii) Aristotle had at least a basic knowledge of the mathematical tradition, on which Euclid later elaborated in his \emph{Elements} (as Aristotle's mathematical examples clearly show \cite{Heath:1949}).  It remains unclear whether  Aristotle's work could influence Euclid. In my view this is unlikely. However Aristotle's logic certainly played an important role in later interpretations and revisions of Euclid's \emph{Elements}. I leave this interesting issue outside of the scope of this book.}.

One may think that the difference between the current meaning of the word ``proof'' in today's mathematics and logic and the meaning of Proclus'  \emph{proof} (Greek ``apodeixis'') is a merely terminological issue, which is due to difficulties of translation from Greek to English. I shall try now to show that this terminological difference points on a deeper problem, which is not merely linguistic. In today's logic the word ``proof'' stands for a logical inference of certain conclusion from some given premises. In fact this is what by and large was meant by proof also by Aristotle and Proclus. Indeed, looking at the \emph{proof} (in Proclus' sense) of Euclids Theorem 1.5 we see that it also qualifies as a proof in the modern sense: we have here a number of premises (which I make explicit in the next Section) and certain conclusions derived from those premises. It is irrelevant now whether or not this particular inference is valid according to today's logical standards; what I want to stress here is only the general setting that involves some premises, an inference (probably invalid) and some conclusions. This core meaning of the word ``proof'' (Greek ``apodeixis'') hardly changed since Proclus' times.

So we get a problem, which is clearly not only terminological: Is it indeed justified to describe the \emph{exposition}, the \emph{specification} and the \emph{construction} as elements of the proof or one should rather follow Proclus and consider these things as independent constituents of a mathematical theorem? 

The question of \emph{logical significance} of the \emph{exposition}, the \emph{specification} and the \emph{construction} in Euclid's geometry has been discussed in the literature; in what follows I shall briefly describe some tentative answers to it. However before doing this I would like to stress that this question may be ill-posed to begin with. As far as one assumes, first, that the theory of Euclid's  \emph{Elements} is  (by and large) sound and, second, that any sound mathematical theory is an axiomatic theory in the modern sense, then, in order to make these two assumptions mutually compatible, one has to describe the \emph{exposition}, the \emph{specification} and the \emph{construction} of each Euclid's theorem as parts of the proof of this theorem and specify their logical role and their logical status. I shall not challenge the usual assumption according to which Euclid's mathematics is by and large sound. (I say ``by and large'' in order to leave some room for possible revisions and corrections of Euclid's arguments and thus avoid controversies about the question whether a given interpretation of Euclid is authentic or not. Although I pay more attention to textual details than it is usually done in modern logical reconstructions of Euclid's reasoning, I am not going to criticize these reconstructions by pointing to their anachronistic character.) However I shall challenge the other assumption according to which any sound mathematical theory is an axiomatic theory in the modern sense. Since I do not take this latter assumption for granted I do not assume from the outset that the problematic elements of Euclid's reasoning  (the \emph{exposition}, the \emph{specification} and the \emph{construction}) play some \emph{logical} role, which only needs to be made explicit and appropriately understood. In what follows I try to describe how these elements work without making about them any additional assumptions and only then decide whether the role of these elements qualifies as logical or not.   

\section{Are Euclid's Proofs Logical?}
Let's look at Euclid's Theorem 1.5 more attentively. I begin its analysis with its \emph{proof}. Among the premisses of this \emph{proof}, one may easily identify Axiom (Common Notion) 3 according to which

\begin{quote}
(Ax.3): If equal things are subtracted from equal things then the remainders are equal
\end{quote}

and the preceding Theorem 1.4 according to which

\begin{quote}
(Prop.1.4): If two triangles have two corresponding sides equal, and have the angles enclosed by the equal sides equal, then they will also have equal bases, and the two triangles will be equal, and the remaining angles subtended by the equal sides will be equal to the corresponding remaining angles.
\end{quote}

 I shall not comment on the role Theorem 1.4 in this \emph{proof} (which seems to be clear) but say few things about the role of the Axiom 3.

Here is how exactly the Axiom (Common Notion) 3 is used in the above Euclid's \emph{proof}. First, \emph{by construction} we have

\begin{quote}
\textbf{Con1}: $BF \equiv AF - AB$ and    
\textbf{Con2}: $CG \equiv AG - AC$
\end {quote}

which is tantamount to saying that point $B$ lays between points $A$, $F$ and point $C$ lays between points $A$, $G$). Second, \emph{by hypothesis} we have

\begin{quote}
\textbf{Hyp}: $AB = AC$ 
\end{quote}

and once again \emph{by construction}

\begin{quote}
\textbf{Con3}: $AF = AG$ 
\end{quote}

Now we see that we have got the situation described in Ax.3: equal things are subtracted from equal things. Using this Axiom we conclude that $BF = CG$.

Notice that Ax.3 applies to all ``things'' (mathematical objects), for which the relation of \emph{equality} and the operation of \emph{subtraction} make sense. In Euclid's mathematics this relation and this operation apply not only to straight segments and numbers but also to geometrical  objects of various sorts including \emph{figures}, angles and solids. Since Euclid's equality is not interchangeable with identity I use for the two relations two different symbols: namely I use the usual symbol for Euclid's equality (even if this equality is not quite usual), and use symbol $\equiv$ for identity. My use of symbols $+$ and $-$ is self-explanatory\footnote{ The \emph{difference} $A - B$ of two figures $A$, $B$ is a figure obtained through ``cutting'' $B$ out of $A$; the  \emph{sum}  $A + B$  is the result of \emph{concatenation} of  $A$ and $B$. These operations are not defined up to \emph{congruence} of figures (for there are, generally speaking, many possible ways, in which one may cut out one figure from another) but, according to Euclid's Axioms, these operations are defined up to Euclid's \emph{equality}. This shows that Euclid's \emph{equality} is weaker than \emph{congruence}: according to Axiom 4 congruent objects are equal but, generally, the converse does not hold. In the case of (plane) figures Euclid's equality is equivalent to the equality (in the modern sense) of their air.     
}.

The other four Euclid's Axioms (not to be confused with Postulates!) have  the same character.
This makes Euclid's Axioms in general, and Ax.3 in particular, very unlike premises like \textbf{Con1-3} and \textbf{Hyp}, so one may wonder whether the very idea of treating these things on equal  footing (as different \emph{premises} of the same inference) makes sense. More precisely we have here the following choice.
 One option is to interpret Ax.3 as the following implication:

\begin{quote}
$\{(a \equiv b - c) \& (d \equiv e - f ) \& (b = d ) \& (c =f)\} \rightarrow (a = b)$
\end{quote} 

and then use it along with \textbf{Con1-3} and \textbf{Hyp} for getting the desired conclusion through \emph{modus ponens} and other appropriate rules. This standard analysis involves a fundamental distinction between premises and conclusion, on the one hand, and rules of inference, on the other hand. It assumes that in spite of the fact that Euclid (as most of other mathematicians of all times) remains silent about logic, his reasoning nevertheless follows some implicit logical rules. The purpose of logical analysis  in this case is to make this ``underlying logic`` (as some philosophers like to call it) explicit.
  
The other option that I have in mind is to interpret Ax.3 itself as a rule rather than as a premiss. Following this rule, which can be pictures as follows:

\begin{quote}
$(a \equiv b - c), (d \equiv e - f ), (b = d ), (c =f)$\\
------------------------------------------------------------------------ (Ax.3)\\
(a = b)
\end{quote}

one derives from \textbf{Con1-3} and \textbf{Hyp}  the desired conclusion. So interpreted Ax.3 hardly qualifies as a  \emph{logical} rule because it applies only to propositions of a particular sort (namely, of the form $X = Y$ where $X, Y$ are some \emph{mathematical} objects of appropriate types). This Axiom cannot help one to prove that Socrates is mortal. Nevertheless the domain of application of this rule is quite vast and covers the whole of Euclid's mathematics. An important advantage of this analysis is that it doesn't require one to make any assumption about hidden features of Euclid's thinking: unlike the distinction between logical rules and instances of applications of these rules the distinction between axioms and premises like \textbf{Con1-3} and \textbf{Hyp} is explicit in Euclid's \emph{Elements}.

There is also a historical reason to prefer the latter reading of Euclid's Common Notions. Aristotle uses the word ``axiom'' interchangeably with the expressions ``common notions'', ``common opinions'' or simply ``commons'' for what we call today logical laws or logical principles but not for what we call today axioms. Moreover in this context he systematically draws an analogy between mathematical common notions and his proposed logical principles (laws of logic). This among other things provides an important historical justification for  calling Euclid's Common Notions by the name of Axioms. It is obvious that mathematics in general and mathematical common notions (axioms) in particular serve for Aristotle as an important source for developing the very idea of logic. Roughly speaking Aristotle's thinking, as I understand it, is this: behind the basic principles of mathematical reasoning spelled out through mathematical common notions (axioms) there are other yet more general principles relevant to reasoning about all sorts of beings and not only about mathematical objects. The fact that Euclid, according to the established chronology, is younger than Aristotle by some 25 years (Euclid's dates unlike Aristotle's are only approximate) shouldn't confuse one. While there is no strong evidence of the influence of Aristotle's work on Euclid, the influence on Aristotle of the same mathematical tradition, on which Euclid elaborated, is clearly documented in Aristotle's writings themselves. In particular, Aristotle quotes Euclid's Ax.3 (which, of course, Aristotle could know from another source) almost verbatim  
\footnote{Here are some quotes: 
 \begin{quote}
 By first principles of proof [as distinguished from first principles in general] I mean the
common opinions on which all men base their demonstrations, e.g. that one of two
contradictories must be true, that it is impossible for the same thing both be and not to be, and
all other propositions of this kind. (Met. 996b27-32, Heath's translation, corrected)
\end{quote}
Here Aristotle refers to a logical principle as ``common opinion''. In the next quote he compares mathematical and logical axioms:
\begin{quote}
We have now to say whether it is up to the same science or to different sciences to inquire
into what in mathematics is called axioms and into [the general issue of] essence. Clearly the
inquiry into these things is up to the same science, namely, to the science of the philosopher.
For axioms hold of everything that [there] is but not of some particular genus apart from
others. Everyone makes use of them because they concern being qua being, and each genus is.
But men use them just so far as is sufficient for their purpose, that is, within the limits of the
genus relevant to their proofs. Since axioms clearly hold for all things qua being (for being is
what all things share in common) one who studies being qua being also inquires into the
axioms. This is why one who observes things partly [=who inquires into a special domain]
like a geometer or a arithmetician never tries to say whether the axioms are true or false.
(Met. 1005a19-28, my translation)
\end{quote}
Here is the last quote where Aristotle refers to Ax.3 explicitly:
\begin{quote}
Since the mathematician too uses common [axioms] only on the case-by-case basis, it must
be the business of the first philosophy to investigate their fundamentals. For that, when equals
are subtracted from equals, the remainders are equal is common to all quantities, but
mathematics singles out and investigates some portion of its proper matter, as e.g. lines or
angles or numbers, or some other sort of quantity, not however qua being, but as [...]
continuous. (Met. 1061b, my translation)
\end{quote}  
The ``science of philosopher'' otherwise called the ``first philosophy'' is Aristotle's logic, which in his understanding is closely related to (if not indistinguishable from) what we call today ontology. After Alexandrian librarians we call today the relevant collection of Aristotle's texts by the name of \emph{metaphysics} and also use this name for a relevant philosophical discipline. 
}.

However important Aristotle's argument in the history of Western thought may be, there is no reason to take it for granted every time when we try today to interpret Euclid's \emph{Elements} or any other old mathematical text. Whatever is one's philosophical stance concerning the place of logical principles in human reasoning one can see what kind of harm can be made if Aristotle's assumption about the primacy of logical and ontological principles is taken straightforwardly and uncritically: one treats Euclid's Axioms on equal footing with premisses like \textbf{Con1-3} and \textbf{Hyp} and so misses the law-like character of the Axioms. Missing this feature doesn't allow one to see the relationships between Greek logic and Greek mathematics, which I just sketched.

Having said that I would like to repeat that  Euclid's  \emph{proof} (apodeixis) is the part of Euclid's theorems, which more resembles what we today call proof  (in logic) than other parts Euclid's theorems. For this reason in what follows I shall call inferences in Euclid's \emph{proofs}, which are based on Axioms, \emph{protological} inferences and distinguish them from inferences of another type that I shall call \emph{geometrical} inferences. This analysis is not incompatible with the idea (going back to Aristotle) that behind Euclid's protological and geometrical inferences there are inferences of a more fundamental sort, that can be called \emph{logical} in the proper sense of the word. However I claim that Euclid's text as it stands provides us with no evidence in favor of this strong assumption. One can learn Euclid's mathematics and fully appreciate its rigor without  knowing anything about logic just like Moliere's  M. Jourdain could well express himself long before he learned anything about prose!

Whether or not the science of logic really helps one to improve on mathematical rigor - or this is rather the mathematical rigor that helps one to do logic rigorously -  is a controversial question that I shall discuss throughout this book and suggest an answer only in the last Chapter. The purpose of my present reading of Euclid is at the same time more modest and more ambitious than the purpose of logical analysis. It is more modest because this reading doesn't purport to assess Euclid's reasoning from the viewpoint of today's mathematics and logic but aims at reconstructing this reasoning in its  authentic archaic form. It is more ambitious because it doesn't take the today's viewpoint  for granted but aims at reconsidering this viewpoint by bringing it into a historical perspective.

\section{Instantiation, Objecthood and Objectivity}

Let us now see where the premises \textbf{Hyp} and \textbf{Con 1-3} come from. As I have already mentioned they actually come from two different sources: \textbf{Hyp} is assumed  \emph{by hypotheis} while \textbf{Con 1-3} are assumed \emph{by construction}. 
I shall consider here these two cases one after the other.

The notion of hypothetic reasoning is an important extension of the core notion of axiomatic theory outlined above; it is well-treated in the literature and I shall not cover it here in full. I shall consider only one particular aspect of hypothetical reasoning as it is present in Euclid. The hypothesis that validates \textbf{Hyp}, informally speaking, amounts to the fact that Theorem 1.5 tells us something about isosceles triangles (rather than about objects of another sort). The corresponding definition (Definition 1.20) tells us that two sides of the isosceles triangle are equal. However to get from here to \textbf{Hyp} one needs yet another step. The \emph{enunciation} of Theorem 1.5 refers to isosceles triangles  \emph{in general}. But \textbf{Hyp} that is involved into the \emph{proof} of this Theorem concerns only \emph{particular} triangle $ABC$. Notice also that the \emph{proof} concludes with the propositions $ABC = ACB$ and $FBC = GCB$ (where $ABC$,  $ACB$, $FBC$ and $GCB$ are angles), which also concern only  \emph{particular} triangle $ABC$. This conclusion differs from the following \emph{conclusion} (of the whole Theorem), which almost verbatim repeats the  \emph{enunciation} and once again refers to isosceles triangles and their angles in general terms.

The wanted step that allows Euclid to proceed from the \emph{enunciation} to  \textbf{Hyp} is made in the \emph{exposition} of this Theorem, which introduces triangle $ABC$ as an ``arbitrary representative'' of isosceles triangles (in general). In terms of modern logic this step can be described as the \emph{universal instantiation}:\\  
$$\forall x P(x) \Longrightarrow P(a/x)$$\\ 
where $P(a/x)$ is the result of the substitution of individual constant $a$ at the place of all free occurrences of variable $x$ in $P(x)$. The same notion of universal instantiation allows for interpreting Euclid's  \emph{specification} in the obvious way. The reciprocal backward step that allows Euclid to obtain the \emph{conclusion} of the Theorem from the conclusion of the \emph{proof} can be similarly described as the  \emph{universal generalization }:\\  
$$P(a) \Longrightarrow \forall x P(x) $$\\ 
(which is a valid rule only under certain conditions that I skip here).

As long as the \emph{exposition} and the \emph{specification} are interpreted in terms of the universal instantiation these operations are understood as logical inferences and, accordingly, as element of proof in the modern sense of the word. A somewhat different - albeit not wholly incompatible - interpretation of Euclid's  \emph{exposition} and \emph{specification} can be straightforwardly given in terms of Kant's \emph{transcendental aesthetics} and \emph{transcendental logic} developed in his \emph{Critique of Pure Reason} \cite{Kant:1999}. Kant thinks of the traditional geometrical \emph{exposition} not as a logical inference of one proposition from another but as a ``general procedure of the imagination for providing a concept with its image''; a representation of such a general procedure Kant calls a \emph{schema} of the given concept (A140). Thus for Kant any individual mathematical object (like triangle $ABC$)  always comes with a specific \emph{rule} that one follows constructing this object in one's imagination and that provides a link between this object and its corresponding concept (the concept of isosceles triangle in our example). According to Kant the representation of general concepts by imaginary individual objects (which Kant for short also describes as ``construction of concepts'') is the principal distinctive feature of mathematical thinking, which distinguishes it from a philosophical speculation.

\begin{quote}
Philosophical cognition is rational cognition from concepts, mathematical cognition is that from the construction of concepts.'' But to construct a concept means to exhibit a priori the intuition corresponding to it. For the construction of a concept, therefore, a non-empirical intuition is required, which consequently, as intuition, is an individual object, but that must nevertheless, as the construction of a concept (of a general representation), express in the representation universal validity for all possible intuitions that belong under the same concept, either through mere imagination, in pure intuition, or on paper, in empirical intuition.... The individual drawn figure is empirical, and nevertheless serves to express the concept without damage to its universality, for in the case of this empirical intuition we have taken account only of the action of constructing the concept, to which many determinations, e.g., those of the magnitude of the sides and the angles, are entirely indifferent, and thus we have abstracted from these differences, which do not alter the concept of the triangle. \\
Philosophical cognition thus considers the particular only in the universal, but mathematical cognition considers the universal in the particular, indeed even in the individual...
 (A713-4/B741-2). 
\end{quote}

Kant's account can be understood as a further explanation of what the instantiation of mathematical concepts amounts to; then one may claim that the Kantian interpretation of Euclid's  \emph{exposition} and \emph{specification} is compatible with its interpretation as the universal instantiation in the modern sense. However the Kantian interpretation doesn't suggest by itself to interpret the instantiation as a logical procedure, i.e., as an inference of a proposition from another proposition. As the above quote makes it clear Kant describes the instantiation as a cognitive procedure of a different sort.
 
Now coming back to Euclid we must first of all admit that the  \emph{exposition} and the \emph{specification} of Theorem 1.5 as they stand are too concise for preferring one philosophical interpretation rather than another. Euclid introduces an isosceles triangle through Definition 1.20 providing no rule for constructing such a thing. (This example may serve as an evidence against the often-repeated claim that every geometrical object considered by Euclid is supposed to be constructed on the basis of Postulates beforehand.) Nevertheless given the important role of constructions in Euclid's geometry, which I explain in the next Section, the idea that every geometrical object in Euclid has an associated construction rule, appears very plausible. There is also another interesting textual feature of Euclid's \emph{specification} that in my view makes the Kantian interpretation more plausible.
 
Notice the use of the first person in the  \emph{specification} of Theorem 1.5 : ``I say that ....''. In \emph{Elements} Euclid uses this expression systematically in the \emph{specification} of every theorem. Interpreting the \emph{specification} in terms of universal instantiation one should, of course, disregard this feature as merely rhetorical. However it may be taken into account through the following consideration. While the  \emph{enunciation} of a theorem is a general proposition that can be best understood \'a la Frege in the abstraction from any human or inhuman thinker, i.e., independently of any thinking \emph{subject}, who might believe this proposition, assert it, refute it, or do anything else about it, the core of Euclid's theorems (beginning with  their  \emph{exposition})  involves an individual thinker (individual subject) that cannot and should not be wholly abstracted away in this context. When Euclid  \emph{enunciates} a theorem this \emph{enunciation} does not involve - or at least is not supposed to involve - any particularities of Euclid's individual thinking; the less this \emph{enunciation} is affected by Euclid's (or anyone else's) individual writing and speaking style the better. However the \emph{exposition} and the\emph{specification} of the given theorem essentially involve an \emph{arbitrary} choice of notation (``Let $ABC$ be an isosceles triangle...''), which is an individual choice made by an individual mathematician (namely, made by Euclid on the occasion of writing his \emph{Elements}). This individual choice of notation goes on par with what we have earlier described as  \emph{instantiation}, i.e. the choice of one individual triangle (triangle $ABC$) of the given type, which serves Euclid for proving the general theorem about \emph{all} triangles of this type. The \emph{exposition} can be also naturally accompanied by drawing a diagram, which in its turn involves the choice of a particular shape (provided this shape is of the appropriate type), to leave alone the choices of its further features like color, etc.

Thus when in the \emph{specification} of Theorem 1.5 we read ``I say that the angle $ABC$ is equal to $ACB$'' we indeed do have good reason to take Euclid's wording seriously. For the sentence ``angle $ABC$ is equal to $ACB$'' unlike the sentence ``for isosceles triangles, the angles at the base are equal to one another'' has a feature that is relevant only to one particular presentation (and to one particular diagram if any), namely the use of letters $A, B, C$ rather than some others \footnote{
Although the choice of letters in Euclid's notation is arbitrary the  \emph{system} of this notation is not. This traditional geometrical notation has a relatively stable and rather sophisticated syntax, which I briefly describe in what follows.}.  The words ``I say that ...'' in the given context stress this situational character of the following sentence ``angle $ABC$ is equal to $ACB$''. What matters in these words is, of course, not Euclid's personality but the reference to a particular act of speech and cognition of an individual mathematician. Proving the same theorem on a different occasion Euclid or anybody else could use other letters and another diagram of the appropriate type.

A competent reader of Euclid is supposed to know that the choice of letters in Euclid's notation is arbitrary and that Euclid's reasoning does not  depend of this choice. The arbitrary character of this notation should be distinguished from the general arbitrariness of linguistic symbols in natural languages. What is specific for the case of \emph{exposition} and \emph{specification} is the fact that here the arbitrary elements of reasoning (like notation) are sharply distinguished from its invariant elements. To use Kant's term we can say that behind the notion according to which the choice of Euclid's notation is arbitrary (at least at the degree that letters used in this notation are permutable) and according to which the same reasoning may work equally well with different diagrams (provided all of them belong to the same appropriate type) there is a certain invariant  \emph{schema} that sharply limits such possible choices. This schema not simply \emph{allows} for making some arbitrary choices but \emph{requires} every possible choice in the given reasoning to be wholly arbitrary. This requirement is tantamount to saying that subjective reasons behind choices made by an individual mathematician for presenting a given mathematical argument are strictly irrelevant to the ``argument itself'' (in spite of the fact that the argument cannot be formulated without making such choices). In general talks in natural languages there is no similar sharp distinction between arbitrary and invariant elements . When I write this paper I can certainly change some wordings without changing the sense of my argument but I am not in a position to describe precisely the scope of such possible changes  and identify the intended ``sense'' of my argument with a mathematical rigor. This is because my present study is philosophical and historical but not purely mathematical.    
 
Thus  Euclid's  \emph{exposition} serves for the formulation of a given universal proposition in terms, which are suitable for a particular act of mathematical cognition made by an individual mathematician. This aspect of the \emph{exposition} is not accounted for by  the modern notion of universal instantiation. It may be argued that this aspect of the \emph{exposition} needs not be addressed in a \emph{logical} analysis of Euclid's mathematics that aims at explication of the  \emph{objective meaning} of Euclid's reasoning and may well leave aside cognitive aspects of this reasoning. I agree that this latter issue lies out of the scope of logical analysis in the usual sense of the term but I disagree that the objective meaning of Euclid's reasoning can be properly understood without addressing this issue. Euclid's mathematical reasoning is \emph{objective} due to a mechanism that allows one to make universally valid inferences through one's individual thinking. Whatever the ``objective meaning'' might consist of this mechanism must be taken into account.     

\section{Proto-Logical Deduction and Geometrical Production} 

Remind that the \emph{proof} of Euclid's Theorem 1.5 uses not only premiss \textbf{Hyp} assumed ``by hypothesis'' but also premisses \textbf{Con 1-3} (as well as a number of other premisses of the same type) assumed ``by construction''. I turn now to the question about the role of Euclid's \emph{constructions} (which, but the way, are ubiquitous not only in geometrical but also in arithmetical Books of the \emph{Elements}) and more specifically consider the question how these \emph{constructions} support certain premisses that are used in following \emph{proofs}.

As it is well-known Euclid's geometrical constructions are supposed to be realized ``by ruler and compass''. In the \emph{Elements} this condition is expressed in the \emph{Elements} through the following three 

\begin{quote}
Postulates: \\
1. Let it have been postulated to draw a straight-line from any point to any point.\\
2. And to produce a finite straight-line continuously in a straight-line.\\
3. And to draw a circle with any center and radius.
\end{quote}

(I leave out of my present discussion two further Euclid's Postulates including the controversial Fifth Postulate.) \\

Before I consider popular interpretations of these Postulates and suggest my own interpretation let me briefly discuss the very term ``postulate'', which is traditionally used in English translations of Euclid's  \emph{Elements}. Fitzpatrick translates Euclid's verb ``aitein'' by English verb ``to postulate'' following the long tradition of Latin translations, where this Greek verb is translated by Latin verb ``postulare''. However according to today's standard dictionaries the modern English verb ``to postulate'' does not translate the Greek verb ``aitein'' and the the Latin verb ``postulare'' in general contexts: the modern dictionaries translate these verbs into ``to demand'' or ``to ask for''. This clearly shows that the meaning of the English verb ``to postulate'' that derives from Latin ``postulare'' changed during its lifetime\footnote{I reproduce here Fitzpatrick's footnote about Euclid's expression ``let it be postulated'':
\begin{quote}
The Greek present perfect tense indicates a past action with present significance. Hence, the 3rd-person present perfect imperative \emph{Hitesthw} could be translated as ``let it be postulated'', in the sense ``let it stand as postulated'', but not ``let the postulate be now brought forward''. The literal translation ``let it have been postulated'' sounds awkward in English, but more accurately captures the meaning of the Greek.
\end{quote}
}.

 Aristotle describes a postulate (aitema) as what ``is assumed when the learner either has no opinion on the subject or is of a contrary opinion'' (\emph{An. Post.} 76b); further he draws a contrast between postulates and  \emph{hypotheses} saying that the latter appear more plausible to the learner than the former (\emph{ibid.}). It is unnecessary for my present purpose to go any further into this semantical analysis trying to reconstruct an epistemic attitude that Euclid might have in mind ``demanding'' the reader to take his Postulates for granted. The purpose of the above philological remark is rather to warn the reader that the modern meaning of the English word ``postulate'' can easily mislead when one tries to interpret Euclid's Postulates adequately. So I suggest to read Euclid's Postulates as they stand and try to reconstruct their meaning from their context forgetting for a while what one has learned about the meaning of the term ``postulate'' from modern sources.

Euclid's Postulates are usually interpreted as propositions of a certain type and on this basis are qualified as axioms in the modern sense of the term. There are at least two different ways of rendering Postulates in a propositional form. I shall demonstrate them at the example of Postulate 1. This Postulate can be interpreted either as the following \emph{modal} proposition:\\
(PM1): given two different points it is always possible to drawing a (segment of) straight-line between these points\\
or as the following \emph{existential} proposition: \\
(PE1): for any two different points there exists a (segment of) straight-line lying between these points.

Propositional interpretations of Euclid's Postulates allow one to present Euclid's geometry as an axiomatic theory in the modern sense of the word and, more specifically, to present Euclid's geometrical constructions as parts of proofs of his theorems. Even before the modern notion of axiomatic theory was strictly defined in formal terms many translators and commentators of Euclid's \emph{Elements} tended to think about his theory in this way and interpreted Euclid's Postulates in the modal sense. Later a number of authors (\cite{Hintikka&Remes:1974}, \cite{Avigard&Dean&Mumma:2009}) proposed different formal reconstructions of Euclid's geometry based on the existential reading of Postulates. 
According to Hintikka and Remes 

\begin{quote}
[R]eliance on auxiliary construction does not take us outside the axiomatic framework of geometry. Auxiliary constructions are in fact little more than ancient counterparts to applications of modern instantiation rules. 
\cite{Hintikka&Remes:1976}, p. 270
\end{quote}

The instantiation rules work in this context as follows. First, through the \emph{universal instantiation} (which under this reading correspond to Euclid's \emph{exposition} and \emph{specification}) one gets some propositions like \textbf{Hyp} about certain particular objects (individuals) like  $AB$ and $AC$. Then one uses Postulates 1-3 reading them as existential axioms according to which the existence of certain geometrical objects implies the existence of certain further geometrical objects, and so proves the (hypothetical) existence of such further objects of interest. Finally one uses another instantiation rule called the rule of \emph{existential instantiation}:\\ 
$$\exists x P(x) \Longrightarrow P(a) $$\\ 
and thus ``gets'' these new objects. Under this interpretation Euclid's \emph{constructions} turn into logical inferences of sort. 
As Hintikka and Remes emphasize in their paper the principal distinctive feature of Euclid's \emph{constructions} (under their interpretation) is that these constructions introduce some \emph{new} individuals; they call such individuals ``new'' in the sense that these individuals are not (and cannot be) introduced through the universal instantiation of hypotheses making part the \emph{enunciation} of the given theorem.

The propositional interpretations of Euclid's Postulates are illuminating because they allow for analyzing traditional geometrical constructions in modern logical terms. However they require a paraphrasing of Euclid's wording, which from a logical point of view is far from being innocent. In order to see this let us leave aside the epistemic attitude expressed by the verb ``postulate'' and focus on the question of \emph{what} Euclid postulates in his Postulates 1-3. Literally, he postulates the following: 

\begin{quote}
P1: to draw a straight-line from any point to any point.\\
P2: to produce a finite straight-line continuously in a straight-line.\\
P3: to draw a circle with any center and radius.
\end{quote}

As they stand expressions P1-3 don't qualify as propositions; they rather describe certain \emph{operations}! And making up a proposition from something which is not a proposition is not a innocent step. My following analysis is based on the idea that Postulates are \emph{not} primitive truths from which one may derive some further truths but are primitive operations that can be combined with each other and so produce into some further operations. In order to make my reading clear I paraphrase P1-3 as follows: 
\begin{quote}
(OP1): drawing a (segment of) straight-line between its given endpoints\\
(OP2): continuing a segment of given straight-line indefinitely (``in a straight-line)''\\
(OP3): drawing a circle by given radius (a segment of straight-line) and center (which is supposed to be one of the two endpoints of the given radius).\\
\end{quote}

Noticeably none of OP1-3 allows for producing geometrical constructions out of nothing; each of these fundamental operation produces a geometrical object out of some other objects, which are supposed to be \emph{given} in advance. The table below specifies inputs (operands) and  outputs (results) of  OP1-3:

\begin {center}
\begin{tabular}{|l|c|r|}
  \hline
  operation & input & output \\
  \hline
  OP1 & two (different) points & straight segment \\
   \hline
  OP2 & straight segment  & (bigger) straight segment \\
  \hline
  OP3 & straight segment and one of its endpoints & circle  \\ 
  \hline
\end{tabular}
\end {center}

PE1 as it stands does not imply that there exists at least one point or at least one line in Euclid's geometrical universe. If there are no points then there are no lines either. Similar remarks can be made about the existential interpretation of other Euclid's Postulates. Thus the existential interpretation of Postulates by itself does not turn these Postulates into existential axioms that guarantee that Euclid's universe is non-empty and contains all geometrical objects constructible by ruler and compass. To meet this purpose one also needs to postulate the existence of at least two different points - and then argue that the absence of any counterpart of such an axiom in Euclid is due to Euclid's logical incompetence. My proposed  interpretation of Postulates 1-3 as operations doesn't require such ad hoc stipulations and thus is more faithful not only to Euclid's text but also to a deeper structure of his reasoning
\footnote{
Remind that the concepts of infinite straight line and infinite half-line (ray) are absent from Euclid's geometry; thus the result of OP2 is always a finite straight segment. However the result of OP2 is obviously not fully determined by its single operand. This shows that OP2 doesn't really fit the today's usual notion of algebraic operation. }.

Hintikka and Remes describe Euclid's geometrical constructions as \emph{auxiliary}. Such a description may be adequate to the role of geometrical constructions in today's practice of teaching the elementary geometry but not to the role of constructions in Euclid's  \emph{Elements}. Remind that Euclid's so-called Propositions are of two types: some of them are Theorems while some other are Problems (see again the above quotation from Proclus' \emph{Commentary} ). In the \emph{Elements} Problems are at least as important as Theorems and arguably even more important: in fact  the \emph{Elements} begin and end with a Problem but not with a Theorem. As we shall now see when a given \emph{construction} makes part of a problem rather than a theorem it cannot be described as auxiliary in any appropriate sense. We shall also see the modern title ``proposition'' is not really appropriate when we talk about Euclid's Problems: while \emph{enunciations} of Theorems do qualify as propositions in the modern logical sense of the term \emph{enunciations} of Problems do not.

I shall demonstrate these features at the well known example of Problem 1.1 that opens Euclid's \emph{Elements}; my notational conventions remain the same as in the example of Theorem 1.5.

\paragraph{[\emph{enunciation:}]}
\begin{quote}
To construct an equilateral triangle on a given finite straight-line.
\end{quote}

\paragraph{[\emph{exposition}:]}

\begin{quote}
Let AB be the given finite straight-line.
\end{quote}

\begin{center}
\includegraphics [scale=0.5]{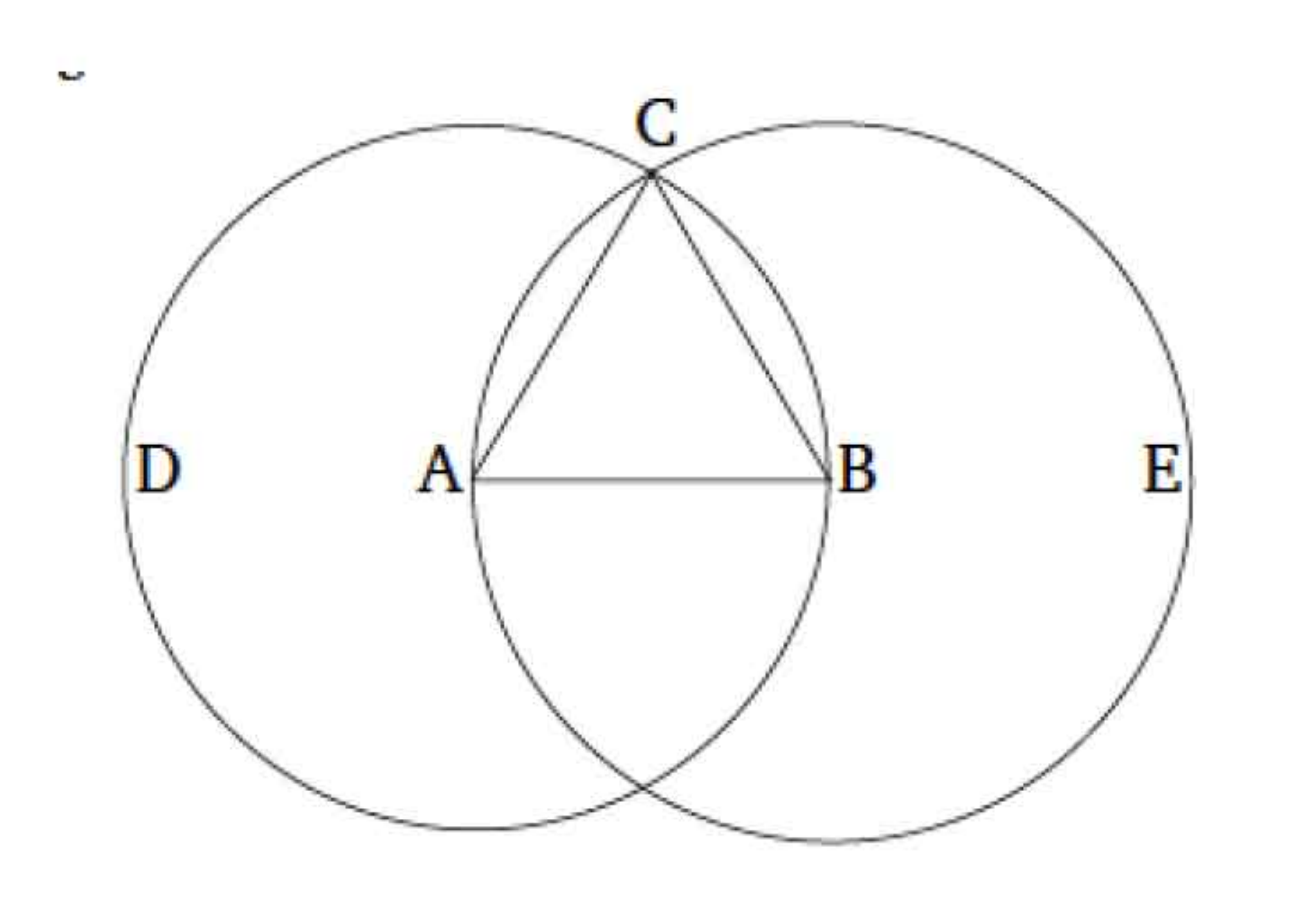}
\end{center}

\begin{center}
Fig. 1.2
\end{center}

\paragraph{[\emph{specification}:]}

\begin{quote}
So it is required to construct an equilateral triangle on the straight-line $AB$.
\end{quote}

\paragraph{[\emph{construction}:]}

\begin{quote}
Let the circle $BCD$ with center $A$ and radius $AB$ have been drawn [Post. 3], and again let the circle $ACE$ with center $B$ and radius $BA$ have been drawn [Post. 3]. And let the straight-lines $CA$ and $CB$ have been joined from the point $C$, where the circles cut one another, to the points $A$ and $B$ [Post. 1].
\end{quote}

\paragraph{[\emph{proof}:]}

\begin{quote}
And since the point $A$ is the center of the circle $CDB$, $AC$ is equal to $AB$ [Def. 1.15]. Again, since the point $B$ is the center of the circle $CAE$, $BC$ is equal to $BA$ [Def. 1.15]. But $CA$ was also shown (to be) equal to $AB$. Thus, $CA$ and $CB$ are each equal to $AB$. But things equal to the same thing are also equal to one another [Axiom 1]. Thus, CA is also equal to CB. Thus, the three (straight-lines) $CA$, $AB$, and $BC$ are equal to one another.
\end{quote}

\paragraph{[\emph{conclusion}:]}

\begin{quote}
Thus, the triangle $ABC$ is equilateral, and has been
constructed on the given finite straight-line $AB$. (Which
is) the very thing it was required to do.
\end{quote}

As one can see at this example \emph{enunciations} of Problems are expressed in the same grammatical form as Postulates 1-3, namely in the form of infinitive verbal expressions. I read these expressions in the same straightforward way, in which I read Postulates: as descriptions of certain geometrical \emph{operations}. The obvious difference between (\emph{enunciations} of) Problems and Postulates is this: while Postulates introduce basic operations taken for granted (drawing by ruler and compass) Problems describe complex operations, which in the last analysis reduce to these basic operations. Such reduction is made through a \emph{construction} of a given Problem: it performs the complex operation described in the \emph{enunciation} of the problem through combining basic operations OP1-3 (and possibly some earlier performed complex operations). The procedure that allows for  performing complex operations by combining a small number of repeatable basic operations I shall call a \emph{geometrical production}. In Problem 1.1 the construction of regular triangle is (geometrically) \emph{produced} from drawing the straight-line between two given points (Postulate 1) and drawing a circle by given center and radius (Postulate 3). This is, of course, just another way of saying that the regular triangle is constructed by ruler and compass; the unusual terminology helps me to describe Euclid's geometrical constructions more precisely.

Let us see in more detail how works Euclid's geometrical production. Basic operations OP1-3 like other (complex) operations need to be \emph{performed}: in order to produce an output they have to be fed by some input. This input is provided through the \emph{exposition} of the given Problem (the straight line $AB$ in the above example). OP1-3 are composed in the usual way well-known from today's algebra: outputs of earlier performed operations are used as inputs for further operations\footnote{
Problem 1.1 involves a difficulty that has been widely discussed in the literature: Euclid does not provide any principle that may allow him to construct a point of intersection of the two circles involved into the \emph{construction} of this Problem. This flaw is usually described as a \emph{logical} flow. In my view it is more appropriate to describe this flow as properly\emph{geometrical} and fill the gap in the reasoning by the following additional postulate (rather than an additional axiom):

\begin{quote}
Let it have been postulated to produce a point of intersection of two circles with a common radius.
\end{quote}  

Even if this additional postulates is introduced here purely ad hoc, the way in which it is introduced gives an idea of how Euclid's Postulates could emerge in the real history.    
}.

Just like Postulates 1-3 \emph{enunciations} of Problems can be read as modal or existential propositions (in the modern logical sense of the term). Then the (modified) \emph{enunciation} of Problem 1.1 reads:
 
 \begin{quote}
 (1.1.M) it is possible to construct a regular triangle on a given finite straight-line:
 \end{quote}
 or
\begin{quote}
(1.1.E) for all finite straight-line there exists a regular triangle on this line.
 \end{quote}

As soon as the \emph{enunciations} of Euclid's Problems are rendered into the propositional form the Problems turn into theorems of a special sort. In the case of existential interpretation Problems turn into \emph{existential} theorems that state (under certain hypotheses) that there exist certain objects having certain desired properties. However this is not what we find in Euclid's text as it stands. Every Euclid's Problem ends with the formula ``the very thing it was required to do'', not ``to show'' or ``to prove''. I can see no evidence in the  \emph{Elements} that justifies the idea that in Euclid's mathematics \emph{doing} is less significant than \emph{showing} and that the former is in some sense reducible to the latter. In the Second Part of this paper I shall argue that \emph{doing} remains as much important in today's mathematics as it was in Greek mathematics, and that the idea of reducing mathematics to  \emph{showing}  or  \emph{proving} (in the precise sense of modern logic) is a unfortunate philosophical misconception.

According to another popular reading Euclid's Problems are tasks or questions of sort. This version of modal propositional interpretation of Euclid's Problems involves a deontic modality rather than a possibility modality: 

\begin{quote}
 (1.1.D) it is required to construct a regular triangle on a given finite straight-line:
 \end{quote}

Indeed geometrical problems similar to Euclid's Problems can be found in today's Elementary Geometry textbooks as exercises. However the analogy between Euclid's Problems and school problems on construction by ruler and compass is not quite precise. \emph{Enunciations} of Euclid's Problems just like the \emph{enunciations} of Euclid's Theorems prima facie express no modality whatsoever. A deontic expression appears only in the \emph{exposition} of the given Problem (``it is required to construct an equilateral triangle on the straight-line $AB$''). I don't think that this fact justifies the deontic reading of the \emph{enunciation} because, as I have already explained above, the \emph{exposition} describes reasoning of an individual mathematician rather than presents this reasoning in an objective form. That every complex construction must be performed through Postulates  and earlier performed constructions is an epistemic requirement, which is on par with the requirement according to which every theorem must be proved rather than simply stated. Remind that the \emph{expositions} of Euclid's Theorems have the form ``I say that...''. This indeed makes an apparent contrast with the \emph{expositions} of Problems that have the form ``it is required to ....''. However this contrast doesn't seem me to be really sharp. Euclid's expression ``I say that...'' in the given context is interchangeable with the expression ``it is required to show that...'', which matches the closing formula of Theorems ``(this is) the very thing it was required to show''. Euclid's expression ``it is required to...'' that he uses in the \emph{expositions} of Problems similarly matches the closing formula of Problems ``(this is) the very thing it was required to do''. The requirement according to which every Theorem must be ``shown'' or ``monstrated'' doesn't imply, of course, that the \emph{enunciation} (statement) of this Theorem has a deontic meaning. The requirement according to which every Problem must be ``done'' doesn't imply either that the \emph{enunciation} of this Problem  has something to do with deontic modalities.

The analogy between axioms and theorems, on the one hand, and postulates and problems, on the other hand, may suggest that Euclid's geometry splits into two independent parts one of which is ruled by (proto)logical deduction while the other is ruled by geometrical production. However this doesn't happen and in fact problems and theorems turn to be mutually dependent elements of the same theory. The above example of Problem 1.1 and Theorem 1.5 show how the intertwining of problems and theorems works. Theorems, generally, involve \emph{constructions} (called in this case auxiliary), which may depend (in the order of geometrical production) on earlier treated problems (as the \emph{construction} of Theorem 1.5 depends on Problem 1.3.) Problems in their turn always involve appropriate \emph{proofs} that prove that the \emph{construction} of the given theorem indeed performs the operation described in the \emph{enunciation} of this theorem (rather than performs some other operation). Such \emph{proofs}, generally, depend (in the order of the protological deduction) on certain earlier treated theorems (just like in the case of \emph{proofs} of theorems).  
Although this mechanism linking problems with theorems may look unproblematic it gives rise to the following puzzle. Geometrical production produces geometrical objects from some other objects. Protological deduction deduces certain propositions from some other propositions. How it then may happen that the geometrical production has an impact on the protological deduction? In particular, how the geometrical production may justify premises assumed ``by construction'', so these premises are used in following \emph{proofs}?

In order to answer this question let's come back to the premise \textbf{Con3} ($AF = AG$) from Theorem 1.5 and see what if anything makes it true. $AF = AG$ because Euclid or anybody else following Euclid's instructions constructs this pair of straight segments in this way. How do we know that by following these instructions one indeed gets the desired result? This is because the \emph{construction} of Problem 1.3 that contains the appropriate instruction is followed by a \emph{proof} that proves that this \emph{construction} does exactly what it is required to do. \emph{Construction} 1.3 in its turn uses \emph{construction} 1.2  while \emph{construction} 1.2 uses \emph{construction} 1.1 quoted above. In other words \emph{construction}  1.1 (geometrically) produces \emph{construction}  1.2 and \emph{construction}  1.2 in its turn produces \emph{construction} 1.3. This  geometrical production produces the relevant part of \emph{construction} 1.5 (the construction of equal straight segments $AF$ and $AG$) from first principles, namely from Postulates 1-3. In order to get the corresponding protological deduction of  premise \textbf{Con3}  from first principles we should now look at \emph{proofs} 1.1, 1.2 and 1.3 and then combine these three proofs into a single chain. For economizing space I leave now details to reader and only report what we get in the end. The result is somewhat surprising from the point of view of the modern logical analysis. The chain of \emph{constructions} leading to \emph{construction} 1.5 involves a number of circles (through Postulate 3). Radii of a given circle are equal by definition (Definition 1.15). Thus by constructing a circle and its two radii, say, $X$ and $Y$ one gets a primitive (not supposed to be proved) premise $X = Y$. Having at hand a number of premises of this form and using Axioms as inference rules (but not as premises!) one gets the desired deduction of \textbf{Con3}. The fact that first principles of the protological deduction of \textbf{Con3} appear to be partly provided by a definition helps to explain why Euclid places his definitions among other first principles such as postulates and axioms.

The above analysis allows for disentangling the protological deduction of \textbf{Con3} from the geometrical production of straight segments  $AF$, $AG$ and so the aforementioned puzzle remains even after we have looked at Euclid's reasoning under a microscope. Even if we can describe in detail the impact of Problems to Theorems and vice versa it remains unclear how the two kind of things can possibly work together. Here is my tentative answer to this question. Every Euclid's \emph{proof}   involves only \underline{concrete} premises like \textbf{Con3} and \textbf{Hyp}, which are related to certain individual objects. It is assumed that such a premise is valid \underline{only if} the corresponding object is effectively constructed. (At least this concerns all premises ``by construction''; as we have seen at the example of Theorem 1.5 hypothetic premises sometimes don't respect this rule.) This fundamental principle links Euclid's geometrical production and protological deduction together.

One may argue that my proposed analysis after all is not significantly different from the standard logical analysis of Euclid's geometrical reasoning according to which Euclid first proves that certain geometrical objects exist and only then prove some further propositions concerning properties of these objects. Is there indeed any significant difference between proving that  such-and-such object exist and producing this object through what I call the geometrical production? There is of course a difference of a metaphysical sort between these two notions: to produce an object is not quite the same thing as to prove that certain object exists. But arguably this difference has no objective significance and so one may simply ignore it. There is however a further difference between the geometrical production and the mathematical existence, which seems me more important. Euclid's \emph{Elements} contain two sets of rules, namely axioms and postulates, supposed to be applied to operations of two different sorts: axioms tell us how to derive equalities from other equalities while postulates tell us how to produce geometrical objects from other geometrical objects. A logical analysis of Euclid's geometry that involves a propositional (in particular existential) reading of postulates aims at replacing these two sets of rules by a single set of rules called \emph{logical}. I would like to stress again that the results my proposed analysis do not exclude the possibility of logical analysis. Such a replacement may be or be not a good idea but in any event logical rules are not made in the Euclid's text explicit and I do not see much point in saying that he uses rules of this sort implicitly. The fact that \emph{we} can use today modern logic for interpreting Euclid is a completely different issue.  An interpretation of Euclid's geometry by means of logical analysis can be indeed illuminating but one should not confuse oneself by believing that Euclid already had a grasp of modern logic even if could not formulate principles of this logic explicitly.

For further references I shall call the 6-part structure of Euclid's problems and theorems \emph{Euclidean structure}. As the above analysis makes it clear the Euclidean structure does not fit into Hilbert's notion of axiomatic theory even when this latter notion is formulated in very general terms as in the above quotes. While Hilbert and his modern followers assume that a mathematical theory is a set of truths, some of which are assumed as axioms and some other are \emph{logically inferred} from axioms, Euclid builds his theory through a combination of two different procedures, which I call protological deduction and geometrical production.  Precipitating what follows I would like to mention here that Hilbert's view on mathematical theory (which is presented more accurately in the next Chapter) is not unique in the 20th century. An influential alternative view has been put forward by Luitzen Egbertus Jan Brouwer; a relevant part of this alternative view is formulated by Brouwer's student Heyting as follows:

\begin{quote}
One of Brouwer's main theses was that mathematics is not based on
logic, but that logic is based on mathematics. [..] If mathematics
consists of mental constructions, then every mathematical theorem is the expression of a result of a successful construction. The proof of the theorem consists in this construction itself, and the steps of the proof
are the same as the steps of the mathematical construction. These are
intuitively clear mental acts, and not applications of logical laws. (quoted by \cite{Scott:1970}, p. 237)
\end{quote}

This general description prima facie better fits Euclid's procedures than the modern axiomatic approach. The problem is that this description does not by itself provide us with an alternative general method of building mathematical theories.  I postpone the discussion on this matter until \textbf{2.2} and conclude the present Chapter with a observation concerning the relevance of Euclid's way of theory-building to today's mathematical practice. 

\section{Euclid and Modern Mathematics}
What has been said above may give one an impression that in Euclid's  \emph{Elements} we deal with an archaic pattern of mathematical thinking that has noting to do with today's mathematics. However this impression is wrong. In fact the Euclidean structure is apparently present in today's mathematics, perhaps in a slightly modified form. Consider the following example taken from a standard mathematical textbook (\cite{Kolmogorov&Fomin:1976}, p. 100, my translation into English):

\paragraph{Theorem 3:}

\begin{quote}
Any closed subset of a compact space is compact
\end{quote}

\paragraph{Proof:}

\begin{quote}
Let $F$ be a closed subset of compact space $T$ and $\{F_{\alpha}\}$ be an arbitrary centered system of closed subsets of subspace $F \subset T$. Then every $F_{\alpha}$ is also closed in $T$, and hence $\{F_{\alpha}\}$ is a centered system of closed sets in $T$. Therefore $\cap F_{\alpha} \neq \emptyset$. By Theorem 1 it follows that $F$ is compact.
\end{quote}

Although the above theorem is presented in the usual for today's mathematics form ``proposition-proof'', its Euclidean structure can be made explicit without re-interpretations and paraphrasing:

\paragraph{[\emph{enunciation:}]}
\begin{quote}
Any closed subset of a compact space is compact
\end{quote}

\paragraph{[\emph{exposition:}]}

\begin{quote}
Let $F$ be a closed subset of compact space $T$
\end{quote}

\paragraph{[\emph{specification}: absent]}

\paragraph{[\emph{construction}:]}

\begin{quote}
[Let] $\{F_{\alpha}\}$ [be] an arbitrary centered system of closed subsets of subspace $F \subset T$. 
\end{quote}

\paragraph{[\emph{proof}:]}

\begin{quote}
[E]very $F_{\alpha}$ is also closed in $T$, and hence $\{F_{\alpha}\}$ is a centered system of closed sets in $T$. Therefore $\cap F_{\alpha} \neq \emptyset$. By Theorem 1 it follows that $F$ is compact.
\end{quote}

\paragraph{[\emph{conclusion}: absent ]}

The absent \emph{specification} can be formulated as follows:

\begin{quote}
I say that $F$ is a compact space
\end{quote}

while the absent \emph{conclusion} is supposed to be a literal repetition of the \emph{enunciation} of this theorem. Clearly these latter elements can be dropped for parsimony reason. In order to better separate the \emph{construction} and the \emph{proof} of the above theorem the authors could first construct set 
 $\cap F_{\alpha}$ and only then prove that it is non-empty. However this variation of the classical Euclidean scheme also seems me negligible. I propose the reader to check it at other modern examples that the Euclidean structure remains today at work.

 Does this mean that the modern notion of axiomatic theory is inadequate to today's mathematical practice just like it is inadequate to Euclid's mathematics? Such a conclusion would be too hasty. Arguably, in spite of the fact that today's mathematics preserves some traditional outlook it is essentially different. So the ``Euclidean appearance'' of today's mathematics cannot be a sufficient evidence for the claim the the Euclidean structure remains significant in it. In order to justify this claim a different argument is needed.

Before I try to provide such an argument I would like to point to the fact that the modern notion of axiomatic theory is used in today's mathematics in two rather different ways. First, it is used as a broad methodological idea that determines the general architecture of a theory but has no impact on details. Such an application of the modern axiomatic method is usually called \emph{informal}. Second, the notion of axiomatic theory is used for building \emph{formal} theories that contain a list of axioms and a set of theorems derived from these axioms according to explicitly specified rules of logical inference. In the next Chapter I shall describe the notion of formal axiomatic theory more precisely and try to explain the precise sense in which it is called formal. 
 
\chapter{Hilbert: Making It Formal}
In the standard textbooks Hilbert's philosophy of mathematics is commonly labelled \emph{formalism} and under this title distinguished from Brouwer's \emph{intuitionism}, on the one hand, and Russell's \emph{logicism}, on the other hand. However, as Hintikka \cite{Hintikka:1997b} rightly remarks, this popular name is very misleading. There are difficulties of two sorts. First of all, Hilbert's work in foundations of mathematics was a long-term project that began in 1890-ies and continued more than 40 years. Although Hilbert unlike Russell never abruptly changed his mind about foundational matters the development of Hilbert's project involved significant shifts in its philosophical underpinning. When one takes this into consideration it becomes impossible to identify Hilbert's views with any particular ``ism''.  Second of all, the meaning of being \emph{formal} is also changing: Hilbert and his contemporaries often use this term not in the same sense in which we use it today, and even today this term is often used in different ways in the mathematical and the philosophical communities. The two difficulties are mutually related because Hilbert's work in foundations strongly affected the changing meaning of being formal. 

I shall not treat here the history of Hilbert's research in foundations systematically\footnote{For the question of historical origins of Hilbert's Axiomatic Method see \cite{Toepell:1986} and \cite{Corry:2006}} but try to reconstruct the core dialectics of Hilbert's ideas, which is crucial for my analysis of today's state of affairs given in the next Chapter. I shall refer in this present Chapter to three Hilbert's texts: first, \emph{Foundations of Geometry} of 1899 \cite{Hilbert:1899}, second, his address \emph{Axiomatic Thought} of 1917 \cite{Hilbert:1918} and, finally, his paper \emph{Foundations of Mathematics} of 1927 \cite{Hilbert:1927}, which makes explicit the philosophical background behind the monumental two-volume work \cite{Hilbert&Bernays:1934-1939} co-authored with Bernays and published in 1934-1939  under the same title. We shall see that although the Axiomatic Method as presented in \emph{Foundations of Geometry} of 1899 and in \emph{Foundations of Mathematics} in both cases qualifies as (or at least is commonly called by Hilbert's contemporaries) formal, the sense of being formal is not the same in the two cases.

\section{\emph{Foundations} of 1899: Logical Form and Mathematical Intuition}
In the last Section I stressed Hilbert's assumption according to which the deduction of mathematical theorems from axioms is purely \emph{logical} and then argued that the geometrical theory of Euclid's \emph{Elements} prima facie falsifies this assumption. However this assumption is not specific for Hilbert's approach. Frege, who sharply criticizes Hilbert's \emph{Foundations} of 1899 on a different ground (that I shortly explain), wholly agrees with Hilbert on this general point. A major difference between Frege's and (early) Hilbert's versions of the Axiomatic Method, which led to a controversy between the two thinkers \cite{Frege:1971}, was the following. Frege assumes as a matter of course that all terms involved into axioms and theorems of a given theory are meaningful and that their meanings are specified in advance and rigidly fixed once and for all (at least within the given theory). Hilbert in his turn allows certain terms of a given to change their meanings and be considered without any fixed meaning at all. A theory of this latter  sort Frege and some other Hilbert's contemporaries call \emph{formal}.  For a mathematically educated reader (let alone logician) this ``informal'' notion of formal theory is, of course, very familiar. Nevertheless for my purpose it is useful to present it here in an explicit form. Then I shall explain the sense in which Frege et al. call such a theory formal.

The first paragraph of the {Foundations} of 1899 reads:

\begin{quote}
Let us consider three distinct systems of things. The things composing the first system, we will call points and designate them by the letters $A$, $B$, $C$,. . . ; those of the second, we will call straight lines and designate them by the letters $a$, $b$, $c$,..; and those of the third system, we will call planes and designate them by the Greek letters  $\alpha$, $\beta$, $\gamma$ . [..]  We think of these points, straight lines, and planes as having certain mutual relations, which we indicate by means of such words as ``are situated'', ``between''; ``parallel'',  ``congruent'', ``continuous'', etc. The complete and exact description of these relations follows as a consequence of the axioms of geometry. These axioms [..] express certain related fundamental facts of our intuition.
\end{quote}

The idea is this. The purpose of foundations of geometry is to develop geometry \emph{ab ovo}. This means that ``fundamental facts of our  [geometrical] intuition'' cannot be here tacitly taken for granted (as this is done in non-foundational geometrical studies) but must be explicitly described and postulated. The proposed method of describing these facts is the following. First, one identifies a list of \emph{types of objects}, which are \emph{primitive} in the sense that they are not defined in terms of some other (types of) objects; they are introduced without any definition. Second, one identifies a list of \emph{primitive relations} between primitive objects; these primitive relations are also introduced without definitions. Finally, one makes up a list of \emph{axioms}, i.e., propositions, which involve only primitive objects and primitive relations between these objects. Every consequence of these axioms qualifies as a geometrical theorem. (I shall specify a relevant notion of consequence in what follows; we shall see that there are in fact two different notions of consequence, which are here in play.)

 Hilbert's Axiomatic Method does \emph{not} assume that primitive objects and primitive relations are given through the usual linguistic meanings of words ``point'', ``between'', etc. Primitive objects are assumed instead to be bare ``things'' (possibly of several different \emph{types}), which are called points, straight lines and the like by a merely linguistic convention having no theoretical significance. Primitive relations are treated similarly. Thus Hilbert's list of types of primitive objects and of primitive relations given in the above quote does not say us anything except that the given axiomatic theory involves three different types of primitive objects and several different relations between these objects. All the relevant information about these objects and these relations is supposed to be captured by axioms, which specify certain facts about these objects and these relations without using any assumption as to \emph{what} are these objects and these relations. 
 
To see how it works consider the First Axiom of Hilbert's \emph{Foundations} of 1899:
 
\begin{quote}
 (A1.0) Two distinct points $A$ and $B$ always completely determine a straight line $a$ (\emph{op.cit.}, p.2).
\end{quote}
     
and remind that words ``points'' and ``straight line'' should not be read here in the usual sense. Notice also a relation between the points and the line, which is expressed by saying that the points determine the line; there is more than one way to translate this expression in terms of relations but Hilbert uses here the binary relation of \emph{incidence} between a given straight line and a given point, which can be also informally expressed by saying that the given point \emph{lies} at the given straight line (or equivalently by saying that the given straight line \emph{goes through} the given point).  This semantic hygiene leaves us with the following \emph{formal} reading of A1.0

\begin{quote}
(A1.1) Given two different primitive objects $A, B$ of basic type $P$ (``points'') there exist a unique primitive object $a$ of another basic type $L$ (``straight lines''),  such that each of $A, B$ and $a$ hold a primitive relation $R$ (``incidence''). 
\end{quote}
 
Although A1.1 may seem to be not very informative it presents what Hilbert's First Axiom ``really says'' more accurately than A1.0. The idea of Hilbert's Axiomatic Method is that a system of propositions like A1.1 provided with an appropriate system of logic may completely determine (in a sense that I try to clarify further in what follows) what the Euclidean (or some other) geometry ``really is''. The same method of theory-building is supposed to apply in various domains of the theoretical inquiry both within and outside the pure mathematics. Whatever is the domain of application of the Axiomatic Method the axioms always involve only abstract objects and abstract relations. What is specific for Euclidean geometry from Hilbert's axiomatic viewpoint is the list of its axioms rather than any particular subject-matter like \emph{space} or \emph{extension}.

Suppose a non-experienced reader looks at A1.1 and asks what this proposition has to do with the Euclidean geometry. An appropriate explanation can be given by translation A1.1 back to A1.0 followed by the ``naive'' reading of A1.0, which turns it into a proposition similar to Euclid's First Postulate. This naive reading of A1.0 refers to a ``fundamental fact of our intuition'', which, by Hilbert's word, this axiom ``expresses''. However in the given context this ``fundamental fact of intuition'' does not \emph{ground} the corresponding axiom A1.0 but merely motivates it. We shall shortly see, however, that in a different version of his Axiomatic Method presented in the \emph{Foundations} of 1927 \cite{Hilbert:1927} Hilbert grants a fundamental role to the geometrical intuition of a \emph{special} sort.   

How a proposition like A1.1 may qualify as an axiom? In his letter to Frege Hilbert says:

\begin{quote}
[A]s soon as I posited an axiom it will exist and be ``true''. [..] If the arbitrarily posited axioms together with all their consequences do not contradict each other, then they are true and the things defined by these axioms exist. For me, this is the criterion of truth and existence.  (\cite{Frege:1971}, p. 12)
\end{quote}

Some comments are here in order. 

(1) Unlike Frege, Hilbert does  \emph{not} think about mathematical axioms as self-evident truths. In the above quote Hilbert speaks of axioms as sheer stipulations, which are  ``true'' in virtue of the fact that they are posited by someone. The only rule restricting the positing of new axioms is the rule according to which each axiom must be self-consistent and any set of such axioms (belonging to the same theory) may contain only mutually consistent axioms. As Hilbert puts this in the above passage ``If the arbitrarily posited axioms [..] do not contradict each other, then they are true''.  One may remark (as did Frege) that given a set of true propositions it is impossible to infer from them a contradiction anyway. However this observation does not make Hilbert's rule redundant because being true does not have its usual meaning. Since being true reduces to being stipulated the question Which stipulations are allowed and which are not? must be treated independently. Thus the consistency condition must be checked before ``axioms become true'', i.e. before one stipulates that a given set of expressions represents a set of mathematical truths. Such a checking requires a special notion of consistency, which applies to linguistic expressions having no definite truth-values. At the time of writing his letter to Frege Hilbert did not formulate yet the appropriate notion of consistency rigorously; we shall shortly see how he tried to solve this problem afterwards.    

(2) Notice the peculiar form of Hilbert's axioms, which involves terms with variable meaning. An expression of this form turns into a proposition only when the meaning of all its terms becomes determined.  So in order to stipulate that a set of axiom-like expressions represents a set of axioms, Hilbert needs to assume that there exist ``things defined by these axioms'', which (a) make all terms in these axioms meaningful and (b) which make these axioms true. In the above quote Hilbert states that the existence of such things is always granted when the corresponding set of axioms is consistent.  (``If the arbitrarily posited axioms [..] do not contradict each other, then [..] the things defined by these axioms exist''.) Notice that the existence of these things has no other prerequisites except consistency. Whence there arise two mutually related questions:  \emph{What} are the things ``defined by axioms''?  and \emph{How} the axioms ``define'' them? Let me consider these two questions in turn.  
 
The former question has at least three different answers. The \underline{first} general answer is this: given an expression like A1.1, which bears on ``bare things'' and ``bare relations'' of multiple types one instantiates these things and these relations in one's mind and so get what Hilbert after Kant calls \emph{objects of thought} or \emph{thought-things} (\emph{Gedankendinge} in German), which are related by corresponding \emph{thought-relations}
\footnote{
Compare with Orwell's \emph{thoughtcrimes}. 
}. 
These thought-things and thought-relations exist merely in virtue of the fact that one thinks of them consistently. They may be or be not supported by some sensual intuitions; the sensual intuition is a separate issue which must not be confused with the capacity to instantiate objects and relations between objects as such. This latter capacity can be also called intuition - not in the sense of Kant's \emph{Transcendental Aesthetics} but exclusively in the sense of Kant's \emph{Doctrine of Method} \cite{Kant:1999}. Hintikka \cite{Hintikka:1997b} quite rightly stresses the fundamental role of this restricted notion of intuition in Hilbert's Axiomatic Method. Even when we think of mathematical objects as ``bare things'' without associating with these things anything over and above the relations stipulated through axioms like A1.1 we think about these objects, by Kant's word, \emph{in concreto} (which shows, by the way, that the usual characterization of such object as \emph{abstract} is somewhat misleading). The mathematical intuition in the relevant restricted sense of the term is the capacity to think concretely about objects and relations between objects without associating to these objects and these relations any additional qualities.     

The \underline{second} answer concerns the role of sensual intuition. Remind that in the introductory part of his \emph{Foundations} of 1899 Hilbert says that his geometrical axioms  ``express certain related fundamental facts of our intuition''. Earlier in 1891 he made the following remark: 

\begin{quote}
Geometry is the science that deals with the properties of space. [..] I can never penetrate the properties of space by pure reflection, much as I can never recognize the basic laws of mechanics, the
law of gravitation or any other physical law in this way. Space is not a product of my reflections. Rather, it is given to me through the senses. (quoted after Corry \cite{Corry:2000}, p. 44)
\end{quote}

In 1894 Hilbert develops this view on geometry: 

\begin{quote}
Among the appearances or facts of experience manifest to us in the observation of nature, there is a peculiar type, namely, those facts concerning the outer shape of things, Geometry deals with these facts [..]. Geometry is a science whose essentials are developed to such a degree, that all its facts can already be logically deduced from earlier ones. Much different is the case with the theory of electricity or with optics, in which still many new facts are being discovered. Nevertheless, with regards to its origins, \underline{geometry is a natural science} (\emph{ib.} p.45)
\end{quote}

\begin{quote}
[A]ll other sciences-above all mechanics, but subsequently also optics, the theory of electricity, etc.- should be treated according to the model set forth in geometry. (\emph{ib.} p.45)
\end{quote}

What Hilbert says here about the empirical character of Geometry prima facie is not compatible with his notion of Geometry as a free creation of mind expressed in his letter to Frege quoted above. It is not impossible, of course, that during this period of time Hilbert had conflicting ideas about the nature of Geometry and could contradict himself. However it seems me suggestive to try to reconcile the two notions of Geometry. As a part of the pure mathematics Geometry is treated as a free creation of mind; the fundamental question here is whether or not the given set of geometrical axioms is consistent while the question where those axioms come from is irrelevant. As a natural science Geometry seeks to express properties of the physical space through an appropriate set of axioms, then ``logically deduce'' from these axioms some further geometrical propositions and finally check these deduced propositions against properties of the physical space. So the two Geometries well fit together: the physical geometry takes care about choosing axioms properly while the mathematical geometry takes care about the consistency of any proposed set of geometrical axioms, about the deduction of new theorems from these axioms and some other relevant problems. This epistemological model is applicable to all natural sciences; what makes geometry ``more mathematical'' than say, the theory of electricity, is the fact that geometry easier allows for an axiomatic treatment because its ``essentials'' are  better developed.       

So we may consider geometry in a larger sense, which combines the axiomatic mathematical geometry, on the one hand, and the empirical physical geometry, on the other hand. Objects of this combined geometry are no longer bare individuals but spatial physical bodies, light rays, etc. Interestingly, the traditional notion according to which geometry presents properties of the physical space in an idealized form is irrelevant to Hilbert's axiomatic setting. Geometrical objects are thought of here either as bare individuals detached from any sensual intuition or as physical bodies as they are perceived by senses; Hilbert's epistemic scheme, which we reconstruct on the basis of the above passages, does not include any intermediate ``ideal'' element between the axiomatic logical reasoning and the sensual perception. We shall see, however, that in his later works Hilbert introduces such ideal elements (\textbf{2.4}).

The \underline{third} answer to the question about Hilbert's mathematical ``things'' and their existence concerns the possibility of interpreting axioms of a given axiomatic theory in terms of another mathematical theory. For example with the help of standard tools of Analytic Geometry A1.0 and other Hilbert's axioms translate into true propositions about real numbers. An interpretation $M$ that translates all axioms of a given axiomatic theory $A$ into true propositions of another theory $T$ is called a \emph{model} of $A$ in $T$; one says also that axioms of $A$ are true in model $M$. Suppose we know which proposition of $T$ is true and which is false. This allows one to reverse the order of ideas about $A$. Observe that in order to check whether or not axioms of $A$ are true in $M$ one does not need to establish consistency of this set of axioms in advance. Moreover, if axioms of $A$ are true in $M$ (i.e., if $M$ is indeed a model of $A$) then one may conclude that $A$ is consistent! Remind Hilbert's remark according to which any consistent set of propositions can be made by a fiat into a system of axioms, which are true and meaningful. Now we proceed the other way round: we first check that our axioms are true and meaningful in some model \footnote{By meaningfulness of a given axiom in a given model I mean the bare fact that this axiom translates into a meaningful proposition. Saying that every axiom of a given theory is true and meaningful in a model of this theory is, of coarse, pleonastic.} and on this basis conclude that the given set of axioms is consistent. However this conclusion is not valid unless $T$, which is the background theory of $M$, is consistent in its turn. So what the above argument really proves is not the absolute but only the \emph{relative} consistency of $A$, i.e., the proposition of the form ``if $T$ is consistent then $A$ is also consistent''.

From a mathematical point of view this third way of interpreting Hilbert-style axioms turns to be the most productive. Already in his \emph{Foundations} of 1899 Hilbert applies this method systematically; in the course of 20th century this method develops  into the modern \emph{model theory}, which remains today an active field of mathematical research still having some philosophical flavor. I would like to stress here that interpreting a Hilbert-style axiomatic theory in terms of another mathematical theory and interpreting such a theory in some intuitive terms \emph{directly} are two very different issues.  Since both procedures go under the same title of ``interpretation'' they are too often confused in the current debates. The idea that Hilbert's axiomatic theory of Euclidean geometry can be either interpreted ``as usual'', i.e., by associating with the terms ``point'', ``straight line'', ``between'', etc. their ``usual'' intuitive meanings, or alternatively, be interpreted arithmetically by  identifying points with pairs of numbers, etc., is plainly misleading because it puts under the same title of interpretation two procedures, which do not belong to the same general type.  

(3) Let us finally discuss Hibert's view according to which axioms of a given mathematical theory ``define'' objects of this theory. Since Hilbert's axioms refer only to bare ``things'' and bare relations and since, according to Hilbert, any consistent set of such axioms allows one to produce a ``system of things'' $S$ satisfying these axioms by a fiat (or more precisely by the very fact that one forms consistent thoughts ``about'' certain things), such $S$ can be thought of as the ``definiendum'' of the axioms. One may ask however whether a given consistent set of axioms defines the corresponding system $S$ \emph{uniquely}. Here is what Hilbert says about this in the same letter to Frege:

\begin{quote}
You say that my concepts, e.g. ``point'', ``between'', are not unequivocally fixed [..].
But surely it is self-evident that every theory is merely a framework or schema of
concepts together with their necessary relations to one another, and that basic
elements can be construed as one pleases. If I think of my points as some system or
other of things, e.g. the system of love, of law, or of chimney sweeps [..] and then
conceive of all my axioms as relations between these things, then my theorems, e.g.
the Pythagorean one, will hold of these things as well. In other words, each and every
theory can always be applied to infinitely many systems of basic elements. For one
merely has to apply a univocal and invertible  one-to-one transformation and stipulate
that the axioms for the transformed things be correspondingly similar. (cit. by \cite{Frege:1971},
p.13).
\end{quote}

There are two important ideas in this passage. First Hilbert stresses here once again that in his  axiomatic setting primitive geometrical terms have no intrinsic meaning: any system of things (i.e., model) satisfying Hilbert's axioms counts as an Euclidean space. This point has been already discussed earlier in this Chapter and I shall not return to it. Then follows this crucial observation: given a model $M$ of a given axiomatic theory one can always get another model $M'$ of the same theory through a one-to-one transformation of elements of $M$ into elements of the new model $M'$ in such a way that relations between elements of $M'$ also satisfy the axioms of the given theory. In the modern language the kind of transformation described here by Hilbert is called \emph{isomorphism}. Apparently Hilbert thinks here about an axiomatic theory that determines its models \emph{up to isomorphism}, i.e., such that all its models are \emph{isomorphic}, i.e., are transformable into each other by some isomorphisms. Such theories are called today  \emph{categorical}. (Beware that \emph{that} sense of being categorical has nothing to do with the category theory!) Isomorphic models can be seen as ``equal'' and representing the same  \emph{structure}, which is invariant under transformations between these model. This leads to a philosophical view on mathematics known as \emph{mathematical structuralism}; according to this view structures are basic mathematical objects. I consider the mathematical structuralism and its significance for  the Axiomatic Method in Chapter \textbf{8}. The idea of the ``replacement of equality by isomorphism'' is also discussed in Chapters \textbf{5, 6}. 

Precipitating this further discussion I would like only to stress here that not every Hilbert-style axiomatic theory is categorical. In fact this is a rather strong property that most of useful axiomatic theories do not enjoy. Apparently Hilbert didn't see this problem before he first published his \emph{Foundations} in 1899; however in his lecture \emph{On the Concept of Number } \cite{Hilbert:1900} delivered in the same year 1899 and published in 1900 Hilbert already introduces an "axiom of completeness" (Vollstandigkeitsaxiom), which requires from any model of a given theory (this time it was arithmetic) this maximal property: any model $M$ of the given theory extended with some new elements is no longer a model. Then he proves that among all models of his theory (without the completeness axiom) there is only one model (up to isomorphism, of course!), which also satisfies the completeness axiom, see \cite{Corry:2004} p. 160 for details. The second edition of Hilbert's \emph{Foundations of Geometry} appeared in 1903 \cite{Hilbert:1903} already contains a geometrical axiom of completeness.  

Let me now return to the question about the sense of being formal. Frege and his contemporaries called Hilbert's axioms for geometry \emph{formal} extending the sense of being formal used by   Trendelenburg \cite{Trendelenburg:1862} when he made popular the expression ``formal logic''. Formal logic in Trendelenburg's sense is by and large what Kant calls \emph{general} logic as distinguished from his transcendental logic: the formal or general logic takes into account only the form of reasoning and is neutral with respect to its content (Kant's transcendental logic is not wholly neutral with respect to the content of reasoning because it takes into account the difference between objects of possible experience and thought-objects of other sorts). The modern usual sense of ``formal logic'' retains this older sense but does not reduce to it because it also includes the idea of symbolic mathematical presentation of logical form, see \textbf{2.3} below. Although Hilbert's axioms are not logical tautologies that hold for all objects and all relations whatsoever, they represent \emph{logical forms} of propositions obtained through the usual contentual reading of the same axioms and allow for alternative instantiations of this logical form (i.e., for alternative models). For example, by reading words ``point'' and ``straight line'' in Hilbert's First Axiom A.1.0 naively, i.e., by associating with these words their usual linguistic meanings, one gets an universal contentual proposition about (all) points and (all) straight lines (in the relevant domain where this axiom applies). However the intended meaning of this axiom expressed more explicitly in A1.1  specifies only a property of abstract relation between abstract objects of two different types. This description is purely formal in the sense that it fixes no domain of objects and no concrete relation. It specifies a specific form of relation between objects but specifies no particular relation and no particular object and no particular type of objects.

Today we would qualify the theory of Hilbert's \emph{Foundations} of 1899 as informal or semiformal at best. This is because this theory is formulated in the natural German with the help of some symbols like any typical introductory mathematical text. Today's paradigmatic examples of formal theories is given by axiomatic theories of sets and of arithmetic like ZF and PA. These latter theories differ from the theory of Hilbert's \emph{Foundations} of 1899 first of all by their symbolic syntax. In \textbf{2.3} we shall see how the idea of using such a symbolic syntax  combines with Hilbert's earlier approach described in this Section. We shall see that the symbolic approach involves some epistemological ideas, which do not make part of the traditional notion of being formal relevant to Hilbert's \emph{Foundations} of 1899. So being formal, semiformal and informal should not be thought of only as a matter of degree. 

\section{\emph{Foundations} of 1899:  Logicality and Logicism}

Consider once again the Hilbert's First  Axiom A1.0         

\begin{quote}
Any two distinct points of a straight line completely determine that line
\end{quote}

and remind that certain words in this sentence including  the words ``points'' and ``straight line'', are \emph{not} supposed  to be understood  in their usual sense. Now remark that some other words like ``any'' and ``two''  \emph{are} supposed to be understood in the usual sense. Clearly this second category of words plays in Hilbert's \emph{Foundations} of 1899 an essential role: unless at least some words in these axioms are meaningful the axioms reduce to an abracadabra! In the last Section I elaborated on words of the former category, now let us look more attentively on words of the latter category. First of all let us see how exactly words are sorted into two sorts here. Words of the first sort refer to primitive geometrical concepts like point, straight line and between (whether these primitive concepts are understood traditionally or in the sophisticated formal way explained above). What about words of the second category? 

In order to answer this question it is helpful to paraphrase A1.0 as follows:

\begin{quote}
If different points $A,B$ belong to straight line $a$ and to straight line $b$ then $a$ is identical to $b$  
\end{quote}

Now leaving out  \emph{geometrical} words and expressions  ``points $A,B$'', ``straight line $a$'', ``straight line $b$'', ``belong to'' we get this list: ``if'', ``different'', ``and'', ``then'', ``is'', and ``identical to''. So the last paraphrase helps us to see that the words belonging to the second list stand for  \emph{logical} notions. 

How to distinguish between logical and non-logical terms more formally? There exist in the literature two main approaches to defining the notion of \emph{logicality}: one develops the idea of logic being \emph{content-free} (so that logical signs are understood as the ``punctuation marks'') and the other that describes itself as \emph{semantic} develops the idea of logic being \emph{content-invariant} \cite{Bonnay:2006}. This later approach dates back to Tarski's proposal \cite{Tarski:1986} to identify logical notions with invariants of all permutations of elements of some given set \footnote{Bonnay \cite{Bonnay:2006} formulates Tarski's Thesis as follows:
begin{quote}
Given a set $M$, an operation $Q_{M}$ acting on $M$ is logical iff it is invariant under all permutations
end{quote}
}. 

This latter approach obviously better squares with Hilbert's Axiomatic Method; the idea here is to make the fixity of meanings of logical terms and the variability of meanings of non-logical terms into a formal criterion allowing one to distinguish between these two sorts of terms.  Tarski accounts for this fixity as invariance under permutations of elements of a given set of individuals (which represents here a certain universe of discourse). This approach to logicality is motivated by Klein's \emph{Erlangen Program} in geometry \cite{Klein:1872};  it establishes a conceptual link between Klein's and Hilbert's works in foundations of geometry, which is both conceptually significant and historically plausible. I discuss it in \textbf{8.3} below. 

Now I would like only to stress Hilbert's fundamental assumption behind his Axiomatic Method (as  presented in his  \emph{Foundations} of 1899) according to which logic is the ground layer foundation of all theories built axiomatically. As Hintikka puts this, for Hilbert

\begin{quote}
The basic clarified form of mathematical theorizing is a purely logical axiom system. (\cite{Hintikka:1997b}, p.20)
\end{quote}

This does not mean, of course, that Hilbert like Russell  in \cite{Russell:1903} tries to \emph{reduce} mathematics to logic. This later version of logicism is certainly not Hilbert's.  In this book I shall use the term ``mathematical logicism'' in a broader sense of epistemic primacy of logic over mathematics. In this broader sense Hilbert's view on mathematics qualifies is another version of logicism.  

As I have already mentioned in \textbf{1.4} the idea of  logic and metaphysics as a foundation of all science dates back to Aristotle. This idea  had apparently little or no influence on Greek mathematics (that followed Euclid rather than Aristotle) but later became quite influential in the medieval Scholasticism. The Early Modern mathematically-laden science that triumphed with Newton's  \emph{Principia} largely rejected the old scholastic pattern of theory-building and developed a very different notion of scientific theory that was described in general terms by Kant in his \emph{Critique of the Pure Reason}. The Kantian philosophy of science and mathematics remained the mainstream until the beginning of the 20th century when the old scholastic pattern of theory-building kicked back under the new name of modern Axiomatic Method.

In Kant's view the \emph{objectivity} of pure mathematics (which underlies the objectivity of the mathematically-laden empirical science) roots in its \emph{objecthood}, i.e., in the universal schemata according to which one constructs mathematical \emph{objects} - but not just in the general character of mathematical concepts. Making difference between the mathematical objectivity and the universal logical validity, according to Kant, is crucial for differentiating between the mathematical reasoning and the philosophical speculation. Here is the famous passage:  

\begin{quote}
Give a philosopher the concept of triangle and let him try to find out in his way how the sum of its angles might be related to a right angle. He has nothing but the concept of figure enclosed by three straight lines, and in it the concept of equally many angles. Now he may reflect on his concept as long as he wants, yet he will never produce anything new. He can analyze and make distinct the concept of a straight line, or of an angle, or of the number three, but he will not come upon any other properties that do not already lie in these concepts.   But now let the geometer take up this question. He begins at once to \underline{construct a triangle}.  Since he knows that two right angles together are exactly equal to all of the adjacent angles that can be drawn at one point on a straight line, he extends one side of his triangle and obtains two adjacent angles that together are equal to the two right ones. [..] In such a way through a chain of inferences that is always \underline{guided by intuition}, he arrives at a fully illuminated and at the same time general solution of the question.'' (\emph{Critique of Pure Reason} \cite{Kant:1999}, A 716 / B 744)
\end{quote}

 Kant's philosophy of mathematics and mathematically-laden science is based upon an analysis of his best contemporary science as represented by Newton's  \emph{Principia} \cite{Friedman:1992}. This does not mean, of course, that Kant derives his philosophical principles from the principles of Newtonian physics; Kant's  \emph{critical} philosophy rather aims at explaining how the type of knowledge best represented by the Newtonian physics is possible (as an objectively valid knowledge). Anyway this method of philosophical work makes Kant's philosophy strongly dependent on the contemporary mathematics and science. Cohen, Natorp and other neo-Kantians who wished to sustain the Kantian project of critical philosophy in the 19th century realized this fact very clearly and made efforts to supply the Kantian philosophy with a historical dimension allowing one to keep track of the progress in sciences and mathematics. \cite{Heis:2007}. It was not quite clear in the 19th century and it still remains a matter of controversy today which (if any) features of Kant's original approach remain sustainable in the context of the current science and mathematics, and which features of this original approach are hopelessly outdated. More radically one may wonder if there is anything at all in Kant's analysis that survives all the dramatic changes in science and pure mathematics that have happened since Kant's own time.   

In spite of a number of interesting attempts of upgrading the Kantian philosophy of mathematics in order to account for new mathematical developments (like the invention of non-Euclidean geometries) at certain point the Kantian line in the philosophy of mathematics has been largely abandoned. Bertrand Russell's intellectual development is representative in this sense: after publishing in 1897 his Kantian \emph{Essay on Foundations of Geometry} \cite{Russell:1897} and a short romance with Hegel \cite{Hylton:1990} Russell learns in 1900 about new works in mathematical logic, publishes during the same year an essay on Leibniz \cite{Russell:1900}, who now becomes the right philosophical ancestor, and  already in 1903 publishes the \emph{Principles of Mathematics} \cite{Russell:1903} where the subject is developed on new logicist grounds. In the \emph{Introduction} to this book Russell explains his attitude to the Kantian line of thought as follows:

\begin{quote}
It seemed plain that mathematics consists of deductions,
and yet the orthodox accounts of deduction were largely or wholly
inapplicable to existing mathematics. Not only the Aristotelian
syllogistic theory, but also the modem doctrines of Symbolic Logic,
were either theoretically inadequate to mathematical reasoning, or at
any rate required such artificial forms of statement that they could not
be practically applied. In this fact lay the strength of the Kantian
view, which asserted that mathematical reasoning is not strictly formal,
but always uses intuitions, i.e. the a priori knowledge of space and
time. Thanks to the progress of Symbolic Logic, especially as treated
by Professor Peano, this part of the Kantian philosophy is now capable
of a final and irrevocable refutation. By the help of ten principles of deduction and ten other premisses of a general logical nature (e.g. implication is a relation"), all mathematics can be strictly and formally deduced. [..]\\
The general doctrine that all mathematics is deduction by
logical principles from logical principles was strongly advocated by
Leibniz... But owing partly to a faulty logic, partly to belief in the logical necessity
of Euclidean Geometry, he was led into hopeless errors in the endeavour
to carry out in detail a view which, in its general outline, is now known
to be correct. The actual propositions of Euclid, for example, do not
follow from the principles of logic alone ; and the perception of this fact
led Kant to his innovations in the theory of knowledge. But since
the growth of non-Euclidean Geometry, it has appeared that pure
mathematics has no concern with the question whether the axioms
and propositions of Euclid hold of actual space or not .....  What pure mathematics asserts is merely that the Euclidean propositions follow from the Euclidean axioms, i.e.,
it asserts an implication. ....  We assert always in mathematics
that if a certain assertion $p$ is true of any entity $x$ or of any set of
entities $x, y, z ...$, then some other assertion $q$ is true of those entities ;
but we do not assert either $p$ or $q$ separately of our entities.
\end{quote}

The above argument, which is supposed to refute Kant, obviously begs the question. From the outset Russell takes it for granted that ``mathematics consists of deductions'' and his following remarks make it clear that by deduction Russell means here a \emph{logical} deduction, i.e. a deduction of propositions from certain other propositions according to some general rules, which are not specific for mathematics. This statement overtly contradicts what Kant says about mathematics, and the following Russell's argument only provides this first statement with some additional details but does not constitute any philosophical objection to Kant.  Kant's own objection to the Leibnizian view on mathematics, to which Russell adheres here, is this. From a \emph{formal} point of view (i.e. as far as only \emph{logical form} of sentences is taken into consideration) mathematics is no different from a mere metaphysical speculation; a speculative metaphysical theory can be developed on an axiomatic basis just like any mathematical theory (think about Spinoza's \emph{Ethics} for example). What makes the crucial difference between mathematics and speculation is the fact that mathematics constructs its \emph{objects} according to certain rules, while speculation proceeds with concepts without being involved in any similar constructive activity. The fact that the speculative thought may also posit some entities falling under these concepts from the Kantian viewpoint does not constitute an objection because such stipulated entities doesn't qualify as \emph{objects} in the strong Kantian sense of the term. Behind an \emph{object} there is a procedure (governed by a certain \emph{rule}) that constructs it; speculative entities are stipulated as mere \emph{thought-things} falling under given descriptions by a fiat. This is the reason why the pure mathematics is \emph{objective} in the sense in which the pure speculation is not. What makes the pure mathematics objective is the rule-like character of object-construction. The formal logical consistency is a necessary but not sufficient condition for claiming that a given axiomatic theory is objectively valid. This objective character of mathematics, according to Kant, allows for application of mathematics in natural sciences (I leave however this further point aside).   Russell's critique of Kant  in the  \emph{Principles of Mathematics} simply does not take into account the Kantian problem of separation of the pure mathematics from the pure speculation. In this respect Russell's Leibnizian approach to mathematics is more traditional than Kant's and in Kantian terms qualifies as dogmatic. Not surprisingly Russell provides his philosophy of mathematics with a metaphysical doctrine that he calls the  \emph{logical atomism}. This is how he describes the relation of this doctrine to logic and mathematics in the  \emph{Introduction} to his \cite{Russell:1918}: 

\begin{quote}
As I have attempted to prove in  \emph{The Principles of Mathematics},
when we analyse mathematics we bring it all back to logic. It all
comes back to logic in the strictest and most formal sense. In the
present lectures, I shall try to set forth in a sort of outline, rather
briefly and rather unsatisfactorily, a kind of logical doctrine
which seems to me to result from the philosophy of mathematics - 
not exactly logically, but as what emerges as one reflects: a
certain kind of logical doctrine, and on the basis of this a certain
kind of metaphysic. 
 \end{quote}
 
As a recent biographer describes Russell's work during this early period of his career
 
 \begin{quote} 
 From August 1900 until the completion of Principia Mathematica in 1910 Russell was both a metaphysician and a working logician. The two are completely intertwined in his work: metaphysics was to provide the basis for logic; logic and logicism were to be the basis for arguments for the metaphysics. (\cite{Hylton:1990}, p. 7-8) 
 \end{quote}

Thus we can see how an older pattern of intellectual work, which many people in the 19th century believed to be definitely sublated by Kant's critical philosophy and other developments, reemerged in the beginning of the 20th century in the context of new mathematics and new symbolic logic. An attempt to describe the general intellectual context of that time would obviously lead me too far but it is my understanding that Russell's case in an extreme form represents a more general intellectual tendency. Even more important is the fact that this tendency towards the revival of the traditional alliance between logical and metaphysical thinking is still very much alive today, and in fact since 1900 this intellectual project has firmly established itself in the philosophical school known as \emph{Analytic Philosophy} (as well as in some other branches of today's philosophy). So my critique of Russell of early 1900ies and, in particular, my attempts to revendicate the Kantian and the Hegelian (see \textbf{4.8} below) lines of philosophical thought in the context of recent mathematics of our own times, aims primarily at the modern proponents of this traditional alliance.    

Russell's interpretation of Kant's work in the philosophy of mathematics as an attempt to fill logical gaps appearing when one tries to reconstruct Euclid's geometry with Aristotle's syllogistic logic hardly correctly describes Kant's intention. However these Russell's words are helpful for a better understanding of his own project. Russell suggests two independent reasons why there are such logical gaps: first, because Euclid's geometry is logically imperfect and, second, because Aristotle's logic is not appropriate for doing mathematics. However the new mathematics (including non-Euclidean geometries) and the new symbolic logic taken together, according to Russell, wholly fix the problem making Russell's Leibnizian dream real. What Russell's \emph{Principles of Mathematics} aim at is made clear by the following lines that I take from the \emph{Preface} to this work:

 \begin{quote}
The second volume, in which I have had the great good fortune
to secure the collaboration of Mr A. N. Whitehead, will be addressed
exclusively to mathematicians; it will contain chains of deductions,
from the premisses of symbolic logic through Arithmetic, finite and
infinite, to Geometry, in an order similar to that adopted in the present
volume ; it will also contain various original developments, in which the
method of Professor Peano, as supplemented by the Logic of Relations,
has shown itself a powerful instrument of mathematical investigation.   
 \end{quote}
 
 (The planned second volume of \emph{The Principles of Mathematics} appeared later as a co-authored independent three-volume work \cite{Russell&Whitehead:1910-1913}.) 
     
Thus we can see that in early 1900ies Hilbert was not alone who thought about logic as the ground layer foundation of mathematics. (This is in spite of the fact that unlike Hilbert in his \emph{Foundations} of 1899 officially sticks to Kant! - notice the Kant's quote used as the epigraph in this book. We shall  shortly see (\textbf{2.4}) that Hilbert's latter version of Axiomatic Method better fit's the Kantian view on mathematics than this earlier version.) However there were also strong opposing voices during the same period of time. Among prominent critics of logical approaches in foundations of mathematics were Poincar\'e \cite{Detlefsen:1992} and Brouwer. Consider, for example, this Brouwer's passage written in 1907:

 \begin{quote}
About mathematical reasoning, I show in the beginning of the chapter
that it is no logical reasoning, that it uses the connectives of logic only
because of the poverty of language, and thus may perhaps keeps alive
the language accompaniment even after the human intellect has already
long ago outgrown the logical argument itself. For, far from the fact that it
would be a ``strange company'' that does not reason logically, I believe that
it is only a matter of inertia, that the words that go with it [i.e., logic] as yet
still exist in modern languages. A pure use of these words hardly occurs,
and [in] impure [form] they are used in daily life, where they have led to
all kinds of misunderstanding and dogmatism [..]. Those misconceptions arose, not because
of insufficient mathematical insight, but because mathematics, lacking a
pure language, makes do with the language of logical reasoning, although
its thoughts reason not logically, but mathematically, which is something
totally different. (quoted after \cite{Dalen:2000}, p. 128-129)
 \end{quote}

I warn the reader that Brouwer's concern expressed in the above quote is not met by using the formal\emph{intuitionistic logic} instead of classical logic in foundations of mathematics because this replacement of logic leaves untouched the assumption about the primacy of logic over mathematics; such a replacement only translates (through Heyting's formalization of the intuitionistic logic \cite{Heyting:1956}) some Brouwer's ideas into a logicist foundational framework. However the very fact that today we have more than one candidate logic for building foundations of mathematics is remarkable. The logicist view on mathematics was particularly appealing in the beginning of the 20th century because at that time the traditional geometry was already split into its Euclidean and multiple non-Euclidean versions but logic still preserved its traditional unity. Today we live in a very different environment.  Here how Gabbay describes it as for 1994:
 
\begin{quote}
In recent years we have witnessed a very strong and fruitful interaction
between traditional logic on the one hand and computer science and 
Artificial intelligence on the other. As a result, there was urgent need for logic
to evolve. New systems were developed to cater for the needs of 
applications. Old concepts were changed and modified and new concepts came
into prominence. The community became divided. Many expressed 
themselves strongly, both for and against, the new ideas. Papers were rejected
or accepted on ideological grounds, as well as technical substance.
In this atmosphere, it seemed necessary to clarify the basic concepts
underlying logic and computation, especially the very notion of a logical
system. [..] The views among members of the community are varied and in many cases, very strongly held. There is at
one extreme the pluralistic view, expressed to me once in a meeting by a
distinguished colleague who said something like ``we use logics like we use
computer languages''. At the other end of the spectrum there is the view
of those who believe there is only one true logic, and all the rest is 
nonsense. Of course there exist several proposals for this true logic with their
respective bands of followers. (\cite{Gabbay:1994}, \emph{Preface},  p. v) 
\end{quote}

So anyone who holds today a logicist view on mathematics (in the broad sense of ``logicist'' explained above) needs, first, to specify  \emph{which} is his or her favorite logic used in foundations, second, to explain  \emph{why} this particular logic is the most appropriate for the purpose \footnote{Alternatively one may consider a variety of ``Non-Classical Mathematics'' each based on its proper logic \cite{Vasyukov:2009}.}. 
and, third (which is perhaps the hardest task), to explain why and in which sense one's favorite logic qualifies as logic. This problematic character of modern logic does not imply that the logicist view on mathematics is no longer tenable but it certainly shows that this view can not and should not be taken for granted. In Chapter \textbf{9} we shall see how the New Axiomatic Methods deals with this new degree of freedom of today's axiomatic thought. 

A more detailed argument against the mathematical logicism has been given in 1907 by Ernest Cassirer \cite{Cassirer:1907}  who continued to push the Kantian line taking into account newest mathematical developments of his time. Referring to Russell \cite{Russell:1903} and new formal logical methods under the name of ``logistics'' Cassirer says: 

 \begin{quote}
 Here rises a problem that lies wholly outside the scope of ``logistics'' [..] All empirical judgements belong to their domain: they must respect the limits of experience. What logistics develops is a system of hypothetical assumptions about which we cannot know, whether they are actually established in experience or whether they allow for some immediate or non-immediate concrete application. According to Russell even the general notion of magnitude does not belong to the domain of pure mathematics and logic but has an empirical element, which can be grasped only through a sensual perception. From the standpoint of logistics the task of thought ends when it manages to establish a strict deductive link between all its constructions and productions. Thus the worry about laws governing the world of objects is left wholly to the direct observation, which alone, within its proper very narrow limits, is supposed to tell us whether we find here certain rules or a pure chaos. [According to Russell]  logic and mathematics deal only with the order of concepts and should not care about the order or disorder of objects.  As long as one follows this line of conceptual analysis the empirical entity always escapes one's rational understanding. The more mathematical deduction demonstrates us its virtue and its power, the less we can understand the crucial role of deduction in the theoretical natural sciences.(\cite{Cassirer:1907}, p. 43)
 \end{quote}
 
So, according to Cassirer, what the formal logical foundations of mathematics can \emph{not} possibly provide (whatever system of formal logic is one's favorite) are the notions of objecthood and objectivity appropriate for doing the modern mathematically-laden empirical science (i.e., the \emph{Galilean} science as I called it \textbf{3.3} above). The popular idea to equate the notion of object with that of logical individual, which stems from Frege \cite{Parsons:2008}, not only leaves this problem open and but also hides it by eliminating a useful terminological distinction, which helps Kant to distinguish objects of possible experience from thought-things of other sorts. Although Cassirer does not provide any concrete solution of this problem he stresses the relevance of Kantian approach to the modern science in the following words (the second phrase I used as an epigraph to this book):

\begin{quote}
 The principle according to which our concepts should be sourced in intuitions means that they should be sourced in the mathematical physics and should prove effective in this field. Logical and mathematical concepts must no longer produce instruments for building a metaphysical ``world of thought'': their proper function and their proper application is only within the empirical science itself. (\cite{Cassirer:1907}, p. 43-44) 
 \end{quote}

The fact that the modern logic indeed tends to create ``metaphysical worlds of thought'' rather than make itself into a part of empirical science, and that today's mainstream philosophy of logic encourages and justifies this overtly metaphysical tendency (usually by presenting it as an innocent intellectual game), appears to me very worrying. Although the hostile attitude towards logic and its mainstream philosophy, which is widely spread in mathematical circles, demonstrates a healthy intellectual reaction, such a negative reaction by itself does not solve the problem. In Chapter \textbf{4} I come back to this Cassirer's argument and after Lawvere point to a way out (\textbf{4.8}).

\section{Axiomatization of Logic: Logical Form versus Symbolic Form}

In his address of 1917 already quoted above  Hilbert says among other things the following:

\begin{quote}
[I]t appears necessary to axiomatize logic itself and to prove that number theory and set theory are only parts of logic. This method was prepared long ago (not least by Frege's profound investigations); it has been most successfully explained by the acute  mathematician and logician Russell. One could regard the completion of this magnificent Russellian enterprise of the axiomatization of logic as the crowning achievement of the work of axiomatization as a whole. (\cite{Hilbert:1918} p. 1113)
\end{quote}

Leaving now aside the purported reduction of number theory (arithmetic) and set theory to logic let us focus on the idea of \emph{axiomatization of logic}. By calling the axiomatization of logic the ``crowning achievement of the work of axiomatization as a whole'' Hilbert suggests that the axiomatization of logic is a continuous extension of the axiomatization of geometry, arithmetic and of any other part of mathematics or natural  science. However the notion of axiomatization, which I have tried to reconstruct above on the basis of Hilbert's  \emph{Foundations} of 1899 does not immediately allow for such an extension. In the nutshell the axiomatization in the sense of  \emph{Foundations} of 1899 works like this: using some fixed logical vocabulary one produces a finite list of axioms, which refer only to abstract objects and abstract relations; an intended ``naive'' interpretation of these axioms and of all theorems derivable from these axioms is supposed to capture the content of the corresponding informal theory in a more precise and ``logically clear'' form. Notice that this whole procedure applies logic as a tool; an axiomatizer needs to have this tool in a ready-made form just like a carpenter needs a ready-made hammer for putting down a nail. So if the above reconstruction of Axiomatic Method is correct in order to axiomatize logic one needs to use logic. How this may possibly work?

Instead of speculating further on this matter let us see how Hilbert axiomatizes logic in his course on \emph{Theoretical Logic} \cite{Hilbert&Ackermann:1928} co-authored with Ackermann and first published in 1928. This work is greatly influenced by earlier works by Frege and Russell; I shall not however trace here these influences but consider Hilbert's and Ackermann's  book on its own rights. The \emph{Introduction} to this book opens with the following words:

\begin{quote}
\emph{Mathematical logic}, also called \emph{symbolic logic} or \emph{logistic}, is
an extension of the formal method of mathematics to the field of logic. It employs for logic a symbolic language like that which has long been in use to express mathematical relations. In  mathematics it would nowadays be considered Utopian to think of using only ordinary language in constructing a mathematical discipline. The great advances in mathematics since antiquity, for instance in algebra, have been dependent to a large extent upon success in finding a usable and efficient symbolism. (quoted after English translation \cite{Hilbert&Ackermann:1950}, p. 1)
\end{quote}

We see that from the very beginning Hilbert and Ackermann introduce here a new kind of logic, which they call mathematical or symbolic\footnote{Saying that symbolic logic is a ``new'' kind of logic I mean that this kind of logic is new with respect to the ``informal'' logic used in Hilbert's \emph{Foundations} of 1899; I don't mean, of course, that symbolic logic first appears in Hilbert and Ackermann's book. In a part of the \emph{Introduction} to this book, which I do not quote here, the authors provide a brief historical sketch of symbolic logic tracing its history back to Leibniz.}. As we shall shortly see Hilbert's notion of axiomatization of logic makes sense \emph{only} in a symbolic setting. The following description of mathematical (symbolic) logic as an ``extension of the formal method of mathematics to the field of logic'' is puzzling. If by formal method one understands the Axiomatic Method in the sense of Hilbert's \emph{Foundations} of 1899 then it is unclear how this application can make logic symbolic. Indeed, Hilbert's \emph{Foundations} of 1899 is written with the usual mixture of informal prose, geometrical diagrams and the traditional algebraic and geometrical symbols; Hilbert's \emph{formal} approach developed in this book is no more symbolic than the approach taken in any other elementary geometry textbook published in the 19th century. 

Notice also that in the above passage Hilbert talks about application of the ``formal method of mathematics'' in logic. So he thinks here about the formal method of mathematics as something established independently from logic and then suggest to ``extend'' this method to the new field of logic. However the formal method of \emph{Foundations} of 1899 is certainly not independent of logic. So talking about ``formal method'' in the above quote Hilbert and Ackermann mean something different. What is then this other formal method? 

Hilbert's reference to \emph{symbolic algebra} provides an important hint. Unlike the notion of logical form that can be understood after Trendelenburg \cite{Trendelenburg:1862} quite independently of any symbolic representation, the notion of \emph{algebraic form} was intimately connected to mathematical symbolism throughout the Modern history of algebra. Descartes' \emph{Geometry} of 1637 \cite{Descartes:1886} which first established algebra as a field of theoretical research at the same time produced what Serfati recently called a ``symbolic revolution'' in mathematics \cite{Serfati:2005}. The key idea of Descartes algebra is the following: the same syntactic operations with symbols may represent geometrical operations with straight segments and arithmetical operations with numbers. In this sense algebraic operations are general \emph{forms} of operations shared by certain arithmetical and geometrical operations \footnote{Beware that this interpretation of algebra is anachronistic. Arnauld in his \emph{New Elements of Geometry} \cite{Arnauld:1683} suggests that algebraic operations are operations with general mathematical magnitudes; his notion of magnitude generalizes upon geometrical magnitudes and arithmetical magnitudes aka numbers. The modern ``abstract'' approach in algebra, which is behind the anachronistic reading of Descartes, has been first systematically developed by van der Waerden in the early 1930-ies \cite{Waerden:1930-31}.}. The role of symbolism is crucial here. Although an abstract notion of operation with abstract non-specified operanda  may make logical sense it can hardly make  \emph{mathematical} sense. The algebraic symbolism invented by Descartes allowed for thinking of abstract operations (applicable both in geometry and arithmetics) in terms of concrete syntactic operations. This idea was used by Boole and other pioneers of symbolic logic. We shall shortly see how Hilbert uses this idea in his new approach to foundations of mathematics. 

So at least one thing that Hilbert most certainly had in mind talking in the above passage about the ``formal method of mathematics'' and suggesting an application of this method in logic is the method of (symbolic) algebra. However in the above passage he describes algebra only as a special case. This is why we cannot  derive the wanted sense of being formal from the notion of algebraic form. A more general notion of form, which turns to be appropriate in this case, is Cassirer's notion of \emph{symbolic form} \cite{Cassirer:1923-25}. I shall not develop it here in its full generality but focus only on its mathematical version relevant to Hilbert's work. 

The passage quoted above continues as follows:

\begin{quote}
The purpose of the symbolic language in mathematical logic
is to achieve in logic what it has achieved in mathematics, namely,
an exact scientific treatment of its subject-matter. The logical
relations which hold with regard to judgments, concepts, etc.,
are represented by formulas whose interpretation is free from
the ambiguities so common in ordinary language. The transition
from statements to their logical consequences, as occurs in the
drawing of conclusions, is analyzed into its primitive elements,
and appears as a formal transformation of the initial formulas
in accordance with certain rules, similar to the rules of algebra;
logical thinking is reflected in a logical calculus. This calculus
makes possible a successful attack on problems whose nature
precludes their solution by purely contentual [\emph{inhaltlische}] logical thinking.
Among these, for instance, is the problem of characterizing those
statements which can be deduced from given premises. (\cite{Hilbert&Ackermann:1950}, p.1)
\end{quote}

The first sentence of this passage clearly shows that Hilbert considers here an application of mathematics to logic as a way to improve on logic with mathematics. Hilbert and Ackermann claim here that by using the symbolic methods mathematics achieves ``an exact scientific treatment of its subject-matter'';  using this evidence the authors suggest that these methods may equally allow for an exact scientific treatment of logic. This project should be certainly distinguished from the idea of improving on mathematics through the clarification of its logical structure purported by Hilbert in his \emph{Foundations} of 1899. Nevertheless Hilbert tends to describe both projects in similar terms, namely in terms of formalization and axiomatization. Notice the expression ``contentual logical thinking'' that appears in the above passage. \emph{Contentual} logical thinking in this context  is opposed to \emph{formal} logical thinking. Interestingly, Hilbert and Ackermann themselves do not use the expression ``formal logic'' in this context and talk instead of mathematical logic and symbolic logic. This can be perhaps explained by the fact that at that time the expression ``formal logic'' was still more commonly used in Trendelenburg's sense, which does not imply that formal logic should always involve symbolic methods. However when Hilbert and Ackermann oppose what they call ``contentual'' logic to symbolic logic the qualification of this symbolic logic as \emph{formal} anyway immediately suggests itself. This is a different sense of being formal, which is more familiar to us today. What Hilbert and Ackermann call ``contentual'' logic would be described today as ``informal'' logic (even if this informal logic qualifies as formal in the weaker Trendelenburg's sense). 

Thus Hilbert and Ackermann's formalization and axiomatization of logic amounts to providing the given system of logic with a new \emph{symbolic} form, not to the specification of the ``logical form of logic''. Even if the authors describe the axiomatization of logic as ``the crowning achievement of the work of axiomatization'' this ``crowing achievement'' is not a simple continuation of the axiomatization in the sense of Hilbert's {Foundations} of 1899. Blurring the distinction between axiomatization of mathematical theories, on the one hand, and the axiomatization of logic, on the other hand, also blurs the distinction between theories and logic, which is fundamental for the 1899 approach. Hilbert and Ackermann's distinction between the ``contentual logic'' and the formal (symbolic) logic opens the new possibility of non-standard interpretation of logical signs on equal footing with non-logical signs. In \textbf{3.4} we shall see how this new possibility is realized in Tarski's topological model of intuitionistic propositional calculus; in  \textbf{4.4} I show how a further exploration of this possibility leads to a significant change of the Axiomatic Method as presented in Hilbert's \emph{Foundations} of 1899. But the conceptual problem that leads to this further development can be seen already at this early stage of the history of Axiomatic Method. 

A fundamental idea behind the Axiomatic Method in the sense of 1899 is a reduction of mathematical reasoning to the form of ``purely logical axiom system'' (Hintikka). But since the distinction between logical systems and mathematical theories is blurred such a reduction becomes senseless because there is no longer any clear sense in which a given axiom system can be called purely logical rather than simply mathematical. One can however clearly distinguish in this new symbolic setting between formal (to wit symbolic) and contentual theories,  and like in 1899 think of formal theories as a foundation for the corresponding contentual theories. Given such a contentual theory $T$ one can design its formal counterpart $F$ and then find another contentual interpretation $T'$ of $F$. Although it is natural to describe this latter procedure as formalization and axiomatization (assuming that $F$ is built axiomatically) this latter kind of formalization and axiomatization is quite unlike the formalization and axiomatization in the sense of 1899. While the formalization in the sense of 1899 involves the notion of logical form the symbolic formalization just described involves, generally, only the notion of symbolic form. 

In Hilbert's thinking both kinds of formalization are merged together, so he hardly distinguishes between them clearly. He apparently assumes that a formal symbolic logical system unlike a formal symbolic mathematical theory does not allow (and does not call for) for multiple alternative contentual interpretations but instead simply clarifies and purifies common vague contentual logical notions expressed in the natural language. This additional assumption apparently allows for systems of formal symbolic logic with a fixed semantics of logical terms. But in fact this assumption produces a tacit shift in the meaning of being formal. If the given symbolic logical system pins down the precise sense of logical notions, which outside the symbolic setting don't have any clear meaning, then the logical symbols used in this logical system are used as proper names of corresponding logical concepts (like the symbol \& conventionally used for denoting the logical conjunction) rather than variables that may acquire different interpretations
\footnote{
I am now talking about variables in the general non-technical sense. Non-logical constants of a given formal theory (in the usual technical sense of the term) also count as variables in this general sense because such constants are differently interpreted in different models of this theory. 
}. 
In mathematics symbols are used in this way when, for example, number $\pi$ is denoted by symbol ``$\pi$'' and number 1 is denoted by symbol ``1''. So the only notion of symbolic form, which is relevant in this case, is very unspecific and applies to any natural or artificial language. A notion of symbolic form, which is more specific and more relevant to logic and mathematics comes into the play \emph{only} when a formal mathematical theory presented by symbolic means allows for alternative contentual interpretations: in this case one can say that the formal symbolic theory grasps the symbolic form shared by all such interpretations. In the case of the system of symbolic logic proposed by Hilbert and Ackermann (see below) this latter specific sense of being formal does not apply. Apparently Hilbert and Ackermann continue to think of their logic as being formal in the traditional Trendelenburg's sense, and do not pay attention to the fact that this traditional sense of being formal does not square with the formal vs. contentual distinction relevant to mathematical theories. 

In \textbf{8.3} we shall see that Tarski's analysis of logicality mentioned in the last Section allows for a reconstruction of Hilbert's thinking, which is more coherent that the above reconstruction based on Hilbert's own words. It is hard to say which of the two reconstructions is more historically correct and I leave this question for a further study. (I would like the Tarski-based structuralist interpretation to be correct but I lack sufficient textual evidences.) It seems me likely that Hilbert was indeed driven by structuralist geometrical ideas described in \textbf{8.3} but since he had to apply a more traditional conceptual apparatus for talking about logic, this led him to incoherences that we notice in the above quote. Leaving aside this historical issues  I would like to stress once again that the idea of interpreting formal logical systems on equal footing with formal mathematical theories suggested by Hilbert's idea of ``contentual logical thinking'' is behind the recent transformations of Hilbert's Axiomatic Method, which I thoroughly discuss later in this book.    
 
Let us now see what kind of axiomatic system of logic Hilbert and Ackermann offer in their book. In fact the authors offer several such systems; for our analysis it is sufficient to consider the simplest one that the authors call \emph{sentential calculus} (Aussagenkalk\"ul) and that is usually called today \emph{propositional calculus}. The given axiomatic presentation of this system of symbolic logic consists of four elements (this analysis into elements is due to myself but not to the authors):

1) Specification of \emph{symbols}, types of symbols, and of the intended interpretations of these symbols. Hilbert and Ackermann distinguish symbols of three types:  (i) capital letters standing for propositional variables, (ii) logical symbols,  (iii) auxiliary symbols like parentheses and commas, which serve for separating and grouping other symbols.

2) Specification of \emph{syntactic rules} according to which well-formed  \emph{formulas} (strings of symbols) are constructed from the specified symbols. These formulas stand for propositions obtained from some given propositions (represented by singular symbols) with a help of logical connectives (also represented by corresponding symbols). So linking propositions by logical connectives is represented here by concatenation of the corresponding symbols into a single string and inserting appropriate auxiliary symbols in a way similar to which this is done with the usual linear alphabetic writing\footnote{Calling the linear alphabetic writing ``usual'' I mean, of course, that it is usual in the geographic area where I live.}  

3)  \emph{Axioms} of this logic, which are distinguished formulas interpreted as logical tautologies, i.e., propositions, which are true for all possible values of variables. 

4) Specification of \emph{logical rules}, which allow one to construct certain formulas (representing propositions) from some other formulas (representing other propositions). This construction is interpreted as logical deduction. The set of formulas deducible from the axioms according to logical rules 3 is interpreted as the set of (all) logical tautologies.

Let us now see whether or not this ``axiomatization'' of propositional logic is similar to the axiomatization of Euclidean geometry in the \emph{Foundations} of 1899. Is this the same notion of axiomatization that is at work in both cases? A common feature is this: in both cases we have a generic set of true propositions called axioms and certain rules that allow one to deduce some other true propositions from these axioms. Here however the analogy between the axiomatic geometry and the axiomatic logic stops. The axiomatic theory of geometry has a logical part and a properly geometrical part. The properly geometrical part is the geometrical axioms and geometrical theorems deduced from these axioms. All the rest (including the rules of deduction) is logic. The considered axiomatic theory of logic consists of logical axioms (generic tautologies), all other tautologies (deducible from the axioms) \emph{and} rules of two sorts plus the specification of symbols. If logical truths (tautologies) are thought of as universal truths about everything and geometrical truths are thought of as truths about some specific kind of things called geometrical objects then the comparison between the axiomatic geometry and the axiomatic logic may make sense. 

It is appropriate to notice here that the specification of certain tautologies as ``logical axioms'' used by Hilbert and Ackerman's in their axiomatization of logic is not necessary for a formal specification of a formal logical system. As has been shown by Gentzen \cite{Gentzen:1934}\cite{Gentzen:1935} a reasonable formal  logical system may have no axiom at the expense of having a bigger number of rules. But there is no logical system that has axioms and has no rules. In that sense rules is a more essential element of a logical system than logical axioms. This fact emphasizes my point that the ``axiomatization of logic'' should not be taken on equal footing with the axiomatization of geometry or axiomatization of any other special theory. What really matters in Hilbert's ``axiomatization of logic'' is developing logic by symbolic means. As I have just mentioned this can be done without using axioms.   

\section{\emph{Foundations} of 1927: Intuition Strikes Back}

As we have seen the formalization of logic (that Hilbert calls it ``axiomatization'') is not an innocent procedure which merely makes logical notions clearer and more explicit; whether it has this effect or not the formalization allows for studying and developing logic by mathematical means. Not surprisingly, the replacement of the traditional non-mathematical ``informal'' logic  by the mathematical symbolic logic has a very significant impact upon Hilbert's ideas about Axiomatic Method and foundations of mathematics. Let us now see how his new foundational project looks like. In the beginning of his paper \emph{Foundations of Mathematics} \cite{Hilbert:1927} that has been delivered in July 1927 at the Hamburg Mathematical Seminar, Hilbert describes this project in the following words:

\begin{quote}
With this new way of providing a foundation for mathematics, which we may appropriately call a proof theory, I pursue a significant goal, for I should like to eliminate once and for all the questions regarding the foundations of mathematics in the form in which they are now posed, by turning every mathematical proposition
into a formula that can be concretely exhibited and strictly derived, thus recasting
mathematical definitions and inferences in such a way that they are unshakable and
yet provide an adequate picture of the whole science. I believe that I can attain this
goal completely with my proof theory, even if a great deal of work must still be done
before it is fully developed.\\
No more than any other science can mathematics be founded by logic alone;
rather, as a condition for the use of logical inferences and the performance of logical
operations, something must already be given to us in our faculty of representation, certain extralogical concrete objects that are intuitively present as immediate experience prior to all thought. If logical inference is to be reliable, it must be possible to survey these objects completely in all their parts, and the fact that they occur, that they differ from one another, and that they follow
each other, or are concatenated, is immediately given intuitively, together with the
objects, as something that neither can be reduced to anything else nor requires
reduction. This is the basic philosophical position that I regard as requisite for mathematics and, in general, for all scientific thinking, understanding, and communication.
And in mathematics, in particular, what we consider is the concrete signs themselves,
whose shape, according to the conception we have adopted, is immediately clear and
recognizable. This is the very least that must be presupposed; no scientific thinker
can dispense with it, and therefore everyone must maintain it, consciously or not. (\cite{Hilbert:1927}, p. 464-465)
\end{quote}

When Hilbert says that mathematics cannot be ``founded by logic alone'' a modern reader acquainted with Hilbert's Axiomatic  Method readily agrees: of course, for doing mathematics one needs in addition to principles of logic some specific mathematical axioms like axioms of set theory! As we shall shortly see Hilbert indeed uses such specific axioms in his  \emph{Foundations} of 1927. But in the above passage he refers to something completely different! He states here that no logical inference is possible without ``certain extralogical concrete objects that are intuitively present as immediate experience prior to all thought'' and then specifies that as far as mathematics is concerned those ``extralogical concrete objects'' are  ``the concrete signs themselves, whose shape, according to the conception we have adopted, is immediately clear and recognizable''. Since mathematical symbolic logic does use concrete signs (symbols) Hilbert's ``logic alone'' cannot be mathematical;  in the given context the \emph{mathematical} logic should be rather understood as the pure logic provided with certain  ``extralogical'' (to wit symbolic) means. According to the new Hilbert's view the immediate intuitive giveness of the ``concrete signs'', which allows one to acknowledge ``the fact that they occur, that they differ from one another, and that they follow each other, or are concatenated'' is an indispensable ingredient of foundations of mathematics. For further references I shall call this specific sort of mathematical intuition, which allows one to manipulate and calculate with mathematical symbols, the \emph{symbolic} intuition. 

Let us compare Hilbert's view on foundations expressed in the above passage with his earlier views expressed in his comments on his  \emph{Foundations} of 1899. In 1899 he founds geometry on the ``pure'' (non-mathematical) logic and some axioms formulated in terms of this logic. The only notion of intuition, which is still indispensable in this setting, is the intuition in the sense of Kant's  \emph{Doctrine of Method}; this  \emph{minimal} intuition can be accounted for by the logical notion of \emph{existential instantiation} and be considered itself as a part of logic (see \textbf{1.3} above). In 1927 Hilbert no longer relies on the ``pure'' informal logic but stresses the foundational impact of symbolic intuition. Hilbert explicitly describes here symbols as ``extralogical''; the following explanation does not allow one to reduce the notion of intuition related to these extralogical objects to the minimal ``logical'' intuition. Indeed, while the mere fact that these objects ``occur'' and ``differ from one another'' does not yet make them extralogical, the fact that ``they follow each other, or are concatenated'' certainly does! Thus Hilbert's new foundational proposal of 1927 unlike that of 1899 essentially involves a non-logical notion of symbolic intuition. 

This does not mean however that by 1927 Hilbert abandon his earlier idea according to which all mathematical theories require a logical background. He rather upgrades this idea as follows: a system of logic, which is appropriate for founding mathematics, is not a system of ``pure'' (non-mathematical) logic but a system of symbolic mathematical logic (which includes an extra-logical symbolic aspect). Here is how Hilbert describes this upgrade himself:

\begin{quote}
[I]n my theory contentual inference is replaced by manipulation of signs [ausseres
Handeln] according to rules; in this way the axiomatic method attains that reliability
and perfection that it can and must reach if it is to become the basic instrument of all
research. (\cite{Hilbert:1927}, p. 467)
\end{quote}

The replacement of the ``contentual inference'' by the manipulation of signs involves two ways of formalization, which work here together but nevertheless can and should be carefully distinguished. The formalization in the sense of 1899 remains here at work, so the manipulation of signs presents here the \emph{logical} form of the given contentual inference. Simultaneously the manipulation with signs presents the \emph{symbolic} form of the same contentual inference. 

Since Hilbert makes it clear that his new method amounts to ``extending the formal point of view of algebra to all of mathematic'' (\emph{ib.} p. 470) we may safely identify the symbolic form with the algebraic form in the given context. However Hilbert's proposal does not reduce to the algebraization of logic the same  sense, in which one can speak of the algebraization of logic in earlier works by Boole, Morgan and others. For Hilbert uses a feature of algebra that plays no special role in earlier works in mathematical logic: I mean the algebraic method of ``ideal elements'' like $-1$ or $\sqrt{-1}$. Leaving a more detailed discussion on this algebraic method until \textbf{2.6} and \textbf{7.4} let me first describe Hilbert's proposal. After the introductory remarks quoted above he first introduces a system of symbolic logic similar to one presented in \cite{Hilbert&Ackermann:1928} and, second, adds two further groups of axioms, which he describes as ``specifically mathematical'', namely ``axioms of equality'' and ``axioms of number''. Then Hilbert shows how this apparatus allows one to do the finitary arithmetic. One may wonder if doing the finitary arithmetics with this heavy logical machinery indeed provides any epistemic advantage over doing it in the traditional way. Hilbert's answer is: No, it does not! As far as the finitary arithmetic is concerned this machinery allows one at best to ``impart information'' \footnote{
\begin{quote}
If we now begin to construct mathematics, we shall first set our sights upon
elementary number theory; we recognize that we can obtain and prove its truths
through contentual intuitive considerations. The formulas that we encounter when
we take this approach are used only to impart information. Letters stand for numerals,
and an equation informs us of the fact that two signs stand for the same thing. (\cite{Hilbert:1927}, p. 469)
\end{quote}
}. If I understand here Hilbert correctly his thinking is this: since usual arithmetical manipulations with natural numbers represented by strings of strokes or by the standard Arabic numerals  are just as intuitively clear as the manipulation of symbols and formulas in the Hilbert's symbolic system, from the foundational viewpoint the difference between the two formalisms is after all not essential (notwithstanding the fact that the former formalism has an advantage of being simpler and more convenient, while the latter formalism has an advantage of making explicit the logical structure of reasoning). However the new proposed formalism is advantageous as soon as one goes beyond the finitary arithmetic. Hilbert suggests thinking about such an extension after the pattern of algebraic extension:

\begin{quote}
Just as, for example, the negative numbers are indispensable in elementary number theory and just as modern
number theory and algebra become possible only through the Kummer-Dedekind ideals, so scientific mathematics becomes possible only through the introduction of ideal propositions. (\cite{Hilbert:1927}, p.471)
\end{quote}

An ``ideal proposition'' is any proposition that is not provable from Hilbert's logical and arithmetical axioms, i.e., any proposition, which is not a proposition of the finitary arithmetic. So any additional axiom and any formal proposition obtained as a formal consequence of the extended axiom system (which includes the same logical and arithmetical axioms plus the new axiom) qualifies as ideal. The only requirement that limits such possible extensions is the requirement according to which the extended system of axioms must be \emph{consistent}. As soon as the consistency is granted one may safely think of ``ideal'' objects and ``ideal'' relations involved into the given ideal proposition as \emph{existent} along the same pattern of thinking, which we have already explained talking about the \emph{Foundations} of 1899 (remind from \textbf{2.1} of thought-things and thought-relations). And in fact one can do more. Since these ideal objects and relations are represented by symbols and strings of symbols, which (unlike the bare thought-things and thought-relations) are bone fide mathematical objects on their own, any further interpretation of these ideal things is an option but not a necessary requirement. In the new symbolic setting these ideal things are concretely represented to begin with, and one may work with them just like in algebra people work with $\sqrt{-1}$. Crucially, working with ideal objects and relations involves the same type of syntactic manipulations as calculating with natural numbers. So even if the Hilbert's \emph{Foundations} of 1927 is an overkill in the case of the finitary arithmetic it's expected advantage is that it allows for an uniform treatment of the whole of mathematics by means similar to those used in the finitely arithmetic. 

The possibility of checking consistency is evidently crucial for Hilbert's project. Although in 1927 Hilbert offers no general solution of this problem he suggests that this problem is relatively easy and ``fundamentally lies within the province of intuition just as much as does in contentual number theory the task, say, of proving the irrationality of $\sqrt{2}$'' (\cite{Hilbert:1927}, p. 471). In the formal symbolic setting where a proof is represented by a string of symbols and formulas are constructed according to precise syntactic rules the proof of consistency of a given set of axioms amounts to a proof showing that there is no string of formulas that ends up with a formula expressing contradiction like $0 \neq 0$ (a simple argument shows that if $0 \neq 0$ cannot be formally proved no other contradiction can be proved either). Hilbert realizes, of course, that such a consistency proof will not itself qualify as \emph{formal} but will belong to his \emph{proof theory}, which in a different place \cite{Hilbert&Bernays:1934-1939} Hilbert calls by the name of \emph{metamathematics}. However since the whole of metamathematics ``fundamentally lies within the province of intuition just as much as does in contentual number theory'' this remark does not lead to the infinite regress in foundations. Thus the intuitive proof theory aka metamathematics (in Hilbert's original sense of this term) in Hilbert's view of 1927 becomes a foundation for the rest of mathematics. Beware that so far this is a declaration of intent, not yet an accomplished project. A more advanced version of the same project is presented in two volumes \cite{Hilbert&Bernays:1934-1939} published by Hilbert with a cooperation with Bernays in 1934-1939. 

As we all well know today in 1927 Hilbert severely underestimated potential difficulties of his proof theory; G\"odel's famous incompleteness theorems and all the following work in the area convinced many people that Hilbert's Program failed \cite{Zach:2003}. However even if Hilbert's foundational project as described in the \emph{Foundations} of 1927 indeed failed, his Axiomatic Method making part of this program certainly survived and until today remains standard. 

\section{Symbolic Logic and Diagrammatic Logic}

What kind of things qualify as symbols? Should they always be the convenient letters written from the left to the right or other systems of writing can be used to the same effect? Is the conceptual structure of the symbolic logic neutral with respect to its symbolic presentation or it is at some degree determined by this presentation? Can one open new logical possibilities by modifying this standard syntax? As every user of \LaTeX  perfectly knows various sorts of complicated non-standard multi-level symbols and diagrams can be encoded into the convenient strings of letters. Tolstoy's printed novels communicate us all sort of things using the same type of linear coding. But this useful feature of the standard symbolic syntax hardly justifies by itself Hilbert's far-reaching idea of the ultimate epistemic reduction of mathematical intuition to the symbolic intuition associated with the symbolic syntax of this particular sort (moreover if one takes into consideration the failure of Hilbert's program in its original form). Instead of trying to answer the above questions directly I propose now to look at some other syntactic possibilities known from the history of symbolic logic. Namely, I shall talk about logical diagrams after Venn's \cite{Venn:1881} of 1881 where this subject is treated systematically.  Then I come back to Hilbert and make some further comments about his idea of formal symbolic mathematics.  

Here is how Venn in 1881 \cite{Venn:1881}  describes the complementary roles of symbols and diagrams in logic talking about ``subdivisions of classes'':   

\begin{quote}
For one thing, we can of course always represent the products of such a subdivision in the language of common [non-symbolic] Logic, or even in that of common life, if we choose to do so. They do not readily offer themselves for this purpose, but when pressed will consent, though failing sadly in the desired symmetry and compactness. The relative cumbousness of such a mode of expression is obviously the real measure of our need for a reformed or symbolic language. [..] The reader will see at once how conveniently and briefly we can thus indicate any desired combination of class terms, and, by consequence, any desired proposition.  [..] \\ That such a scheme is complete there can be no doubt. But unfortunately, owing to this very completeness, it is apt to prove terribly lengthy.  [..] 
This then is the state of thing which a reformed scheme of diagrammatic notation has to meet. It must correspond in all essential respects to that regular system of class subdivision which has just been referred to under its verbal and its literal or symbolic aspect. Theoretically, as we shall see, this is perfectly attainable. (\cite{Venn:1881},  p. 102-103)
\end{quote}

As we can see Venn discusses here the use of symbols and diagrams in logic from a practical rather than theoretical viewpoint. Unlike Hilbert Venn does not try to distinguish his symbolic logic from the ``pure'' logic;  instead he distinguishes between the symbolic logic and  the ``common'' logic meaning the ``informal'' logic that applies only the natural language without helping itself with special symbols. Then Venn stresses that diagrams become helpful when the corresponding symbolic expressions turn to be ``terribly lengthy''.  Does the length of symbolic expression play any role from a theoretical or foundational viewpoint? Hilbert could argue that the length of symbolic expressions is irrelevant to the \emph{theoretical}  logic  and to the foundations of mathematics (as far as this length remains finite) even if it matters practically. I am not convinced by this argument. Hilbert's project aims at reduction of the ``ideal'' objects and propositions to their ``real'' counterparts, i.e., to their corresponding symbolic expressions. The idea that a long symbolic calculation is just as reliable as a short one is obviously a simplifying idealization. Even if this idealization is acceptable for certain purposes there is no reason to disregard a more realistic picture in a theoretical study. If Venn is right that diagrams help to tackle with the complexity of symbolic logical operations the diagrammatic notation must be taken as seriously as the symbolic notation.

It may be argued that symbols unlike diagrams do not involve the idea of \emph{resemblance} to what these things stand for, and that this fact makes symbols appropriate and diagrams inappropriate in foundations of mathematics. In old good times, so the argument goes, when mathematics dealt with circles, triangles and the like people could picture these things with diagrams and use these diagrams in their proofs. However since the modern mathematics involves highly abstract concepts, which do not allow for such a straightforward intuitive representation, diagrams become irrelevant. If such abstract mathematical concepts allow for any intuitive representation at all such representations are purely \emph{symbolic} and do not involve any relation of resemblance between symbolic constructions their corresponding mathematical objects.  

This argument is based on several misunderstandings. The idea that traditional geometrical diagrams in some sense \emph{resemble} certain ideal objects is a Platonic interpretation of the traditional geometrical practice that can and must be distinguished from this practice itself. One does not need this interpretation for doing traditional geometry with traditional geometrical diagrams. Notice that Venn's logical diagrams just like Venn's and Hilbert's logical symbols are not supposed to resemble anything at all. 

The idea that symbolic constructions do \emph{not} resemble anything should be also taken critically. In my view the relation of resemblance between mathematical objects and their material representations is ill-construed to begin with. We can rather establish such a relations between different material representations of the same object.  For example an Euclidean circle can be represented both with the traditional diagram (found in Section \textbf{6.2} in this book) and also symbolically with a single letter. Then one may observe that each of these two representations is unlike the other and further observe that a letter used for denoting the circle is unlike any other letter that can be used for the same purpose. Here we can see the difference. The capacity to read diagrams and the capacity read letters equally require the capacity to distinguish between graphical shapes and to identify tokens of the same shape. But in the case of diagrams this graphical typing reflects the typing of corresponding mathematical objects while in the case of letters it does not. In that particular respect diagrams are more informative than letter symbols (although letters, of course, have other epistemic advantages). However I cannot see any general philosophical justification of the idea that using graphical types in one way rather than in another way is more appropriate in mathematics and logic. Which one is more appropriate in a given context is rather a technical question. 

Let us read again more attentively what Hilbert says in the passage quoted in the last Section: 

\begin{quote}
[A]s a condition for the use of logical inferences and the performance of logical
operations, something must already be given to us in our faculty of representation, certain extralogical concrete objects that are intuitively present as immediate experience prior to all thought. If logical inference is to be reliable, it must be possible to survey these objects completely in all their parts, and the fact that they occur, that they differ from one another, and that they follow each other, or are concatenated, is immediately given intuitively, together with the objects, as something that neither can be reduced to anything else nor requires
reduction. (\cite{Hilbert:1927}, p. 464-465)
\end{quote}

Hilbert describes the above claim as his  ``basic philosophical position that I regard as requisite for mathematics and, in general, for all scientific thinking, understanding, and communication''. Yet one may remark that in this claim Hilbert mentions objects that ``follow each other'' and ``are concatenated'' which is hardly appropriate in a claim aiming at the full philosophical generality because the notions of order (``following each other'') and concatenation are applicable to objects of some sorts but may be not applicable to objects of some other sorts. If Hilbert would like to claim here that the order and the concatenation of the ``extralogical concrete objects'' are indeed \emph{necessary} in ``all scientific thinking, understanding, and communication'' he would need to provide an additional argument justifying this latter claim, which is actually missing. The following lines make it clear that Hilbert has here in mind nothing else but the familiar symbolic algebraic notation:

\begin{quote}
And in mathematics, in particular, what we consider is the concrete signs themselves,
whose shape, according to the conception we have adopted, is immediately clear and
recognizable. This is the very least that must be presupposed; no scientific thinker
can dispense with it, and therefore everyone must maintain it, consciously or not.  (\cite{Hilbert:1927}, p. 465)
\end{quote}
   
Why signs rather than diagrams? Why the concatenations of signs rather than geometric constructions of other sorts? Well, Hilbert's project aims at showing that the manipulation with signs (symbols) is sufficient for all mathematical purposes. Even if this foundational project works out it does not close the possibility of a different foundation with uses a different intuitive background and a different basic geometry. Notice also that Hilbert's foundational project involves his understanding of purposes of doing mathematics. Without trying to reconstruct the exact Hilbert's view on this issue I shall show in \textbf{3.3} that Hilbert's formal mathematics is not quite appropriate for applications in natural sciences and try to explain why. As soon as the application of mathematics in science counts as a purpose of doing mathematics this argument shows that the manipulation with signs in the way suggested by Hilbert is not sufficient for all mathematical purposes.      
  
In order to see a possibility of Axiomatic Method, which involves other mathematical intuitions than symbolic, it is instructive to look back at Euclid's \emph{Elements} from a particular point of view suggested by Friedman \cite{Friedman:1992}. Friedman suggest this point of view as his reconstruction of Kant's view but this does not matter in the given context: 

\begin{quote}
Euclidean geometry [..] is not to be compared with Hilbert's axiomatization [of Euclidean geometry in his \cite{Hilbert:1899}], say, but rather with Frege's \emph{Begriffsschrift}. It is not a substantive doctrine, but a form of rational representation: a form of rational argument and inference. Accordingly, its propositions are established, not by quasi-perceptual acquaintance with some particular subject matter, but, as far as possible, by the most rigorous methods of proof - by the proof-procedures of Euclid, Book I, for example. There remains a serious question about Euclid's axioms, of course; when pressed, Kant would most likely claim that they represent the most general conditions under which alone a concept of extended magnitude - and therefore a rigorous conception of an external world - is possible (see A163/B204). And, of course, we now know that Kant is
fundamentally mistaken here. (\cite{Friedman:1992}, p. 94-95)
\end{quote}

Euclidean geometry can be equally compared in this sense with Hilbert's symbolic logic: ``when pressed'' Hilbert like Kant could say that unless the symbolic syntax of his logic is taken for granted one, strictly speaking, cannot reason mathematically. Or perhaps Hilbert's mathematical genius and the common sense would win in this case over his philosophical ideas and he would take a more flexible attitude claiming only a particular way of building foundations for mathematics and asking his opponent to do this better. 

Although Euclidean geometry cannot any longer serve us as an organon of scientific reasoning I can see no reason why a theory capable to play this role today should be anything like a symbolic logical calculus. Moreover, I believe that if such a scientific organon is possible at all it must reflect objective features of the physical world learned through experience (as Euclidean geometry does it at a limited degree) rather than be based on speculative ideas about the correct thinking underpinned by some metaphysical theories. In Chapter \textbf{4} I present a tentative organon of this sort, which is Lawvere's categorical logic. (I shall tell more about Venn's diagrammatic logic during this discussion.) In the end of Chapter \textbf{6} I consider Voevodsky's Univalent Foundations as another candidate organon of the same sort. Now let us return back to Hilbert and see how his foundational project of 1927 develops in his later work.  

\section{\emph{Foundations} of 1934-1939: Doing is Showing?}

\emph{Foundations of Mathematics} \cite{Hilbert&Bernays:1934-1939} published by Hilbert and Bernays in two volumes in 1934 (first volume) and 1939 (second volume) is a systematic technical development of the project outlined in Hilbert's paper of 1927 \cite{Hilbert:1927}. Leaving this technical development wholly aside I shall now focus only on the authors' discussion about ``formal axiomatics'' in the beginning of the first volume. Here is the passage of my interest:

\begin{quote}
The term axiomatic will be used partly in a broader and partly in a narrower sense.We will call the development of a theory axiomatic in the broadest sense if the basic notions and presuppositions are stated first, and then the further content of the theory is logically derived with the help of definitions and proofs. In this sense, Euclid provided an axiomatic grounding for geometry, Newton for mechanics, and Clausius for thermodynamics.\\
In Hilbert's \emph{Foundations of Geometry} [of 1899] the axiomatic standpoint received a
sharpening regarding the axiomatic development of a theory: From the factual and
conceptual subject matter that gives rise to the basic notions of the theory, we retain only
the essence that is formulated in the axioms, and ignore all other content. Finally, for
axiomatics in the narrowest sense, the \emph{existential form} comes in as an additional factor.
This marks the difference between the \emph{axiomatic method} and the \emph{constructive} or \emph{genetic}
method of grounding a theory. While the constructive method introduces the objects of a theory only as a \emph{genus} of things, an axiomatic theory refers to a fixed system of things (or several such systems),
and for all predicates of the propositions of the theory, this fixed system of things constitutes
a \emph{delimited domain of subjects}, about which hold propositions of the given theory.\\
There is the assumption that the domain of individuals is given as a whole. Except for
the trivial cases where the theory deals only with a finite and fixed set of things, this is an
idealizing assumption that properly augments the assumptions formulated in the axioms.\\
We will call this sharpened form of axiomatics (where the subject matter is ignored
and the existential form comes in) \emph{formal axiomatics} for short. (quoted after bilingual edition \cite{Hilbert&Bernays:2010}, p.1a-2a) 
\end{quote}

We have already discussed the distinction between formal and contentual axioms and I shall not return to it now but comment on the authors' distinction between the constructive (genetic) and the axiomatic methods of theory-building. My first comment concerns Euclid. As we can see in the above passage Hilbert and Bernays qualify Euclid's method as axiomatic in a ``broader'' sense. This ``broader'' sense of being axiomatic includes what Hilbert and Bernays call constructive or genetic method. This is clear from another description of Euclid's method that Hilbert and Bernays provide later in their book: 

\begin{quote}
Euclid's axiomatics was intended to be contentual and intuitive, and the intuitive meaning of the figures is not ignored in it. Furthermore, its axioms are not in existential form either: Euclid does not presuppose that points or lines constitute any fixed domain of individuals. Therefore, he does not state any existence axioms either, but only construction postulates. (\emph{ib.}, p. 20a)
\end{quote}

The above quote clearly shows that Hilbert is well aware about the difference between his and Euclid's approaches to theory-building, which I have emphasized in \textbf{1.2}
\footnote{Nevertheless Hilbert and Bernays in a different place describe Euclid's (as well as Newton's and Clausius') method in this way: ``the basic notions and presuppositions are stated first, and then the further content of the theory is \emph{logically} derived with the help of definitions and proofs'' (\emph{ib.}, p. 1, my emphasis). As we have shown in Chapter \textbf{1} this description of Euclid's ``genetic'' method is incorrect. The fact that this method amounts to building certain non-logical objects like triangles from given basic objects according to certain rules is not a sufficient reason for calling it logical.}. 
The reference to \emph{Foundations} of 1899 in the work of 1934 suggests that in 1934 Hilbert retains his earlier notion of Axiomatic Method. But as a matter of fact in 1934 Hilbert and Bernays use the \emph{symbolic} version of this method. As we shall now see this fact has important consequences for the author's distinction between constructive (genetic) and axiomatic (in the narrow sense) ways of building mathematical theories. 

Remind that in the symbolic setting mathematical proofs are chains of formulas constructed according to certain fixed rules (with a help of some basic construction including formulas representing axioms). In this setting to \emph{prove} (a proposition expressed by formula $F$) is to \emph{construct} (a chain of formulas that ends up with $F$). Thus doing mathematics in Hilbert's formal axiomatic setting does not reduce contentual constructions in mathematics altogether but reduces all such constructions to constructions of a special sort, namely, to the finite symbolic constructions. Notice however that this reduction does not concern the \emph{metamathematical} proofs including the proofs of consistency (of a given set of formal axiom). In the metamathematics everything works like in the traditional contentual mathematics developed by ``genetic'' methods. One not only performs here certain constructions (chains of formulas) but also makes certain judgements about these constructions, including judgements of the form ``such-and-such construction is impossible'' (remind that in order to prove the consistency of a set of axioms it is sufficient to show that formula $0 \neq 0$ is not derivable from these axioms). I would like to stress that this new turn of the dialectics of ``doing'' and ``showing'' is relevant only to the symbolic version of Hilbert's formal Axiomatic Method but not to the ``informal'' version of this method presented in Hilbert's \emph{Foundations} 1899. Although all main features of \emph{Foundations} of 1899 are present in the new \emph{Foundations} of 1934, these new \emph{Foundations} have some new features which significantly change the whole picture.

In this context I would like to consider more attentively Hilbert's algebraic motivation explained in his  \emph{Foundations} of 1927. Remind that Hilbert suggests thinking of his \emph{ideal propositions} (represented by formulas) after the pattern of algebraic ``ideal objects'' like $\sqrt{-1}$. In \textbf{2.4} I have tried to explain the analogy (after Hilbert) but now I would like to stress a point where it fails: while in algebra manipulations with symbols represent manipulations with ideal objects in Hilbert's formal setting manipulations with symbols represent logical inferences and other logical operation, which allows one to \emph{say} different things about ideal objects but not manipulate with them. 

Let me demonstrate this point with two historical examples. Consider, first, the following interesting passage from Arnauld's  \emph{New Elements} of 1683 :

\begin{quote}
 What cannot be multiplied by its nature can be multiplied through a mental fiction where the truth presents itself as certainly as in a real multiplication. In order to learn the distance covered during 10 hours by one who covers 24 lieu per 8 hours I multiply through a mental fiction 10 hours by 24 lieu that gives me the imaginary product of 240 hour times lieu, which I divide then by 8 hours and get 30 lieu.  By the same mental fiction one multiplies a surface by another surface  even if the product has 4 dimensions and cannot exist in nature. One may discover many truths through multiplications of this sort. \\ 
 People say that this is because the imaginary products can be reduced to lines. [..]  But there is no evidence that [relevant] proofs  depend on those lines, which are in fact wholly alien to them.  (\cite{Arnauld:1683} p. 38-39, my translation from French) 
\end{quote}

 Traditionally (in particular, in the early Arab algebra) the product of two straight lines is construed as a rectangle having these given lines as its sides; the product of three lines is a solid but in order to form products of four and more linear factors in this geometrical way one needs higher dimensions, which according to Arnauld ``cannot exist in nature''. Nevertheless he is ready to consider such higher-dimensional geometrical constructions as useful fictions on equal footing with products of distances by time intervals, which don't have any immediate physical interpretation either but are demonstratively useful for calculations. In the last quoted paragraph Arnauld refers to Descartes' proposal to construct the multiplication of geometrical magnitudes differently, so the product of straight lines is again a straight line; in modern words Descartes' definition of multiplication of straight lines makes this operation algebraically closed. (Descartes uses an auxiliary line 1 as a unit and then considers similar triangles with sides $1, a, b, c$ such that $\frac{1}{a} = \frac{b}{c}$, which gives him the wanted definition of $c = ab$, see \cite{Descartes:1886}) Arnauld finds this trick artificial and unnecessary. What makes him confident about higher-dimensional geometrical products and quasi-physical units like hour times lieu is the symbolic algebraic calculus supporting these otherwise problematics notions. Using this symbolic calculus one forms the product of four factors $p = abcd$ as easily as any product of two or three factors. One cannot easily imagine the product of two surfaces (just like one cannot give a physical sense to hours times lieu) but one can easily concatenate the string $ab$ and the string $cd$ and think of this operation as multiplication of two surfaces. Descartes' alternative definition of the geometrical product aims at providing a ``clear and distinct''  intuitive underpinning of this operation that avoids the talk of higher dimensions. Arnauld finds Descartes' construction of product unnecessary because in his eyes the symbolic calculus provides such an intuitive underpinning by itself. The ``proofs'' that Arnauld mentions in this context are nothing but symbolic calculations.  

Now consider this passage from MacLaurin's \emph{Treatise of Fluxions} of 1742 where the author shows how the symbolic algebra allows for operating with Newton's \emph{fluxions} in a precise way in spite of the ``obstruse'' intuitive nature of those things.
        
\begin{quote}
The improvement that have been made by it [the doctrine of fluxions] [..] are in a great measure owing to a facility, conciseness, and great extend of the method of computation, or algebraic part. It is for the sake of these advantages that so many symbols are employed in algebra.  [..]  [Algebra] may have been employed to cover, under a complication of symbols,  obstruse doctrines, that could not bear the light so well in a plain geometrical form; but, without doubt, obscurity may be avoided in this art as well as in geometry, by defining clearly the import and use of the symbols, and proceeding with care afterwards. \\
(quoted by \cite{Cajori:1929},  v2, p. 330) 
\end{quote}

The above passages from Arnauld and MacLaurin well illustrate Hilbert's remarks in his  \emph{Foundations} of 1927 where he stresses an epistemic impact of the ``formal method of albera'' (see above). Let us however compare the above historical examples with Hilbert and Bernay's new symbolic axioms for plane Euclidean geometry found in the \emph{Foundations} of 1934. This new set of axioms involves a single type of primitive objects called  \emph{points} (instead of two types specified in the \emph{Foundations} of 1899) and two primitive ternary relations $Gr(x,y,z)$ and $Zw(x,y,z)$, which are informally interpreted as  \emph{points $x,y,z$ lie on the same straight line} ($Gr$ stands for German \emph{Gerade}) and  \emph{point $y$ lies between point $x$ and point $z$} ($Zw$ stands for German \emph{zwischen}). The First Axiom reads:

\begin{quote}
$(x)(y) Gr(x, x, y)$
\end{quote}  

that translates into the prose as follows: 
\begin{quote}
 \emph{all points $x,y$ lie on the same line} 
 \end{quote}  
 
 (prefix $(x)$ in Hilbert's notation reads ``for all $x$''). Using this axiom, other geometrical axioms expressed similarly, logical axioms and, finally, the appropriate deductive rules  (which in the given symbolic setting are syntactic rules allowing for building new formulas from some given formulas) one may formally derive certain \emph{theorems} (which are some new formulas generated from the axioms by the given rules). However this formal deduction does not allow for constructing any new  \emph{geometrical} object: remind that ``an axiomatic theory refers to a fixed system of things (or several such systems), and for all predicates of the propositions of the theory, this fixed system of things constitutes a \emph{delimited domain of subjects}, about which hold propositions of the given theory''\footnote{One may suggest that this method allows for introducing new objects defined as classes (or sets) of points satisfying holding certain relations definable in terms of the primitive relations; for example a straight line can be tentatively defined as a class of points such that any three points $x,y,z$ of this class hold relation $Gr(x, y, z)$. However in fact the formal deduction cannot be interpreted in this ay because in the given setting a \emph{class} of things does not count as a new thing. So one needs a background set theory for following this suggestion.} The only way of producing new \emph{ideal objects} within the formal axiomatic framework is to suggest a new system of axioms and assure its consistency. This is very unlike the way, in which the \emph{ideal objects} are produced in the traditional symbolic algebra. Arnauld first takes it for granted that the product of three straight lines is a parallelepiped; next he represents the straight lines by symbols $a,b,c$ and represents their product (i.e., the parallelepiped) by the string $abc$ (I modernize his notation slightly but this is not essential here). Then he considers string $abcd$, thinks how to interpret it geometrically and gets a vague idea of a 4-dimensional object that ``cannot exist in nature''. Since the algebraic rules of multiplication do not limit the number of factors 4-element strings turn to be just as well-manageable as 3-elements strings. Thus one may safely operate with 4-dimensional and $n$-dimensional parallelepipeds symbolically, no matter how ``ideal'' such things appear to the geometrical intuition. A symbolic construction, namely the concatenating of four primitive symbols (letters), represents here an ideal geometrical construction (the construction of 4-dimensional parallelepiped). 
 
Hilbert's symbolic constructions unlike Arnauld's and MacLaurin's symbolic constructions represent not ideal constructions themselves but propositions and systems of propositions ``about'' ideal objects. This makes Hilbert's ideal objects (i.e., all mathematical objects except syntactic objects) fundamentally non-constructive like Plato's ``ideal numbers'' (see \textbf{5.4} below). One may stipulate their existence by adopting appropriate axioms (keeping in mind the consistency requirement) but one cannot construct them from simpler elements.  In \textbf{3.3} I argue that this feature of Hilbert's formalism makes it inappropriate for building mathematical theories useful in sciences.

For further references I shall call the Axiomatic Method as it is described in Hilbert's \emph{Foundations} of 1927 \cite{Hilbert:1927} and implemented in the \emph{Foundations} of 1934-39  \cite{Hilbert&Bernays:1934-1939} by the name of \emph{Formal} Axiomatic Method and call the procedure of reconstruction of mathematics with this method the \emph{formalization} of mathematics. It is essential for what follows not to confuse this specific notion of formalization with any other notion that can be occasionally called by the same name.

A  modernized version of Formal Axiomatic Method, which includes basics of the modern model theory, is presented (together with a philosophical underpinning) in Tarski's textbook  \cite{Tarski:1941}. A more detailed historical study of the Axiomatic Method would require an analysis of this and many other similar works but this task is out of the scope of the present book. I focus my attention on Hilbert's work because for my purposes it is essential to analyze how ideas emerge and less essential to study how they solidify and become an orthodoxy.

\chapter{Formal Axiomatic Method and the 20th Century Mathematics}
The Formal Axiomatic Method has been proposed by Hilbert about a century ago and it is appropriate to ask how it performed during the past century. It appears to me that its impact is somewhat controversial. On the one hand, during this time period the Formal Axiomatic Method was and still remains \emph{the} standard method of theory-building in eyes of logicians and logically-minded mathematicians, physicists, biologists and philosophers. On the same side of the scale I put the progress in the logico-mathematical investigations (some of which use the title of \emph{foundations of mathematics}), which apply this method in some form. Finally, in some form this method is widely used in the current educational practice. (I specify the relevant forms of the Axiomatic Method below in this Chapter.) But on the other hand, one can also observe that by this date the Formal Axiomatic Method had no significant effect on either in the mainstream mathematics or in natural sciences, which remained largely ``informal''; none of recent significant advances in mathematics (like the recent proof of Poincar\'e conjecture) used formal logical methods. 

So the today's situation is somewhat schizophrenic. When Hilbert in his \emph{Foundations} of 1927 suggested his Formal Axiomatic Method as the ``basic instrument of all research'' (\cite{Hilbert:1927}, p. 467) he really meant it, and this clearly did not happen. At the same time most mathematicians and logicians (as well as many philosophers and a few physicists) agree that the Axiomatic Method is useful; some of these people also believe that this method is indispensable in their discipline. When they are asked what do they really mean by Axiomatic Method they likely refer to the modern axiomatic set theory or to another example of axiomatic theory built by Hilbert's method. If we now ask mathematicians and physicists how it is possible that the Axiomatic Method plays an important role in their discipline without having any clearly visible effect on it the most popular answer will be likely this: the Axiomatic Method matters only in the \emph{foundations} of science while the mainstream science cares very little about its own foundations. At that point philosophers and logicians join the discussion and argue like this. It is, of course, not normal that the foundations of your science are not properly taken care of. But this is understandable because you, guys, have a lot of other things to do. And perhaps you are not quite qualified for the job because the foundations is a somewhat philosophical and logical subject. So forget about foundations and leave this subject to us and give us in return some jobs at mathematical and physical departments. We shall take care about the foundations and you will do the rest. 

I am not satisfied by this division of labor (even if I have a philosophical interest to mathematics and don't mind a job at a mathematical department) because, in my view, it produces a wrong notion of foundation of a discipline, which allows a foundation to be wholly detached from the discipline itself. A notion of foundation that seems me satisfactory is described by Lawvere and Rosebrugh in the following words:

\begin{quote}
A foundation makes explicit the essential general features, ingredients, and operations of a science, as well as its origins and generals laws of development. The purpose of making these explicit is to provide a guide to the learning, use, and further development of the science. A ``pure" foundation that forgets this purpose and pursues a speculative ``foundations" for its own sake is clearly a nonfoundation.
(\cite{Lawvere&Rosebrugh:2003}, p.235) 
\end{quote}  

Having this Lawvere's notion of foundations of mathematics in mind I cannot reserve for Hilbert's Axiomatic Method a place in foundations of today's mathematics until I can see more clearly the role of this method in a broader mathematical context. In the following Section I consider the case of ``speculative foundations'', i.e., the 20th century research of formal axiomatic set theories, and in the next Section the case of Bourbaki's mathematics, which is not wholly ``speculative'' in that sense because it has an important continuing impact on mathematical education and reflects some important features (even if ignores some other important features) of the mainstream mathematics of the second half of the 20th century. We shall see that although in both cases the Formal Axiomatic Method plays a central role, in neither of these case  this method is used as intended by Hilbert. In the second half of this Chapter I shall do two other things: first, put forward an argument explaining why the Formal Axiomatic Method is not quite appropriate for natural sciences and, second, discuss an unusual application of this method by Tarski, which I consider as a step towards the rethinking of Hilbert's Axiomatic Method in Lawvere's work.

\section{Set Theory}

Remind the notion of \emph{system of things} required by the Formal Axiomatic Method. Can one provide a formal axiomatic theory of such \emph{systems}? Unless one assumes that a system $U$ of systems of things is an element of itself this project does not involve a circularity but provides a somewhat \emph{restricted} notion of system of things (as an element of $U$) that can serve for developing various formal axiomatic theories on the top of the formal theory of $U$. This anachronistic description of Zermelo's idea to axiomatize set theory explains why and how the later development of Hilbert's formal approach to foundations of mathematics involved not only the Formal Axiomatic Method itself but also the axiomatic theory of sets\footnote{Zermelo's principal motivation for axiomatizing set theory was saving Cantor's so-called ``naive'' set theory from paradoxes \cite{Peckhaus:2008}}.   
 
Since Zermelo's pioneering works in the axiomatic set theory the mainstream research in set theory focused on studies of various formal theories of sets of models of such theories. This makes set theory a rare and arguably the most important example of a modern mathematical theory developed wholly within a formal axiomatic setting. So in order to see how the Formal Axiomatic Method works in today's mathematics it is useful to consider the case of set theory quite independently from any foundational claims made about this theory. For being more concrete let us consider Cantor's  Continuum Hypothesis (CH), which in 1900 has been listed by Hilbert \cite{Hilbert:1902} as the number one among 23 open mathematical problems that Hilbert at that time judged to be the most important\footnote{The Continuum Hypothesis conjectured by Cantor states that there is no cardinal number strictly bigger than the minimal infinite cardinal number $\aleph_{0}$ (which can be described as the ``number of all natural numbers'') and strictly smaller that the cardinal number $2^{\aleph_{0}}$ of the set of all subsets of some set having the cardinal number $\omega$ (for example, the set of all series of natural numbers, including infinite series). Number $2^{\aleph_{0}}$ has been identified by Cantor with the number of points on a given continuous line or surface; hence the name of this conjecture}. 

CH is a very peculiar example of mathematical problem because today there is still no common opinion as to whether this problem is solved or still remains open! And this peculiar situation is obviously due to the fact that the modern set theory unlike (almost) the rest of mathematics is developed in a formal axiomatic setting. The story in brief is the following. In 1938  G\"odel \cite{Godel:1938} discovered that ZF (which is an improved version of Zermelo's axiomatic theory of sets so called after the names of Zermelo and Fraenkel \cite{Fraenkel&Bar-Hillel&Levy:1973}) is consistent with CH by building a model of ZF in which CH holds. In 1963 Cohen \cite{Cohen:1963}  discovered that ZF is also consistent with the negation of CH by building a model of ZF in which CH does not hold. So it is well established today that neither CH nor its negation can be derived from the axioms of ZF \cite{Kunen:1980}. 
What remains controversial is whether or not this independence result provides a definite answer to the original question by allowing one to claim that the original question is ill-posed. An additional axiom - or some wholly new system of axioms for set theory - may eventually help, of course, to settle the problem in the sense that CH or its negation can be deduced from the new system of axioms. There are obvious trivial ``solutions'' of this sort like considering CH itself as an axiom. Then, however, it remains to show that the system of axiom for set theory solving the CH problem is a ``right'' one, and so the proposed solution is ``genuine''. I cannot see how this can be done on purely mathematical grounds; any possible argument to the effect that one system of axioms for set theory is ``more natural'' than some other has a speculative nature and lacks any objective validity. Even if one gets some non-trivial proof of CH from some system of axioms that appear to be in some sense natural one can hardly claim that this system of axioms is the ``right one'' solely because it solves the CH problem and because such a proposed solution is smart and elegant. Although this situation is not unprecedented and may be compared, in particular, with the fate of the  \emph{Problem of Parallels} in geometry of the 19th century (see \textbf{7.3} below) it makes a sharp contrast with the mainstream mathematics that still manages to provide yes-no answers to many well-posed questions.    

It may be argued that the formal axiomatic framework makes explicit a relativistic nature of mathematics, which we should learn to live with; according to this viewpoint it is pointless to ask whether CH is true or false without further qualifications, and all that mathematicians can do is to study which axioms do imply CH (modulo some specified rules of inference), which axioms imply its negation, and which do neither (like the axioms of ZF). More generally, the only thing that mathematics can do according to this point of view is to provide true propositions of the  \emph{if - then} form:  \emph{if} such-and-such propositions are true  \emph{then} certain other propositions are also true. I cannot see how such a deductive relativism (or ``if-thenism'') about mathematics can be sustainable. It is incompatible not only with the common mathematical practice but also, more specifically, with the current practice of studying formal axiomatic systems. Denote $S$ the proposition saying  that CH is independent from the axioms of ZF (in the sense that neither $S$ nor its negation can be derived from these axioms). $S$ is commonly seen as an established theorem on a par with any other firmly established mathematical theorem. However $S$ is not expressed in the  \emph{if - then} form; it is expressed as an ``absolute'' mathematical truth about ZF and CH, which does not refer to any particular formal framework. The proof of $S$ (that comprises the construction of G\"odel's  model $L$ verifying CH and Cohen's forcing construction falsifying CH) is a piece of rather sophisticated ``usual'' or ``informal'' mathematics but not a formal inference within certain axiomatic theory. So a consistent if-thenist would not hold without further qualifications that CH is independent from the axioms of ZF  but rather say that it depends of one's assumptions.

Remind that Hilbert's foundational project of 1927 was supported by the hypothesis according to which all metamathematical questions concerning the consistency and the independency of axioms expressed in a formal language were trivial or, to put it more precisely, treatable by finitary means. And we know today this hypothesis is false. When people suggest today ZF as a foundation of mathematics they no longer hope to prove the consistency of this formal theory (since such a proof requires a stronger metatheory) but rather take a pragmatic attitude according to which this theory can be used unless one eventually discovers that it is contradictory; if this happens there are always ways of modifying the axioms of ZF that may block the eventual contradiction. As I have already said the known proofs of $S$ are far from being trivial or finitary.

Thus in spite of the fact that that the modern set theory no longer considers sets naively but works instead with various formal axiomatic theories of sets this modern theory like any other modern mathematical theory relies on non-formalized proofs. What is specific for the modern set theory is its  \emph{object} rather than its method. Instead of studying sets ``directly'' in the same way in which, say, group-theorists study groups, set-theorists study formal axiomatic theories of sets. However the \emph{methods} used by modern set-theorists are not essentially different from methods used in other parts of today's mathematics. It remains in my sense an open question whether or not such a roundabout way of studying sets has indeed proved effective. True, at the present there is no clear alternative to it. However it is not inconceivable that in the near future the mathematical community may bring about an improved ``naive'' concept of set that would allow one to study sets like groups. It is not inconceivable that such an old-fashioned way of thinking about sets could after all allow for a real progress in the CH problem. In any event it seems me important to keep such a possibility open and not try to take it out of the table using philosophical arguments.

The distinction between a \emph{theory} and \emph{metatheory} (and, more generally, between \emph{mathematics} and \emph{metamathematics}, which dates back to Hilbert, is helpful for making things clearer. In modern set theory a \emph{theory} is ZF or another formal axiomatic theory while proofs of independence of CH from the axioms of ZF and similar results belong to a \emph{meta-theory} that tell us important things about formal axiomatic theories. I would like, however, to stress here the fact that this terminology, which remains standard in the community of people working in \emph{foundations} of mathematics, is heavily philosophically-laden and overtly clashes with the language in which the wider mathematical community usually describes its own activities. Namely, the distinction between mathematics and metamathematics implies the view according to which formal axiomatic theories are ``usual'' mathematical theories while metatheories belong to a special domain of metamathematics that lays somehow ``beyond'' the usual mathematics and has some philosophical flavor. But if we leave  now philosophy aside and describe the same subject-matter in the language of mathematicians we come to a different view: formal axiomatic theories are \emph{not} mathematical theories in the usual sense of the word while their corresponding metatheories are comparable with ``usual'' theories from any other area of mathematics! This linguistic confusion reflects a gap between what today's mathematics \emph{is} and what in the opinion of certain people it \emph{must} be. In any event it seems me essential to fix it by suggesting a more neutral language. By the analogy with the distinction between an \emph{object language} and a \emph{metalanguage} in formal semantics I suggest to use the term ``object-theory'' for what in formal axiomatic studies (but not in the rest of mathematics) is usually called simply a ``theory''
\footnote{Remind from \textbf{2.4} that Hilbert's distinction between real and ideal mathematical objects translates into Hilbert's distinction between (contentual) mathematics and \emph{metamathematics} as follows: mathematics studies ideal objects with a help real syntactic constructions; metamathematics studies real syntactic constructions without using anything ideal. This, remind, was Hilbert's original idea supposed to help mathematics to ``get real'' without leaving the ideal ``Cantor's Paradise''. When in the light of G\"ode's incompleteness theorems and other developments it became clear that metamathematics cannot do solely with finitary means, some limited ``ideal content'' - i.e., some more advanced \emph{mathematical} content - was allowed in it. The title of ``The Mathematics of Metamathematics'' appeared in 1969 \cite{Rasiowa&Sikorsky:1963} perfectly illustrates this shift, which has relaxed the boundary between mathematics and metamathematics.}.

One may argue as follows. True, any formal object-theory requires some supporting informal metatheory. True, metatheories are, generally, sophisticated. But metatheories can be formalized in their turn and studied with the metametatheories. Although this regress is potentially infinite and is not going to lead us to any ultimate self-evident ground it nevertheless deepens our understanding of foundations of mathematics. Although informal instruments cannot be then wholly taken away like the Wittgenstein's ladder they can be viewed as a part of general philosophical underpinning of mathematics rather than as a part of mathematics proper\footnote{A similar view is developed by Shapiro in \cite{Shapiro:1991}}.   
 
Once again I claim that the above view is based on an a priori idea about mathematics and its foundations that is not justified neither by the old nor by the recent mathematical practice. First of all it confuses firm (meta)mathematical results like the independence of CH with a general philosophical discussion. It reflects an existing trend in the Analytic philosophy of mixing mathematical and speculative arguments indiscriminately, which is rejected by the majority of mathematicians who do not want to allow for philosophical speculation in mathematical papers. I support this rejection from the philosophers' side and insist that although philosophical speculation and mathematics may indeed fertilize each other they must be carefully distinguished and kept apart.  The fact that metamathematical results are obtained ``informally'' does not preclude them to be firm and mathematically valid; such results must be sharply distinguished from their philosophical interpretations, which do not have and cannot possibly have any objective validity. Second of all the above view once again takes it for granted that a formalized object-theory is a self-standing mathematical theory. This seems me very dubious. Metamathematical results concerning ZF and and its likes make the core content of the modern set theory rather than only a foundation of this theory, which can be left aside unless one has a special interest in studying foundations (while \emph{philosophical} interpretations can and should be left aside when one works in set theory as a mathematician). 

Thus we can see that even in set theory where formal methods are applied systematically the (metamathematical) informal methods remain essential. Let us now see how the Axiomatic Method applies in the modern mathematics outside set theory and mathematical logic.  

\section{Bourbaki}
The multi-volumed \emph{Elements of Mathematics} \cite{Bourbaki:1939-1988} produced in 1939-1998 by a group of (mainly French) mathematicians using the pseudoname \emph{Nicolas Bourbaki} is the most recent serious attempt to write a self-contained compendium of the core contemporary mathematics after the Euclid's example (interpreted liberally)\footnote{The original French title is \emph{El\'ements de math\'ematique}, which uses the unusual singular form `` math\'ematique'' (while the usual French word for mathematics is ``math\'ematiques''). So a more accurate English translation of the title is \emph{Elements of Mathematic}. This unusual singular form of the word is supposed to stressed Bourbaki's aim of the unification of mathematics.}. Although Hilbert's \emph{Foundations} of 1899 fall under the same description the two works differ in their purpose. Remind that Hilbert's work of 1899 is focused on the Euclidean geometry, which in the end of the 19th century was already only a relatively small part of what was commonly known under the name of geometry in the mathematical community. Hilbert rebuilt here an old theory with a new Axiomatic Method, clarified the logical structure of this old theory,  and left it to other people to do a similar job for more recent theories like the Riemanian geometry (see for example Veblen and Whitehead \cite{Veblen&Whitehead:1932}). So in his \emph{Foundations} of 1899 Hilbert presented a  \emph{method} allegedly applicable everywhere in mathematics and beyond, but unlike Euclid he did not try in this work to account for basic mathematical \emph{concepts} sufficient for developing the rest of his contemporary mathematics.  Bourbaki in his turn like Euclid aims at providing a genuine self-contained  introduction into the contemporary mathematics, which systematically presents not only its method but also its basic content. A concise general description of this project, which makes explicit some grounding ideas behind it, has been published in 1950 as a separate article \cite{Bourbaki:1950}

The section of this article named ``Logical Formalism and the Axiomatic Method'' begins as follows:

\begin{quote}
After more or less evident bankruptcy of the different systems [..] it looked, at the beginning of the present [20th] century as if the attempt had just about been abandoned to conceive of mathematics as a science characterized by a definitely specified purpose and method; instead there was a tendency to look upon mathematics as a ``collection of disciplines based on particular, exactly specified concepts'', interrelated by ``a thousand roads of communications'' [..] [quoted by the author from \cite{Brunschvicg:1912}, p.447] Today, we believe however that the internal evolution of mathematical science has, in spite of appearance, brought about a closer unity among its different parts, so as to create something like a central nucleus that is more coherent than it has ever been. The essential aspect of this evolution has been the systematic study of the relation existing between different mathematical theories, and which has led to what is generally known as the ``axiomatic method.'' (\cite{Bourbaki:1950}, p.222) 
\end{quote}

After this recognition of the unifying power of the Axiomatic Method Bourbaki makes an interesting move by distinguishing between the \emph{logical} aspect of the Axiomatic Method from another aspect, which can be called \emph{structural} (see Chapter \textbf{8}); in Bourbaki's view this is the latter rather than former aspect that makes the Axiomatic Method a powerful instrument of the unification; as we shall now see Bourbaki points here to his proper version of Axiomatic Method rather than Hilbert's Formal Axiomatic Method in its original form as described in the last Chapter: 

\begin{quote}
[E]very mathematical theory is a concatenation of propositions, each one derived from the preceding ones in conformity with the rules of a logical system [..] It is therefore a meaningless truism to say that this ``deductive reasoning'' is a unifying principle for mathematics. [..] [I]t is the external form which the mathematician gives to his thought, the vehicle which makes it accessible to others, in short, the language suited to mathematicians; this is all, no further significance should be attached to it. \\
What the axiomatic method sets as its essential aim, is exactly that which logical formalism by itself cannot supply, namely the profound intelligibility of mathematics. [..] Where the superficial observer sees only two, or several, quite distinct theories, lending one another ``unexpected support'' [quoted by the author from \cite{Brunschvicg:1912}, p.446] through the intervention of a mathematician of genius, the axiomatic method teaches us to look for the deep-lying reasons for such a discovery, to find the common ideas of these theories, buried under the accumulation of details properly belonging to each of them, to bring these ideas forward and to put them in their proper light. (\cite{Bourbaki:1950}, p.223)
\end{quote}

In order to illustrate his point Bourbaki uses the example of the ``abstract'' group theory; the author describes here this theory as an axiomatic theory construed after the pattern of Hilbert's \emph{Foundations} of 1899 \cite{Hilbert:1899} as a  \emph{system of things}, which are subject to the following three axioms (modulo a slight change Bourbaki's original notion).  

\textbf{G1}: $x\circ (y\circ z) = (x\circ y)\circ z$  (associativity of $\circ$)

\textbf{G2}: there exists an item \emph{1} (called \emph{unit}) such that for all $x$  $x\circ 1 = 1\circ x = x$ 

\textbf{G3}: for all $x$ there exists $x^{-1}$ (called \emph{inverse} of  $x$) such that $x\circ x^{-1} = x^{-1}\circ x =  1$.

A system of things satisfying these axioms as called a \emph{group}. Expression $x\circ y = z$ stands here for an abstract binary \emph{algebraic operation}, which in the given context is to be understood as a (uninterpreted) logical ternary \emph{relation} $R(x, y, z)$ having this special property: if $x = x'$ and $y = y'$ then $R(x, y, z) \leftrightarrow R(x', y', z)$ (which means informally that the output of the operation is uniquely determined by its input).   

 The above axiomatic theory of groups (which I shall denote \textbf{GT} for further references) have various interpretations, which were known and studied before \textbf{GT} the rise of Hilbert's Axiomatic Method: by interpreting variables $x,y,z,$ as invertible  geometrical transformations (like motion) and interpreting the operation $\circ$ as composition of these transformation one gets the notion of group of geometrical transformation; by interpreting variables $x,y,z,$ as whole numbers and interpreting $\circ$ as $+$ (addition of whole numbers) one gets the additive group of whole numbers, etc. Such examples belong to different domains of mathematics and many of them play some significant role in their proper domains. But until \textbf{GT} was axiomatically formulated as above \footnote{As Bourbaki notices here \textbf{GT} can be defined through different axioms. What determines the identity of this theory is its set of true propositions (including both axioms and theorems inferred from these axioms) but not a given set of axioms.} and until it brought about the precise \emph{general} notion of group those examples could not be understood as instances and special cases of one and the same thing (for no such thing was known yet!) and the links between these different groups, which were eventually guessed by some smart mathematicians, looked as unsystematic and sometimes even mysterious. Thus in this case the Axiomatic Method helps to bring about the new powerful mathematical concept of a group and develop the corresponding general theory, which unifies a large spectrum of significant mathematical results from different areas of mathematics.

As we see Bourbaki points here on a general phenomenon, which is not specific to the Axiomatic Method (and moreover to Hilbert's Formal Axiomatic Method), namely to the mathematical concept-building and its unifying role. Such basic mathematical concepts as \emph{number}, \emph{figure} and the like can be similarly seen as abstractions generalizing upon more specific examples like sets of dots or sets of strokes (for number) and circles, polygons, etc. (for figures). Although the generalization upon and abstraction from specific features of previously known examples is not the only way in which emerge new mathematical concepts this way of emergence is a major one. And in Bourbaki's example it works through the Axiomatic Method. Let us see how it works more attentively.

First of all we need to distinguish between two ways in which an axiomatic theory \emph{unifies} its content. When a set of contentual propositions is logically deduced from certain propositions belonging to the same set and chosen as axioms this unifies all these propositions into a single (contentual) theory. As we have seen Bourbaki recognizes this fact but in the last quote he clearly points to a different way of unification, which is equally made possible by the Axiomatic Method. This different way of unification is made possible by the \emph{formal} character of Axiomatic Method, where ``formal'' is to be understood in the sense of Hilbert's \emph{Foundations} of 1899 rather than in the sense of his \emph{Foundations} of 1927 \cite{Hilbert:1927} (which is referred to in the above quote in the expression ``logical formalism''). This second axiomatic unification amounts to the following: a  \emph{formal} (as opposed to contentual) axiomatic theory unifies its interpretations (models) by identifying certain common features of these interpretations and abstracting from all other specific features of those interpretations. The usual talk of interpretation of a given formal theory takes it for granted that the formal theory is given first and interpreted next. Now we reverse the perspective and consider the (would-be) interpretations as given (as contentual mathematical theories and fragments of such theories) and then think how to make up a formal theory, which captures common features of these things and thus unifies them. The example of \textbf{GT} is used by Bourbaki to illustrate the latter but not the former way of unification. 

Thus Bourbaki shows - in my view quite correctly - that the Formal Axiomatic Method has a unifying capacity, which is absent from contentual versions of the Axiomatic Method. However Bourbaki's version of Axiomatic Method is not identical to Hilbert's! Let me now describe the difference.

As an example of a \emph{theorem} of \textbf{GT} Bourbaki mentions this proposition \textbf{P}: 

For all $x,y, z$ if  $x\circ y = x\circ z$ then $y =z$

which follows from \textbf{G1} - \textbf{G3} almost immediately. I claim that the simplicity of this example does not allow it to represent correctly Bourbaki's Axiomatic Method. Notice that among \emph{objects} of \textbf{GT} there is no (abstract) groups just like among objects of Euclidean (3D) geometry developed in Hilbert's  \emph{Foundations} of 1899 there is no such thing as (3D) Euclidean space. As a  \emph{system of things} (model) of \textbf{GT} any group is a domain where all axioms and theorems of  \textbf{GT} hold; objects of this theory are \emph{elements} of the given group but not this group itself. But Bourbaki's theory of (abstract) groups (see \cite{Bourbaki:1939-1988} vol. 2, Chapter 1, Section 6) like any other presentation of this theory does treat groups as its objects, distinguishes between different groups, classify them and makes various constructions with them. Notice that \textbf{GT} by itself does not allow one even to formulate the notion of \emph{subgroup}! Take also in consideration that \textbf{GT} is not \emph{categorical} (in the usual model-theoretic sense of the term), which simply amounts to saying that not all groups are isomorphic. So axioms \textbf{G1} - \textbf{G3} provide nothing but the general notion of abstract group and in this sense can be compared with a  \emph{definition} of some traditional mathematical object like triangle; theorems of  \textbf{GT} like \textbf{P} are to be compared with propositions like ``all triangles have three angles'' implied by the definition of triangle. Remind the famous passage from Kant's \emph{Critique of Pure Reason} that has been already quoted above       

\begin{quote}
``Give a philosopher the concept of triangle and let him try to find out in his way how the sum of its angles might be related to a right angle. He has nothing but the concept of figure enclosed by three straight lines, and in it the concept of equally many angles. Now he may reflect on his concept as long as he wants, yet he will never produce anything new. He can analyze and make distinct the concept of a straight line, or of an angle, or of the number three, but he will not come upon any other properties that do not already lie in these concepts.   But now let the geometer take up this question. He begins at once to \underline{construct a triangle} .... In such a way through a chain of inferences that is always \underline{guided by intuition}, he arrives at a fully illuminated and at the same time general solution of the question.'' (\emph{Critique of Pure Reason} \cite{Kant:1999}, A 716 / B 744)
\end{quote}
    
Kant's point apparently applies to the case of group theory too. Although there are some general truths like \textbf{P} which logically follow from the concept of group (i.e., from axioms \textbf{G1} - \textbf{G3}) those truths are of little mathematical interest. The genuine mathematical work in the (abstract) group theory begins when groups are conceived of (or ``constructed'') as individual objects. In this way mathematicians prove non-trivial facts about groups that do not follow from \textbf{G1} - \textbf{G3} alone (the mathematical reader may think of Lagrange's theorem for a simple example). I do not mean to suggest here that Kant's philosophy of mathematics in its original form fully applies to the modern mathematics in general and to the modern group theory in particular. So we need to study more precisely \emph{how} groups and other modern mathematical concepts are constructed as objects; this issue will be thoroughly discussed in Chapter \textbf{8}. However already at this point  it is clear that \textbf{GT} by itself no more deserves the name of group theory than a definition of triangle deserves the name of a \emph{theory} of triangle. Kant's point about the proper geometrical study of triangles and a ``philsophical'' study of triangles that seeks to get the relevant knowledge directly from definitions (see \textbf{1.3} above) is perfectly relevant in the case of group theory! 

Before I make explicit the way in which groups and other mathematical concepts are constructed as (mathematical) objects I need to clarify a widespread terminological ambiguity. Bourbaki call abstract groups \emph{abstract} by the contrast with such ``concrete'' examples of groups as the additive group of whole numbers, the group of Euclidean motions and the like. This way of distinguishing between the abstract and the concrete in mathematics is specific for the given mathematical context and should not be confused to the general distinction between abstract concepts and concrete instants of those concepts, which is applicable everywhere in mathematics and beyond. The former distinction describes a way of obtaining some mathematical concepts from some others: given the concept of the additive group, the concept of group of Euclidean motions, the concept of group of permutations, etc., one forms the general concept of group through abstraction from certain specific features of each of these examples. The obtained concept of so-called \emph{abstract} group like any other mathematical concept may be though of both \emph{in abstracto}, i.e., as a self-standing non-represented concept, and in  \emph{in concreto} as represented by an individual object or a set of such objects. The confusion between the two senses of ``abstract'' and ``concrete'' in mathematics is not without a reason; it can be explained by the fact that at certain points of history certain concepts obtained through abstraction from earlier known mathematical concepts cannot be immediately represented  \emph{in concreto}, and so remain speculative rather than genuinely mathematical until they get supported by newly developed corresponding modes of representation. In Chapter \textbf{8} below I show more precisely how this dialectics of abstract and concrete works in the history of mathematics. Notice that the title of ``abstract group theory'' used by Bourbaki in 1950 today is an anachronism: what Bourbaki in 1950 calls an ``abstract group'' today is commonly called simply a group.  

Let's now see what kind of representation of the abstract group concept is used by Bourbaki in their volume on algebra, which includes the group theory. Like in other similar cases they use for it a \emph{set-theoretic} representation; the relevant set theory is developed in the first volume of Bourbaki's  \emph{Elements} (I shall tell more about it shortly.) Crucially, the set-theoretic representation does not reduce to interpreting \textbf{GT} in set theory. One gets a set-theoretic model of \textbf{GT} by interpreting variables $x, y, z,$ as elements of some set $G$ and interpreting the group operation $\circ$ in terms of Cartesian product of sets which reduces it to the primitive set-theoretic relation of membership. However such a model of \textbf{GT} is nothing but one particular group $G$ but not a domain where live \emph{all} groups accounted for by group theory! So the group theory as developed in Bourbaki's \emph{Elements} is not just an interpreted version of \textbf{GT} but a theory of  \emph{set-theoretic models} of \textbf{GT} developed on the basis of (Bourbaki's) set theory. Thus all Bourbaki's groups live in an universe of sets (that interprets their set theory). The same is the case for objects (structures) of other sorts treated in Bourbaki's  \emph{Elements}. In this sense set theory qualifies as a foundation of all Bourbaki's mathematics.  

We see that Hilbert's Axiomatic Method is used by Bourbaki in a very peculiar way, which is not made clear by Bourbaki's remark about the Axiomatic Method in his paper of 1950 \cite{Bourbaki:1950}. Namely, the Hilbert-style axiomatic theory of groups \textbf{GT} is used by Bourbaki for defining \emph{objects} of their theory  of groups, which is not an informal version of \textbf{GT} but something quite different! Objects of this sort Bourbaki call  \emph{structures}. Formal Axiomatic Method involves the idea that  \emph{isomorphic} structures are essentially the same. Thus Bourbaki's structures have a double identity criterion: in one sense the identity of a given structure is given by the identity of its underlying set (which is formally treated in the corresponding set theory) and in a different sense the identity of structures is their isomorphism. In Chapter  \textbf{8} we shall see how this conceptual tension translates into a philosophical controversy between \emph{structuralism} and \emph{set-theoretic substantialism}.  The issue of identity of structures will be also discussed in Part \textbf{2} of this book. Here I would like to stress that the peculiar way, in which Bourbaki applies Hilbert's Axiomatic Method, allows Bourbaki's mathematics to fit (by and large) the pattern of traditional mathematics, where mathematical concepts are represented by individual objects. 

Let us now look at Bourbaki's set theory presented in the first volume of his \emph{Elements} \cite{Bourbaki:1939-1988}. This theory follows the formal method of Hilbert's \emph{Foundations} of 1927 closer than any other known to me axiomatic exposition of set theory.   The First Chapter of the volume, which has the title \emph{Description of Formal Mathematics}, begins with an account of \emph{signs} and \emph{assemblies} (strings) of signs provided with a definition of mathematical theory according to which a theory

\begin{quote}
... contains rules which allow us to assert that certain assemblies of signs are \emph{terms} or \emph{relations} of the theory, and other rules which allow us to assert that certain assemblies are \emph{theorems} of the theory. (quoted by English translation \cite{Bourbaki:1968}, p. 16)  
\end{quote}

Then follows a description of operations that allow for constructing new assemblies of signs from some given assemblies; the simplest operation of this sort is the \emph{concatenation} of two given assemblies $A, B$ into a new assembly $AB$.  On such a purely syntactic basis Bourbaki introduces some logical concepts necessary for a formal axiomatic treatment of set theory. Although  Bourbaki's version of axiomatic set theory is not identical to ZF it is similar in its character; the differences between the two ways of formalizing set theory are not relevant to the present discussion and I leave them aside. The general Bourbaki's strategy in set theory is to develop formally all set-theoretic concepts used in various specific theories treated in other volumes after  the pattern of group theory, so that all these specific theories appear as (or can be translated into) fragments of the formal Set-theory. 

As a matter of fact the formal notation used in the first volume never reappears in the following volumes, so the possibility of translation of specific theories into the formal set theory remains theoretical and is never explicitly explored. Instead the formal set theory presented in the first volume Bourbaki everywhere else uses its informal version, which is systematically exposed in an unpublished draft of the first volume of the \emph{Elements} \cite{Bourbaki:19??}. After a philosophical introduction  the author introduces (informally) the concept of \emph{fundamental set} and the relation of membership between sets and their elements. (The author calls a set \emph{fundamental} in order to distinguish a well-defined set concept from a more general notion of collection.) Then the author introduces (with the usual informal notation) the concept of \emph{subset}, which is the subject to the following axiom:  

 \begin{quote}
Any predicate of type $A$ defines a subset of $A$; any subset of $A$ can be defined through a predicate of type $A$. (\cite{Bourbaki:19??}, p.7, hereafter my translation from French)
 \end{quote}

(Predicate of type $A$ is a predicate $P$ such that for every element $a$ of set $A$ $P(a)$ has a definite truth-value. The subset $S$ of set $A$ defined by $P$ consists of such and only of such $a$ for which $P(a)$ is true.)  

Next Bourbaki introduces the concept of \emph{complement} of a given subset, of \emph{powerset} $P(A)$ of a given set $A$ (i.e., the set of all subsets of $A$); of \emph{union}, \emph{intersection} and \emph{cartesian product} of sets (described as \emph{operations} on sets), of \emph{relation} and \emph{function} between sets. Having these basic concepts in his disposal the author says:

\begin{quote}
In any mathematical theory one begins with a number of fundamental sets, each of which consists of elements of a certain type that needs to be considered. Then on the basis of types that are already known one introduces new types of elements (for example, the subsets of a set of elements, pairs of elements) and for each of those new types of elements one introduces sets of elements of those types. 
\\  
So one forms a family of sets constructed from the fundamental sets. Those constructions are the following: \\
1) given set $A$, which is already constructed, take the set $P(A)$ of the subsets of $A$;\\
2) given sets $A$, $B$, which are already constructed, take the cartesian product $AxB$ of these sets.\\
The sets of objects, which are constructed in this way, are introduced into a theory step by step when it is needed. Each proof involves only a finite number of sets. We call such sets \emph{types} of the given theory; their infinite hierarchy constitutes a \emph{scale of types}.   (\cite{Bourbaki:19??}, p.43-44)
\end{quote}
  
On this basis the author describes the (informal) concept of  \emph{structure} as follows:  
\begin{quote}
We begin with a number of fundamental sets: $A, B, C, ... , L$ that we call \emph{base sets}. To be given a \emph{structure} on this base amounts to this: \\
1) be given properties of elements of these sets;
2) be given relations between elements of these sets;
3) be given a number of types making part of the scale of types constructed on this base;
4) be given relations between elements of certain types constructed on this base;
5) assume as true a number of mutually consistent propositions about these properties and these relations. 
 (\cite{Bourbaki:19??}, p.44-45)
\end{quote}

What I want to stress is the fact that principles of building mathematical theory described in the Bourbaki's draft are not so different from Euclid's: Bourbaki like Euclid begins with principles of building mathematical objects but not with axioms. Axioms (in the modern sense of the term) appear only in the very end of the above list (the 5th item) and as I have already argued they play the role of definitions. While for Euclid the basic data is a finite family of \emph{points}  and the rest of the geometrical universe is constructed from these points by Postulates for Bourbaki the basic data is a finite family of \emph{sets} and the rest is constructed as just described. While for Euclid the basic type of geometrical object is a \emph{figure} for Bourbaki the basic type of mathematical object is a \emph{structure}. In both cases the constructed objects come with certain propositions that can be asserted about these objects without proofs because they immediately follow from corresponding definitions. In both cases the construction of objects is a subject of certain rules but not the matter of a mere stipulation. In order to continue this analogy one may compare the notion of \emph{isomorphism} of structures in Bourbaki with Euclid's notion of \emph{congruence} of figures. 

We see that Bourbaki's mathematics has two different layers. On the ground layer it has a formal set theory built with Formal Axiomatic Method. But the rest of this mathematics is developed with an informal notion of set described through constructive postulates like ``take the cartesian product of given sets'', which are similar to Euclid's Postulates. Calling such postulates constructive I mean not that they involve only finitary operations (they are obviously infinitary) but that they are not existential propositions. However the formal theory of sets provides a theoretical possibility of translating all these informal set-theoretic postulates into formal existential axioms and ultimately - of translating all of Bourbaki's mathematics into strings of symbols operated according to purely syntactic rules, which can be also described as constructive postulates (about operating with strings of  symbols). So this translation allows (in principle) for replacement of informal set-theoretic infinitary postulates into finitary syntactic postulates. 

It is interesting to observe that the ground set-theoretic layer of Bourbaki's mathematics is seen by the author as undesirable; it looks like Bourbaki would be happy to get rid of it but does not know how to do this exactly. Here is what he says about it in the earlier quoted draft (my translation from French):

\begin{quote}
The reader will see that the nature of elements of fundamental sets can be always easily left undetermined and that this point of view is often useful. From here there is only one step to thinking that only structure matters and that the true aim of mathematical theory is a study of structure independently from sets that may represent it. Perhaps it is indeed possible to study structures themselves and forbid oneself to consider fundamental sets. However because of the commodity of language and the invincible habit of mind we take the ``ontological'' approach, i.e., stipulate fundamental sets for each theory.      
\end{quote}

and in his manifesto of 1950 Bourbaki says in a footnote:

\begin{quote}
We take here a naive point of view and do not deal with the thorny questions, half philosophical, half mathematical, raised by the problem of the ``nature'' of the mathematical ``beings'' or ``objects''. Suffice it to say that the axiomatic studies of the nineteenth and twentieth centuries have gradually replaced the initial pluralism of the mental representation of these ``beings''  thought of at first as ideal ``abstractions'' of sense experiences and retaining all their heterogeneity  by an unitary concept, gradually reducing all the mathematical notions, first to the concept of the natural number and then, in a second stage, to the notion of set. This latter concept, considered for a long time as ``primitive'' and ``undefinable'', has been the object of endless polemics, as a result of its extremely general character and on account of the very vague type of mental representation which it calls forth; the difficulties did not disappear until the notion of set itself disappeared (and with it all the metaphysical pseudo-problems concerning mathematical ``beings'') in the light of the recent work on logical formalism. From this new point of view, mathematical structures become, properly speaking, the only ``objects'' of mathematics. (\cite{Bourbaki:1950}, p. 225-226)
\end{quote}

Since ``the recent work on logical formalism'' referred to by Bourbaki in the last quote does not contain anything that goes beyond Hilbert's Formal Axiomatic Method these passages are rather puzzling. A way to understand Bourbaki's misgivings about set theory is this. A strictly formal version of the Axiomatic Method like one applied in Bourbaki's volume on set theory does not require set theory developed in advance; however except this volume the set theory (namely its informal version) is used by Bourbaki essentially. So the misgiving may be due to the fact that, in particular, Bourbaki's group theory is not a formal theory like \textbf{GT} (which by itself does not require the notion of set) but a theory of set-theoretic models of \textbf{GT}. Although this fact may indeed explain Bourbaki's reference to the ``logical formalism'' the idea of making structures into objects is hardly compatible with the Formal Axiomatic Method (for structures determined by formal axiomatic theories are not objects of their theories). So what Bourbaki really aim at is not a strictly formal axiomatic treatment of group theory and the rest of mathematics but rather finding an appropriate replacement for set theory in its role of universal vehicle carrying all specific mathematical structures. In \textbf{8.5} I shall come back to this issue and consider category theory as such an alternative vehicle.

Coming now back to the question about the role of the Formal Axiomatic Method in the 20th century mathematics we observe the following. Bourbaki formulated an important part of the contemporary mathematics in set-theoretic terms and thus showed a theoretical possibility to formalize this part of mathematics by reducing it to formal axiomatic set theory. However the Formal Axiomatic Method did not become in Bourbaki's hands the ``basic instrument of all research'' as Hilbert suggested in 1927 (\cite{Hilbert:1927}, p. 467); it is used instead only for foundational purposes (notwithstanding the fact that the syntax of Bourbaki's formal set theory  lacks the intuitive clearness through which Hilbert hoped to ground the infinitary reasoning in mathematics). As an instrument of research and more importantly as an instrument of presenting a ready-made knowledge in a systematic form Bourbaki uses a different version of Axiomatic Method, which I have tried to describe in this paragraph. Describing this other method as \emph{informal} one should keep in mind that the relevant difference between being formal or informal is not a matter of degree. As we have seen Bourbaki's method involves features, which are not present in Hilbert's method, in particular the idea of using formal axioms for defining objects of a given theory rather than developing this theory itself. Another important difference concerns the doing-showing dilemma. While Formal Axiomatic Method makes mathematical (as distinguished from metamathematical) reasoning into a (syntactic) construction Bourbaki uses informal set-theoretic constructions and proves (informally) theorems about these constructions and in this sense more closely follows Euclid's example. The fact that this new way of doing-and-showing can be in principle represented formally should not be, of course, neglected. However it does not provide any substance to the claim that the Axiomatic Method used by Bourbaki throughout his \emph{Elements} qualifies as formal in Hilbert's sense.

\section{Galilean Science and Set-Theoretic Foundations of Mathematics}

Remind from the last Section Bourbaki's claim according to which ``every mathematical theory is a concatenation of propositions, each one derived from the preceding ones in conformity with the rules of a logical system'', which the author suggests as self-evident. This is a statement of what I have called above a mathematical logicism in the large sense and argued that it is not self-evident at all.  In Chapter \textbf{1} I have shown that Euclid's mathematics does not fit the logicist description of mathematics as a ``concatenation of propositions'' because in addition to  theorems it contains \emph{problems}, which aim at constructing some objects rather than prove some propositions. In \textbf{2.4 - 2.6}) we have seen that Hilbert's Formal Axiomatic Method (in the sense of his \emph{Foundations} of 1927 and later works) employs a very special form of mathematical logicism, which reduces a ``concatenation of propositions'' to a symbolic construction and in this way assures a fundamental role of mathematical intuition of a special  sort. Finally we have observed that Bourbaki's mathematics, in fact, is not quite unlike Euclid's: it cannot be described as a concatenation of propositions either because it similarly involves  constructions of objects of special sort called structures. The idea of reduction of all constructive postulates to existential propositions plays a role \emph{only} at the ground level of the formal set theory, where informal set-theoretic constructions indeed translate into symbolic constructions (or more precisely, a \emph{possibility} of such translation is \emph{shown} albeit no such translation is actually \emph{done}). Now I would like to return to the issue of mathematical logicism and argue that mathematics has a non-propositional constructive aspect, which makes possible application of mathematics in empirical sciences and technology. This view is obviously not original and dates back at least to Kant. Nevertheless it is appropriate to spell it here again in modern terms.

It is hardly controversial that mathematics deals with forms of possible human experience;  in its simplest and most general form this claim is simply tantamount to saying that mathematics applies across a wide range of human practices. Today this is even more true than it was in Kant's time: crucial technologies, on which depend our well-being, in many ways depend on mathematical considerations and cannot be sustained and further developed without mathematical expertise; mathematics today makes part of any engineering education. In Kant's time the only properly mathematized science was (Newtonian) mechanics; the following progress of science in the 19th century has brought us to the point when every physical theory deserving the name has a mathematical aspect. Today physics and chemistry are mathematized and the mathematization of biology is in progress. Using mathematical models also becomes an usual practice in social sciences.

Let me now be more specific and ask which general forms of experience are relevant to today's science and technology. This question is obviously yet too large and too general to be answered in any detail here. I would like however to stress only the following point. At least since Galileo's times science practices an active intervention of humans into the nature through experiments rather than a passive observation and description of the observed phenomena. So we are talking now about the mathematically-laden science where mathematics serves for guiding human interactions with the environment rather than simply for describing how this environment appears to our senses. As van Fraassen puts this

\begin{quote}
The real importance of theory, to the working scientist, is that it is a factor in experimental design. (\cite{Fraassen:1980}, p. 73)
\end{quote}

Thus mathematical forms of possible experience relevant to the modern science are forms of such possible \emph{interactions} with the environment rather than only linguistic and logical forms that allow for spelling out some plausible hypotheses about the world and deriving from them some consequences according to certain rules. The forms of the latter kind may be sufficient for developing a speculative science along the older scholastic pattern but they are certainly not sufficient for developing the modern mathematically-laden science and the modern mathematically-laden technology.           
 
As far as the pure mathematics is conceived as a domain of abstract logical possibilities the fact that mathematics proves ``unreasonably effective'' \cite{Wigner:1960}  in its applications to empirical sciences and technology remains a complete mystery. The mystery is dissolved as soon as one observes that mathematics explores not everything that can possibly \emph{be} the case (which is a hardly observable domain unless one delimit the sense of ``possibly'' in one way or another) but rather what we can possibly \emph{do} within the limits of our human capacities (which are steadily growing with the progress of science and technology). What are these limits is a tricky question. On the one hand, mathematics systematically ignores certain apparent limits by exaggerating relevant capacities: this is usually called the \emph{mathematical idealization}. For example, mathematicians pretend that they can count up to $10^{10^10}$ just as easily as up to 10 or that they can draw a straight line between two stars just as easily as they can draw a line between two points marked on a sheet of paper. This strategy usually works until the point where the empirical constraints become pressing and people invent new mathematics that takes these constraints into account as this, for example, happened when people realized that the old good Euclidean geometry is not appropriate for describing the physical space at large astronomical scales (in spite of the fact that it still works amazingly well at the scale of a planetary system like ours). One the other hand, it also happens that in a real experiment people observe what in terms of the assumed mathematical description of this experiment qualifies as impossible as this happened in the Michelson-Morley experiment supposed to measure parameters of the ether flow around the Earth. In such cases people say that the assumed mathematical description (and hence the corresponding physical theory) is wrong and look for a new one. Sometimes the suitable mathematics can be found in a nearly ready-made form and only used for building a new physical theory but sometimes in order to fix the problem one needs to develop the appropriate mathematics from the outset as this happened in the history of the electro-magnetism, for example. This picture suggests the view on the pure mathematics as a proper part of the modern Galilean science.

Russell's neo-Leibnizian logicism about mathematics promises nothing more and nothing less than that: to make mathematics a part of logic, so that any mathematical form of possible experience turns into the form of a \emph{proposition} (and forms of logical inference of propositions from some other propositions), which may eventually refer to some experience. Russell's view is quite radical in this respect, and many people including Hilbert who were directly involved into reforming mathematics on the basis of new logic in the beginning of the 20th century didn't share Russell's philosophical views. Anyway, as I have already argued, a weaker form of the neo-Leibnizian approach is intrinsic to Hilbert's Formal Axiomatic Method.  Even if forms of possible experience delivered by a formal axiomatic theory do not qualify as \emph{logical} forms in the precise sense of the term they nevertheless are forms of possible empirical \emph{propositions} rather than forms of empirical interactions or anything else. This, in my view, explains the very little success that Hilbert's Axiomatic Method has had so far in physics and other natural sciences.

As far as we want to continue to develop the Galilean science (and the technology connected to this type of science) our mathematics must provide for it forms of possible empirical interaction rather than just forms of propositions. In other words it must provide forms appropriate for \emph{doing} various things in the world but not only forms for talking about this world and \emph{showing} how the world looks like in a mathematical representation. Since formal axiomatic theories are not appropriate for this job we need to learn how to build mathematical theories differently.

As Kant shows in great detail the traditional geometry and his contemporary algebra are useful in the Galilean science because these mathematical theories are \emph{constructive} in the sense that they involve rules for constructing their objects (explicitly or implicitly). Today we cannot hope, of course, to get a new mathematical theory that would allow for identifying a physical object with a mathematical object in the same way, in which one may identify (modulo the mathematical idealization), say, a planet with an Euclidean sphere. Today's physicists describe particles using the mathematical group theory and manipulate with these particles in experiments using a special hi-tech equipment;  they don't expect that mathematical manipulations with groups would map their experimental manipulations in a direct way. Nevertheless  the constructive character of the mainstream informal mathematical practice, which I have stressed earlier in this Chapter, still helps physicists and other scientists to design their experiments and their equipments.  Scientists make up mathematical models of their experimental systems, manipulate both with the models (theoretically) and with the experimental systems (in real experiments) and see whether the manipulations of the two sorts work coherently. This is, of course, an oversimplified picture of the scientific experiment (for more details see \cite{Fraassen:1980}) but it is sufficient for seeing that the possibility to establish a correlation between mathematical manipulations, on the one hand, and experimental manipulations, on the other hand, remains essential in today's mathematically-laden experimental science.

I cannot see how such a correlation cannot be possibly established when the only type of mathematical objects available for manipulation are syntactic objects, which represent (logical forms of) propositions related to certain experimental settings but do not represent these experimental settings themselves. Given a proposition expressed in a formal language one may interpret it in terms of physical data and evaluate whether under this interpretation the given proposition is true or false. From a number of so interpreted true propositions one may deduce some further propositions, consider them as physical predictions and finally check these predictions against the available physical data. So far so good but notice that from a methodological viewpoint this way of doing natural science resembles Ptolemean astronomy saving phenomena rather than the Galilean physics intervening into the nature with experiments! In order to design an experiment one needs a mathematical model (representation) of a given physical environment itself but not of formal propositions interpretable in terms of this environment; mathematical manipulations (i.e. some further constructions) with a model of the former sort may serve as a guide for real experimental manipulations within the given environment, mathematical manipulations with a model of the latter sort cannot be directly used for this purpose.

Bourbaki's \emph{Elements} present a possible compromise between the formal axiomatic approach and a more traditional constructive approach of doing mathematics. (Beware that I use here the term ``constructive'' in the sense of Hilbert's distinction between ``construction postulates'' and propositional ``existential form'' (\cite{Hilbert&Bernays:2010}, p. 20). This sense of being constructive does not imply anything about the admissibility of mathematical constructions of a given type. The admissibility of constructions depends on concrete constructive postulates but I am now taking about constructive postulates in general.) So we have in Bourbaki's \emph{Elements} a formal axiomatic theory (of sets) on the ground level, informal set-theoretic constructions on upper levels and the (theoretical) possibility to translate the higher-level mathematics into the formal ground-level terms. Arguably such a translation (however impractical it may be) plays an important  justificatory role since it allows for reduction of problematic set-theoretic constructions (like the construction of the powerset of a given infinite set) to formal deductions from axioms of set-theory (including existential axioms) to (long and tedious but anyway finitary) symbolic constructions. So the question of admissibility of problematic set-theoretic constrictions reduces to the question of admissibility of the corresponding symbolic constructions. And in this latter case we have a clear criterion of admissibility: a given formal deduction is admissible if and only if it, first, conforms the explicit rules of deduction and, second, the given formal theory of sets is \emph{consistent}, which is tantamount to saying that from the given axioms the proposition $0\neq 0$ can \emph{not} be deduced by the aforementioned rules (i.e. that these rules do not allow to build a symbolic construction with the axioms on the one end and $0\neq 0$ on the other end). 

Although having such a device for checking suspicious mathematical constructions and suspicious mathematical proofs seems to be a good idea Bourbaki's attempt to apply this device in practice brings a rather controversial outcome. First of all it does not allow for an effective proof-checking: formal versions of Bourbaki's proofs are too long and cumbersome to be survey able by a human mathematician and at the same they are not adopted for an automatic proof-checking with a computer. This is why  the theoretical possibility of translating set-theoretic constructions into symbolic one may serve only as an \emph{epistemic ground} justifying certain basic set-theoretic constructions like the powerset. 

It must be however stressed that in Bourbaki's \emph{Elements} the Formal Axiomatic Method does not perform this foundational task as perfectly as Hilbert has imagined it in 1927. There are two fundamental difficulties here that Hilbert did not acknowledge. First, the criterion of formal consistency cannot be used straightforwardly because the proof of consistency of Bourbaki's set theory or any other formal set theory requires using some stronger theory. So we don't really know whether or not the background set theory is consistent. Second, what can and what cannot be proved in a formal theory depends on the chosen background logic and such a choice, as I have already stressed, is a matter of controversy. Given these fundamental difficulties I doubt that formal Set theories can play the foundational role they are supposed to play in Bourbaki or elsewhere in the modern mathematics. The powerset operation is admitted by Bourbaki not \emph{because} its formal existential version (i.e., the formal powerset axiom) is formally consistent with other axioms of set theory but rather on an intuitive ground by an appropriate modification of the intuition relevant to the finite case. The finite intuition underpinning Euclid's First Postulate is similarly modified in the case of the Second Postulate, which allows for an infinite extension of a given straight segment. Although the infinite powerset operation requires a deeper modification of the finite intuition this is nevertheless a matter of degree rather than matter of principle. This is why any sharp distinction between ``real'' and ``ideal'' objects in mathematics is in my view unjustified. In Chapter \textbf{7} I show how in the history of mathematics ``ideal'' objects become ``real'' through the development of appropriate intuitions. The power set operation can be, in my view, stipulated directly as a constructive postulate (``constructive'' in the above sense!) rather than in the roundabout way through the existential powerset axiom. Having said that I do not deny that the formal logical analysis has helped to clarify some important issues in foundations of set theory and in foundations of mathematics in general. We have learned about certain set-theoretic paradoxes and ways of avoiding them. In particular we learned to distinguish sets from \emph{proper classes} (like the class of all sets, see \textbf{5.8}). So even formal axiomatic set theories cannot provide ultimate foundations they show how us how to avoid known contradiction by forbidding certain constructions like the set of all sets.

Using a formal set theory as a foundation and using a set theory informally for building structures are related but still different issues, which we need to distinguish carefully. We have seen that the idea of building mathematical objects as \emph{structured sets} comes from the Formal Axiomatic Method but then lives an independent life within theories, which are nor formal neither axiomatic in Hilbert's sense of these terms. It is appropriate to ask whether the informal ``set-theoretical language'' is indeed sufficient for doing modern mathematics. We shall see in Chapter \textbf{8} that the answer is rather in negative. Here I shall only make a remark concerning possible applications of Bourbaki's structural mathematics in natural sciences. Although the informal Bourbaki's structural mathematics just like Euclid's mathematics involves constructions (but not only logical deductions) the set-theoretic nature of these constructions apparently makes an obstacle for interpreting these constructions in physical terms and applying them in physics and other natural sciences. Consider the case of group theory: in spite of the fact that this theory is widely used in the 20th century physics physicists normally think of groups as \emph{groups of transformations} (see \textbf{6.1} ) rather than construe them as structured sets. This fact is hardly surprising since Cantor's notion of infinite set  is of metaphysical rather than physical origin and has no physical sense; naively this can be expressed by saying that infinite sets do not exist in nature. The later development of set theory, which eventually led to modern axiomatic theories of sets, clarified the logical aspect of set theory but did not have anything to do with physics either. This makes Bourbaki's structural mathematics quite unlike Euclid's mathematics where the basic mathematical concepts such as number and (geometrical) magnitude are directly relevant to the material practices of counting and measuring, which are even today indispensable in any empirical science. The question, which I want to stress now is the following: is the present detachment of foundations of mathematics from the foundations of natural science indeed an epistemic necessity or rather an outcome of bad epistemic strategy? 

It is true that the contemporary mathematics is so sophisticated and the experience relevant to the contemporary fundamental physics is so unlike the everyday human experience  that we cannot hope to find a pre-established cognitive mechanism providing a natural link between mathematics and physics at the level of their foundations. However constructing such a link artificially can be an epistemic strategy on its own. As I have already mentioned in \textbf{2.2} the idea to use metaphysics as a guide and formal logic as a tool for building foundations of mathematics is related to the anti-Kantian turn in the philosophy in the very beginning of the 20th century, which gave a new credit to the traditional pre-Modern and pre-Galilean patterns of doing science. I claim that the epistemic success of this type of science in the past is a sufficient reason to continue to keep its basic epistemic strategy untouched and think how to conform new mathematical and physical results with it rather than give up this strategy and replace it by a new form of Scholasticism. This may sound like knocking into an open door in the large scientific and mathematical community but not in the community of people professionally working in the foundations of mathematics, most of which still follows the logicist agenda established by Frege and Russell a century ago and works in an isolation from the mainstream mathematics and physics. I would like to stress that linking foundations of mathematics to foundations of physics is in \emph{not} in my view  only a pragmatic question of having better mathematical tools for physics but a fundamental question concerning foundations of the modern mathematized natural science. Although in this book I cannot present a systematic defense of the Galilean Science I claim that the irrelevance of set-theoretic foundations of mathematics to the contemporary physics is a strong reason for  rejecting these attempted foundations and looking for a replacement. 

\section{Towards the New Axiomatic Method: Interpreting Logic} 

In the 20th century the part Symbolic Logic, which describes itself as ``philosophical'', went through a booming development; for an overview I refer the reader to the last continuing edition of the \emph{Handbook of Philosophical Logic} edited by Gabbay and Guenthner \cite{Gabbay&Guenthner:2001-} \footnote{If one asks what is specifically ``philosophical'' about the multiple formal systems presented in this Handbook then, I think, the answer is twofold: all these formal systems are designed with certain philosophical motivations and/or used for treating some philosophical problems. The relevant notion of being philosophical derives from a particular notion of philosophy, which can be roughly identified with the Analytic Philosophy.}. Without trying to survey here the recent history of philosophical logic I shall try only to answer this question: Whether or not the development of logic after Hilbert brought about any new notion of Axiomatic Method? 

In order to answer this question let me first of all stress that Hilbert's Formal Axiomatic Method is extremely flexible and leaves one the freedom not only for building various axiomatic theories (including incompatible ones) but also for making a choice of the background logic. The mere replacement of (formalized) Classical logic by the (formalized) Intuitionistic logic or any other system of formal logic does require the replacement of the Formal Axiomatic Method by a new method of theory-building\footnote{Although this second degree of freedom (which adds to the free choice of axioms) has been not previewed by Hilbert himself I don't qualify its discovery as a modification of Hilbert's Axiomatics Method. However this new degree of freedom undermines Hilbert's suggested epistemic justification of his method and thus creates new epistemological problems. I shall discuss these problems and suggest a solution in Chapter \textbf{9}.}. However there is, in my view, at least one continuing development in logic, which indeed changes the sense of the Formal Axiomatic Method but not only presents a new application of the same method. 

Remind from \textbf{2.3} that in their logical textbook of 1928 \cite{Hilbert&Ackermann:1928} Hilbert and Ackermann distinguish between ``contentual'' and formal (i.e., symbolic) logic along with distinguishing between contentual and formal axiomatic mathematical theories. Do Hilbert and Ackermann also mean here a possibility of providing their system of symbolic logic with alternative interpretations (alternative models) along with alternative interpretations of formal axiomatic theories based on this system of logic?  As I have already argued this is not the case:  the authors rather think of formalization of ``contentual'', i.e., informal, logical concepts as a way of making these concepts sharper and better manageable with a help of symbolic means. A similar attitude is expressed in 1932 by Gentzen  in a footnote: 
\begin{quote}
If the words 'sentence' ('theorem') and 'proof are used informally as constituents of
our language they are of course intended to mean something quite different from the
purely formally introduced concepts of 'sentence' ('theorem') and 'proof (and even
under an intuitive interpretation the latter concepts are still considerably narrower than
the former); the context should make it clear in each case how these concepts are
intended. (\cite{Gentzen:1932}, p. 312)   
\end{quote}
So in Gentzen's view a formal language allows for a narrower and more precise formulation of logical concepts than the natural language. It is an essential part of this view that each  formal logical concepts has certain \emph{intended} interpretation, which can be expressed by formal means more precisely than by the natural language. Tarski is his \emph{Introduction}  of 1941 \cite{Tarski:1941} treats logic as a formal theory such that it necessarily makes part of any other (formal) theory and which can be called, in this sense, the \emph{minimal} theory. (He does not discuss in this textbook the issue of multiplicity of logics and talks about \emph{the} logic in singular.) After explaining the notion of model of a given formal theory and providing some examples Tarski remarks:

\begin{quote}
For precision it may be added, that the considerations which we sketched here are applicable to any deductive theory in whose construction logic is presupposed, but their application to logic itself brings about certain complications which we would rather not discuss here. (\cite{Tarski:1941} p. 119)
\end{quote}

Let me now point to such a complication, which is of philosophical rather than mathematical character. In his earlier paper \emph{Sentential Calculus and Topology} \cite{Tarski:1956} first published in 1938 Tarski develops topological interpretations of Classical and Intuitionistic propositional calculi. Under these interpretations the syntax of propositional logic is interpreted in terms of elements (open subsets) of a given topological space and operations with these elements, so a given well-formed formula $F$  designates (not a proposition but) a certain element $\phi$ of the given space construed from other elements. Then Tarski proves that $F$ is derivable in the given calculus if and only if $\phi$ has certain property $P$ (which is not the same in the Classical and the Intuitionistic cases). 

 Let me briefly explain for a non-mathematical reader what is  \emph{topology} and an \emph{open set}. Topology is a ``gummy geometry'', which studies those properties of geometrical objects, which remain invariant under invertible continuous transformations. Think about a circle (more precisely a circumference of a circle) and allow it to change its form however you want but without cutting it and without gluing it to itself. Such a ``gummy'' object qualifies as a circle (or ``one-dimensional sphere'') in the \emph{topological} sense of the term (even if it may look like an oval and so not qualify as circle in the usual Euclidean sense of the term). A standard way to define a topological object (usually called a topological \emph{space}) using set theory is the following. For a given set $T$ (of ``points'') one specifies which subsets $U \subseteq T$  count as \emph{open}; the complement $T\backslash U$ of an open set is called \emph{closed}. The choice of open (sub)sets is restricted by several axioms, which I shall not list here. Consider a circle (or any closed curve) on the Euclidean plane and think of these two things as sets of points. Then the ``natural'' topology is one where the interior area of the circle counts as open while the rest (the rest of the plane together with the curve) counts as closed. In this setting a continuous transformation is defined as one that reflects opens, i.e., one, which never transforms closed sets into open.

I shall not reproduce here details of this Tarski's construction and discuss only its problematic significance for the Formal Axiomatic Method. Remind that a pillar of Formal Axiomatic Method is the symbolic logic and that this pillar is twofold: it comprises a symbolic syntax, on the one hand, and a logical content, on the other hand. Each of the two aspects of symbolic logic has its specific epistemic impact: the logical content comprises basic laws of reasoning; the symbolic syntax makes these laws explicit through the symbolic intuition (\textbf{2.4}). This is why the  \emph{logical} content of formal theories (built by Hilbrt's Formal Axiomatic Method) unlike their non-logical content is indispensable. Indeed, given a list of mutually consistent axioms one may after Hilbert stipulate \emph{thought-things} and \emph{thought-relations} satisfying these axioms and consider them as mathematical objects; since the axioms are written with a symbolic language this language  provides these ``ideal'' objects with a symbolic representation, which makes it easy to reason about these objects as if they were ``real''. Further interpretations of the given formal theory can be interesting and useful (and may even constitute the pragmatic raison d'\^etre of the given theory) but, strictly speaking they are not necessary. However this way of construing mathematical objects does not work unless one has the notion of logical consistency in one's disposal. On the syntactic level the logical consistency translates into the impossibility of building a string of formulae, which begins with the axioms and ends up with a distinguished formula like $0 \neq 0$ (provided one follows fixed rules about building strings of formulas). However this syntactic counterpart of the logical consistency is not sufficient by itself to make the Formal Axiomatic Method work. In order to make the method work one must assume, in particular, that the symbolic expression $0 \neq 0$ expresses a contradiction. So the logical content unlike the non-logical content cannot be dispensed with in formal theories. It is always possible, of course, to study systems of symbols on their own rights but the Formal Axiomatic Method does not fall out of such a study unless one brings logical concepts into it. So interpreting  logical and non-logical terms in a formal theory are indeed two rather different issues.  Understandably Tarski did not want to enter into these details in his introductory chapter on model theory  \cite{Tarski:1941}.

With Tarski's topological interpretation of the Classical and Intuitionistic propositional calculi we get in addition to the \emph{intended} logical interpretation of the given symbolic calculi a non-intended one.  Tarski  \cite{Tarski:1956} uses for this interpretation Kuratowski's semi-formal set-theoretic approach to topology \cite{Kuratowski:1948-50}; a fully formalized version of this latter theory would comprise  (i) a background logic, (ii) a formal set theory, and (iii) an axiomatic theory of topological spaces interpreted in the underlying set theory. What sense one can make of the fact that a theory, which involves these three foundational levels (i), (ii), (iii) interprets a logical calculus, which can be a fragment of its own background logic (i)? 

The standard notion of model  and its usual epistemological underpinning (as developed in  \cite{Tarski:1941}) hardly allows for thinking about Tarski's topological interpretation as epistemically significant. However, if one sees this construction from a  different perspective, which brings one back to Boole and Venn, the same construction appears as a proper part of propositional logic (rather than its interpretation), which accounts (in topological terms) for the \emph{universe of discourse} relevant to this logic . We shall discuss this alternative view on logic in the next Chapter ( see \textbf{4.3} below). Let me only mention here that the issue of universe of discourse in Venn's work is closely related to his use of logical diagrams (see \textbf{2.5} above), and that the logical diagrams are topological objects in the sense that their sizes and shapes do not matter but their topological properties (and only such properties) do.  

Tarski himself makes the following remark about his proposed topological interpretation of the two propositional calculi (which the translator calls \emph{sentential}):

\begin{quote}
The present discussion seems to me to have a certain interest not only from the purely formal point of view; it also throws an interesting light on the content relations between the two systems of the sentential calculus and the intuitions underlying these systems [zugrundeliegende Intuitionen]. ( \cite{Tarski:1956}, p. 421)
\end{quote}
      
Does Tarski claims here seriously that certain topological intuitions underly logical calculi or he rather talks about intuitions here in a vague sense, which does not imply any epistemological commitment on his part?  Be it as it may Tarski's topological interpretation of propositional logic points to a connection between geometry and logic, which is very unlike that assumed in Hilbert's \emph{Foundations} of 1899 \cite{Hilbert:1899} and other works following the same line of thought. Remind from \textbf{2.1} that in 1899 Hilbert thinks of logic as an ultimate foundation of geometrical theories, which are treated on equal footing  with physical theories and theories of other sorts. This way of founding geometry on logic allows Hilbert to avoid in foundations of geometry any appeal to geometrical intuition (or at least this is his intention). When the foundations of Euclidean geometry are remade by Hilbert with a symbolic logical calculus the intuition comes back and its foundational role becomes explicit but this concerns only the \emph{symbolic} rather than the properly geometrical intuition. Now under Tarski's topological interpretation the propositional fragment of logic itself appears to be geometrical in a sense. (The same observation can be made already on the basis of using diagrams in logic.) This suggests that some geometrical (and in particular topological) intuitions may after all play a role in foundations too. 
\newpage

\chapter{Lawvere: Pursuit of Objectivity} 
In a paper of 2003 Lawvere makes the following general remark about foundations of mathematics and Axiomatic Method:
\begin{quote}
In my own education I was fortunate to have two teachers who used the
term ``foundations'' in a common-sense way (rather than in the speculative
way of the Bolzano-Frege-Peano-Russell tradition). [..] The orientation
of these works seemed to be ``concentrate the essence of practice and in turn
use the result to guide practice''. I propose to apply the tool of categorical
logic to further develop that inspiration.\\
Foundations is derived from applications by unification and concentration,
in other words, by the \emph{axiomatic method}. Applications are guided by
foundations which have been learned through education.
( \cite{Lawvere:2003}, p. 213, italic is author's) 
\end{quote}

The author's attitude to foundations of mathematics, which he describes as commonsensical, assumes a permanent interaction between the foundations and the current mathematical practice (``applications''). It is opposed in this sense to the ``Bolzano-Frege-Peano-Russell tradition'' (making part of today's Analytic philosophy), which attempts to provide Formal Axiomatic Method with a philosophical underpinning disregarding its problematic status in the context of the current mathematical research. Lawvere's notion of Axiomatic Method is very general; his more specific proposal concerns the \emph{categorical logic}, which is an area of logic invented by Lawvere himself \cite{Bell:2005}, \cite{Marquis&Reyes:2012}.

Here the first time in this book we meet the mathematical notion of \emph{category}, which is going to play an important role in what follows. I shall introduce this notion twice in this book: in this current Section and then more geometrically in  \textbf{6.3} below; the original motivation of category theory is (partly) explained and in  \textbf{8.5}. Here I allow myself to reproduce the introduction of the notion of category, which Lawvere designed for a philosophical reader in his \cite{Lawvere:1969}. This particular introduction is of special interest for us because it appears to me that in this philosophical publication Lawvere makes an attempt to use a new original axiomatic approach, which at least at this time (1960-70ies) he did not use in his properly mathematical papers. 

\begin{quote}     
The formalism of category theory is itself often presented in ``geometric'' terms. In fact,
to give a category is to give a meaning to the word \emph{morphism} and to the {commutativity}
of diagrams like

$$\xymatrix{A\ar[r]^f & B} \xymatrix{&B\ar[dr]^g\\ A\ar[ur]^f\ar[rr]_h && C} \xymatrix{A \ar[r]^f \ar[d]_a & B\ar[d]^b \\ A' \ar[r]_g  & B'}$$

which involve morphisms, in such a way that the obvious associativity and identity conditions
hold, as well as the condition that whenever

$$\xymatrix{A\ar[r]^f & B}, \xymatrix{B\ar[r]^f & C}$$

are commutative then there is just one $h$ such that

$$\xymatrix{&B\ar[dr]^g\\ A\ar[ur]^f\ar[rr]_h && C}$$

is commutative. \\ To save printing space, one also says that $A$ is the \emph{domain}, and $B$ the \emph{codomain} of $f$
when
$$\xymatrix{A\ar[r]^f & B}$$

is commutative, and in particular that $h$ is the composition $f.g$ if

$$\xymatrix{&B\ar[dr]^g\\ A\ar[ur]^f\ar[rr]_h && C}$$

is commutative. We regard \emph{objects} as co-extensive with identity morphisms, or equivalently
with those morphisms which appear as domains or codomains. As usual we call a
morphism which has a two-sided inverse an \emph{isomorphism}. (\cite{Lawvere:1969}, p. 283)
\end{quote}  
\footnote{Given morphism $\xymatrix{A\ar[r]^f & B}$  its two-sided inverse is morphism $\xymatrix{B\ar[r]^g & A}$ such that $f.g = A$ and $g.f = B$.}

Although Lawvere in the first phrase of the above quote uses inverted commas around the word ``geometrical'', I suggest that he thinks about the geometric aspect of categorical diagrams seriously. Notice that the concept of category is introduced here right from the stretch as in the case of a ``purely formal'' axiomatic introduction of this concept by usual symbolic means (like one found in the beginning of Lawvere's \cite{Lawvere:1966a}, see also \cite{Eilenberg&MacLane:1945} of 1945 where the mathematical notion of category are first introduced by Eilenberg and MacLane). So this ``geometrical'' introduction of the categorical formalism is to be contrasted with the standard introduction by means of symbolic sogic (i.e., by the Formal Axiomatic Method in the precise sense of this expression specified above). Lawvere's remark that with a help of the above diagrams one is supposed to ``give a meaning'' not only to the term ``morphism''  but also to the term ``commutativity'' suggests that Lawvere treats this latter term also as primitive (albeit ``primitive'' should not be understood here in the formal logical sense).  The usual definition of commutative diagram is this: a given categorical diagram is commutative when for all pairs of its nodes $A, B$ all morphisms between these nodes obtained by composition of composable morphisms  shown at this diagram are equal; for example the commutativity of the square diagram shown above amounts to the equality  $f.b = a.g$. Lawvere's idea, if I understand him correctly, is rather to postulate the commutativity for simple cases like $\xymatrix{A\ar[r]^f & B}$ and then transport it geometrically to more complicated cases like the square diagram and beyond. 

Although this introduction of the concept of category addressed to philosophers is not sufficient for a rigorous treatment of the concept of \emph{adjoint functors} discussed further in the same paper, I believe that this is not a merely informal description of the notion of category either. Even if it does not qualify yet as a new form of rigorous Axiomatic Method it gives an idea how such a new method may look like.  A crucial characteristic feature of the new method is the involvement of \emph{geometry} at the foundational level and its dialectical interplay with logic. In \textbf{4.9} we shall see how this dialectics between logic and geometry develops in a more involved and more rigorous context of topos theory. 

The axiomatic introduction of the notion category found in the very beginning of Lawvere's Ph.D. thesis \cite{Lawvere:1963} has the same geometrical flavour: 
\begin{quote}
[In a category w]e identify objects with their identity maps and we regard a diagram
 $$\xymatrix{A \ar[r]^f & B}$$
as a formula which asserts that $A$ is the (identity map of the) domain of $f$ and that $B$ is
the (identity map of the) codomain of $f$. Thus, for example, the following is a universally
valid formula
 $$\xymatrix{A \ar[r]^f & B}\Rightarrow\xymatrix{A \ar[r]^A & A}\wedge\xymatrix{A \ar[r]^f & B}\wedge\xymatrix{B \ar[r]^B & B}\wedge Af = f = fB$$
\end{quote}

What is remarkable here is the way in which Lawvere combines the diagrammatic and the logical symbolic notation. The diagram 
$$\xymatrix{A \ar[r]^f & B}$$
in the category theory usually represents a mathematical \emph{object}, namely a particular morphism $f$ in some category. I assume that Lawvere as anybody else uses this diagram in this sense too (even if he may understand the notion of being a particular object in his proper way). However he also reads it as an \emph{assertion} and combines it with standard propositional connectives. ``Officially'' this is only an unusual symbolic convention, which does not change the sense of the matter. However in fact it touches upon the core of Formal Axiomatic Method. Remind from \textbf{2.4} that this method assumes a distinction between ``real'' and ``ideal'' mathematical constructions: only symbolic constructions qualify in this sense as real while all their interpretations qualify as ideal. As I have stressed earlier this distinction creates a gap between the formalized and the ``real'' (i.e., informal) mathematics because manipulations with ``real'' symbols in formal theories represent some logical operations but not operations with the ``ideal'' objects themselves (\textbf{3.3}). I would like to stress now that Hilbert's distinction between real and ideal mathematical objects is built into the Formal Axiomatic Method technically but is not a matter of philosophical interpretation of this method. A user of this method cannot ignore this distinction even if he or she is not inclined to describe it in Hilbert's original terms. More commonly this distinction is described today as the distinction between the formal syntax and its informal semantics. 

Let me repeat a part of my earlier argument using a simple example. Hilbert in \cite{Hilbert&Bernays:1934-1939} uses symbolic expression $Zw(x, y, z)$ for denoting a predicate saying that given point $y$ lies \emph{between} given points $x, z$; here the symbolic expression itself makes part of the formal syntax of the given theory (formalized Euclidean geometry) and what this symbolic expression stands for makes part of the informal semantics of this theory. As soon as values of $x, y, z$ are fixed $Zw(x, y, z)$ expresses a proposition. Now let us tentatively identify $Zw(x, y, z)$ with a geometrical \emph{object} (construction), which makes the corresponding proposition true, namely with a triple of points $<x, y, z>$ such that $y$ lies \emph{between} $x$ and $z$. So we get in the same parcel, first, an object, and, second, a true proposition ``about'' this object - just like in the case of expression  $f: A \rightarrow B$ read both as a particular morphism $f$ and an (asserted) proposition saying that $A$ is the domain and $B$ is the codomain of $f$. 

In some special cases such a constructive interpretation of Hilbert's formalism seems to work. Consider formula $Zw(x, y, z) \rightarrow Zw(z, y, x)$ and read it, first, as intended (i.e., as a logical implication) and second, as a description of the geometrical operation, which turns this geometrical construction
$$\xymatrix{X\ar@{-}[r]&Y \ar@{-}[r]&Z}$$
into that
$$\xymatrix{Z\ar@{-}[r]&Y \ar@{-}[r]&X}$$
by permuting endpoints $X, Y$. In this particular case indeed there is a structural similarity (which can be described as a precise isomorphism if one likes)  between operations with symbols $x, y, z$ in Hilbert's formulas and operations with symbols $X, Y, Z$ making part of the traditional geometrical notation used together with traditional geometrical diagrams.

However such a geometrical interpretation of formulas obviously does not extend to the whole of Hilbert's formalized Euclidean geometry. Notice that formula $Zw(x, y, z)$ is meaningful even if it expresses a false proposition; in such cases we still have a symbolic construction but have no corresponding geometrical construction. (Mutatis mutandis this remark applies to  Lawvere's notation: unless the diagram $A \rightarrow B$ is \emph{commutative} it does not represent any actual morphism. This is why Lawvere in the above quote reads this diagram as an \emph{assertion} but not as a mere proposition.) Further, if we consider a bit more complex formulas like this one

 $$\forall x \forall y \forall z (Zw(x, y, z) \rightarrow Zw(z, y, x))$$
 
 (which under the intended interpretation says that $Zw(x, y, z) \rightarrow Zw(z, y, x)$ is universally valid) we find no obvious geometrical interpretation for symbol $\forall$, and hence no geometrical counterpart of the extension of formula $Zw(x, y, z) \rightarrow Zw(z, y, x)$ with prefix $\forall x \forall y \forall z$ . Clearly if we consider formal deductions we lose any structural similarity between symbolic constructions, which represent these deductions, and geometrical constructions making part of the informal intended interpretation of the given theory. So in order to make sense of this theory we should do exactly what Hilbert asks us to do: consider only the formulas as real mathematical constructions and treat geometrical constructions and geometrical operations appearing in the intended informal interpretation of this theory as a metaphorical \emph{fa\c{c}on de parler} about ideal mathematical objects and their relations. 

Since, as I have already argued, mathematicians in their actual practice generally do not tend to reduce mathematical constructions to symbolic constructions they develop a mathematical notation called by some philosophers (but rarely by mathematicians themselves) \emph{informal} or \emph{semi-formal}. This latter sort of notation just like the traditional geometrical notation helps one to describe mathematical constructions in terms of certain symbolic and diagrammatic constructions and does not, generally, require making difference between ``real'' and ``ideal'' mathematical objects (although such a difference may eventually appear in some more specific contexts, for example, if one treats real numbers as real and imaginary numbers as ideal). It is often assumed that the semi-formal presentation of mathematical theories is essentially a useful shorthand to a corresponding formal presentation, which can be obtained through some tedious routine procedure of formalization; in particular this is the official position of Bourbaki in his \emph{Elements}. Using this latter example I have argued in \textbf{3.2} that this position is not tenable (even if the formalization is workable) because the semi-formal presentation allows for manipulating with mathematical objects other than symbols while the purely formal presentation does not; as a part of the same argument I have also argued in \textbf{3.3} that manipulations with mathematical objects other than symbols is an epistemically (but not only pragmatically) significant aspect of mathematics. So even if a purely formal presentation of mathematical theories is possible this presentation does not presents an essential aspect of these theories. Beware that talking about the formal presentation I mean here a presentation made by the Formal Axiomatic Method in the precise Hilbert's sense. My argument does not concern the issue of symbolic presentation as such. I do not deny that any mathematical content can be presented in a symbolic form - even if I do not think that such a presentation is always better than a presentation that combines symbols with diagrams and some prose.

Commutative diagrams first appeared in mathematics as a part of semi-formal notation\footnote{
According to MacLane \cite{MacLane:1998}, p. 29, commutative diagrams were first used in early 1940ies by W. Hurewicz in topology some time before the official birth of category theory in 1945 \cite{Eilenberg&MacLane:1945}} and are commonly used today (both within and outside the category theory proper) as such. Lawvere's attempt to use commutative diagrams as a part of \emph{formal} syntax is a bold attempt to bridge the gap between the formal and the semi-formal (aka ``usual'') notation in  mathematics. Instead paying the usual lip service according to which the semi-formal notation can be replaced by a purely formal notation ``in principle'' Lawvere does not want to tolerate the gap between the ``official'' foundations and the practice.  Clearly the symbolic convention according to which commutative diagrams stand for assertions, does not solve the problem. However as we shall shortly see Lawvere's contribution into Axiomatic Method involves much more; in particular, we shall see that the idea of geometrical interpretation of logical quantifiers, which sounds absurd in the context of Hilbert's approach, becomes a part of Lawvere's novel approach. A central role in Lawvere's approach is played by \emph{categorical logic} which he designs as a tool of axiomatic thinking that helps to ``concentrate the essence of practice and in turn use the result to guide practice'' (\cite{Lawvere:2003}, p. 213)

Before we discuss Lawvere's notion of categorical logic I would like to mention two Lawvere's achievements, which do not involve the categorical logic proper but can be described as unusual applications of the usual Formal Axiomatic Method. This ``classical'' aspect of Lawvere's work is equally important for our analysis.

\section{Elementary Theory of the Category of Sets}
An important Lawvere's achievement made wholly within the standard Formal Axiomatic Method is his Elementary Theory of the Category of Sets (ETCS) first presented in his thesis in 1963 \cite{Lawvere:1963}, \cite{Lawvere:2004} and published the next year as a separate paper \cite{Lawvere:1964} \footnote{A longer version \cite{Lawvere:2005} of this paper dating back to 1964 has been recently republished with new author's commentaries.}. The category of sets $S$ is category having (all) sets as \emph{objects} and (all) functions between these sets as \emph{morphisms}. If a universe of sets is given in advance (say, through the axiomatic theory of sets $ZF$) category $S$ can be understood as a specific way of thinking about this universe. (Beware that the totality of \emph{all sets} is a proper class but not a set, see \textbf{5.8}.) Lawvere's ETCS reverses this order of ideas and introduces $S$ axiomatically using, first, MacLane and Eilenberg's axioms for a general category and, second, some additional axioms which distinguish $S$ among other categories up to categorical equivalence (which is an equivalence relation weaker than that of isomorphism, see \textbf{6.7}). The major difference between ZF and other similar axiomatic theories of sets, on the one hand, and ETCS, on the other hand, lies in the choice of primitive non-logical constants. While ZF and other similar theories take the binary relation of  \emph{membership} between sets as primitive (usually denoted by symbol $\in$) ETCS takes as primitive the binary operation of \emph{composition} of functions (which is a ternary relation between functions). 

As McLarty and Lawvere rightly stress \cite{Lawvere:2005} this way of thinking about sets has important advantages. Surely ETCS qualifies as an important example of what is colloquially called the \emph{categorical thinking}. I would like to stress, however, that ETCS is built by exactly the same method as ZF, namely by the standard Formal Axiomatic Method, so in this particular case the new categorical thinking wholly complies with the earlier established method of building axiomatic mathematical theories. No \emph{categorical logic} is involved in the ETCS in its original form. The theorem according to which the axioms of ETCS specify the category of sets up to categorical equivalence (see \textbf{6.7}) is explicitly described by Lawvere as a ``metatheorem'' - once again in full accordance with the Formal Axiomatic Method in its canonical form.  

\section{Category of Categories As a Foundation}

Lawere's axiomatic theory of Category of Categories As a Foundation (CCAF), which equally stems from his thesis and which has been presented in a separate paper in 1966 \cite{Lawvere:1966a}, like ETCS rests on a formal axiomatic theory of general categories; this time Lawvere writes down the appropriate axioms explicitly once again using the Formal Axiomatic Method without trying to modify it. This background part of CCAF Lawvere calls the  \emph{elementary theory} (ET).  After the introduction of the axioms of ET and providing some definitions Lawvere says:

\begin{quote}By a category we of course understand (intuitively) any structure which is an interpretation of the elementary theory of abstract categories, and by a functor we understand (intuitively) any triple consisting of two categories and a rule $T$ which assigns, to each morphism $x$ of the first category, a unique morphism $xT$ of the second category in such a way that ... \end{quote}

(follows the usual definition of \emph{functor} as a structure-preserving map, which I explain and criticize in \textbf{8.6} below). 

ET is a preparatory step towards an extended theory, which Lawvere calls \emph{basic theory} (BT). BT begins with  the introduction of the  \emph{category of categories} and a re-introduction of the notion of functor:

\begin{quote}Of course, now that we are in the category of categories, the things denoted by the capitals will be called categories rather than objects, and we shall speak of functors rather than morphisms.\end{quote}

Remind that in his ETCS paper \cite{Lawvere:1964} Lawvere distinguishes between his theory and the relevant metatheory explicitly.  In the CCAF paper this distinction is blurred. The definition of functor through the structure-preserving rule $T$ is clearly metatheoretical. But the following re-definition of functor as a morphism in a category of categories brings us back to the theoretical (i.e. \emph{elementary}) level. The usual way of handling the difference between a theory and its metatheory in this case is the following. First of all one needs to specify what is meant by ``all'' categories. A natural candidate is the class of all \emph{small} categories, i.e., all categories in which morphisms (including identity morphisms aka objects) form \emph{sets}. Then having a set theory in one's disposal one may consider the class of all set-theoretic models of ET and then using the meta-theoretical notion of functor conceive of a  \emph{large} category \textbf{\emph{C}} of all small categories\footnote{Since the totality of all sets is not a set but a proper class (see \textbf{5.8}) the category of all small categories is not small. In the axiomatic theory of sets called after Quine \emph{New Foundations} (NF) the totality of all sets is a set. However NF turns to be not appropriate for modeling categories, see \cite{McLarty:1991}.}. Since \textbf{\emph{C}} is large and its objects are small \textbf{\emph{C}} is not its own object. All this reasoning is clearly metatheoretical with respect to ET. What we get at the end is a class of set-theoretic models of ET (i.e., the class of small categories), which with a help of the metatheoretical notion of functor (as morphism between small categories) is made into another (not set-theoretical) model of ET. Whether or not the obtained large category of all small categories is legitimate is an instance of the general problem about  legitimacy of large collections like the putative collection of \emph{all sets} (\textbf{8.7}). 

This line of thought relies on a set theory and for this reason is not appropriate for Lawvere's purpose, which is to construe his category of categories \emph{as foundation}. A rationale behind this project is that small categories are certainly not all categories that one typically encounters in the mathematical practice (see \textbf{8.5}). So Lawvere's idea is different: to think of the hypothetical category $CAT$ of all categories as an intended model of ET and then add to ET new axioms which distinguish $CAT$ between other categories; then pick up from $CAT$ an arbitrary object $A$ (i.e., an arbitrary category) and finally specify $A$ as a category by internal means of $CAT$  (stipulating additional properties of $CAT$ when needed). So Lawvere replaces the set-theoretic bottom-up approach outlined in the last paragraph by an original top-down approach. 

More precisely it goes as follows (I omit details and streamline the argument). Stipulate the existence of \emph{terminal object}  \emph{1} in $CAT$, i.e., the object with exactly one incoming functor from each object of $CAT$. Then identify objects (= identity functors) of $A$ as functors in $CAT$ of the form \emph{1}$\rightarrow A$. Stipulate also the existence of initial object \emph{0}, i.e. the object with exactly one outgoing functor into each object of  $CAT$.  Then consider in  $CAT$ object  \emph{2} of the form \emph{0}$\rightarrow$\emph{1}  and stipulate for it some additional properties among which is the following: \emph{2} is a universal generator which means that:

\textbf{G} (generator):  for all $f$, $g$ of the form:
$$\def\dar[#1]{\ar@<2pt>[#1]^f \ar@<-2pt>[#1]_g} \xymatrix{ A \dar[r] & B }$$
and such that $f \neq g$ there exist $x$ such that:    
$$\def\dar[#1]{\ar@<2pt>[#1]^f\ar@<-2pt>[#1]_g} \xymatrix{ \emph{2}\ar[r]^x & A \dar[r] & B }$$

and $xf \neq xg$.

\textbf{U} (universal): if any other category $N$ has the same property, then there are $y$, $z$ such that:

$$\xymatrix{ A \ar@<2pt>[r]^y & B \ar@<2pt>[l]^z }$$

and  $yz = \emph{2}$. 

This allows Lawvere to identify functors (morphisms) of $A$ as functors of the form \emph{2}$\rightarrow A$ in $CAT$. The fact that \emph{2} is the universal generator (it is unique up to isomorphism as follows from the above definition) assures that categories are determined ``arrow-wise": two categories coincide if and only if they coincide on all their arrows. This new definition of functor also allows one to make sense of the notion of a component of a given functor of the form $h$: $A\rightarrow B$ , which in ET is understood as a map $m$ sending a particular morphism $f$ of $A$ into a particular morphism $g$ of $B$ . In BT, $m$ turns into this commutative  triangle:

$$\xymatrix{&\emph{2}\ar[dl]\ar[dr]\\ A\ar[rr] && B}$$ 

Thus categories and functors are no longer built ``from their elements" but rather ``split into" their elements when appropriate. Although the notion of functor as a structure-preserving map can be recovered in this new context it no longer serves for defining the very notion of functor. Rule $T$ used by Lawvere for defining functors in the \emph{elementary theory} disappears in BT without leaving any trace.

Further consider this triangle which Lawvere denotes \emph{3}:
  
$$\xymatrix{&\emph{0}\ar[dr]\ar[dl]\\ \emph{1}\ar[rr] &&\emph{2}}$$

(It should satisfy a universal property which I omit). \emph{3} serves for defining composition of morphisms in our ``test-category" $A$ as a functor of the form \emph{3}$\rightarrow A$ in $C$. Finally, in order to assure the associativity of the composition Lawvere introduces category \emph{4}, which is pictured as follows:

$$\xymatrix{&\emph{3}&\\ \emph{0}\ar[ur]\ar[rr]\ar[dr] &&\emph{2}\ar[ul]\\&\emph{1}\ar[ur]\ar[uu]&}$$ 

(The associativity concerns here the path  $\emph{0}\rightarrow\emph{1}\rightarrow\emph{2}\rightarrow\emph{3}$.)

This construction provided with appropriate axioms makes A into an ``internal model" of ET in the following precise sense: If $F$ is any theorem of ET, then ``for all $A$, $A$ satisfies $F$" is a theorem of BT.

It must be mentioned that this Lawvere's work contains a technical flaw that has been noticed by Isbell in his review  \cite{Isbell:1967}. This flaw has been later fixed, in particular, by McLarty \cite{McLarty:1991} who also provides some additional clarifications on CCAF, which I use in the following discussion. 

Let us see once again what we get here. We have a formal first-order theory BT and its hypothetical intended model $CAT$. Lawvere and McLarty describe $CAT$ as a ``metacategory'' , and McLarty quite rightly, in my view, stresses that ``whether there are useful axioms on $CAT$ making $CAT$ an object in itself'' still remains an ``open question'' (\cite{McLarty:1991}, p. 1259). That BT does not make $CAT$ into an \emph{object} is hardly surprising. Similarly, Hilbert's  \emph{Foundations of Geometry} \cite{Hilbert:1899} does not make Euclidean \emph{space} into an object although it introduces points, straight lines and planes as (primitive) objects. With respect to the given axiomatic theory the Euclidean space is the ``system of things'', i.e., an \emph{universe of discourse} and in this sense it has a \emph{metatheoretical} status; $CAT$ is metatheoretical with respect to BT in the same sense. So far BT works just like any other formal axiomatic theory and so fully complies with the Formal Axiomatic Method.

An unusual feature of BT concerns the fact that it contains ET as its part. ET is sufficient for the axiomatic introduction of the notion of category: any ``system of things'' that satisfy ET is a category by definition. Why one needs more axioms in category theory? An obvious answer is this: ET is not categorical, i.e., it does not define a category up to isomorphism; there are many different categories that one wishes to study but not just one. We have already seen in \textbf{3.2} that the usual axiomatic group theory has the same property and that in order to make this latter axiomatic theory useful one needs to use it together with set theory, which serves  for building and handling its various models. CCAF assumes no such background metatheory. Instead it upgrades ET to BT with additional axioms and so determines a category of a special sort, namely $CAT$ - in a sense of ``determines'' that does not implies that $CAT$ is an object as just explained \footnote{BT does not determine $CAT$ up to isomorphism (or  in any other relevant sense) \emph{uniquely} \cite{McLarty:1991} but let me not bother about this problem now}. Although $CAT$ is not itself an object but a ``system of things'' one must expect that those ``things'', i.e., elements of $CAT$, namely categories and functors, are determined through the theory of $CAT$ (i.e., through BT) \emph{as objects}. This is exactly what BT claims to achieve. This means, however, that BT significantly strengthens the concept of category with respect to its usual definition through ET. To be a category in the sense of being an object of $CAT$ assumes more axioms than to be a category in the sense of being a model of ET. Let us see what this strengthening exactly amounts to. 

There are two ways of thinking about the additional axioms of BT, which complement each other. First, these additional axioms  pick some \emph{specific} notion of category. They strengthen the notion of category in the same sense in which the definition of isosceles triangle strengthen the definition of triangle. Clearly that such a strengthening by itself has no bearing on the foundations of the category theory or the foundations of mathematics in general. Second, these additional axioms allow for a \emph{re-introduction} of basic category-theoretic concepts like \emph{object} (as a functor of the form $1 \rightarrow A$) , \emph{morphism} (as a functor of the form $2 \rightarrow A$) and the rest. This second sense of ``strengthening'' certainly does have a bearing on foundations and is of our special interest in our discussion of the Axiomatic Method. Indeed although BT is construed as a formal theory in full accordance with the Formal Axiomatic Method we see that within this theory (given its intended interpretation $CAT$) Lawvere develops a  \emph{different} way of introduction primitive concepts, which essentially involves $CAT$ and hence some previously assumed category-theoretic constructions. So we have here a foundation involving two levels: the level of elementary theory ET and the level of description of the same primitive categorical concepts by categorical means made available through BT. I would like to stress that the distinction between these two foundational levels does not coincide with the standard distinction between a formal theory and its metatheory. Although Lawvere describes $CAT$ as metacategory he does not provide it with a metatheory. McLarty \cite{McLarty:1991} provides some elements of such a metatheory by considering some candidates for $CAT$ but it is not relevant to my present point. The second-level axiomatic foundation I am talking about from a formal viewpoint is fully accounted for by the extension of ET to BT. So from a formal point of view it amounts to a particular interpretation of this extension. Still this particular interpretation is explicit in Lawvere's paper and it amounts to determining basic categorical concepts by categorical (rather than standard formal ) means. I would like to stress however that these categorical means do not qualify as elements of  \emph{categorical logic}, which I consider below. 

It is appropriate to ask whether or not this latter ``purely categorical'' way of laying out axiomatic foundations may work independently from the Formal Axiomatic Method. Surely it does not work in such an independent way in Lawvere's CCAF paper \cite{Lawvere:1966a}; McLarty's development of Lawvere's idea \cite{McLarty:1991} does not go in this direction but rather cleans it up from the point of view of the standard Formal Axiomatic Method. (McLarty  after Lawvere challenges the set-theoretic foundations built with this method but does not challenge the method itself.) In \textbf{4.9} I argue that in another work Lawvere does use a new axiomatic method that I shall call the New Axiomatic Method; in \textbf{9.3} I describe this New Method in general terms.

\section{Conceptual Theories and their Presentations}
In two abstracts \cite{Lawvere:1966b}, \cite{Lawvere:1967} published in 1966 and 1967 correspondingly Lawvere presents an improved version of yet another idea first developed in his thesis \cite{Lawvere:1963}, \cite{Lawvere:2004}, namely, the idea of presenting a formal theory (more precisely every first-order theory with equality) as a category $T$ of a special sort and presenting models of $T$ as functors from this theory to the category of sets (i.e., functors of the form $T \rightarrow S$). Terms and formulas in this setting are specific morphisms in $T$, truth-values form a distinguished object $L$ in $T$, and  (every instance of) the existential quantifier $\exists$ is also a morphism in $T$. In this context Lawvere makes a distinction between a non-interpreted theory $T$ and its syntactic presentation $\mathcal{L}$:

\begin{quote}
Given a first-order language-with-axioms $\mathcal{L}$, the associated theory $T$ may be thought of as the ``\emph{Sinn}'' [= meaning] of $\mathcal{L}$ and the category $S^{[T]}$ as the ``\emph{Bedeutung}'' [= reference] of $\mathcal{L}$. (\cite{Lawvere:1966b}, p.295)
\end{quote}
 
referring to Frege's \cite{Frege:1892}, which we discuss in \textbf{5.10 - 5.11}. $S^{[T]}$ is the category of models of $T$, i.e., of specific functors of the form $T \rightarrow S$. No details are provided about the ``association'' with a given language $\mathcal{L}$ an appropriate theory $T$. 

The following discussion aims at placing Lawvere's ideas into the historical context discussed in earlier Chapters. For a more technical exposition I refer the reader to \cite{Marquis&Reyes:2012}.

Lawvere's ``language-with-axioms'' $\mathcal{L}$ is a standard formal theory like ZF formulated with a list of axioms written down with usual symbolic means. So we have here a standard Formal Axiomatic setting with formal theories, on the one hand, and set-theoretic models of these theories, on the other hand. Now Lawvere introduces into this general framework a third intermediate element, which he calls after Frege \cite{Frege:1892} the \emph{meaning} of $\mathcal{L}$ and identifies this  \emph{meaning} (rather than $\mathcal{L}$ itself) with a  non-interpreted theory.  As a later Lawvere's comment makes it clear $\mathcal{L}$ in this context is to be thought of as a particular \emph{presentation} of the corresponding theory $T$; $T$ itself is thought of as a \emph{conceptual} (rather than formal) theory, which remains invariant under changes of its various presentations:  

\begin{quote}
Since in practice many abstract concepts [...] arise by means other than presentations, it is more accurate to apply the term ``theory'', not to the presentations as had become traditional in formalist logic, but rather to the more invariant abstract concepts themselves which serve a pivotal role, both in their connection with the syntax of presentations, as well as with the semantics of representations. ( \cite{Lawvere:2004}, p. 8)
\end{quote}

Lawvere's conceptual theories like usual formal theories include a core part, which can be identified as a pure \emph{logic}. Although in \cite{Lawvere:1966b} Lawvere does not provide any explicit criterion for distinguishing between logical and non-logical elements of a theory it is clear that the notions of existential quantifier and truth-value belong to the former category. Saying that Lawvere in this paper \emph{interprets} logic in category-theoretic terms is actually in odds with Lawvere's way of thinking about this matter: Lawvere would rather say that the category-theoretic terms allow for formulating logic concepts in an invariant form, which does not depend on this or that symbolic presentation of these concepts. Nevertheless I shall use this external language because it helps me to compare Lawvere's categorical logic with Tarski's topological interpretation of propositional logic mentioned above (\textbf{3.4}). 

Remind that Tarski \cite{Tarski:1956} provides a topological interpretation of Classical and Intuitionistic propositional logic, which puts into one-one correspondence logical operations, on the one hand, and some operations with elements of a topological space, on the other hand. Lawvere's categorical interpretation of logic works similarly but it is more powerful because it interprets more logic. First, it works for the first-order logic but not only for propositional logic. (Here the crucial role plays Lawvere's interpretation of logical quantifiers as functors, which we shall discuss shortly.) Second, it also interprets truth-values (\emph{true}, \emph{false}), which in Tarski's topological interpretation remain uninterpreted. These features allow Lawvere to think about his categorical interpretation of logic (also called a \emph{categorical semantics} of logic by some authors \cite{Awodey&Bauer:2009}) as a genuine \emph{conceptual formulation} of logic and to call it \emph{categorical logic} rather than by any other name. The following quote from Lawvere's textbook \cite{Lawvere&Rosebrugh:2003} co-authored with Rosebrugh shows that Lawvere's categorical logic has a philosophical but not only technical aspect; the same passage provides an informal explanation of Lawvere's interpretation of logical quantifiers as functors:   

\begin{quote}
The term ``logic'' has always had two meanings - a broader one and a narrower one:\\
(1) All the general laws about the movement of human thinking should ultimately
be made explicit so that thinking can be a reliable instrument, but\\
(2) already Aristotle realized that one must start on that vast program with a more
sharply defined subcase.\\
The achievements of this subprogram include the recognition of the necessity of
making explicit\\
(a) a limited universe of discourse, as well as\\
(b) the correspondence assigning, to each adjective that is meaningful over a whole
universe, the part of that universe where the adjective applies. This correspondence
necessarily involves\\
(c) an attendant homomorphic relation between connectives (like and and or) that
apply to the adjectives and corresponding operations (like intersection and
union) that apply to the parts ``named'' by the adjectives.\\
When thinking is temporarily limited to only one universe, the universe as such
need not be mentioned; however, thinking actually involves relationships between
several universes. [..] Each suitable passage from one universe of discourse to another induces\\
(0) an operation of substitution in the inverse direction, applying to the adjectives
meaningful over the second universe and yielding new adjectives meaningful
over the first universe.
The same passage also induces two operations in the forward direction:\\
(1) one operation corresponds to the idea of the direct image of a part but is called
``existential quantification'' as it applies to the adjectives that name the parts;\\
(2) the other forward operation is called ``universal quantification'' on the adjectives
and corresponds to a different geometrical operation on the parts of the first
universe.\\
It is the study of the resulting algebra of parts of a universe of discourse and
of these three transformations of parts between universes that we sometimes call
``logic in the narrow sense''. Presentations of algebraic structures for the purpose
of calculation are always needed, but it is a serious mistake to confuse the arbitrary
formulations of such presentations with the objective structure itself or to arbitrarily
enshrine one choice of presentation as the notion of logical theory, thereby obscuring
even the existence of the invariant mathematical content. In the long run it is best
to try to bring the form of the subjective presentation paradigm as much as possible
into harmony with the objective content of the objects to be presented; with the
help of the categorical method we will be able to approach that goal. ( \cite{Lawvere&Rosebrugh:2003}, p. 193 - 194)
\end{quote}

Notice that Lawvere and Rosebrugh put at the first place of their description of ``logic in the narrow sense'' the stipulation of ``limited universe of discourse''.  So the authors definitely associate themselves with the logical tradition, which originates in Leibniz and then develops through Boole's pioneering works in symbolic logic, through De Morgan, Jevons, Venn and Peirce. Answering Schr\"oder's criticism Frege  calls logic developed in this tradition ``\emph{calculus ratiocinator}'' (calculus or reasoning) and opposes it to another tradition (to which he adheres himself), where logic is understood as ``\emph{lingua characteristica}'' aka \emph{characteristica universalis} (universal language of thought); this latter tradition also dates back to Leibniz but then develops relatively independently by Frege, Peano, Russell, Quine and Church \cite{Heijenoort:1967},\cite{Hintikka:1997a}, \cite{Hilplinen:2004}. While the \emph{calculus ratiocinator} tradition applies logic ``locally'' leaving it up to the user to determine the universe of discourse in every concrete application, the \emph{characteristica universalis} tradition tends to apply logic to the fixed metaphysical universe that is supposed to include \emph{all} that there is
\footnote{
\begin{quote}
For Frege there is [..] only one possible \emph{Begriffsschrift}, for there is only one kind of human thinking it must reflect. Frege's \emph{Formelsprache} is not a particular development beyond our ordinary language; it is a purified and streamlined version of the entire ordinary language itself. It is calculated to replace ordinary language, at least in its mathematical uses, not to extend it. [..]
This syndrome of ideas characterizes what I have called \emph{language as the universal medium}. Since the meanings (references) of the expressions of our language cannot be expressed in that language we cannot rationally consider varying them, either, at least not in a way that could be specified in language and theorized about. In this sense, our language cannot be reinterpreted. Hence all model theory of our actual language is impossible, for the basic idea of all model theory is precisely to let the interpretation of the language vary. [..] And hence the meanings of our language cannot be changed, it can be used for one purpose only, viz. to speak of this one actual world of ours. Hence a kind of one-world assumption is implicit in the idea of language as the universal medium.  (\cite{Hintikka:1997a}, p. x-xi)
\end{quote} 
}. Thus Lawvere's categorical logic belongs to the \emph{calculus ratiocinator} tradition rather than the \emph{characteristica universalis} tradition. 

In his book \cite{Venn:1881} already quoted above Venn includes a chapter titled \emph{The Universe of Discourse, and its symbolic representation} where he, in particular, remarks (referring to a famous example of syllogism, which includes a premises ``All men are mortal''): 
 
\begin{quote}
Hence we constantly make assertions about all men
without the slightest intention of being bound by our words
beyond a reference to a comparatively small selection of
mankind. (\cite{Venn:1881}, p. 182)
\end{quote}

and further 

\begin{quote}
All and nothing therefore, in any application of our
formulae, are to be interpreted in accordance with the limits
which we may decide to lay down at the outset of the
particular logical processes in question.(\cite{Venn:1881}, p. 186)
\end{quote}

Since the choice of the universe of discourse is for Venn a matter of ``application'' of logic rather than matter of ``symbolic statement'' of logic, Venn unlike Lawvere describes the  notion of universe as ``extra-logical'' (\emph{ib.}, p. 184). In this context Venn makes an interesting link between the issue of universe of discourse and logical diagrams:

\begin{quote}
It has been said above that this question of the Universe
only arises when we apply our formulas. Now diagrams are
strictly speaking a form of application, and therefore such
considerations at once meet us when we come to make use
of diagrams. I draw a circle to represent X, then what
is outside of that circle represents not-X, but the limits of
that outside are whatever I choose to consider them. (\cite{Venn:1881}, p. 186)
\end{quote}
 
Even if laying down logical diagrams according to given symbolic formulas may indeed qualify as ``a form of application'' of these formulas, the logical diagrams certainly make part of Venn's logic and in this sense are quite unlike \emph{men} and any other external objects, which one may wish to include in one's universe of discourse at a given occasion. So the ``application'' of logical formulas to diagrams is a form of self-application. In other words, the use of logical diagrams makes the notion of universe of discourse in some sense \emph{internal} after all. Thus Venn's use of diagrams in logic can be seen as a step toward Lawvere's position, who considers ``making explicit a limited universe of discourse'' as an essential determination of logic rather than a matter of its application. Tarski's topological interpretation of Classical and Intuitionistic propositional calculi (see\textbf{3.4}) can be also most naturally understood as a specification of the universe of discourse of a given calculus. Beware however that Tarski himself does not consistently think of this construction as a proper part of logic. 

The next necessary determination of logic according to Lawvere and Rosebrugh is classical and can be found already in Boole: given a predicate (adjective) $P(x)$ and an universe $U$ consisting of individuals $x$ such that $P(x)$ is meaningful (i.e. has a definite truth-value) for all $x$, one is also given a part (subclass) of $U$ consisting of such individuals $x$ for which $P(x)$ is $true$. I would like to stress that in the standard set-theoretic foundations of mathematics using formal axiomatic theories of sets like ZF this principle is rendered into the form of an axiom schema of set theory (called \emph{Separation} or \emph{Restricted Comprehension} Axiom Schema) \cite{Fraenkel&Bar-Hillel&Levy:1973} but not as a principle of logic.  This technical point is a part of Lawvere's disagreement with the ``speculative'' ``Bolzano-Frege-Peano-Russell tradition'' in foundations of mathematics, which the standard set-theoretic foundations inherit (\cite{Lawvere:2003}, p. 213). Now we see that this disagreement is neither purely technical nor purely ideological but that it also concerns Lawvere's philosophical understanding of logic as a conceptual \emph{tool} (\emph{calculus ratiocinator}) rather than the \emph{characteristica universalis}, i.e.,  the ``universal medium''  (\cite{Hintikka:1997a}, p. xi). 

The last point of determination of logic according to Lawvere and Rosebrugh concerns the ``attendant homomorphic relation'' between logical connectives, on the one hand, and the ``mereology'' of the corresponding universe of discourse; as along as this universe is thought of as a class, this mereology reduces to a number of operations with subclasses of this class like unions and intersection. The ``attendant relation'' is, of course, well-known since Boole; the structure shared by propositions, on the one hand, and subclasses of a given class, on the other hand, is called after Boole the \emph{Boolean algebra}. Tarski's topological interpretation of propositional calculus (\textbf{3.4}) allows for a similar ``homomorphic relation'' for the case of \emph{intuitionistic} propositional calculus; the corresponding structure is called \emph{Heyting algebra}. I would like once again to stress that Lawvere and Rosebrugh talk here about the ``homomorphic relation'' between logical connectives and the mereology of the universe of discourse as a proper part of logic but not as a model-theoretic issue. 

Whether there exists only one metaphysical universe of everything or one has a freedom to assume as many universes of discourse as one likes may seem to be a philosophical issue that can and arguably must be separated from the logical formalism as such. However as the further explanations of Lawvere and Rosebrugh make it clear such a complete separation of ``application'' of logic from the logic itself is not compatible  with their approach. In fact the assumption about the plurality of universes of discourse is necessary for the most original technical aspect of Lawvere's proposal, namely for his interpretation of logical quantifiers \emph{there exist} and \emph{for all} as functors \emph{adjoint} to \emph{substitution} (also understood as a functor).

An \emph{adjoint situation} (called also an \emph{adjunction}) is a pair of categories $A, B$ with two functors $f, g$ going in opposite directions:

 $$\xymatrix{ A \ar@<2pt>[r]^f & B \ar@<2pt>[l]^g}$$

provided with certain \emph{natural transformations} (i.e., morphisms between functors) and
satisfying certain conditions, which make these functors in some special sense mutually inverse. Given an adjoint situation as above functor $g$ is called \emph{left adjoint} to functor $f$ and functor $f$ is called \emph{right adjoint} to functor $g$, in symbols $g \dashv f$. A given functor has at most one (up to unique isomorphism) left adjoint and one right adjoint. Here is the precise definition
\footnote{The concept of adjunction has been first introduced in 1958 by Kahn \cite{Kan:1958} in the context of algebraic geometry.}.
In addition to functors $f, g$ adjunction $g \dashv f$ comprises natural transformations $\alpha : A \rightarrow fg$ and $\beta : gf \rightarrow B$ such that $(g\alpha)(\beta f) = g$ and $(\alpha f)(g \beta) = f$. (As above I write here the composition in the geometrical order and following Lawvere in \cite{Lawvere:1966a} do not distinguish between objects and their identity morphisms.) The required condition is that the following two triangles commute: 

$$\xymatrix{g \ar[r]^{g\alpha}\ar[dr]_{1_{g}} & gfg\ar[d]^{\beta g}\\ &g}$$

$$\xymatrix{ f \ar[r]^{\alpha f}\ar[dr]_{1_{f}} & fgf\ar[d]^{f \beta}\\ &f}$$

Once again I would like to stress that saying that quantifiers are \emph{interpreted} as adjoint functors provides an external description of Lawvere's discovery; as the above quote from Lawvere and Rosbrugh makes it clear Lawvere himself thinks of this functorial interpretation of quantifiers as the ``objective content'' of the notion of quantification. Let me now  discuss some of this conceptual objective content beginning with the notion of \emph{substitution}; in this discussion I shall continue to trace a historical link between the categorical logic and the earlier \emph{calculus ratiocinator} tradition in logic as represented by Venn \cite{Venn:1881}. 

While Venn thinks of multiple universes of discourse as certain independent domains of application of his logic, Lawvere and Rosebrugh claim that ``thinking actually involves \emph{relationships} between several universes'' (my emphasis) and then talk about the ``passage'' from one universe of discourse to another. In a category-theoretic setting universes are presented as objects $X, Y, ..$ of some category and passages between universes are presented as morphisms like   
 
 $$\xymatrix{X \ar[r]^f & Y}$$
   
(Notice that passages between two given universes are, generally, many and they are, generally, not inversible.)  Thus we have here a \emph{category} of universes of discourse but not just a class of these things. Let us for simplicity think of these universes as \emph{sets} but have in mind that the functorial construction of logical quantifiers does not require this assumption. (We shall see in \textbf{4.5} how this construction works in a more general setting and how in such a more general setting it reveals its distinctive geometrical aspect.) Suppose now that we have a one-place predicate (a property) $P$, which is meaningful on set $Y$, so that there is a subset $P_{Y}$ of $Y$ (in symbols $P_{Y} \subseteq Y$) such that  for all $y \in Y$ $P(y)$ is true just in case $y \in P_{Y}$. Now using these data (together with morphism $f$ as above) we can define a new predicate $R$ on $X$ as follows: we say that for all $x \in X$ $R(x)$ is true when $f(x) \in P_{Y}$ and false otherwise. So we get subset $R_{X} \subseteq X$ such that for all $x \in X$ $R(x)$ is true just in case $x \in R_{X}$. Let us assume in addition that every subset $P_{Y}$ of $Y$ is determined by some predicate $P$ meaningful on $Y$. Then given morphism (``passage'') $f$ from ``universe'' $X$ to ``universe'' $Y$ we get a way to associate with every subset  $P_{Y}$ (every part of universe $Y$) a subset $R_{X}$ and, correspondingly, a way to associate with every predicate $P$ meaningful on  $Y$ a certain predicate $R$ meaningful on $X$. Since subsets of given set $Y$ form Boolean algebra $B(Y)$ we get a map between Boolean algebras (notice the change of direction!): 

$$f^{*}: \xymatrix{B(Y) \ar[r] & B(X)}$$

Since Boolean algebras themselves are categories (with objects subsets and maps inclusions of subsets) $f^{*}$ is a functor. For every proposition of form $P(y)$ where $y \in Y$ functor $f^{*}$ takes some $x \in X$ such that $y = f(x)$ and produces a new proposition $P(f(x)) = R(x)$ (for a single given $y$ it may produce a set of different propositions of this form). Since it replaces $y$ in $P(y)$ by $f(x) = y$ it is appropriate to call $f^{*}$ \emph{substitution} functor. 

The \emph{left} adjoint to the substitution functor $f^{*}$ is functor 

$$\exists_{f} :  \xymatrix{B(X) \ar[r] & B(Y)}$$

which sends every $R \in B(X)$ (i.e. every subset of $X$) into $P \in B(Y)$ (subset of $Y$) consisting of elements $y \in Y$, such that \emph{there exists} some $x \in R$ such that $y = f(x)$; in (some more) symbols 

$$\exists_{f} (R) = \{y | \exists x (y = f(x)  \wedge   x \in R)\}$$

In other words $\exists_{f}$ sends $R$ into its \emph{image} $P$ under $f$. Now if (as above) we think of $R$  as a property $R(x)$ meaningful on $X$ and think of $P$ as a property $P(y)$ meaningful on $Y$ we can describe $\exists_{f}$ by saying that it transformes $R(x)$ into $P(y) = \exists_{f}xP'(x,y)$ and interpret $\exists_{f}$ as the usual  existential quantifier. 

The \emph{right} adjoint to the substitution functor $f^{*}$ is functor
$$\forall_{f} :  \xymatrix{B(X) \ar[r] & B(Y)}$$ 

which sends every  subset $R$  of $X$ into subset $P$ of $Y$ defined as follows: 

$$\forall_{f} (R) = \{y | \forall x ( y = f(x)  \Rightarrow   x \in R)\}$$
 
and thus transforms $R(X)$ into  $P(y) = \forall_{f}xP'(x,y)$.  

Notice that functors  $\exists_{f}$ and $\forall_{f}$ are defined here as adjoints to functor $f^{*}$, i.e., quite independently from their interpretation as logical quantifiers explained above. 

The very fact that that in this setting quantifiers arise ``naturally'' through the functorial adjunction is remarkable from a mathematical point of view. According to Marquis and Reyes ``[t]his was a key observation that convinced many mathematicians that this was the right analysis of quantifiers'' (\cite{Marquis&Reyes:2012}, p.710). 

Let me now focus on philosophical aspects of this Lawvere's achievement.

Lawvere's approach to logicality like Tarski's approach discussed in \textbf{2.2} is semantical because it essentially uses the notion of universe of discourse. Nevertheless it is very unlike Tarski's. While Tarski thinks about logic in a traditional vein as an invariant structure over a given universe of discourse, Lawvere's categorical logic is a device that allows one to \emph{translate} a proposition meaningful in a given universe $Y$ into another proposition meaningful in another given universe $X$ taking into account the relationship between the two universes expressed by morphism $f: X \rightarrow Y$; as we have seen the corresponding ``translation functor'' $f^{*}$ expresses the simple idea of substitution of a given variable $y$ by term $f(x)$. The fact that logical quantifiers are defined in this setting through the substitution shows that the operation of translation between different universes is not just a useful extra feature of logic but its very conceptual basis. The only logical structure that is shared by all universes is the \emph{propositional} (Boolean, Heyting or other) structure. However for the first-order and higher-order logic the translational nature of logic becomes essential. Logic is understood here no longer as a system of universal \emph{forms} of thought, which are not sensitive to differences between various domains of its application, but rather as a universal  \emph{translational protocol}, which allows one to navigate between different domains. This makes Lawvere's conception of logic fundamentally relational: notice that quantifiers applied in domain $X$ are relativized to a particular ``passage'' $f: X \rightarrow Y$ from this domain to certain codomain $Y$; one needs to specify such a passage and the corresponding codomain every time when one talks about ``all $x$ from $X$'', for example, when one says that  \emph{all men are mortal}. Does this make sense?  Let us see. The content of the proposition \emph{all men are mortal} can be expressed by saying that being a man implies being mortal. So this content can be expressed by means of propositional logic without using quantifiers. But when  quantifiers are really needed - as for example in the case of Hilbert's First Axiom (see \textbf{2.6} above):

$$\forall x\forall y Gr(x, x, y)$$

(saying that any two points are aligned) - one should have in mind that talking about ``all points'' does not make sense \emph{within} the universe of points $Pt$ but requires a ``higher'' or, better to say, simply some \emph{external} viewpoint, which allows one to look at $Pt$ from outside and see it as a whole. In the given example such an external viewpoint can be specified rather straightforwardly: since the First Axiom characterizes a geometrical space, one should stipulate an appropriate category $G$ of spaces and make a new independent variable in the above expression to range over these spaces. Then the universal quantifiers in the above axiom can be understood as Lawvere quantifiers; they depend not only on category $G$ but on a chosen functor $f: Pt \rightarrow G$ from points to spaces. From the traditional Hilbertian viewpoint functor $G$ and $f$ belong to the metatheory and so must be treated separately but Lawvere shows us that without these things the first-order logic used in the given formal theory cannot work. Since ``a quantifier is an operation in logic that moves a statement from one context to a related context'' \footnote{http://www.ncatlab.org/nlab/show/quantifier as for April 2012} one cannot use quantifiers working within a single context.   

The informal explanation of the notion of logical quantifier given by  Lawvere and Rosebrugh in the above quote in terms of universes of discourse makes it clear that these universes must \underline{not} be understood as external domains of interpretation of a given theory; as objects of a given logical category these things are proper elements of this category and so they do not provide by themselves an interpretation (semantics) in the model-theoretic sense of the term. In \cite{Lawvere:1969} where Lawvere presents an advanced version of his categorical notion of elementary theory he refers to objects of the given logical category as  \emph{types}. This terminology points to Lawvere's ``1963 observation [..], that cartesian closed categories serve as a common abstraction of type theory and propositional logic'' (\cite{Lawvere:2006}, p.1), which deserves our special discussion. For a systematic study of this subject (that extends far beyond Lawvere's original observation of 1963) see  \cite{Lambek&Scott:1986} and \cite{Jacobs:1999}.

\section{Curry-Howard Correspondence and Cartesian Closed Categories}
The idea of logical calculus that not simply applies to different domains of individuals but explicitly distinguishes between different \emph{types} of individuals dates back to Russell who has coined the term ``theory of types'' (\cite{Russell:1903}, Appendix B); the idea can even be traced further back to Aristotle's distinction between different \emph{genus} of things and his principle according to which switching between different genus in a reasoning (\emph{metabasis}) is not allowed. 
One may remark that the type distinction reflects a feature of our natural languages:
\begin{quote}
Types are inherent in everyday language, for example, when we   distinguish between ``who'' and ``what'' or between ``somebody'' and ``something''. ( \cite{Lambek&Scott:1986}, p. 125
\end{quote}

and further remark that distinguishing between different types of objects is tacitly made in the mathematical practice:

\begin{quote}
In our mathematical practice we have learned to keep things apart. If we have a rational number and set of points in the Euclidean plane, we cannot even imagine what it means to form the intersection.  [..] If we think of a set of objects, we usually think of collecting things of a certain type, and set-theoretical operations are to be carried out inside that type. Some types might be considered as subtypes of  some other types, but  in other cases two different types have nothing to do with each other. That does not mean that their intersection is empty, but that it would be insane to even \emph{talk} about the intersection. (\cite{Bruijn:1995}, p. 31)
\end{quote}

Indeed the standard Set-theoretic foundations of mathematics do not allow for distinguishing between types of objects (at least at the foundational level) and, formally, do allow for crazy set-theoretic operations mentioned in the above quote (since \emph{every} mathematical object is a set one may always form an intersection of two objects). Notice that the type distinction between \emph{points} and \emph{straight lines}, which is made explicit in Hilbert's \emph{Foundations} of 1899, disappears in his \emph{Foundations} of 1934, where the Formal Axiomatic Method takes its mature symbolic form (remind that in the latter case Hilbert treats only points as primitive objects). This is not surprising because the system of (symbolic) logic used by Hilbert is not typed. One could expect that the replacement of the underlying logical calculus by a typed logic may solve the problem without effecting the foundations of the Formal Axiomatic Method itself; however at least one particular development in type theory, namely one that has eventually led to the discovery of the so-called \emph{Curry-Howard correspondence} aka \emph{Curry-Howard isomorphism}, as a matter of fact does touch upon the foundations of the Axiomatic Method and connects the type theory to categorical logic, as we shall now see. 

As Lawvere (\cite{Lawvere:2006}, p.2) and some other people notice the name of \emph{Curry-Howard isomorphism} is misleading because the term ``isomorphism'' is used in it loosely. I would like to stress that it is also misleading in a different sense, namely in the sense that it does not properly reflect the original philosophical motivation behind this discovery. In 1924 Sch\"onfinkel published a paper  \cite{Schonfinkel:1924}, \cite{Schonfinkel:1924a} aiming at deepening Hilbert's formalization of logic, which, as I have already stressed in Chapter 2, indeed does not provide a \emph{purely formal} treatment of logical concepts like \emph{proposition} and \emph{variable}; the sense in which Hilbert treats logical concepts ``formally'' is significantly weaker than the sense in which he treats formally non-logical concepts (\textbf{2.3}). Sch\"onfinkel's idea was to reduce the logical concepts, which so far were generally seen as basic, to a small number of allegedly more fundamental syntactic operations like  substitution and permutation of signs\footnote{ 
\begin{quote}
The successes that we have encountered thus far on the road taken encourage us to
attempt further progress. We are led to the idea, which at first glance certainly
appears extremely bold, of attempting to eliminate by suitable reduction the remaining fundamental notions, those of proposition, propositional function, and variable,
from those contexts in which we are dealing with completely arbitrary, logically
general propositions (for others the attempt would obviously be pointless). To examine
this possibility more closely and to pursue it would be valuable not only from the
methodological point of view that enjoins us to strive for the greatest possible conceptual uniformity but also from a certain philosophic, or, if you wish, aesthetic point
of view. For a variable in a proposition of logic is, after all, nothing but a token
 that characterizes certain argument places and operators as belonging
together; thus it has the status of a mere auxiliary notion that is really inappropriate
to the constant, ``eternal'' essence of the propositions of logic.
It seems to me remarkable in the extreme that the goal we have just set can be
realized also; as it happens, it can be done by a reduction to three fundamental signs. ( \cite{Schonfinkel:1924a} , p. 358-359)
\end{quote}
}
Independently similar ideas inspired Huskell Curry in 1920ies (before he first came across Sch\"onfinkel's paper during the academic year of 1927-1928) who gave to this field of study its current name of  \emph{combinatory logic} \footnote{For more detailed historical accounts see \cite{Cardone&Hindley:2009}, \cite{Seldin:2009}.} Here is how Curry describes the aim and the scope of combinatory logic in a later co-authored monograph:

\begin{quote}
Combinatory logic is a branch of mathematical logic which concerns itself with the ultimate foundations. Its purpose is the analysis of certain notions of such basic character that they are ordinarly taken for granted. These include [(i)] the process of substitution, usually indicated by the use of variables; and also [(ii)] the classification of the entities constructed by these processes into types or categories, which in many systems has to be done intuitively before the theory can be applied.  It has been observed that these notions, although generally presupposed, are not simple; they constitute a prelogic, so to speak, whose analysis is by no means trivial.  (\cite{Curry&Feys&Craig:1958}, p. 2)
\end{quote} 

 Purposes (i) and (ii) mentioned by Curry in the above quotes are mutually dependent
 \footnote{In the above quote Curry uses the term ``category'' interchangeably with the the term ``type'' - and further in his book he uses the former term more often than the latter. Although in the category theory the term ``category'' is used in a different sense the two different uses of this term are compatible as long as one assumes in the spirit of CCAF that an object of a category is, generally, itself a category.}. 
Since a formal logical calculus is seen as a bare symbolic calculus where signs do not have any previously assumed meaning one needs to make explicitly certain distinctions between different types of symbolic constructions without which this calculus cannot qualify as logical - including, in particular, the distinction between individuals, propositions and logical connectives. The idea of Combinatory logic as Curry describes it requires making all such distinctions \emph{formally} without appealing to the usual meaning of words ``individual'', ``proposition'', etc. Thus pushing the formal approach to logic to the extreme shows the necessity of typing, so one may argue that the type distinction is always present in logic whether one describes explicitly or not. As Jacobs puts this 
 
 \begin{quote}
 A logic is always a logic over a type theory. (\cite{Jacobs:1999}, p. 1)
 \end{quote}

 In 1969 William Howard  reformulated and extended Curry's results in a note \cite{Horward:1969-1980} that has been first published only in 1980. Instead of using Combinatorial logic Howard used the formalism of (simply) typed \emph{lambda calculus} invented by Alonzo Church in late 1920s \cite{Cardone&Hindley:2009} and first published in 1933. Curry was of, of course, aware about the fact that the formalism of lambda-calculus comes close to that of combinatory logic but he claimed that his formalism provides a deeper foundational analysis (\cite{Curry&Feys&Craig:1958}, p.6 - 9). Unlike Curry Howard did not stress the philosophical motivation and the foundational significance of this result but formulated it in terms of structural correspondence between two families of formal calculi, namely, the simply typed lambda calculi, on the one hand, and the formal deductive systems, on the other hand. Although such a presentation may have certain advantages for a mathematical reader not interested in foundational issues, it leaves behind the philosophical content of Curry's work and gives a wrong impression that we are dealing here with some unexplained structural similarity rather than with a conceptually clear mathematical fact concerning the foundations of logic.    
 
Lawvere's 1963 insight that the relevant structure shared by type theories and deductive systems can be described as a category of special sort (that is called \emph{cartesian closed}) has been systematically developed by Joachim Lambek in late 1960s and early 1970s, see \cite{Lambek:1968},  \cite{Lambek:1969},  \cite{Lambek:1972}; a systematic analysis of relationships between combinatory logic, lambda-calculus and cartesian closed categories, which also contains some historical notes, is found in Lambek's and Scott's monograph \cite{Lambek&Scott:1986}. For my present purpose it is sufficient to stress that the notion of cartesian closed category (CCC) used by Lawvere for logical purposes connects his work to the earlier attempts to deepen foundations of logic and mathematics made by Sch\"onfinkel, Church, Curry and their collaborators. It is appropriate also to mention that the notion of CCC first appears (without the name) in Lawvere's work as a property of the category of sets \cite{Lawvere:1964} but not as a part of logic. Lawvere's thinking about CCC is made explicit in the following Introduction to his paper \cite{Lawvere:1969a} of 1969 where CCC first appears in press under that name:

\begin{quote} 
 Cartesian closed categories [..] serve as a common abstraction of type theory and propositional logic [..] (\cite{Lawvere:1969a}, p. 134)
 \end{quote} 
 
For Lawvere CCC is not an exclusively \emph{logical} category but rather a category that plays a (central) role in logic in particular. In that sense there is a clear difference between Curry's strategy who pursues ``the ultimate foundations'' of logic by studying general features of syntactic structures and Lawvere's strategy, who uses the category theory for describing ``abstract structures'' found everywhere in mathematics \emph{including} mathematical logic and who does not grant to syntactic structures any special importance. 

The abstraction of CCC, which lies behind the so-called Curry-Howard correspondence, sheds a new light on the \emph{doing versus showing} dilemma, which I have stressed earlier. We have seen that in Euclid these two aspects of mathematical reasoning are intrinsically intertwined - in spite of the fact that Euclid  explicitly qualifies each ``Proposition'' of his \emph{Elements} either as a \emph{problem} or as a \emph{theorem} (using  appropriate endings, see \textbf{1.4}). In Hilbert's Axiomatic Method (in its mature form of 1927 and later) the relationships between \emph{doing} and \emph{showing} are arranged as follows: what we \emph{do} are syntactic constructions and what we \emph{show} by doing these ``real'' syntactic constructions are certain facts concerning ``ideal'' mathematical objects in terms of which we interpret these constructions (\textbf{2.6}). (A different sort of \emph{showing} is reserved for metamathematics; I leave it now aside.) The setting of Curry-Howard Correspondence does not apply Hilbert's rigid distinction between ``real'' and ``ideal'' mathematical objects and allows for \emph{doing} something else than only building formulas: the simply typed lambda calculus  is not a formal deductive machinery but a kit for building mathematical objects of sort like Euclid's geometrical ``calculus'' of constructions with ruler and compass. The Curry-Howard correspondence amounts to the observation that the rules of simply typed lambda calculus can be used (mutatis mutandis) for making formal deductions; the notion of CCC makes the mathematical sense of this correspondence precise. Thus in the given case the relationships between \emph{doing} and \emph{showing} are arranged differently. We no longer assume that we \emph{do} things of one sort (build strings of symbols) and then use these things for \emph{showing} something about things of another sort (by interpreting these strings of symbols as propositions about ``ideal'' mathematical objects); instead we explore the same structure of CCC by \emph{doing} allowed calculations and  \emph{showing} with these calculations certain facts \underline{about this very structure}. 
As we can see this arrangement is similar to Euclid's and dissimilar to Hilbert's. It is appropriate to notice here that CCC is not the only kind of category that allows for this new form of synthetic mathematical reasoning (albeit historically it was invented first). A similar synthetic approach is supported, in particular, by \emph{locally} cartesian closed categories (LCCC), which provide a categorical framework for Martin-L\"of's intuitionistic type theory with dependent types \cite{Seely:1984}, see \textbf{5.11} and \textbf{6.9}.

\section{Hyperdoctrines}
The idea of quantifiers as adjoints to substitution, which I explained in \textbf{4.3} pretending that we were talking about sets, was first mentioned by Lawvere in \cite{Lawvere:1967} and then fully elaborated in the \emph{Dialectica} paper \cite{Lawvere:1969} with a help of the notion of CCC; the categorical construction, which supports the quantification (and in fact the full first-order logic) Lawvere calls a  \emph{hyperdoctrine}. A hyperdoctrine consists of a CCC $T$ of ``types'' and functor $h$ that associates (i) with every object $A$ of $T$ - a category $P(A)$ of ``parts'' of $A$, which in the given context are also thought of as ``predicates'' or ``attributes''\footnote{A given part $P_{A}[\phi]$ consists of those and only those ``elements'' of $A$, which have attribute $\phi$; now this idea is expressed in a form, which doesn't require the reference to elements. This does \emph{not} mean that we get rid of elements altogether but means only that we don't use the notion of element as primitive; in fact this notion can be recovered in the given setting and, moreover, in a generalized form.}, and (ii) with every morphism $f: A \rightarrow B$ - a functor $s_{f}: P(A) \leftarrow P(B)$ (beware the reversal of the arrow!) thought of as ``substitution'' (in the sense explained in \textbf{4.3} but this time in a more general setting). It is assumed that every $P(A)$ is a partial order (parts are partially ordered by inclusion), which is tantamount to saying that it is a category having at most one morphism $P_{A}[\phi] \rightarrow P_{A}[\phi]$ with for given domain and codomain. It is also assumed that $P(A)$ is CCC. This allows for thinking about this partial order as a \emph{deductive} order and denote it $\vdash_{A}$. 

Notice once again how closely this setting reproduces the traditional ideas of British symbolic and diagrammatic logic as presented by Venn \cite{Venn:1881}! An important difference, however, is that in Lawvere's case the order $\vdash_{A}$ is construed \emph{concretely} as a structure on given type $A$ rather than as an abstract structure (like Boolean algebra) responsible for logical deduction, which merely \emph{applies} to $A$ and to any other ``universe of discourse''. This change of viewpoint is the key idea that allowed Lawvere to capture algebraically the first-order logic after the example of the traditional algebraic treatment of propositional logic by Boole, Venn and others. 

As we already know from \textbf{4.3} the universal and the existential quantifiers are recovered as the right and the left adjoints to functor $s_{f}$, so the quantification also turns into a \emph{concrete} structure associated with the given morphism $f$, which maps one given type (the domain of quantification) to another given type (the codomain). The fact that a sound notion of quantification requires not only a specification of its domain $A$ but also a specification of (i) its codomain $B$ and (ii) a particular map $f: A \rightarrow B$ between the domain and the codomain, is a discovery of great philosophical significance, which I have already tried to explain. The reader is advised to re-read now \textbf{4.3} and see that the above discussion concerns the general hyperdoctrine but not only to the hyperdoctrine of sets used then as a convenient example\footnote{
Here is how the category of sets presents itself as a hyperdoctrine. First consider this category as the category $T$ of types (think of each set $A$ as a particular type). Next associate (i) with every set $A$ its powerset (i.e. set of all subsets of $A$) $P(A)$ and (ii) with every function $f: A \rightarrow B$ a new function between powersets $s_{f}: P(B) \rightarrow P(A)$ (substitution) which sends given subset $B'\subseteq B$ to subset $A'\subseteq A$ formed as follows: $a \in A$ is a member of $A'$ if and only if $b = f(a)$ is a member of $B'$. Elements of powerset $P(A)$ (i.e., the subsets of set $A$)
are partially ordered by inclusion, this partial order has the structure of Boolean algebra and the resulting logic is classical.  
}.

Observe that functor $h: T \rightarrow P$ that associates to the category $T$ of types a further structure, is ``something more'' than a ``usual'' contravariant functor, which maps our given category $T$ to another one.  Objects of category $P$ (which I did not describe so far) are categories $P(A)$ (so it is a category of categories) and its morphisms are functors of the form $s_{f}: P(B) \rightarrow P(A)$. Crucially, we cannot disregard here the internal structure of each $P(A)$ and thus treat $P$ as an abstract category of some appropriate kind. As we shall see in \textbf{6.9} $P$ can be described in a more compact way as a \emph{2-category} but now it is more appropriate for my purpose to describe functor $h$ as a \emph{fibration} (in particular, because this is how Lawvere describes it in \cite{Lawvere:1970b}). The notion of fibration comes from geometry: it is a way of ``thickening'' a given geometrical object $G$ by associating with every point $p$ of $G$ a new object $f(p)$ called a \emph{fiber} in a way, that reflects the structure of the base object $G$. For a suggestive example think of a hairy head: if the hair is normally cut (i.e., is neither too short nor too long) one can see how the form of the hair reflects the form of the head without being reduced to it (if the hair is too short one can see only the form of the head, and if it is too long the form of the head becomes invisible). Now we have a similar situation: the hyperdoctrine $h$ transforms every ``thin'' type $A$ into a ``thick'' category $P(A)$ of attributes over this type and every ``thin''  morphism $f: A \rightarrow B$ into a ``thick'' functor $s_{f}: P(B) \rightarrow P(A)$ between the categories of attributes. 

Let me say a few more words about the geometrical notion of fibration, which will be used in Chapter \textbf{6}. 
In accordance to the direction of $h$ I informally described the geometrical fibration as a process of ``thickening'', which is non-trivial in the sense that it adds more structure. In the case of hyperdoctrines it adds a logical structure. However in order to give a geometrical definition of fibration it is more appropriate to describe fibration in terms of the opposite process that maps the ``hairy head'' onto its base ``bold head'' (which is the simplest way of cutting one's hair). Hence the definition of \emph{fiber bundle} (which is a special case of fibration corresponding to intuitive examples like that of hairy head) as a continuous surjective (i.e., \emph{onto}) map $E \rightarrow B$, which \emph{locally} looks like a projection of a product space $B \times F$ onto one of its components $F$. A trivial example of fiber bundle is given by a cylinder $C$ seen as the product $O \times L$ of its base circle over its side (the ``hair'') and the projection $p: O \times L \rightarrow O$ of this cylinder onto its base. (In order to reverse the process think of this cylinder as formed by hairs growing out of its base.) A non-trivial example is obtained from the former one by twisting the side $L$ of the cylinder: this twist produces  a M\"obius strip $M$. Observe that $M$ and $C$ grow from the same base $O$ with the same hair $L$;  however $M$ cannot be (globally) described as a product space and the map $f: M \rightarrow O$ is no longer a projection (Fig. 4.1)

\begin{center}
\includegraphics [scale=0.5]{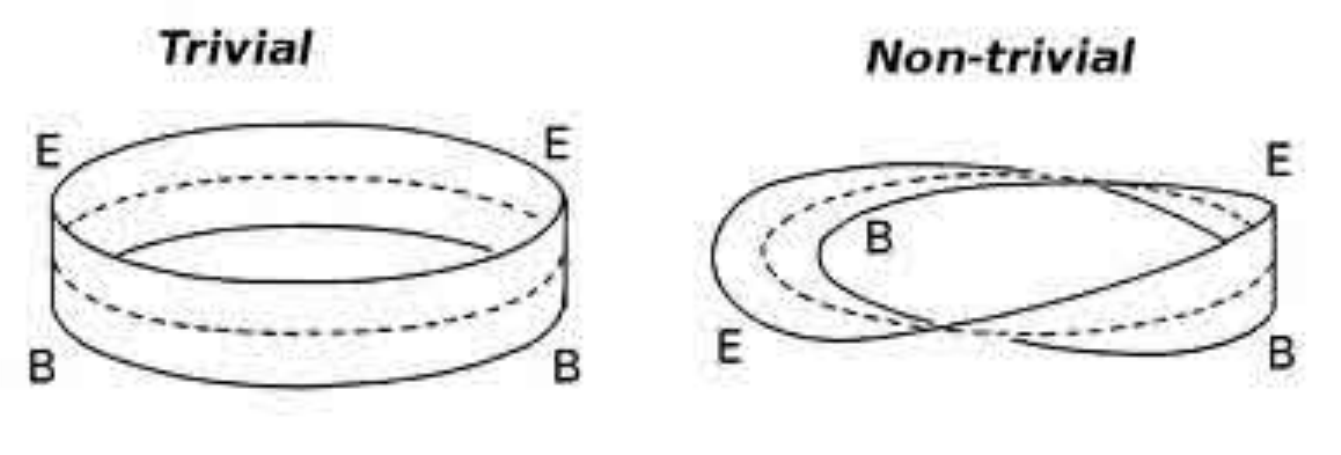}
\end{center}

\begin{center}
Fig. 4.1
\end{center}

The notion of fiber bundle allows for a generalization called \emph{Hurewicz fibration} or simply fibration, which turns to be particularly appropriate in a categorical context as we shall see in \textbf{6.8}. Hurewicz fibration must not be confused with \emph{Grothenideck fibration} discussed in \textbf{6.5}. How these two categorical notions of fibration relate to each other will be also explained in \textbf{6.8}.

The fact the geometrical aspect of the notion of hyperdoctrine is something more than just a useful intuitive way of thinking about it has been demonstrated by Lawvere's work in topos theory \cite{Lawvere:1970a} almost simultaneously with the first public appearance of hyperdoctrines in  \cite{Lawvere:1969} and the second appearance of hyper doctrines in \cite{Lawvere:1970b}; we shall better see how works this dialectical interplay between geometry and logic in  in \textbf{4.9}. 

Although in the late 1960ies and early 1970ies the ``internalization'' of quantifiers with hyperdoctrines was justly commonly perceived as the most impressive result of Lawvere's work, I would like to mention now another Lawvere's achievement, which attracted the common interest more recently. This is the internalization of \emph{equality} (aka \emph{identity} in logic). The idea first appears in \cite{Lawvere:1967} (before the introduction of hyperdoctrines) and then becomes central in \cite{Lawvere:1970b}. Technically it comes in the same parcel with the internalization of quantifiers. Consider in the category $T$ of types (making part of hyperdoctrine $h$ as above) morphism $\delta_{A}: A \rightarrow  A \times A$ which we want to play the role of the diagonal map that sends each ``element'' of type $A$ to its cartesian square (which always exists because $T$ is cartesian closed). Then consider in $P(A)$ the ``maximal'' attribute $1_{A}$, which we interpret as the predicate ``identically true over $A$'' (if $A$ is a set then $1_{A}$ is the maximal subset of $A$, which however we must formally distinguish from $A$ itself); any deduction in $P(A)$ that has the form $1_{A} \rightarrow X$ we interpret as a \emph{proof over A}. Finally consider the adjunction $\exists_{\delta_{A}} \dashv s_{\delta_{A}}$ and the image $id_{A} = \exists_{\delta_{A}}(1_{A})$, which is an attribute in $P(A \times A)$; we take this latter attribute to be the internal equality relation for terms of type $A$. The adjunction brings the canonical map $1_{A} \rightarrow s_{\delta_{A}}(id_{A})$, which we read as a proof (over $A$) of reflexivity of $id_{A}$.

As we shall see in \textbf{6.9} the geometrical view on hyperdoctrines as fibrations sheds a new light on the notion of identity and reveals its highly non-trivial character. The link between the proof theory and the homotopy theory pointed to by Lawvere in 1970 \cite{Lawvere:1970b}
\footnote{
\begin{quote}
For deductions over $X$, one may take provable entailments (so that the category $P(X)$ reduces to a preodered set) or one may take suitable ``homotopy classes'' of deductions in the usual sense. One can write down an inductive definition of the ``homotopy'' relation, but the author does not understand well what results. (\cite{Lawvere:1970b}, p. 3-4) 
\end{quote}
}
recently became a field of active research, which led Voevodsky to the idea of the new \emph{Univalent Foundations} of mathematics discussed in \textbf{6.10} below.   
       
\section{Functorial Semantics}
The model-theoretical notion of interpretation (model) is treated by Lawvere in \cite{Lawvere:1966b}, \cite{Lawvere:1967} under the title of \emph{Functorial Semantics}. The idea is to identify a model of a given theory $T$ (presented as a category as just explained) with a functor $T \rightarrow S$ to the category of sets, which preserves the additional structure (i.e., the features of theory $T$ distinguishing it from a general category). In \cite{Lawvere:1969} Lawvere uses the same notion of model but this time describes the category of sets itself as a hyperdoctrine.  Thus the Functorial Semantics is a rather straightforward translation into a category-theoretic setting of the standard Tarski's set-theoretic semantics. However this translation is far from being trivial: as we shall now see it brings about (or ``reveals'' if one likes) some genuinely new features. 

Remind that the older notion of interpretation involves the idea of \emph{substitution} of non-logical terms of formal theory $T$ by certain objects and relations making part of another (usually informal) theory $T'$ like in the case of the arithmetical model of Euclidean geometry described by Hilbert in his \emph{Foundations} of 1899; what is then usually called a model is a particular construction in $T'$. What this standard setting lacks is a properly mathematical theory of substitution. Tarski's model theory describes concrete cases of substitution precisely but still has no general mathematicaldescription of this procedure. Lawvere's category-theoretic approach construes the relevant notion of substitution as a mathematical object, namely a functor. 
In additional to technical advantages of this approach, which I cannot analyze here, it also changes the usual way of thinking about relationships between theories and their interpretations (models) discussed in the last Chapters. With the functorial approach to model theory the notion of \emph{categoricity} of a formal theory (i.e., the uniqueness of its model up to isomorphism) loses its usual appeal: in the new setting one expects to obtain a \emph{category} $M_{T}$ of models (i.e., of certain functors) of given theory $T$, which has algebraic properties making it manageable, but not necessarily a category consisting of a single object (up to isomorphism). Further, this setting allows for considering given theory $T$ itself as a special object of $M_{T}$, so by Lawvere's word

\begin {quote} The theory appears itself as a generic model \cite{Lawvere:2004}, p.19.\end {quote} 

It is appropriate to ask whether or not Lawvere's theories are \emph{formal} in anything like Hilbert's sense of the word. As we have already seen the answer is rather in negative because Lawvere distinguishes between ``invariant'' conceptual theories from their formal presentations, and by default  uses the word ``theory'' for invariant theories. Let us now see more precisely how Lawvere describes the relationships between the formal and the conceptual, and how in this context he thinks about theories and their semantics. I come back to Lawvere's functorial semantics in \textbf{9.2}.

\section{Formal and Conceptual}
Lawvere's \emph{Dialectica} paper \cite{Lawvere:1969} begins as follows:

\begin{quote}
That pursuit of exact knowledge which we call mathematics seems to involve in an essential
way two dual aspects, which we may call the Formal and the Conceptual. For
example, we manipulate algebraically a polynomial equation and visualize geometrically
the corresponding curve. Or we concentrate in one moment on the deduction of theorems
from the axioms of group theory, and in the next consider the classes of actual groups to
which the theorems refer. Thus the Conceptual is in a certain sense the subject matter
of the Formal.
\end{quote}

The above passage can be interpreted coherently if one reminds that ``formal'' in the contemporary mathematical parlance implies ``symbolic''. So Lawvere talks here about (possibly uninterpreted) symbolic calculi, on the one hand, and their non-symbolic interpretations, on the other hand. I would like to stress that although such a distinction makes perfect sense with respect to our contemporary mathematics it does not apply to mathematics throughout all of its history. It certainly does not apply to mathematics practiced before the 17th century when symbolic methods were first introduced. 

Later in the same paper Lawvere considers the ``Theories'' on equal footing with the Formal and the Conceptual and suggests (without providing details) to think about these three things as a pair of adjoint situations: 

 $$\xymatrix{ Theories^{op} \ar@<2pt>[r]^{semantics} & Conceptual \ar@<2pt>[l]^{structure}}$$
 $$\xymatrix{ Formal \ar@<2pt>[r] & Theories \ar@<2pt>[l]}$$

which compose into  
 
  $$\xymatrix{ Formal^{op} \ar@<2pt>[r] & Conceptual \ar@<2pt>[l]}$$

I shall not try to reconstruct missing mathematical details but try to interpret these suggested adjoint situations heuristically. The sense of the first one is clear: an interpretation of a given \emph{theory} in terms of appropriate \emph{concepts} provides this theory with an appropriate \emph{semantics} while a reciprocal operation amounts to extraction of certain theoretical  \emph{structure} from a given conceptual domain. This resembles the usual notions of semantics and structure but one should keep in mind that Lawvere's \emph{Theories} are not formal theories in the usual syntactic sense but are connected with the \emph{Formal} through another pair of adjoint functors. Since in a different place Lawvere describes \emph{theories}  as ``invariant abstract concepts'' it makes sense to qualify \emph{concepts} in 
 $$\xymatrix{ Theories^{op} \ar@<2pt>[r] & Conceptual \ar@<2pt>[l]}$$
as (variable) \emph{concrete} concepts and include under this latter title what Kant would call ``representation'' (like the curve produced as a graph of a given polynomial in Lawvere's example).

As far as I can see the main purpose of this tripartition is to keep after Hilbert and Tarski the distinction between a theory and its (various) semantics and at the same time distinguish between the given theory and its (various) \emph{formal} syntactic presentations. Crucially Lawvere's ``invariant'' notion of theory unlike the ``semi-formal'' notion of theory developed in Hilbert's \emph{Foundations} of 1899 (which in a sense is also invariant with respect to the choice of its possible formalization) is itself a mathematical object, namely, a category.

Recently Joyal summarized Lawvere's tripartition as follows:
\begin{quote}
Ideally, every formal system of logic should exhibit three layers: a conceptual layer which specifies a
certain class $C$ of categories and functors, a semantic layer which specifies natural
examples of categories in $C$ (the semantic domains) and a formal layer which
specifies a language and a deduction system for constructing algebraically the categories
in $C$. The layers are not independent of each other and each clarifies the
others. But the conceptual layer has the central role as a kind of middle-man: \\
\emph{Formal layer} $\leftarrow$ \emph{Conceptual layer} $\rightarrow$ \emph{Semantic layer}
(\cite{Joyal:2011}, p. 19)
\end{quote}  

\section{Categorical Logic and Hegelian Dialectics}
As we have seen Lawvere critisizes the standard Hilbert-style Axiomatic Method by arguing that this method fixes attention upon particular syntactic presentations of logical structures and ignores the invariant content of logical and mathematical theories. Lawvere's novel categorical logic aims at revealing this invariant content and putting the syntactic presentations at their appropriate (modest) place. As Lawvere makes it explicit at many instances (some of which are referred to below) his project is underpinned by Hegel's philosophy and, in particular, by Hegel's distinction between the \emph{objective} and the \emph{subjective} logic \cite{Hegel:2010}. I shall try now to persuade the reader that Lawvere's Hegelian perspective on science, mathematics and category theory is indeed crucial for a philosophical understanding of his work. I would like to stress that in order to appreciate the relevance of Hegel's dialectical philosophy in categorical logic one does not need to be a ``Hegelian'' - although this must go without saying it seems me appropriate to make this remark in the present intellectual context of  continuing ideological and cultural (and very rarely properly philosophical) battles between the so-called ``Analytic'' and ``Continental'' philosophical traditions, which make an echo of the Cold War and other more tragic episodes of the passed century.  It is also important to have in mind that in \cite{Lawvere:1970a} Lawvere states explicitly that he himself is not  a Hegelian (p. 74, footnote 2). In a later work Lawvere describes his understanding of the mutual impact of mathematics and philosophy as follows: 

\begin{quote}
It is my belief that in the next decade and in the next century the technical advances forged by category theorists will be of value to dialectical philosophy, lending precise form with disputable mathematical models to ancient philosophical distinctions such as general vs. particular, objective vs. subjective, being vs. becoming, space vs. 
quantity, equality vs. difference, quantitative vs. qualitative etc. In turn the explicit attention by mathematicians to such philosophical questions is necessary to achieve the goal of making mathematics (and hence other sciences) more widely learnable and useable. Of course this will require that philosophers learn mathematics and that mathematicians learn philosophy. (\cite{Lawvere:1992}, p. 16)
\end{quote} 

Interestingly, Lawvere shares with his philosophical adversaries including Russell and other Analytic philosophers the idea that mathematics allows to put philosophical distinctions and arguments in a sharper and better disputable form, which allegedly allows for ``doing philosophy mathematically'' -  the view, which is denied by most of so-called Continental thinkers. (I don't want to make myself a subject to the Analytic/Contintental distinction, which is my view is very ill-construed
\footnote{While the expression ``Analytic philosophy'' is a self-description the expression ``Continental philosophy'' is invented by Analytic philosophers in the second half of the 20th century and to the best of my knowledge has been never used as a self-description.}
 but I would like to remind that I am also not enthusiastic about the idea of using mathematical methods in philosophy. I explained some reasons for it the above Introduction.), 

At the same time, as we shall shortly see, Lawvere's disagreements with the basics of Analytic philosophy are too profound for qualifying his work as a mere ``formalization of Hegel's dialectics by means of modern logic''. This makes Lawvere's philosophical thinking and its impact onto mathematics quite unique in its kind. As a philosopher Lawvere inherits from different schools of thought and remains remarkably unaffected by the unfortunate Analytic versus Continental distinction. I believe that today such a position is a \emph{sine qua non} for doing philosophy seriously.   

Let me now turn to Hegel and, more specifically, to his distinction between the \emph{objective} and {subjective} logic. Here is a relevant passage: 

\begin{quote}
What is to be considered is the whole Notion, firstly as the Notion in the form of being, secondly, as the Notion; in the first case, the Notion is only in itself, the Notion of reality or being; in the second case, it is the Notion as such, the Notion existing for itself (as it is, to name concrete forms, in thinking man, and even in the sentient animal and in organic individuality generally [..]). Accordingly, logic should be divided primarily into the logic of the Notion as being and of the Notion as Notion - or, by employing the usual terms (although these as least definite are most ambiguous) into 'objective' and 'subjective' logic. (\cite{Hegel:2010}, 79)
\end{quote}

Hegel's \emph{Notion} is a category that comprise both  (i) the reality (aka being) and human (and more generally also animal, as Hegel hints) (ii) thinking about reality. This notion of notion, which is crucial for Hegel's so-called ``objective idealism'', is not of our primary concern here because Lawvere rejects Hegel's objective idealism (\cite{Lawvere:1970a}, p. 74) and does not use it. So the relevant part of the content of the above quote is this: logic is divided into two parts one of which is (i) the logic of being aka the objective logic, while the other is (ii) the logic of thinking aka the subjective logic (it is called by this latter name because thinking requires a thinking subject). Hegel's subjective logic is what today (and also in Hegel's times) is commonly called logic.  Let us first focus on the objective logic. 

For explaining Hilbert's notion of objective logic I cannot help but use a bit of Hilbert's dialectical method of definition - even if I'm trying now to reduce the use of this method to minimum in order to make my explanation more accessible to people who are not used to it. This simply amounts to saying that I am going to define this notion through several  approximations rather than give immediately a formal definition. 

In the first approximation I identify Hegel's objective logic (=logic of being) with \emph{ontology}. In this first approximation the above quote simply tells us about the familiar distinction between logic (in the usual sense of logic of thinking) and ontologyr. This a useful approximation but it is by no means sufficient for our purpose.

At the next dialectical step we should distinguish between ontology as a part of metaphysics (meaning the metaphysics of Aristotle, of Schoolmen, Wolf's metaphysics or the contemporary Analytic metaphysics) and the objective logic proper. As Hegel tells us a way to pursue a dialectical reasoning consists in following the history of relevant notions. This is what we need at this second step.  The relevant historical character is Kant:

\begin{quote}
Recently Kant has opposed to what has usually been called logic another, namely, a transcendental logic. What has here been called objective logic would correspond in part to what with him is transcendental logic. He [Kant] distinguishes it from what he calls general logic in this way, [a] that it treats of the notions which refer a priori to objects, and consequently does not abstract from the whole content of objective cognition, or, in other words, it contains the rules of the pure thinking of an object, and [b] at the same time it treats of the origin of our cognition [..]. It is to this second aspect that Kant's philosophical interest is exclusively directed. (\cite{Hegel:2010}, 81)
\end{quote}

So our second approximation amounts to saying after Hegel that the objective logic ``corresponds in part'' to Kant's transcendental logic. I find it actually useful in this second approximation to boldly identify Hegel's objective logic with Kant's transcendental logic. As we shall see at the final third step, the difference between the two is much neater than the difference between the objective logic and the traditional ontology.

We have already encountered Kant's notion of transcendental logic earlier in this book (see \textbf{1.3, 2.2}). As Hegel reminds us in the above quote Kant distinguishes his transcendental logic from the ``general logic'', which is most often called today simply logic, by specifying that the transcendental logic is not wholly topic-neutral: while the general logic deals with abstract logical individuals the transcendental logic deals only with objects of possible experience, i.e., with would-be empirical objects. This is why the general logic applies everywhere but the transcendental logic applies only in mathematics and (mathematized) empirical sciences. Notice now that there is a sense in which the transcendental logic provides a replacement for the traditional metaphysical ontology. 

Let me explain this with an example. Doing the traditional metaphysics \footnote{I keep talking about the ``traditional metaphysics'' but not simply about ``metaphysics'' not for saying that the metaphysics is doomed to be swept away by the future scientific progress (it is very hard to make today such an optimistic prediction) but because Kant has a different and more specific notion of metaphysics, which is irrelevant in the given context.} one may stipulate that all entities (or alternatively - only entities of a special sort called \emph{physical entities}) live in space and time, and then ask whether space and time are absolute or relational, elaborate on formal properties of spatio-temporal relations and so on. Kant's transcendental logic reformulates these questions into questions about our \emph{thinking} about space and time, i.e., into \emph{logical} questions (in the special ``transcendental'' sense of ``logical''). Although some questions about space and time can be asked and answered in both settings similarly, the two settings, generally, give rise to different questions and different answers. This is why the difference between the two approaches does not reduce to two different answers to the  question Does the world (as we usually think about it) really exist or it is a creation of our brains? This latter question belongs to the metaphysical ontology rather than to the transcendental logic, and so from the point of view of transcendental logic it is ill-posed.

The third approximation (that will be my last) is more difficult because here we touch upon the heart of Hegel's project: somehow to get rid of the subjective bent of Kant's transcendental logic without bringing back the traditional metaphysics. The main idea is to keep doing something like logic (rather than the sort of speculative physics called metaphysics) but make it objective. 
\footnote{The claim that Kant's transcendental logic is ``subjectively bent'' is less evident and less simple than it may appear. Natorp \cite{Natorp:1912} denounces such a claim as a misinterpretation - although his argument can be  perhaps better understood as an attempt to improve on Kant rather than simply reconstruct his thought. However Hegel's understanding of Kant's alleged subjectivism is far from being naive; here is a relevant passage, which also provides a useful exposition of Kant's understanding of objecthood:
\begin{quote}
An object, says Kant, is that in the notion of which the manifold of a given intuition is unified. But all unifying of representations demands a unity of consciousness in the synthesis of them. Consequently it is this unity of consciousness which alone constitutes the connection of the representations with the object and therewith their objective validity and on which rests even the possibility of the understanding. Kant distinguishes this unity from the subjective unity of consciousness, the unity of representation whereby I am conscious of a manifold as either simultaneous or successive, this being dependent on empirical conditions. On the other hand, the principles of the objective determination of notions are, he says, to be derived solely from the principle of the transcendental unity of apperception. Through the categories which are these objective determinations, the manifold of given
representations is so determined as to be brought into the unity of consciousness. According to this
exposition, the unity of the notion is that whereby something is not a mere mode of feeling, an intuition, or even a mere representation, but is an object, and this objective unity is the unity of the ego with itself. In point of fact, the comprehension of an object consists in nothing else than that the ego makes it its own, pervades it and brings it into its own form, that is, into the universality that is immediately a determinateness, or a determinateness that is immediately universality. As intuited or even in ordinary conception, the object is still something external and alien. When it is comprehended, the being-in-and-for-self which it possesses in intuition and pictorial thought is transformed into a positedness; the I in thinking it pervades it. (\cite{Hegel:2010}, 1293)
\end{quote}
}. In the following passage Hegel compares his objective logic with the traditional metaphysics and at the same time describes it as an improved version of Kant's philosophical ``critique'': 

\begin{quote}
The objective logic, then, takes the place rather of the former metaphysics which was intended to be the scientific construction of the world in terms of thoughts alone. If we have regard to the final shape of this science, then it is first and immediately ontology whose place is taken by objective logic [..] But further, objective logic also comprises the rest of metaphysics in so far as this attempted to comprehend with the forms of pure thought particular substrata taken primarily from figurate conception, namely the soul, the world and God [..]. [Objective l]ogic, however, considers these forms free from those substrata, from the subjects of figurate conception; it considers them, their nature and worth, in their own proper character. Former metaphysics omitted to do this and consequently incurred the just reproach of having employed these forms uncritically without a preliminary investigation as to whether and how they were capable of being determinations of the thing-in-itself, to use the Kantian expression - or rather of the Reasonable. Objective logic is therefore the genuine critique of them - a critique which does not consider them as contrasted under the abstract forms of the a priori and the a posteriori, but considers the determinations themselves according to their specific content. (\cite{Hegel:2010}, 85)
\end{quote}
 
Observe that the price paid by Hegel for fixing Kant's critical philosophy is high: along with the subjective bend Hegel gets rid of Kant's a priori versus a posteriori distinction;  arguably this step diminishes the role of empirical data in sciences and after all allows for the ``scientific construction of the world in terms of thoughts alone'' in a new dialectical form. As far as I am concerned this price is too high. However this critical remark does rule out Hegel's dialectical method of reasoning as such but rather rises the problem of how to use this method in the modern science without compromising against the empirical character of this science. 

Let us now turn to  \emph{subjective} logic. It requires for its explication a dialectical procedure of a different sort because this kind of logic is most commonly simply identified with logic.  Apparently this common opinion did not change much since Hegel's times - albeit the booming development of logic during the last century perhaps made it a bit less ``ossified'' and more ``fluid'': 

\begin{quote}
This part of the logic which contains the Doctrine of the Notion [..] is issued under the particular title System of Subjective Logic, for the convenience of those friends of this science who are accustomed to take a greater interest in the matters here treated and included in the scope of logic commonly so called, than in the further logical topics treated in the first two parts [which cover the Objective Logic]. For these earlier parts I could claim the indulgence of fair-minded critics on account of the scant preliminary studies in this field which could have afforded me a support, material, and a guiding thread. In the case of the present part, I may claim their indulgence rather for the opposite reason; for the logic of the Notion, a completely ready-made and solidified, one may say, ossified material is already to hand, and the problem is to render this material fluid and to re-kindle the spontaneity of the Notion in such dead matter. If the building of a new city in a waste land is attended with difficulties, yet there is no shortage of materials; but the abundance of materials presents all the more obstacles of another kind when the task is to remodel an ancient city, solidly built, and maintained in continuous possession and occupation. Among other things one must resolve to make no use at all of much material that has hitherto been highly esteemed.  (\cite{Hegel:2010}, 1277) 
\end{quote}

So Hegel's treatment of subjective logic has three basic aspects: (i) revealing the subjective character of what is commonly called logic (i.e., the formal logic),  (ii) ``re-kindling the spontaneity of Notion'' in the common ``ossified'' logical categories and finally (iii) ``passing out of subjectivity into objectivity \cite{Hegel:2010}, 1441''; this latter dialectical procedure can be called the \emph{objectification} of logic. 

As Hegel suggests in the last quote and elsewhere (see, for example, \cite{Hegel:2010}, 86) the very distinction between the objective and the subjective logic should be thought of as dialectical to wit polemical. What he really aims at is a transformation of the common notion of logic into what he calls the objective logic or simply logic. A better name for subjective logic is the ``logic of Notion'', which is the third concluding part of Hegel's (objective) logic. (Beware that the three parts of Hegel's logic - the logic of Being, the logic of Essence and the logic of Notion -  are to be thought of as consecutive stages of a single process of dialectical reasoning rather than three parts of the same whole, which co-exist side by side in some intellectual space.) However since Hegel's notion of logic is very far from being common Hegel struggles (dialectically!) against the common view and use such labels as ``objective'' and ``subjective'' in this struggle. Let me however leave this very preliminary elucidation of Hegel's logic at this point and turn back to Lawvere. (The dialectical character of the above explanation implies its open-endedness, and the interested reader is advised to pursue a further study of Hegel's logic independently.) 
  
As we already know from the above quote (\cite{Lawvere:1992}, p. 16) Lawvere's project involving Hegel's logic is twofold: Lawvere wants to (i) recast Hegel's dialectical logic in mathematical terms of categorical logic and (ii) use Hegel's dialectical logic as a guide in his mathematical research including his research in categorical logic. Lawvere also makes it clear that his strategy is to merge (i) and (ii)  into a single mathematico-philosophical project, so when I distinguish between (i) and (ii) and keep them apart, I'm clearly doing an opposite move. The reason I am doing it is the following: just as in Chapter \textbf{2} I discussed various philosophical influences on and philosophical interpretations of Hilbert's Formal Axiomatic Method but did \emph{not} discuss the impact of this method on philosophy (which was indeed very significant), talking in this present Chapter about Lawvere's work I'm going to discuss only (ii) and leave (i) aside. Since in Lawvere's work (i) and (ii) are very closely merged my choice is a choice of perspective rather than a choice of material.  One and the same piece of Lawvere's work can be usually seen both (i) as a piece of Hegelian dialectics mathematized with the category theory and (ii) as a piece of mathematics inspired and guided by Hegelian dialectics - and the author's intention, as I understand it, is to synthesize these two aspects into a single whole. Thus my choice to focus on (i) in this book makes my perspective on Lawvere's work special and deliberately incomplete.  

Let us now see how Lawvere distinguishes between the objective and the subjective logic in the context of categorical logic: 
\begin{quote}
Arising [..] from the needs of geometry, category theory has developed such notions as adjoint functor, topos, fibration, closed category, 2-category, etc. in order to provide

(i) a guide to the complex, but very non-arbitrary constructions of the concepts and their interactions which grow out of the study of space and quantity.

It was only the relentless adherence to the needs of that basic subject that made category theory so well-determined yet powerful. [..] If we replace ``space and quantity'' in (i) above by ``any serious object of study'', then (i) becomes my working definition of \emph{objective logic}. Of course, when taken in a philosophically proper sense, space and quantity do pervade any serious field of study. Category theory has also objectified as a special case

(ii) the subjective logic of inference between statements. Here statements are of interest only for their potential to describe the objects which concretize the concepts.  
 (\cite{Lawvere:1994}, p. 16)
\end{quote} 

If we judge the above Lawvere's definition of objective logic by Hegel's standard we immediately notice a gap (or rather a leap) in it: while Lawvere relates his objective logic to categories of space and quantity \emph{directly}, Hegel arrives to these categories only after some dialectical development that begins with categories of Being, Nothingness and Becoming. However instead of trying to find a place for these and some other missing Hegelian categories\footnote{For such an attempt see \cite{Meles:2012}.} in the categorical logic I shall rather change the perspective (in accordance with my choice of perspective explained above) and comment on Lawvere's definition of logic as it stands using some key ideas coming from Hegel's logic but without assuming that Lawvere's logic is supposed to represent Hegel's logic fully and faithfully. 

Notice now how Lawvere relates the ``study of space and quantity'' and ``any serious field of study''. He, first, (a) tells us about a category-theoretic guide to conceptual constructions growing out of the study of space and quantity, second, (b) tentatively defines the objective logic as a category-theoretic guide to conceptual constructions growing from any serious field of study and, finally, (c) states that ``in a philosophically proper sense, space and quantity do pervade any serious field of study'', so that (a) and (b) are in fact equivalent and one may think of the objective logic as the logic of space and quantity without restricting the generality. This micro-dialectics seems me quite remarkable and illuminating. Why Lawvere does not define the objective logic as the logic of space and quantity to begin with? I suppose, because according to the common opinion (formal) logic must apply in any serious field of study independently of its subject-matter; so if Lawvere would say from the outset that his objective logic applies only to the study of space and quantity one could think that this special logic applies only in a specific field and for this reason perhaps even does not deserve the name of logic. This is why Lawvere says first that the objective logic applies in any serious field of study and only then he makes it clear that in the ``philosophically proper sense''  the logic of any serious field of study \emph{is} the logic of space and quantity.  

What is then the philosophically proper sense in which ``space and quantity do pervade any serious field of study''? Since Lawvere does not explain this I give my own explanation. It is useful for my explanation to begin with Kant and Neo-Kantians and only then turn to Hegel. Remind from \textbf{2.2} (a) Cassirer's critique of Russell's mathematical logicism, which stresses the fact that Russell's approach throws the general theory of magnitude away from the pure mathematics, and (b) the related claim according to which  logic and mathematics must not be ``instruments for building a metaphysical 'world of thought''' but must apply  only within the mathematized empirical science. Both Cassirer's claims are closely related because he assumes after Kant that  an adequate general mathematical theory of magnitude is necessary and sufficient for distinguishing objects of (possible) experience from metaphysical ``individuals'', which one may stipulate for free. Talking about logic Cassirer has in mind Kant's idea of transcendental logic, which uses such a general theory of magnitude for this very purpose. Now by exchanging the general theory of magnitude for the theory of space and quantity (which in the given context seems me unproblematic ) we can distinguish after Cassirer a philosophical sense in which this later theory indeed pervade any serious field of study.  If by a serious field of study we mean one that, first, apply mathematics and, second, deal with empirical data then we need a general theory of space and quantity in order to distinguish between those conceptual constructions, which can possibly represent empirical objects, and those which can not.\footnote{Mathematics and logic themselves are serious fields of study insofar as they develop tools for dealing with empirical data; thus they also deal with empirical data albeit indirectly.}. As long as an appropriate theory of space and quantity is built into a system of logic this system of logic can be described as ``logic of space and quantity'' and at the same time apply in any serious field of study (but not elsewhere). 

The above Kantian reconstruction of Lawvere's ``philosophically proper sense'' in which  ``space and quantity pervade any serious field of study'' does not involve any assumption according to which the space and quantity are in some sense subjective, so this can hardly be a controversial point here. It is unsatisfactory rather in a different respect: the very idea of a theory allowing for the sharp distinction between useful and ``metaphysical'' (in the sense of ``merely speculative'') mathematical constructions is doubtful. Even if I strongly disagree with Wigner's view according to which  the effectiveness of mathematics in natural sciences is ``unreasonable'' and must be thought of as  a ``miracle'' and the ``wonderful gift'', I realize that we are not living an the Kantian paradise where all sound mathematics applies in natural sciences and technologies immediately. So we need a more flexible and more dynamic way of thinking about these issues than Kant's philosophy may provide. Here the model of Hegel's dialectical logic becomes relevant. Rather then be a device for philosophical critique of existing science like Kant's transcendental logic the categorical logic is conceived by Lawvere after Hegel as a dialectical tool for re-configuring the foundations of modern science, which involve, in particular, the categories of space and quantity and the way in which these categories are applied in empirical studies (in particular, in measurement). Thus there is a further philosophically proper sense in which space and quantity pervade any serious field of study: a progress in the general theory of space and quantity is a necessary condition of the substantial progress in the fundamental empirical research. Since the objective logic is supposed to be a guide for further progress in any serious field of study it must qualify as a dialectical logic of space and quantity. 
 
Lawvere's theory of space and quantity is found in \cite{Lawvere:1992}. The idea is to think of presheaf $X: C^{op}\rightarrow Set$ as a rule that assigns to every object $c$ of $C$ a set of maps $X$ going from $c$ to \emph{space} $X$ where this given object $c$ is ``placed'', so $c$ plays here the role of ``test space'' or ``shape'' and the bigger space $X$ is determined in terms placements of this basic shape. A \emph{quantity} is a co-presheaf $Y: C \rightarrow Set$ thought of as a rule that assigns to every object $c$ (a test space) the set of incoming maps, each of which in its turn assigns certain ``value'' in this given space to other objects. In \cite{Lawvere:2005a} develops this theory further and, in particular, provides the categories of space and quantity with metrics. A lot still remains to be done in order to make this conceptual apparatus into a working tool for physicists and other scientists. During the last decade the search for applications of category-theoretic methods in natural sciences became a field of active research which I can not overview here.

As we have seen earlier Lawvere achieves an \emph{internalization of logic} with respect to appropriate categories in the sense that basic logical concepts such as propositional functions, logical connectives, quantifiers, truth-values, theories and models, are understood as category-theoretic constructions. In the above quote Lawvere describes this move in Hegelian terms as an ``objectification'' of logic. The key idea here is that the  ``subjective logic of inference between statements'' (which, remind, is commonly identified with logic as such) must not be thought of as a self-standing systems of laws and rules, which provides ultimate foundations of mathematics and natural science, but must emerge as an aspect of basic conceptual constructions of science, which involve categories of space and quantity, and perhaps some other. In my view, this is the most important impact of categorical logic on the Axiomatic Method. 

Let me to be clear at this point: the very idea that laws of logic are not self-standing is not new. What we call today logical laws like the law of non-contradiction, Aristotle describes as metaphysical laws. Aristotle's syllogistics is equally grounded in his metaphysics and ontology. Russell and other people who attempted to revive metaphysics in the 20th century grounded their metaphysics onto the new mathematical logic rather than the other way round (this is why Russell calls his metaphysical atomism ``logical''). While Aristotle's metaphysical views were closely related to his own research in natural sciences the new logical metaphysics emerged in the 20th century under the name of Analytic metaphysics stands wholly apart the contemporary fundamental research in sciences but survives either in the form of self-contained intellectual game, which combines playing with linguistic intuitions with some amount of formal rigor, or in the form of applied discipline (applied ontologies in computer sciences). Thus the problem of \emph{grounding} logic remains today largely open. It becomes only more acute in the present situation when we have got hundreds formal calculi, which are offered under the brand of ``logic''. Hegel's notion of objectification of logic as a modern replacement for the traditional idea of grounding logic in metaphysics suggests a strategy for solving this problem.

\section{Toposes and their Internal Logic}

Now I shall demonstrate Lawvere's dialectical method at work with an example of his axiomatic treatment of the geometrical notion of \emph{topos}. This is in fact something more than just an example because it amounts to a mathematical discovery, which is important for our study on its own rights.  Remind Lawvere's remark that ``[t]he formalism of category theory is itself often presented in ``geometric'' terms''  (\cite{Lawvere:1969}, p. 283, already quoted in the introductory part of Chapter \textbf{4}). As we shall now see the link between categorical logic and geometry is not limited to the geometric way of presentation but is more profound. It has been first made explicit in Lawvere's paper \emph{Quantifiers and Sheafs} of 1970 which begins as follows:

\begin{quote}
The unity of opposites in the title is essentially that between logic and geometry, and there are compelling reasons for maintaining that geometry is the leading aspect. At the same time, in the present joint work with Myles Tierney there are important influences in the other direction: a Grothendieck ``topology'' appears most naturally as a modal operator, of the nature ``it is locally the case that'', the usual logical operators, such as $\forall$, $\exists$, $\Rightarrow$ have natural analogues which apply to families of geometrical objects rather than to propositional functions, and an important technique is to lift constructions first understood for ``the'' category $\underline{S}$ of abstract sets to an arbitrary topos. We first sum up the principle contradictions of the Grothendieck-Giraud-Verdier theory of topos in terms of four or five adjoint functors [..] enabling one to claim that in a sense logic is a special case of geometry. ( \cite{Lawvere:1970a}, p. 329)
\end{quote}

The unity of opposites in the title is that between logic and geometry because the term ``quantifier'' refers to logic while the term ``sheaf'' refers to geometry. Since the geometrical background, to which Lawvere refers here, is not generally known I'll try to present it briefly for the non-mathematical reader.

The notion of \emph{topos} first appeared in the circle of Alexandre Grothendieck about 1960 as a twofold upgrade of the notion of topological space (for an informal explanation of the notion of topological space see \textbf{3.4} above). The first upgrade amounts to considering a given topological space $T$ together with the \emph{sheafs} of functions from open subsets of this space to some target sets; a sheaf respects conditions, which allow for seeing the target sets as ``momentary images'' of the same \emph{continually variable} set varying over $T$. (If the target sets are provided with an extra structure, say, with the group structure, one may similarly think of \emph{groups} continuously varying over a given topological space.) In order to get from the notion of sheaf to that of topos we need first to render the former notion into category-theoretic terms. Think of $T$ as a category with objects open subsets of $T$ (\emph{opens} for short) and morphisms set-theoretic inclusions of these subsets, so in the resulting category there is at most one morphism going from one given object to another  (for any pair of opens $U, V$ we either do or do not have $U \subseteq V$; categories with at most one morphism with a fixed domain and a fixed codomain are called partial orders). Then a sheaf can be defined as a functor $T^{op} \rightarrow S$ from the category $T^{op}$ obtained from $T$ by the ``reversal of arrows'' to the category of sets $S$, which satisfies certain conditions assuring that the target variable set varies (with respect to $T$) \emph{continuously}. 
\footnote{ 
The fact that the arrows must be reversed in this case was difficult to understand without using the category-theoretical notion of functor; this was a major difficulty for earlier attempts to develop a ``topology without points'' \cite{Johnstone:1983}.  
} 
One gets now a basic example of \emph{topos} by considering  the category of all sheafs on a given topological space together with maps between those things (since sheafs are functors the maps are ``functors between functors'' aka \emph{natural transformations}). This topos can be naturally thought of as a space (or rather spacetime) of sets continuously varying over $T$.   
 
 The second upgrade amounts to a generalization of the usual notion of topology. Given (usual) topological space $T$ one may always associate with a given open $U$ its \emph{covering family} $C_{U}$ which is a collection of opens $V_{i} \subseteq U$ such that their union equals $U$ (i.e., each point of $U$ belongs to at least one of $V_{i}$); in particular, $T$  itself is always covered by at least one collection of its opens. Grothendieck observed that the notion of covering family makes sense not only for partial orders but also for categories of more general sort and defined a covering family of a given object to be a collection of incoming morphisms (not necessarily monomorphisms) closed under certain operations. This led him to a more general notion of topology called \emph{Grothendieck topology} defined by distinguishing among all collections of morphisms sharing a codomain those, which count as covering families of this given object. A category $C$ provided with a Grothendieck topology $J$ is called a \emph{site} $(C, J)$. A sheaf over a site is defined just like in the case of topological space. The Grothendieck topos is a category of sheaves over some given site. For a systematic introduction I refer the reader to  \cite{MacLane&Moerdijk:1992}, ch. 2-3.   
 
The notion of topos invented by Grothendieck and developed by his collaborated mentioned by Lawvere in the above quote was not originally supposed to have any special relevance to logic; the discovery of such a special relevance is wholly due to Lawvere. This is what Lawvere, modestly calls the ``influences in the other direction'' meaning the impact of logical considerations. Let us see what this impact amounts to.   

In the beginning of his seminal  paper \cite{Lawvere:1970a} Lawvere provides his definition of topos usually called today the definition of \emph{elementary} topos ; the title ``elementary'' reflects the fact that Lawvere's definition unlike Grothendieck's original construction almost straightforwardly translates into the standard first-order formal language \cite{McLarty:1992}. According to this definition an (elementary) topos $T$ is CCC with a \emph{subobject classifier}, which plays in a general topos the role similar to that played by 2 (the two-point set) in the category of sets (which also qualifies as a topos in the sense of Lawvere's definition). 2 classifies subsets of a given set $S$ in the sense that if one asks whether a given element $p \in S$ belongs to subset $U \subseteq S$ there are just 2 possible answers: yes and no; this allows for identifying every subset  $U$ with a particular function $u: S \rightarrow 2$, which sends every $p$ belonging to $U$ to ``yes'' and every $p$ not belonging to $U$ to ``no''. Correspondingly the set $2^{S}$ of all such functions is identified with the set of all subsets (the powerset) of set $S$. Given two objects $A, B$ of CCC the exponential object $A^{B}$ always exists but in order to get a distinguished object $\Omega$ playing the role of ``object of truth values'', so that for all $A$ $\Omega^{A}$ represents the space of subobjects of $A$, one needs an additional postulate. By a subobject of $A$ one means here any incoming \emph{monomorphism} $f$, i.e.,  such $f$ that for all $g, h$ $g\circ f = h \circ f$ implies $g = h$ (given the composition is written in the geometrical order). Given two subobjects $f_{1}, f_{2}$ of the same object $A$ consider morphism $h$ such that $f_{1} = h\circ f_{2}$; according to the definition of subobject there is no more than one such morphism. This shows that subobjects of a given object are partially ordered. In the case of the category of sets, which qualifies as a particular elementary topos, the partial order of subobjects is the (complete) Boolean lattice while in the general case the lattice is Heyting. For a systematic treatment see  \cite{McLarty:1992}, and for the most complete compendium of topos theory existing to the date and developed from the elementary viewpoint see \cite{Johnstone:2002}. Evidently Lawvere's earlier work on categorical axioms for set theory, which we reviewed in \textbf{4.1}, helped Lawvere to formulate his axioms for (i.e., the definition of) an elementary topos. It was Lawvere but not Grothendieck who first thought of sheaves as \emph{continuously variable sets} and observed that the category of such things shares with the category of usual ``static'' sets a number of basic properties.      

As Lawvere notices in the above quote the concept of elementary topos is more general than that of Grothendieck topos: there is a large class of elementary toposes, which are not Grothendieck. In particular, every category $(C^{op}, S)$ of all \emph{presheaves}, i.e., of contravariant functors from a small category $C$ (without topology) to the category of sets, is an elementary topos but not Grothendieck topos. However the topological aspect of topos is treated in the elementary setting too: Grothendieck topology is recaptured as a modal operator satisfying simple axioms. 

Lawvere's axioms for elementary topos helped many people outside the community of specialists in algebraic geometry to enter into this field  and make a fruitful research in it. Everyone who learns today the topos theory begins with Lawvere's axioms for elementary topos. This makes Lawvere's axiomatization of topos theory a true  success story of Axiomatic Method in the 20th century mathematics. Although Lawvere's axiom for toposes are not given in a wholly formalized form they allow for a natural formalization (without going the roundabout way through a membership-based set theory) as it has been shown (or better to say \emph{done}) by McLarty  \cite{McLarty:1992}. However I insist that this just one aspect of the story but not the whole story. 

Remind from the above quote Lawvere's claim that ``in a sense logic is a special case of geometry''. It is stunning if we think about it against the background of Hilbert's Axiomatic Method. For every modification of this method (either exemplified by Hilbert's \emph{Foundations} of 1899 or by any modern textbook in the axiomatic set theory or even by \cite{McLarty:1992}) assumes what I have called above (\textbf{2.2}) the mathematical logicism in the large sense of the term, i.e., the view according to which (some system of) logic provides a foundation of mathematical theories by helping formulate axioms of this theory and then derive from these axioms some further theorems. Now Lawvere tells us that ``geometry [rather than logic!] is the leading aspect'' and that there is a sense in which logic is a special sort of geometry. I shall purport now to explain this claim both from a mathematical and philosophical viewpoints, and then show that Lawvere's dialectical Axiomatic method indeed does not reduce to Hilbert's.

When Lawvere talks in the above quote about logic as a special case of geometry he refers to the \emph{internal} logic of a given topos. We have already discussed the idea of internalization of logic in category theory meaning the possibility to construe basic logical concepts including logical connectives and truth-values by category-theoretic means. This sense of being internal remains wholly relevant when we talk about the internal logic of a topos but what is special about the topos logic is the fact that toposes (unlike CCCs in general) have also a reach geometrical content (this is particularly true in the case of Grothendieck toposes) and thus provide a dialectical interplay between geometry and logic. So in his \cite{Lawvere:1970a} Lawvere not only axiomatized the topos theory by bringing into this difficult geometrical theory the logical clarity of first-order logic: he also discovered that toposes, which were earlier invented by Grothendieck as geometrical categories, are also appropriate for ``doing logic \emph{in} them'', i.e., that toposes are not only geometrical but also logical categories! In particular, Lawvere shows how his earlier invented notion of quantifier as adjoint to substitution (see \textbf{4.3}) is geometrically realized in any topos. (This particular result gave the paper its title.) I would like to mention here that the internal logic of a given topos can be presented in the usual way with a symbolic syntax (called the Mitchell-B\'enabou language or \emph{internal language}) and the corresponding semantics of that language (called the Kripke-Joyal semantics) associated with that topos, see \cite{MacLane&Moerdijk:1992}, p. 296 - 318).

How the two aspects of Lawvere's axiomatization of topos theory - (i) the logical clarification of the geometrical notion of topos with first-order axioms \emph{and} (ii) the internalization (and hence the ``geometrization'') of logic in a topos) - relate to each other? It is fair to say that they relate dialectically but one may want to be more specific here.  McLarty's book \cite{McLarty:1992} shows how this relationships looks like when one develops Lawvere's theory of elementary categories and toposes with the standard Formal Axiomatic Method
\footnote{Arguably, unless this theory is treated in this way it does not deserve the title of ``elementary''. I would not like to make an issue of this terminological point; I use the expression ``elementary topos'' interchangeably with ``Lawvere topos''.}. In the Preface to his book McLarty describes his general logical framework as follows:  
\begin{quote}
Our metatheory avoids excluded middle and choice so it is sound in any topos, except when we are explicitly concerned with constructions in $Set$. ( \cite{McLarty:1992}, p. vii)
\end{quote}  

After introducing categories and toposes (in a semi-formal way) McLarty comes to the internal logic of toposes  (Chapter 14)  and then in Chapter 16 titled ``From the Internal Language to the Topos'' he shows how a given topos can be described \emph{internally} in terms of its own internal language. The internal view on a topos, generally, does not fully coincide with the external view - in spite of the fact that the ``external logic'' (i.e., the metatheory) is of the same sort that the internal logic of the given topos, as this is indicated in the above quote. It is suggestive to compare this effect with a similar situation in geometry where one may describe one and the same space $S$ either extrinsically as embedded into some outer spaces $T$ or intrinsically in terms of embedding of a ``test spaces'' $P$ in $S$. Draw a straight line $L$ one a sheet of paper and then fold the paper. Extrinsically $L$ is no longer straight but intrinsically it has not changed. So the extrinsic and the intrinsic views on the situation are no longer the same (see \textbf{8.8} below).  

Coming back to Lawvere I suggest that the geometrisation of logic in topos theory can be philosophically interpreted as a way of \emph{objectification} of logic (in the Hegelian sense of the term explained above). This may work, of course, only if the geometrical notion of topos itself proves objective in the appropriate sense. And this requires a topos to be not just an abstract mathematical concept but a concept providing a connection between the geometrical intuition and the world of experience, between the pure mathematics and the natural sciences. Lawvere clearly states that establishing of such a connection is indeed his aim:

\begin{quote}
[E]xperience with sheaves, permutation representations, algebraic spaces, etc. shows that a ``set theory'' for geometry should apply not only to \emph{abstract} divorced from time, space, ring of definition, etc., but also to more general sets, which do in fact develop along such parameters.  ( \cite{Lawvere:1970a}, p. 329)
\end{quote}  

(By categories of ``more general sets'' Lawvere means here toposes: remind the metaphor of sheaf as a continuously variable set.)\footnote{It is not immediately clear how a purely mathematical study in sheaf theory and algebra may have a bearing onto the physical time and the physical space. However the context of Lawvere's work makes it clear that Lawvere speaks here about time and space as objective physical categories rather than metaphors or subjective intuitions. So I can see two possible answers. Either Lawvere believes that a serious mathematical study brings about objective - and hence physically relevant - results even if this study is pursued quite independently of any physical considerations. Or he intends to bring himself this objective meaning into the field by establishing a connection with natural sciences through such basic categories as time and space (and quantity).}

It is appropriate to ask: How it is possible to combine Hegel's dialectical logic with the Axiomatic Method? How Lawvere manages to do this? As far as one thinks about the Axiomatic Method in Hilbert's vein any combination of this method with Hegel's dialectical logic seems to be impossible for the following simple reason: while on Hegel's account the ``subjective'' logic of inference is the concluding element of a dialectical process, which involves categories of Magnitude (Lawvere's Space and Quantity) at an earlier stage, on Hilbert's account such a subjective logic of inference is supposed to be fixed from the outset and then used for an axiomatic introduction of further categories including that of Magnitude. Remind however that in  (\cite{Lawvere:2003}, p. 213, see the quote in the introductory part of Chapter \textbf{4}) Lawvere describes the Axiomatic Method as  the ``unification and concentration'' of mathematical practice; obviously this description is much broader than Hilbert's notion of Axiomatic Method, which we thoroughly studied in Chapter \textbf{2}. So we have no reason to assume that Lawvere's Axiomatic Method is the same as Hilbert's.

Lawvere's co-authored textbooks \cite{Lawvere&Schanuel:1997} and \cite{Lawvere&Rosebrugh:2003} give some more hints about Lawvere's Axiomatic Method. Both books begin with an introduction of basic category-theoretic concepts at the example of the category of finite sets and their mappings (i.e., functions) and end up with the elementary topos. Logic in the usual sense of  logic of inference appears (or rather emerges) only at this final stage in the form of the internal logic of elementary topos. Throughout the exposition mathematical concepts are given appropriate physical meanings, in particular, spatial and temporal meanings.  Thus the structure of this exposition fits the idea of objective logic: first one builds an empirically meaningful system of objective categories, and only then on this objective basis introduces a system of subjective logic, which involves truth-values, types, connectives, quantifiers and the rest of the usual formal logical machinery. This subjective logic not only reflects the general features of the corresponding objective categories but is also fine-tuned by specific objective features of each given example: in particular, the set (or more generally the \emph{object}) of truth-values of the internal logic $L$ of given topos $T$ is determined by this very topos, which is an objective and physically meaningful conceptual construction. Thus in the accordance with Hegel the subjective logic emerges here on the top of the objective one.

One may argue that in order to cover all this material ``more accurately'' one should replace the informal introduction of categories and toposes given in \cite{Lawvere&Schanuel:1997} and \cite{Lawvere&Rosebrugh:2003} by a more formal treatment - or at least to make sure that such a routine formalization is unproblematic and then skip it. After all, so the argument goes, the basic category and topos theory requires proving a number of theorems, and one cannot possibly prove anything without using one system of ``underlying'' logic or another. If one doesn't specify any such system of logic this means that one takes it for granted, perhaps even without knowing about it. To make this underlying logic explicit is a purpose of logical analysis; such an analysis is necessary for clarification and justification of proofs.              

The argument sounds convincing but we may further ask: What justifies one's choice of the underlying logic? and How the principles of one's favorite logic are grounded?  One may reply that these are Big Philosophical Questions that must be not meshed with mathematics; if these questions can by reasonably answered at all the answers cannot possibly be mathematical answers. What mathematician can and must do is to fix some reasonable system of logic on pragmatic grounds and do mathematics with it. Then he may try the same with a different system of logic. 

Lawvere's ambition, as I understand it, \emph{is} to give to the aforementioned Big Questions not only philosophical but also mathematical answers. On the philosophical side he relies onto Hegel's idea of objective logic. And on the mathematical side he uses the category theory for building objective categories, which ground logic (i.e., the subjective logic of inference) and give logic a particular shape (which turns to be variable at a certain degree). True, a \emph{presentation} of objective categorical constructions requires following some basic logical rules. This is particularly true when we are talking about an \emph{axiomatic} presentation. However this only reflects the fact that every presentation is subjective - including the case when the presented thing is itself objective. So the above argument does not really demonstrate the foundational significance of the subjective logic. It only shows that some subjective logic $L$ is needed for presenting mathematical theories. Lawvere's Hegelian requirement is that $L$ must be grounded in mathematical structures, which reflect objective empirical features of our world, and thus count as an aspect of objective logic.   

A possible way out of this dialectical circle is McLarty's bootstrap: in his presentation he uses from the outset a logic without the rule of excluded middle, which is sound in any topos; using this basic logic he introduces categories and toposes and the notion of internal logic of a topos. Thus it turns out that the basic logic used from the beginning of the exposition can be an internal logic of any topos. This allows for developing the general theory of categories and toposes  \emph{internally} in any topos: every theorem of the general ``external'' theory still holds true in every internal version of the theory - albeit the internal version typically comprises theorems, which do not hold externally, i.e., in the general case.  From a logical point of view this is perfectly consistent (notwithstanding my earlier arguments against equating the``external'' with the ``general''). So McLarty's book \cite{McLarty:1992} shows how one can meet formal requirements of Hilbert's Axiomatic Method without buying its usual philosophical underpinning. I would like to stress however that a formal treatment of topos theory like one presented in McLarty's book by itself does not reveal the objective character of this theory and of this logic. In Chapter \textbf{9} I shall interpret Lawvere's work differently and describe the New Axiomatic Method that better captures the original features of Lawvere's axiomatic reasoning.  

\addcontentsline{toc}{chapter}{Conclusion of Part 1}
\chapter*{Conclusion of Part 1}

In his textbooks \cite{Lawvere&Schanuel:1997} and \cite{Lawvere&Rosebrugh:2003} Lawvere and his co-authors stress the objective character of the category theory in general and of the categorical logic in particular by suggestive real-life examples relevant to the practice of modern science; in other places Lawvere also stresses the fact that the category theory has emerged from a wide mathematical practice of ``structural'' mathematics (see Chapter \textbf{8}) and helped to ``unify and concentrate'' this practice in textbooks. All of this qualifies as an evidence for the claim of objectivity of the categorical logic. However much more, in my view, remains to be done in order to re-configure the standard notion of Axiomatic Method, which has been designed on very different philosophical principles. Here Euclid's example can be helpful. Euclid's geometry as presented in Euclid's \emph{Elements} is a systematically organized discipline, which played the role of paradigmatic example of the systematic organization of knowledge for quite a while after it was first created. This remarkable organization was achieved by Euclid not with some ``background logic'', about which he was not aware, or which he decided to hide from us, but rather by internal mathematical means, which I tried to describe in Chapter \textbf{1}. During the time of its flourishing Euclid's geometry also justly qualified as the only empirically justified theory of physical space; as such it served both the fundamental physics and technology. It is only relatively recently that it proved insufficient on both accounts. Hilbert's Formal Axiomatic Method designed as a replacement for Euclid's method of organizing mathematical knowledge (which after Hilbert's suggestion we also conveniently call by the name of Axiomatic Method) does not qualify as an adequate replacement because of its formal character: it makes no difference between the case when a given axiomatic theory contains a valuable piece of knowledge and the case of a purely fictitious theory having to epistemic value at all. True, Hilbert's studies in foundations of mathematics allowed for a great progress in logic through the application of mathematics to systems of formal symbolic logic and formal mathematical theories built with these logical systems. So we have learned a great deal about such formal systems -  including the fact that some key questions about them, which in Hilbert's view could be easily answered, in fact can not. But however valuable this new knowledge might be for its own sake it can not constitute a sufficient basis for designing a general method of building scientific theories and, in particular, of anything deserving the name of foundations of mathematics. Since the Formal Axiomatic Method does not take into account any other features of theories except formal, building a meaningful theory by this method  turns to be a matter of pragmatic, aesthetic or some other prejudged choice - or perhaps just a matter of good luck and the so-called ``scientific intuition''. After Cassirer and Lawvere I believe that this formal turn makes science basically irrational, and that the right choice can be made only theoretically, albeit not merely speculatively.

Thus the problem of replacement of  Euclid's Axiomatic Method remains wide open. Lawvere's work is a unique attempt to solve this problem with the novel mathematical technique of categorical logic. A further progress in this field requires a systematic study of the relevance of the categorical logic in the natural sciences and in the modern technologies, which I cannot include this book. 

\part{Identity and Categorification}
\chapter{Identity in Classical and Constructive Mathematics}

\section {Paradoxes of Identity and Mathematical Doubles}
Changing objects (of any nature) pose a difficulty for the metaphysically-
minded logician known as the \emph{Paradox of Change}. Suppose a green
apple becomes red. If A denotes the apple when green, and B when it
is red then $A = B$ (it is the same thing) but the properties of $A$ and $B$
are different : they have a different color. This is at odds with the Indiscernibility
of Identicals thesis according to which identical things have
identical properties. A radical solution - to explain away and/or dispense
with the notion of change altogether was first proposed by Zeno around 500 BC and remains popular among philosophers (who often appeal to the relativistic spacetime to justify Eleatic arguments). Unlike physics,
mathematics appeared to provide support for the Eleatic position : for
some reason people were more readily brought to accept the idea that
mathematical objects did not change than to accept a similar claim about
physical objects - in spite of the fact that mathematicians had always
talked about variations, motions, transformations, operations and other
process-like notions just as much as physicists.

The Paradox of Change is the common ancestor of a family of paradoxes of identity which might
be called \emph{temporal} because all of them involve objects changing in time. \emph{Chrisippus' Paradox}, \emph{Stature}, \emph{The Ship of Theseus} belong to this family \cite{Deutsch:2002}. However time is not the only cause of troubles about identity: space is another. The \emph{Identity of Indiscernibles} (the thesis dual to that of the \emph{Indiscernibility of Identicals}) says that perfectly like things are identical. According to legend in order to demonstrate this latter thesis, Leibniz challenged a friend during a walk to find a counter-example among the leaves of a tree. Although there are apparently no perfect doubles among material objects, mathematics appears to provide clear instances immediately: think about two (different) points. But the example of geometrical space brings another problem: either the \emph{Identity of Indiscernibles} thesis is false or our idea of perfect doubles like points is incoherent. In what follows I shall refer to this latter problem as the \emph{Paradox of Doubles}. Mathematics looks more susceptible to this paradox than physics. However were she living today, Leibniz's friend might meet his challenge by mentioning the indiscernibility of particles in Quantum Physics \cite{French:1988}

The example of two distinct points $A, B$ (Fig. 5.1) does not, it is usually argued, refute the Identity of
Indiscernibles because the two points have different relational properties: in Fig. 5.1 $A$ lies to the left of $B$ but $B$ does not lie to the left of itself
\footnote{These relational properties of the two points depend on their shared space: the argument
doesn't go through for points living on circle. I owe this remark to John Stachel.}:

\begin{center}
\includegraphics[scale=0.2]{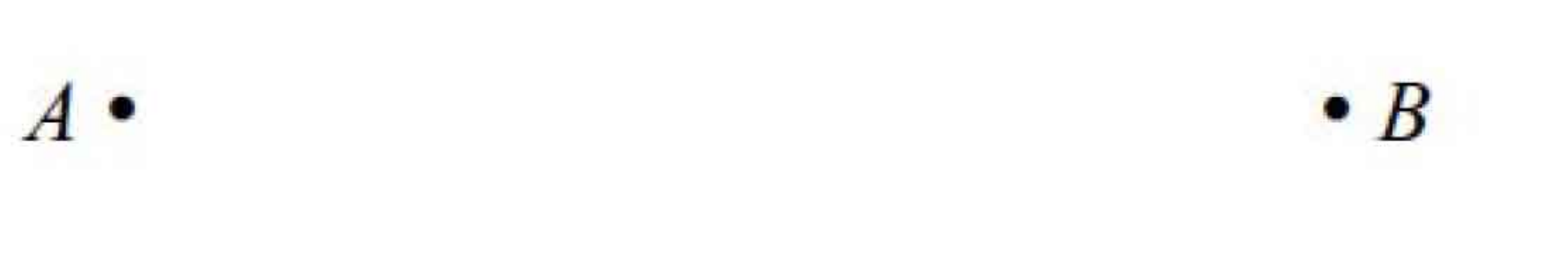}
\end{center}

\begin{center}
Fig. 5.1
\end{center}

The difference in the relational properties of $A$ and $B$ amounts to saying that the two points have
different positions.  However the example can be easily modified to meet the argument. Consider two coincident
points (Fig. 5.2): now $A$ and $B$ have the \emph{same} position.

\begin{center}
\includegraphics[scale=0.33]{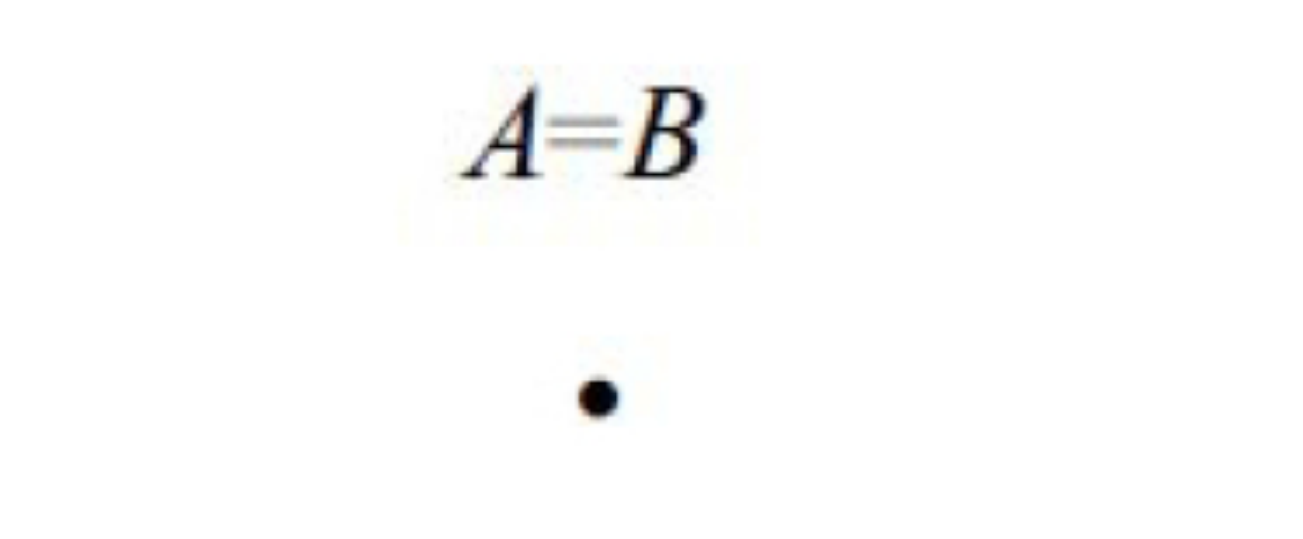}
\end{center}

\begin{center}
Fig. 5.2
\end{center}

It might be argued that coincident points are an exotic case, one which can and should be excluded
from mathematics via its logical regimentation. But this is far from evident - at least if we are
talking about classical Euclidean geometry. For one of the basic concepts of Euclidean geometry is
\emph{congruence}, and this notion (classically understood) presumes coincidence of points: figures $F, G$
are congruent iff by moving $G$ (without changing its shape and its size) one can make $F$ and $G$
coincide point by point.

The fact that geometrical objects may coincide differentiates them significantly from material
solids like chairs or Democritean atoms. The supposed \emph{impenetrability} of material solids counts
essentially in providing their identity conditions \cite{Lucas:1973}. Thus, identity works differently
for material atoms and geometrical points.

This fact shows that Euclidean geometrical space cannot be viewed as a realistic model of the
space of everyday experience as is often assumed. One needs the third dimension of physical space to establish in practice the relation of congruence between (quasi-) 1- and 2-dimensional material objects through the application of a measuring rod or its equivalent.

We see that the alleged contradiction with the Identity of Indiscernibles is not the only difficulty
involved here. Indeed the whole question of identity of points becomes unclear insofar as they are
allowed to coincide. Looking at Fig. 5.2 we have a surprising freedom in interpreting ``='' sign.
Reading ``='' as \underline{identity} we assume that $A$ and $B$ are two different names for the same thing.
Otherwise we may read ``=''  as specifying a \underline{coincidence} relation between the (different) points $A$
and $B$. It is up to us to decide whether we have only one point here or a family of superposed
points. The choice apparently has little or no mathematical sense. One may confuse coincidence
with identity here without any risk of error in proofs. However this does not mean that one can
just assimilate the notions of identity and coincidence. For identity so conceived would be very
ill-behaved, allowing for the merger of different things into one and the splitting of one into many.
(Consider the fact that Euclidean space allows for the coincidence of any point with any other
through a suitable motion.) Perhaps it would be more natural to say instead that the relations of
coincidence and identity while not identical in general, coincide in this context?

For an example from another branch of elementary mathematics consider this equation: $3 = 3$. Just
as in the previous case there are different possible interpretations of the sign ``=''  here. One may
read ``=''  either as identity, assuming that 3 is a unique object, or as a specific relation of equality
which holds between different "doubles" (copies) of 3. Which option is preferable depends on a
given context. There is a unique natural number $x$ such that $2 < x < 4$; $x = 3$. Here ``=''  stands for
identity. But when one thinks about the sum $3 + 3$ or the sequence $3,3,3,...$ it is convenient to think
of the 3s as many. In this latter case $3 = 3$ still holds but now ``=''  is being read as equality rather than identity. Again the choice looks like a matter of convenience rather than of theoretical importance.

Similarly, in one sense \emph{cube} is a particular geometrical object, while in a different sense there exist
(in some suitable sense of ``exist'') many cubes. When one proves that there exist exactly 5
different regular polyhedrons, and says that the cube is one of them, one speaks about the cube in
the first sense. When one considers a geometrical construction, which comprises several cubes,
one thinks about the cubes in the second sense. However no distinction between the two
meanings of the term "cube" can be found in standard textbooks, and it is not even clear whether
such distinction can be sharply made.

The above examples might make one think that the notion of identity simply plays no significant
role in mathematics. $2 \times 2 = 4$ remains true independently of whether the sign "=" is read as equality,
or as identity, whether equality is treated as identity, or identity is weakened to equality. It looks
as if here one may choose one's interpretation according to personal taste or preferred
philosophical position.

However such a liberal attitude to identity in mathematics looks suspicious from the logical point
of view. Claims of existence and uniqueness of mathematical objects satisfying given descriptions
(definitions) play an important role in mathematics. Such a claim means that a given description
indeed picks out (identifies) an object, not just a property. The standard definition of the \emph{unit} of a
given group $G$ is an example. Obviously a claim that such-and-
such an object is unique makes sense only if its identity conditions are fixed. But as we have
seen they may in fact be very loose. It is clear that 3 is the only natural number bigger than 2 and
smaller than 4 but it is not clear that 3 indeed refers to an unique object. But how can mathematics
hang together as a body of knowledge if it apparently does not meet Quine's ``no entity without
identity'' requirement?

The unit of a group $G$ is defined as the element $1\in G$ such that for any element $x\in G$ (including $1$
itself) $1 \otimes x = x \otimes 1 = x$, where $\otimes$ is the group operation. The existence of 1 is guaranteed by
definition but its uniqueness is proved. Suppose $1'$ is another element of the group satisfying the
same condition: $1' \otimes x = x \otimes 1' = x$. Then taking first $x = 1$, and then $x = 1'$ we have $1' \otimes 1=
1 \otimes 1 '= 1 = 1'$.

There are several ways to approach this problem. I now explore them.

\section{Types and Tokens}
The remedy, which readily comes to mind on the part of anyone familiar with contemporary
Analytic metaphysics, is that of the type/token distinction. Consider another example, which
prima facie looks very like the above mathematical cases. There are 26 letters in the English
alphabet, and the letter $a$ is one of them. In the last phrase the letter a is referred to as a particular
thing, namely a particular letter of the alphabet. But in this phrase itself there are five such things.
Hence the letter $a$ is not a particular thing. The standard way of dissolving this puzzle is to say
that here we have one $a$-type and five $a$-tokens.

In explaining the distinction, one starts from tokens: an $a$-token is a piece of paper with
typographic pigment, or another material object (e.g. a piece of printer's type) representing the
letter $a$. Obviously a-tokens are many. The second step is to explain what the $a$-type is. Let me however show instead that the type/token distinction does not fix the problem of identity of mathematics anyway: whatever mathematical types might be they do not correspond to well-distinguishable tokens.

\emph{The} natural number \textbf{3} (which I write in bold for further references) indeed looks like a type but
the 3s, which we find in the series 3,3,.. or in the formula 3+3 do not look like tokens from the
viewpoint of standard examples (like particular chairs). For formula 3+3 may be applied to many
different situations: one might add 3 chairs to 3 chairs, 3 points to 3 points, or even (taking a
liberal attitude) 3 chairs to 3 points. Arguably such application amounts to instantiation of both
3s (in formula 3+3) by certain sets of objects. That is certainly not how good tokens behave: the
fact that types can be instantiated but tokens cannot is essential; if we allow for the instantiation
of tokens by other tokens we either lose the type/token distinction or must provide it with a new
relational sense (which looks like an interesting project but I cannot pursue it here).

The case of points (or more structured geometrical figures like triangles) at first sight looks more
promising. Apparently points are well-distinguishable tokens of the same type. Unlike the case
of numbers it is common in mathematics to denote different point-tokens by different labels such
as $A$ and $B$. However this works only until coincident points are taken into consideration. For in
the case of coincident points we cannot distinguish a singular point-token from a stock of
point-tokens. It is tempting in this case to think of the stock of points as a ``place'' occupied by a
family of singular point-tokens. But this again involves a reiteration of the type/token distinction
on another level as in the case of 3-tokens. Point-locations initially considered as tokens can
themselves be instantiated by second-order tokens stocked there. Once again this destroys the
usual distinction between point-tokens and the point-type. It is a condition of acting as a
(classical) token that the object so acting have determinate identity conditions - as concrete
symbols like printed numerals do. But our hypothetical number- and point-tokens do not meet
this condition. So the type/token distinction (at least in its usual form) does not help us to handle
the identity issue in mathematics. (This also makes me doubt how well it works outside
mathematics.)

\section{Frege and Russell on The Identity of Natural Numbers}
Frege considered it a principal task of his logical reform of arithmetic to provide absolutely
determinate identity conditions for the objects of that science, i.e. for numbers. Referring to the
contemporary situation in this discipline he writes in the Introduction to his \cite{Frege:1884}:
\begin{quote}
How I propose to improve upon it can be no more than indicated in the present work. With numbers
[..] it is a matter of fixing the sense of an identity. (English translation \cite{Frege:1964}, p.Xe)
\end{quote}
Frege makes the following critically important assumption : identity is a general logical concept,
which is not specific to mathematics:
\begin{quote}
It is not only among numbers that the relationship of identity is found. From which it seems to
follow that we ought not to define it specially for the case of numbers. We should expect the concept
of identity to have been fixed first, and that then from it together with the concept of number it
must be possible to deduce when numbers are identical with one another, without there being need for this purpose of a special definition of numerical identity as well. (\emph{Ibid.}, p.74e)
\end{quote}
In a different place \cite{Frege:1903}, \cite{Frege:1962} Frege says clearly that the concept of identity is absolutely stable across all possible domains and contexts:
\begin{quote}
Identity is a relation given to us in such a specific form that it is inconceivable that various forms of
it should occur (p.254 in the edition \cite{Frege:1962}, my translation from German)
\end{quote}

Frege's definition of natural number, as modified by Russell \cite{Russell:1903} later became standard. I present it here informally in Russell's simplified version. Intuitively the number 3 is what all collections consisting of three members (trios) share in common. Now instead of looking for a common form, essence or type of trios let us simply consider all such things together. According to Frege and Russell the collection (class, set) of all trios just \emph{is} the number 3. Similarly for other numbers\footnote{Following Russell \cite{Russell:1903} I use here words \emph{class}, \emph{collection}, and \emph{set} interchangeably ignoring their technical meanings if any. This terminological freedom is helpful for rethinking the concept of set (or class etc.) without smuggling in ready-made solutions through the existing terminology.}.

Isn't this construction circular? Frege and Russell provide the following argument which they
claim allows us to avoid circularity here: given two different collections we may learn whether or
not they have the same number of members without knowing this number and even without the
notion of number itself. It is sufficient to find a one-one correspondence between members of two
given collections. If there is such a correspondence, the two collections comprise the same number
of members, or to avoid any reference to numbers we can say that the two collections are
equivalent. I shall follow current usage in calling this equivalence Humean (see \cite{Hume:1978}, book 1, part 3, sect. 1). Now we check that this relation is indeed an equivalence in the usual sense, and define natural numbers as equivalence classes under this relation.

This definition reduces the question of identity of numbers to that of identity of classes. This
latter question is settled through the axiomatization of set theory in a logical calculus with
identity. Thus Frege's project is realized: it has been seen how the logical concept of identity
applies to numbers. (In fact this does not work that smoothly as I show in \textbf{5.8} below.)
In an axiomatic setting ``identities'' in Quine's sense (that is, identity conditions) of mathematical
objects are provided by an axiom schema of the form
$\forall x \forall y (x = y \leftrightarrow  \ldots)$
called in \cite{Keranen:2001} the Identity Schema (IS).

This does not resolve the identity problem though because any given system of axioms, generally speaking, has multiple models \cite{Benacerraf:1965}. The case of isomorphic models is similar to that of equal numbers or coincident points (naively construed): there are good reasons to think of isomorphic models as one and there is also good reason to think of them as many. So the paradox of mathematical doubles reappears. Thus the logical analysis \`a la Frege-Russell certainly clarifies the mathematical concepts involved but it does not settle the identity issue as Frege believed it did.

In the recent philosophy of mathematics literature the problem of the identity of mathematical
objects is usually considered in the logical setting just mentioned: either as the problem of the
non-uniqueness of the models of a given axiomatic system or as the problem of how to fill in the
Identity Schema. For my present purposes it is important, however, to return to the problem in
its original informal version, which inspired Frege and Russell 100 years ago. Such a return to
the starting point is, in my view, helpful and perhaps necessary if one wishes (as I do) to
consider the Category-theoretic approach to identity discussed in this paper as a viable
alternative to the approach taken by Frege, Russell and their followers. At the first glance the
Frege-Russell proposal concerning the identity issue in mathematics seems judicious and innocent
(and it certainly does not depend upon the rest of their logicist project): to stick to a certain
logical discipline in speaking about identity (everywhere and in particular in mathematics). The
following historical remark shows that this proposal is not so innocent as it might seem.

\section{Plato}
Given a sequence like 3,3,3,.. mathematicians conveniently talk about multiple \emph{copies} of the
same number (similarly about copies of a given set, or space) . Such talk about copies carries
echoes from Plato. A glance at Plato's philosophy of mathematics shows some features which
might be attractive for a mathematician resistant to the logical regimentation of talk of identity in
different contexts proposed by Frege and Russell. If I understand Plato correctly, according to
him identity applies only to the immutable ideas, and only ideas exist. (So Plato's view in this
respect is in accord with Quine's dictum about  \emph{no entity without identity} \cite{Quine:1969}, p. 23.) Material things don't exist but  \emph{become} ( they change, come into and go out of being ) and hence have no proper identities: this is another possible way out of the Paradox of Change. Mathematical things occupy an intermediate position between material stuff and ideas: they involve a weaker sort of becoming and a softer form of identity. In the case of numbers such ``soft identity'' is equality. Things in the three layers of Plato's ontology are partially ordered by ``distorted copying'' where ideas are the maximal elements, mathematical objects are distorted copies of ideas, and material objects are distorted copies of mathematical objects (and hence also of ideas). The distortion of self-identical ideal numbers amounts to their replacement by families of equal mathematical
numbers. For example, there is a unique ideal number \textbf{3} and an indefinite number of its equal
mathematical copies. In other words numbers in mathematics are defined up to equality but not up to identity
\footnote{Plato's philosophy of mathematics should not be confused with the \emph{Mathematical Platonism} in the sense of \cite{Balaguer:1998}, which has little if anything to do with historical Plato. For an introduction to Plato's philosophy of mathematics see \cite{Pritchard:1995}.}.

There are multiple passages where Plato speaks of ``$X$ itself'', ``$X$ (thought of) through itself'' and ``Idea of $X$'' interchangeably or explains the latter through the former. For example in his \emph{Symposium} (210-211) Plato does this with the notion of Beauty, and in \emph{Phedon} (96-103) with number 2. (In this latter dialog Socrates rejects the view that 2 could be thought of as sum of two units pointing to the fact that 2 can be equally obtained through division of given unit into two halves. Since each of the two operations is the reverse of the other none of them can be viewed as bringing 2 about. So one needs to think of the idea of 2 independently of operations of this sort.) I interpret these passages in the sense, which seems me straightforward: the ``identity to itself'' applies to ideas but neither to material things, nor to mathematical things (as they are usually thought of). To see that Plato's ``idea of 2'' is indeed something else than mathematical number see last chapters of Aristotle's Metaphysics where the author criticizes the \emph{Unwritten Doctrine} developed by Plato in the later period of his life \cite{Findlay:1974}. Here the distinction between ideal and mathematical numbers is made explicit. Aristotle stresses the fact that each ideal number is unique while their mathematical copies are many (Met. 987b) and the fact that ideal numbers are not a subject of arithmetical operations (Met. 1081a-1082b).

Thus unlike Frege Plato does not suppose that the notion of identity applies to whatever there is
(or whatever occurs) indiscriminately. Instead Plato thinks of identity as a specific property of
things he calls \emph{ideas} and notices the fact that in mathematics the identity requirement is relaxed (through the talk of an ``intermediate'' character of mathematical objects). In what follows I shall show that this Platonic insight is particularly appealing in the context of our contemporary Category-theoretic mathematics.

Plato hints at the following division of labor: mathematicians work on equalities whilst
philosophers take care of identities. In the case of arithmetic this is exactly what mathematicians
(and philosophers like Frege) have been doing for centuries. In geometry however the situation is
more complicated because equality in this discipline may mean - and historically did mean -
different things.

Euclid uses the term ``equality'' (Greek \emph{ison}) in the sense of equicompositionality (of plane geometrical
figures ) but there are other equivalencies in geometry, which may be considered as better
``working substitutes for identity'': for example congruence, (geometric) similarity, and affinity.
For there is a sense in which the ``same figure'' means a figure of the same shape and the same
size, and there is another sense in which it means only a figure of the same shape, and the notion
of "same shape" can itself also be specified in different ways. In addition geometry unlike
arithmetic allows for the identification of its objects (of geometrical figures) by directly naming
them, usually through naming of their most important points. This allows us to distinguish two
different triangles $ABC$ and $A'B'C'$ which are the ``same'' in any of above senses. There is
apparently no clear argument, which would allow us to choose one of these senses of \emph{the same} as basic and eliminate the others as an abuse of the language. In particular, as I have shown in \textbf{5.1}, the pointwise specification of figures cannot do this job. So the situation in geometry
(even classical geometry!) is exactly that which \cite{Manin:2002} describes for a different purpose:
\begin{quote}
There is no equality in mathematical objects, only equivalences! (p. 8)
\end{quote}

\section{Definitions by Abstraction}
To pursue his project of reducing the various informal meanings of "the same" in mathematics to a
standard notion of identity captured in a universal logic Frege proposed the method of "definition
by abstraction". In his \cite{Frege:1884} Frege gives the following example of such definition:
\begin{quote}
The judgment ``line $a$ is parallel to line $b$'', or, using symbols $a//b$ , can be taken as identity. If we do this, we obtain the concept of direction, and say: ``the direction of line $a$ is identical with the
direction of line $b$''. Thus we replace the symbol // by the more generic symbol =, through removing
what is specific in the content of the former and dividing it between $a$ and $b$. (translation \cite{Frege:1964}, p. 74e)
\end{quote}

Notice that the procedure as described here by Frege involves a change of notation: in the formula
$a = b$ the symbols $a,b$ no longer stand for lines but denote the same direction. Calling this formal
procedure \emph{definition by abstraction} Frege suggests its interpretation. The idea is that the procedure  picks out a property common to all members of a given equivalence class. In \textbf{6.1} I shall show that this procedure can be interpreted differently.

As our earlier quotations from Frege's \cite{Frege:1884} clearly show, in treating an equivalence $E$ ``as identity'' Frege does not mean to replace identity by something else. He aims at the exact
opposite: to introduce identity where mathematicians usually use only equivalencies.
Definition by abstraction is problematic from the logical point of view \cite{Scholz&Schweitzer:1935}. But I want to stress a different point. Even if definition by abstraction were justified logically it would not provide what a mathematician normally looks for. Frege's direction (not to be confused with
orientation!) is hardly an interesting mathematical notion; this concept might play at most an
auxiliary role in geometry and can easily be dispensed with. The idea of a \emph{family} of parallel lines does the same job as Frege's abstract direction but is more convenient and more intuitive.
Similarly it is more convenient to think of a natural number as a family of equal doubles rather
than a unique abstract object. Such abstract numbers would be much like Plato's ideal numbers.
Plato certainly had a point in arguing that such things do not belong to mathematics!
Frege would most likely answer that the question of convenience does not matter because his
proposal is logically justified and the more traditional mathematical practice and parlance is not.

\section{Relative Identity}

The Theory of Relative Identity is a logical innovation due to Geach (\cite{Geach:1972}, ch.7) motivated by
the same sort of mathematical examples as Frege's definition by abstraction. Like Frege Geach
seeks to give a logical sense to mathematical talk ``up to'' a given equivalence $E$ through replacing $E$ by identity but unlike Frege he purports, in doing so, to avoid the introduction of new abstract
objects (which in his view causes unnecessary ontological inflation). The price for the ontological
parsimony is Geach's repudiation of Frege's principle of a unique and absolute identity for the
objects in the domain over which quantified variables range. According to Geach things can be
same in one way while differing in others. For example two printed letters aa are same as a \emph{type} but different as \emph{tokens}. In Geach's view this distinction does not commit us to a-tokens and a-types as entities but presents two different ways of describing the same reality. The unspecified (or \emph{absolute} in Geach's terminology) notion of identity so important for Frege is in Geach's view is incoherent \footnote{For recent discussion see \cite{Deutsch:2002}.}.

Geach's proposal appears to account better for the way the notion of identity is employed in
mathematics since it does not invoke ``directions'' or other mathematically redundant concepts. It
captures particularly well the way the notion of identity is understood in category theory.
According to Baez and Dolan \cite{Baez&Dolan:1998}

\begin{quote}
In a category, two objects can be ``the same in a way'' while still being different (p.7)
\end{quote}

so in category theory the notion of identity is relative in Geach's sense. But from the
logical point of view the notion of relative identity remains highly controversial. Let $x,y$ be
identical in one way but not in another, or in symbols: $Id(x,y) \& \neg Id'(x,y)$. The intended
interpretation assumes that $x$ in the left part of the formula and $x$ in the right part have the same
referent, where this last (italicized) same apparently expresses absolute not relative identity. So
talk of relative identity arguably smuggles in the usual absolute notion of identity anyway. If so,
there seems good reason to take a standard line and reserve the term ``identity'' for absolute
identity.

We see that Plato, Frege and Geach propose three different views of identity in mathematics.
Plato notes that the sense of the ``same'' as applied to mathematical objects and to the ideas is
different: properly speaking, sameness (identity) applies only to ideas while in mathematics
sameness means equality or some other equivalence relation. Although Plato certainly recognizes
essential links between mathematical objects and Ideas (recall the \emph{ideal numbers}) he keeps the two domains apart. Unlike Plato Frege supposes that identity is a purely logical and domain-independent notion, which mathematicians must rely upon in order to talk about the sameness or
difference of mathematical objects, or any other kind at all. Geach's proposal has the opposite
aim: to provide a logical justification for the way of thinking about the (relativized) notions of 
sameness and difference which he takes to be usual in mathematical contexts and then extend it to
contexts outside mathematics.

As Geach puts it

\begin{quote}
Any equivalence relation ... can be used to specify a criterion of relative identity. The procedure is
common enough in mathematics: e.g. there is a certain equivalence relation between ordered pairs of
integers by virtue of which we may say that $x$ and $y$ though distinct ordered pairs, are one and the same rational number. The absolute identity theorist regards this procedure as unrigorous but on a relative identity view it is fully rigorous. (\cite{Geach:1972}, p.249)
\end{quote}

\section{Internal Relations}
In his paper ``The classification of relations'' of 1899 \cite{Russell:1990} Russell says:
\begin{quote}
Mr. Bradley has argued much and hotly against the view that relations are ever purely ``external''. I
am not certain whether I understand what he means by this expression but I think I should be
retaining his phraseology if I described my view as the view that all relations are external. (p.143)
\end{quote}

In arguing that relations are, generally, internal Bradley \cite{Bradley:1922} means roughly the
following: the relata of a given relation cannot, generally, be thought of independently of
each other and of the relation in question. (So relations, if any, such that their relata can be
thought of independently are external.) Bradley makes indeed a stronger claim:

\begin{quote}
Relations exist only in and through a whole, which cannot in the end be resolved into relations and
terms. [..] The opposite view is maintained (as I understand) by Mr. Russell. But for myself, I am
unable to find that Mr. Russell has ever really faced this question (\emph{ib.}, p.127)
\end{quote}
\footnote{
Notice that in Bradley's view there is no duality between external and internal relations since internal relations are not supposed to be defined independently of their relata which would be an absurdity. (See \cite{Hylton:1990} for further discussion.)}

As we can see each of the two authors admits that he hardly understands arguments of the other.
Since Russell's outright rejection of internal relations they have been under great suspicion
amongst Analytic philosophers. Today the neglect of internal relations is not only the
consequence of underlying inclinations in systematic metaphysics but also a matter of available
logical means. For the main tradition of (modern) logical systems is developed in keeping with
Russell's rejection of internal relations, so one may ask whether or not the standard modern
notion of $n$-placed relation as $n$-placed predicate can be possibly understood as internal relation.
Let us see. Consider the standard procedure of interpretation of given relation $R(x, y)$ in given
domain $D$. Here $x, y$ are logical variables which take their values among members of $D$ that is
usually thought of as a class of individuals. When $x, y$ take values $a, b$ from $D$ $R(a, b)$ takes a certain truth value (usually $true$ or $false$) just like function $f(x, y) = x + y$ takes value 3 when $x$ takes value 1 and $y$ takes value 2 (and + is interpreted in the usual way). Noticeably the substitution of $a, b$ for $x, y$ is proceeded uniformly for any binary relation and in this sense it doesn't depend on $R$. To put it in other words the substitution is formal: one first substitutes $a, b$ for $x, y$ and then looks for the true-value of $R(a, b)$. So relata $a, b$ are assumed here independently of $R$. This meshes well with Russell's view according to which all relations are external.

But can $R(x, y)$ be possibly understood as internal? Consider relation $NEXT(m, n)$ between natural numbers which says that number $m$ is followed by number $n$. Arguably natural numbers cannot be correctly thought of without $NEXT$. This means that this relation is internal. But this claim apparently has nothing to do with the order in which $NEXT$ is interpreted: nothing prevents one to pick up numbers 1, 2, substitute them for $m, n$ in $NEXT(m, n)$ and then see that $NEXT (1, 2)$ is true. Thus the logical machinery involved here allegedly has no bearing on the metaphysical controversy between the external and internal understanding of relations. So given relation $R(x, y)$ might be internal as well as external.

However the above argument is not convincing. For it involves an interplay between the formal
analysis of the concept in question and implicit assumptions made about this concept. As far as
we already know what are natural numbers then we can claim, of course, that 1 is followed by 2.
We can also write down this truth in a more fashionable way as $NEXT(1, 2)$. Formal logic is used
in this case for description of a ready-made concept. In such a case logic has no bearing on how
the concept in question is built, and so it is metaphysically neutral. But when logic is used for
concept-building like in the case of foundations of mathematics then specific features of logical
apparatus get directly involved into emerging concepts. In practice the distinction between the
two ways of applying logic can hardly be ever made rigidly: the major application of logic is a
logical \emph{reconstruction} of given background (in particular of common mathematical practice) but not an external description of ready-made concepts nor creation of new concepts from nothing.
Logical reconstruction is making of new concepts from something. I will not elaborate this point
here and only notice that foundations of mathematics obviously involve concept-building even if
it has a descriptive function too (with respect to the common mathematical practice).

Is it possible to stipulate the relation $NEXT$ between natural numbers without assuming a fulfledged
notion of natural number in advance? A positive answer is given with Hilbert's axiomatic
method. One assumes some class of individuals $N$ as domain of binary relation (two-placed
predicate) $NEXT(m, n)$, and stipulates certain formal properties of $NEXT$ as axioms. The idea is
that a system of axioms of this kind will turn abstract individuals into numbers. Or to put it more
accurately, elements of given class $N$ will be thought of as natural numbers as far as they verify
some properly chosen axioms. Think about Peano's arithmetic.

Is this procedure indeed compatible with the internalist account of relations? The answer is not
trivial, in my view. On the one hand, there is obvious reason to think of $NEXT$ introduced
axiomatically as internal: unless $NEXT$ (with its formal properties) is taken into account elements
of $N$ are thought of as abstract individuals but not as numbers. But on the other hand, the
stipulation of relata of $NEXT$ as individuals is incompatible with a strong version of internalism about relations according to which these relata cannot be thought of without its relation at all, not even as abstract ``things'' without properties. So the standard logical apparatus is indeed incapable to represent relations which are internal in this strong sense.

Apparently Bradley defends such a strong version of internalism about relations when he says
that ``a whole [..] cannot [..] be resolved into relations and terms''. True, this radical position
undermines the very notion of relation, so after all Russell's account of relations should be
probably preferred. However Bradley's remark points to a real problem which shows that the
notion of relation (or at least in its Russell's restricted version) is far less powerful than it seems.
Notice that any axiomatic theory built by the standard Formal Axiomatic Method  assumes its objects (for example sets) to be individuals. However we have seen that the identity of basic mathematical objects like points, circles or natural numbers is highly problematic. The blunt
stipulation of such things as individuals doesn't resolve the problem but turns it into a new form:
given two classes $N$ and $N'$ ($N'$ might be a ``copy'' of $N$) both satisfying axioms of arithmetic which of the two classes is the class of natural numbers  \cite{Benacerraf:1965} ? ( See also Chapter \textbf{8} below.) I suppose that in order to get a satisfactory solution of the identity problem in mathematics we should give up the idea that mathematical objects always form classes and look for different ways of getting multiple objects into a whole. In the next Section I analyze the notion of class and show its limits.

\section{Classes}
Sets of chairs or crowds of people are usually considered as a paradigm cases for our thinking
about the notion of \emph{many}. There are different examples though. Think about clouds in the sky or waves at see surface. One can always count persons or chairs or at least in principle so. But one can hardly count clouds and waves. The problem is not that they are too many but that there is no definite criterion for distinguishing one from another. Clouds and waves are certainly many but this kind of many is in general not countable. For a mathematical example think about families of equal numbers or of coincident points: the question of the cardinality of such multiplicities (to choose a term for many with the broadest meaning) is apparently senseless. In what follows I shall specify the sense of ``countable'' relevant to this context. I shall term the wanted concept \emph{weak} countability in order to avoid confusion with countability in the usual set-theoretic sense.

In his \cite{Russell:1903}, Chapter VI , Russell distinguishes between extensional and intensional ``genesis of classes'': the former proceeds through the ``enumeration of terms'' while the former proceeds as follows: one takes a predicate $P(x)$ and considers class $\{x \mid P(x)\}$ consisting of all such $x$ that $P(x)$ is true
\footnote{
The principle according to which for each given predicate $P(x)$ there exists class $\{x \mid P(x)\}$ is called the \emph{Comprehension Principle}.} 
. For example class $\{1,2,3\}$ can be defined either through the direct enumeration of its elements $1, 2, 3$ (extensional genesis) or as the class of natural numbers smaller than $4$ (intensional genesis). According to Russell the extensional genesis of classes through enumeration is possible only when the number of elements (terms) is finite. However Russell claims that this constraint is only ``practical'' and ``psychological'' but not logical and theoretical. In particular he says:

\begin{quote}
[L]ogically, the extensional definition appears to be equally applicable to infinite classes, but
practically, if we were to attempt it, Death would cut short our laudable endeavor before it had
attained its goal. Logically, therefore, extension and intension seem to be on a par. (\emph{ib.}, p 69)
\end{quote}
After claiming the essential equivalence of extensional and intensional viewpoints Russell goes
ahead and claims the priority of extension:
\begin{quote}
A class [..] is essentially to be interpreted in extension. [..] But practically, though not theoretically,
this purely extensional method can only be applied to finite classes. [..] [A]lthough any symbolic
treatment must work largely with class-concepts and intension, classes and extension are logically
more fundamental for the principles of Mathematics. (\emph{ib.} p.81)
\end{quote}

These arguments are not convincing. True, theories often extend domains of possible application
of available practical means through relaxing certain constraints. For example, since Ancient times
people tend to think about distances between celestial bodies and between pebbles on sand on
equal footing. In many cases such theoretical extension works and allows for improvement of
existing practical means; in other cases taking practical constrains into theoretical consideration
allows for improvement of theories (think about Gauss' work in geodesy which motivated his
geometrical discoveries). However I cannot see how this might help to settle the issue of
extension and intension. What Russell says about Death is irrelevant: an immortal god would no
better succeed to accomplish the task of finishing enumeration of an infinite series than a mortal
human because the enumeration of an infinite series has no end. So to the contrary of Russell's opinion, the difficulty of the infinite enumeration is not practical nor psychological but certainly theoretical and logical. Russell refers to the mathematical (Cantorian) notion of infinite set but he misses an essential point of Cantor's invention. In his \cite{Cantor:1883} Cantor says roughly the following. Count 1,2,3,..This counting never ends - not practically nor theoretically - but we may stipulate a new ideal object $\omega$ as the limit of this process just like we stipulate an irrational number $r$ as a limit of a series of its rational approximations. Then $\omega$ can be understood as a number of all (finite) natural numbers, and so the talk about the set of all natural numbers becomes reasonable. Cantor proposes here a specific extension of the usual finitary enumeration, and I don't think that the philosophical distinction between the theory and the practice much clarifies this Cantor's proposal. 

Observe that Cantor's invention has no immediate bearing on the issue of predication,
so Russell's idea that a predicate may bring about anything like Cantorian set (remind that Russell
doesn't distinguish between sets and classes) is a very strong independent hypothesis. The
following development of logic and set theory imposed well-known constraints upon the use of
classes but these commonly accepted constraints, in my view, are not sufficient. 
Set-theoretic ``antinomies'' including Russell's paradox forced Zermelo \cite{Zermelo:1908}, \cite{Heinzmann:1986} to restrict Russell's ``intensive genesis'' through the \emph{Aussonderungsaxiom}, which allows for ``genesis'' of set $S$ with property $P$ only when one is given another set $M$ such that $P$ is \emph{definite} on $M$ (which is tantamount to saying that $P$ has definite truth-value for every element of $M$ ; then $S$ comprises all those elements $x$ of $M$ for which $P(x)$ is \emph{true}). Russell changes his mind about classes already in 1906 \cite{Russell:1906a}, \cite{Russell:1906b}, \cite{Heinzmann:1986} by putting forward \emph{No Class Theory} according to which what one needs in logic is only a \emph{domain} $U$ of individuals but not any longer classes constructed out of $U$. In 1908 \cite{Russell:1908}, \cite{Heinzmann:1986} Russell changes his mind again and puts forward his type theory.

Bernays in his \cite{Bernays:1958} purports to save Russell's early liberal notion of class through a formal distinction between classes and sets. According to Bernays sets are classes having a specific property of being individuals, that is, capable of being elements of other classes. For sets Bernays accepts an improved version of the \emph{Aussonderungsaxiom} (which he proves as a theorem). However classes in Bernays' view are formed by properties ``automatically'', so one even doesn't need a quantifier for it and can simply write $\{x \mid P(x)\}$ to denote the class of all $x$ such that $P(x)$ is true (this class can be a set or a proper class dependently on predicate $P$). Moreover Bernays doesn't exclude the possibility that classes can be produced in other ways not mentioned in his theory \footnote{The distinction between proper classes and sets has been introduced earlier by von Neumann. See Fraenkel 's \emph{Historical Introduction} to \cite{Bernays:1958}, p.32-33}:

\begin{quote}
This point of view suggests also to regard the realm of classes not as fixed domain of individuals but as an open universe, and the rules we shall give for class formation need not to be regarded as limiting the possible formations but as fixing a minimum of admitted processes for class formation. ( \cite{Bernays:1958}, p. 57)
\end{quote}

Bernays' liberal notion of class remains very popular among mathematicians. People have learnt that the notion of set shouldn't be applied without caution but thanks to Bernays they feel free to talk about classes of anything. This has changed the way of thinking even about elementary mathematical concepts. The idea that the Euclidean plane contains the class of all circles would sound completely weird in the 19th century but today's mathematical students usually don't feel any inconvenience about it. In the eyes of many this freedom of thinking about infinite collections (``Cantorian Paradise'') is a very important achievement of mathematics of 20th century.

A usual worry about such extensional representation of mathematical concepts concerns the issue of infinity: why we need such huge collections where we can do well with only few examples? Now I want to stress a different point. It concerns the fact that thinking of, say, circles on the Euclidean plane, as forming a class we are obliged to take circles as full-fledged individuals with definite identity criteria. But as I have tried to show in the beginning of this Chapter such criteria are hardly available.

Before I elaborate on this crucial point let me make a methodological remark. When I criticize the
class-based representation of mathematical concepts I do not assume that there exist the only
right way to represent mathematical concepts. I believe that mathematical concepts are exactly
what we think about them, and that there is a sense in which the same concepts can be
represented differently. The class-based representation is a way of thinking about mathematical
concepts which proved to be in many ways successful. My critical efforts directed against this
approach aim at revealing its hidden assumptions and constraints and at giving place for
alternative approaches, which look more promising. When I say ``circles are not individuals'' I
mean that the class-based representation of circles clashes with what people usually think about
circles in many standard contexts. I recognize that this clash alone provides no strong argument
against the class-based representation: perhaps we should fix the traditional way of thinking
about mathematical objects rather than modern formal methods. However in \textbf{6.2} I shall show that these traditional intuitions support some important contemporary mathematical developments, so in order to promote these developments we need to elaborate on these intuitions rather than rule them out.

\section{Individuals}
Bernays understands the notion of individual in the logical sense as an element of a domain of
quantification, that is, an element of some class. The extensionality property of classes (which
Russell rightly stresses as indispensable) implies that individuals so understood (elements of
classes) must have unproblematic identity criteria. To see this remind how the Axiom of
Extensionality is written in $ZF$:
$$EXT: \forall x \forall y(\forall z(z \in x \leftrightarrow z \in y) \rightarrow x = y) $$
Informally this axiom says that sets are wholly determined by their elements. Although the
identity of sets is introduced in $ZF$ independently of $EXT$ the intuitive appeal of this axiom
certainly depends of the fact that it can be used for ``checking identity'': given two sets one can
check whether or not they are the same through checking their elements
\footnote{
To make $EXT$ into an instance of Keranen's identity schema (\textbf{5.3}) we need to replace the implication
by the biconditional: 
$$EXT': \forall x \forall y(\forall z(z \in x \leftrightarrow z \in y) \leftrightarrow x = y) $$
$EXT'$ is true in $ZF$ but is not used neither as a definition nor as an axiom for the reason of logical parsimony. In fact $ZF$ allows for another instance of the identity schema obtained from $EXT'$ by the reversal of $\in$: 
$$INT : \forall x \forall y(\forall z(x \in z \leftrightarrow y \in z) \leftrightarrow x = y) $$ 
Taking $INT$ as giving the sense of identity brings about rather unusual way of thinking about sets, which I developed in my \cite{Rodin:2003b}.}.

Remark however that in $ZF$ there is no distinction between sets and elements as different types: elements of sets are also sets. So $EXT$ reduces the question about identity of given pair of sets to the question about identity of some other pairs of sets. If given sets $x, y$ are infinite then checking the identity $x = y$ through $EXT$ reduces the problem to checking an infinite number of identities. Prima facie this doesn't look helpful. In fact $EXT$ is helpful for checking identity $x = y$ only when the questions about identity of elements of $x, y$ have obvious answers or at least are easier to answer. If identity of elements of $x, y$ is just as problematic as identity of $x, y$ themselves then $EXT$ looses all of its appeal.

Bernays assumes the extensionality of classes but in order to avoid quantification over classes he
modifies EXT into this open formula
$$EXTCl : \forall z(z \in x \leftrightarrow z \in y) \leftrightarrow x = y)$$
which he uses as the definition of identity (equality) of classes (here $x, y$ are classes while $z$ ranges over sets). So the extensionality of classes in Bernays' account becomes also automatic and doesn't require a special axiom. Anyway $EXTCl$ provides classes by definite identity conditions just like sets. However according to Bernays certain classes (proper classes) cannot be elements of other classes. Why not? Because it is known that making classes elements of other classes in certain cases leads to contradiction. But this is a mere recognition of the fact but not an
explanation of the phenomenon. The colloquial explanation according to which proper classes are
``too big'' or ``over-comprehensive''  \cite{Bernays:1958}) for being elements of something bigger (because there is nothing bigger?) certainly cannot be viewed as satisfactory.

Here is my explanation, which implies a substantial revision of Bernays' point of view. I suppose
that multiplicities like ``all sets'' cannot be viewed as individuals because their elements are not
individuals either and hence have no definite identity conditions. Such multiplicities cannot be
thought as classes (or as elements of other classes) on the pain of loosing the sense of extensionality. Although we can think about all sets in a way we cannot think of all sets as individuals.

Indeed, in the traditional (pre-Cantorian) mathematics the individuation is always finitary and associated with naming: one stipulates, for example, points $A, B, C,..$ (some of which might appear to be identical) but not an infinite set of individual points. This doesn't, of course, preclude one of speaking, say, about ``any point of given line''; the difference with the modern point of view is that this expression doesn't commit one to an infinite set of points. Cantor's notion of infinite set is based on the assumption that individuals
can form not only finite but also infinite collections. In other words he assumes that thinking
about all points of given line we can still think of these points as individuals like $A, B, C$. Cantor
provides the following justification of this view. He shows that a properly generalized procedure
of counting (enumeration) of elements of given set works in the infinite case too. (This applies to
all infinite sets but not only for sets which are countable in the usual technical sense, see about
Bernays' Numeration Theorem below.) This doesn't really prove that elements of infinite sets are
individuals in precisely the same sense in which elements of finite sets like $\{A, B, C\}$ count as individuals but this shows that at least one essential feature of finite sets is preserved in the infinite case, namely, the fact that elements of infinite sets may be brought into one-to-one correspondence just like elements of finite sets $\{A, B, C\}$ and $\{D, E, F\}$. This gives indeed a reason to think of elements of infinite sets as individuals by analogy with the finite case.

I shall call multiplicities having a cardinality \emph{weakly countable} and require classes to be weakly countable. Given this additional requirement for classes I shall call elements of given class \emph{individuals}.  Thus my hypothesis is that weak countability implies (at least a weak form of) individuation. Equating the weak countability with having certain cardinality I take the most liberal attitude possible intended to preserve the whole of Cantorian set theory. More constructively-minded people might prefer to equate the weak countability with the usual countability, or even to insist that infinite enumeration is impossible.

This hypothesis is in accord with Russell's point that all classes are in a certain sense
``denumerable''. Unlike Russell Bernays says nothing about enumeration of classes but 
proves for sets his \emph{Numeration Theorem} ( \cite{Bernays:1958}, p.138) which improves upon Cantor's infinitary enumeration in terms of formal rigor and states that every set has a certain cardinality. The theorem doesn't hold for proper classes. Nevertheless Bernays assumes that proper classes consist of well-distinguishable elements, and that the extensionality property holds for proper classes. In my view this assumption is ungrounded. Just like Russell in \cite{Russell:1903} Bernays apparently thinks that a mere predication brings about some sort rudimentary enumeration. I don't think that this view is tenable.

Consider predicate \emph{human} for example. The collective term \emph{humans} unlike the term \emph{all presently living humans} is not associated with any particular group of people. The expression \emph{all humans} does not make much sense unless it is further specified (we cannot count all future generations). Nevertheless we can speak about humans as a multiplicity. When we talk about sets in mathematics the situation is not different. Multiplicities of all sets or of all singletons don't deserve the name of classes because such multiplicities have no definite cardinalities and hence there is no reason to think of their elements as individuals. (Remind that according to Bernays every set as an element of the class of all sets is an individual.) 

Bernays disqualifies Russell's aforementioned definition of cardinal numbers as classes of
equivalent sets because he wants to define cardinal numbers as sets. Hence the idea to identify cardinal numbers with certain ordinals. This technical solution causes Benacerraf's problem already mentioned: why we should call cardinal number one particular set of given cardinality rather than another? Such a definition of cardinal number differs drastically from Frege's and Russell's earlier proposals discussed in \textbf{5.3} above.

Thus my point is that weak countability required by classes shouldn't be always taken for
granted and expected to be found everywhere in mathematics. As the phenomenon of
``mathematical doubles'' suggests many mathematical objects might be accountable in terms of
internal relations (in particular internal equivalences) which don't allow for considering these
object as full-fledged independent individuals. Moreover the unique multitask notion of
individual (and hence the unique notion of identity) should be likely given up in favor of various
specific structures. Some structures of this sort appear in the Intuitionistic mathematics as we shall now see. 

\section{Extension and Intension}

Consider after Frege \cite{Frege:1892},  \cite{Frege:1952} expressions ``Morning Star'' and ``Evening Star'': they have different \emph{meanings} but \emph{refer} to the same object, namely to planet Venus. Now if we think about Morning Star and Evening Star as \emph{predicates} then the previous remark translates as follows: the two predicates have different \emph{intensions} but one and the same \emph{extension}. For a simple mathematical example think of predicates Equilateral Triangle and Isogonal Triangle (meaning the usual Euclidean figures). Clearly the two predicates have different intensions. However they have the same extension: a given triangle is isogonal if and only if it is equilateral. Importantly, this latter claim is not an immediate consequence of the two definitions but a theorem based on some further geometrical assumptions (cf. Theorems 5, 6 of of the First Book of Euclid's \emph{Elements}).    

The intension/extension distinction has a long history that I shall not try to overview here. However I shall make few remarks about the fate of this traditional distinction in the 20th century logic. The most characteristic feature of logic developed during this period is its \emph{formal} character described in  \textbf{2.3} above.  The formalization of logic in the beginning of the 20th century in works of Frege, Russell and their followers put extension-related and intension-related logical notions in very unequal positions. Fitting \cite{Fitting:2006-2011} describes the situation as follows:

\begin{quote}
In classical first-order logic intension plays no role. It is extensional by design
since primarily it evolved to model the reasoning needed in mathematics.
\end{quote} 

By classical first-order logic Fitting means, of course, the formal predicate logic of Frege and Russell. The second sentence should be understood in the context of another remark of the author, according to which

\begin{quote}
Mathematics is typically extensional throughout - we happily write ``$3 + 2 = 2 + 3$'' even
though the two terms involved may differ in meaning.
\end{quote}

Thus since mathematics is extensional throughout and the modern (classical) logic is designed to model the mathematical reasoning intension has no place in this logic. So when intension shows up in a non-mathematical context the classical logic doesn't help to account for it and as a result the notion of intension appears to be problematic or even mysterious. 

I agree with Fitting, of course, that Frege-Russell's logic gives no place for intension but I don't quite agree with his explanation of this fact. First of all I cannot see that ``mathematics is extensional throughout''. The equality  $3 + 2 = 2 + 3$, as well as the above example of equilateral and isogonal triangles, demonstrates, in my view, that the intension/extension distinction \emph{is} relevant to mathematics: expressions $3 + 2$ and $2 + 3$ have different meanings but refer to the same object, namely to the number 5. So by writing  ``$3 + 2 = 2 + 3$'' one does \emph{not} necessarily ignore the difference between meanings of the two terms. Similarly, the theorem according to which a given triangle is equilateral if and only if it is isogonal does not involve and does not require the ignorance of the fact that concepts of being equilateral and being isogonal are different (in their intensions)! On the contrary, if being equilateral would simply \emph{mean} being isogonal then there were no theorem. Similarly, if the expression ``Morning Star'' would simply \emph{mean} the same thing as the expression ``Evening Star'' then the claim that Morning Star and Evening Star is one and the same planet would be trivial while in fact it is not. Thus I cannot find any significant difference between mathematical and real life examples in this respect.  In my view the principle reason why early systems of formal logic (like Frege-Russell logic) are purely extensional (in the sense that they don't treat intension explicitly) is rather the following.

In the real life
\footnote{ By the ``real life'' I understand here the domain of intended applications of logic. I assume that this domain is sufficiently large to include usual linguistic examples as well as simple mathematical examples.}
some objects are used as \emph{signs}, which \emph{refer} to some other objects. This setting is straightforwardly modeled by mathematical means with taking some mathematical objects to be sings and some other mathematical objects to be referents of those signs. In this way one gets both a formal \emph{syntax} and formal \emph{semantics}.  Mathematics simply reflects in this case a broader ``real life'' situation providing us with important epistemic and practical advantages similar to those obtained through application of mathematics in natural sciences and technology.

Let us now see how intension can be taken into a formal account. A way to do this is to assume after Frege that in addition to signs and their referents there exist such things as \emph{meanings} and then find an appropriate symbolic representation of meanings. Let $M_{P}$  and $R_{P}$ be the meaning and the reference of a given predicate $P$. Symbolic expressions ``$M_{P}$'' and ``$R_{P}$'' stand for (i.e., \emph{refer to}) some mathematical objects  ($M_{P}$ and  $R_{P}$) which represent the meaning and the referent of $P$ correspondingly: these mathematical objects are specified by \emph{semantics} of our logic. Remind that the symbolic expressions ``$M_{P}$''  and ``$R_{P}$'' not only \emph{refer} but also \emph{mean} something. If, for example, the meaning of $P$ is represented by some specific function then the meaning of $M_{P}$ includes a definition of this function. The meaning of $P$ and the meaning of $M_{P}$ are, generally, different.  

Now observe that a purely extensional formal logic is \emph{natural} in a sense, in which any intensional logic construed along the above lines is not. While a purely extensional formal logic simply reflects a broader ``real life'' context an intensional logic construed as above does not do this because it maps real life meanings to mathematical objects (i.e., referents of certain concepts) but not to mathematical meanings. In other words, an intensional logic so construed treats intension by extensional means, and in a sense reduces the former to the latter. In my view this lack of naturalness is a reason why in a formal mathematized setting the logical intension always looks more problematic than the logical extension.      

Although the above argument does not rule out the possibility of representing meanings by certain mathematical objects (examples of intensional logics designed in this way are given in \cite{Fitting:2006-2011}) it points to a general difficulty of this approach and suggests this obvious alternative: to represent real life meanings directly by mathematical meanings rather then by mathematical objects of some sort. The problem is, of course, that we don't really know what is a mathematical meaning: it appears to be no less elusive than the real life meaning. As we shall now see this straightforward way of formalizing intension can nevertheless work with the price of a substantial reconsideration of Frege's meaning/reference distinction. 

Frege (see \cite{Frege:1892}, \cite{Frege:1952}) has introduced his distinction between meaning and reference trying to clarify questions like this: Is the Morning Star the same thing as the Evening Star or not? Frege's answer is roughly the following: as it stands the question is ambiguous and has no definite answer; in order to get a yes-know answer one should specify whether one asks this question about the meaning of ``Morning Star'' (resp. ``Evening Star'') or about the reference of  ``Morning Star'' (resp. ``Evening Star''). Thus the initial ambiguous question splits into two questions each of which has a definite answer: (i) Do expressions ``Morning Star'' and ``Evening Star'' have the same meaning? (No); (ii) Do expressions ``Morning Star'' and ``Evening Star'' have the same referent? (Yes).  Notice that this solution is in line with Frege's idea that the identity relation is unique and applies indiscriminately to meanings, referents and what not. Now consider the following alternative solution. One agrees that the question ``Is Morning Star the same thing as Evening Star?'' has no yes-no answer but instead of distinguishing between meaning and referents one distinguishes between two different identity relations by saying that there is a sense of identity in which Morning Star and Evening Star are the same and there is another sense of identity (i.e., a different identity relation) in which Morning Star and Evening Star are different.  One wants to speak here about the identity of meaning (intensional identity) and the identity of reference (extensional identity) without positing meanings and references as separate entities. This sounds like  Geach's relative identity (\textbf{5.6}) but now I am talking about a very different approach. In the next Section I show how this idea is realized formally in Martin-L\"of's type theory. As we shall see this theory allows for representing intensional aspects of real world examples by intensional mathematical notions. 

\section{Identity in the Intuitionistic Type Theory}

The intuitionistic type theory with dependent types developed by Martin-L\"of \cite{Martin-Lof:1984}  involves two kinds of identity relations
\footnote{The standard version of this theory involves four different kinds of identity (\cite{Martin-Lof:1984}, page 59). Following Awodey and Warren \cite{Awodey&Warren:2009}, \cite{Awodey:2010} I simplify the original account by deliberately confusing some syntactic and semantical aspects. Then we are left with the following two forms of identity described above in the main text. The version of type theory presented in \cite{Warren:2008} applies the definitional equality also to \emph{contexts}. For simplicity I don't consider the type-theoretic notion of context in this paper. For a more recent exposition of Martin-L\"of's theory and discussion of related philosophical issues see \cite{Granstrom:2011} and \cite{Sommaruga:2000}}.
First, we have here a notion of \emph{definitional equality} of types (written $A=B$) and of terms belonging to the same given type ($x = y : A$). Corresponding rules assure that the definitional equality is an equivalence relation and that definitionally equal types and terms are mutually interchangeable through substitution in the usual way. Second, we have here a notion of \emph{propositional equality} $Id_{A}(x, y)$ that reads as a proposition saying that objects $x,y$ of type $A$ are equal. This second kind of identity (equality) does not apply to types. It does not apply to terms belonging to different types either. 

If two terms $x,y$ of given type $A$ are definitionally equal they are interchangeable through substitution and hence also propositionally equal. If the converse is also the case (i.e. if any propositionally equal terms are definitionally equal) the corresponding version of the theory is called \emph{extensional}; otherwise it is called \emph{intensional}. In the extensional theory the difference between the two kinds of identity is trivial: even if it can be formally maintained it is wholly redundant from a pragmatic viewpoint. However in the intensional theory this difference turns to be fruitful  and mathematically non-trivial as we shall later see.     

Frege's logic accounts for forms of reasoning about things as they (supposedly) \emph{are} without accounting for how we come to \emph{know} these things and without accounting  (at least in the special case of identity statements) for the \emph{justification} of one's claims. Martin-L\"of's constructive logic in its turn is designed as an instrument of inquiry, which never appeals to an entity without specifying explicitly the way, in which one comes to know this entity. The name \emph{constructive} refers here to the fundamental epistemic assumption behind this logic according to which the best way to know a thing is to construct (or perhaps \emph{reconstruct}) it. Since the Intuitionistic Type theory is designed for dealing primarily with purely mathematical reasoning rather than with reasoning in natural sciences, this epistemic assumption reduces to the maxim according to which any mathematical object (that one may want to consider for purely mathematical or some other purposes) must be explicitly constructed rather than simply found somewhere in Nature or on the Platonic Heaven. However a similar constructive approach in natural sciences is well-known too; it dates back at least to Kant's First Critique \cite{Kant:1998} and has a continuing history afterwards. For a version of the constructive approach in natural sciences, which is developed against the background of the 20th century physics, I refer the reader to Fraassen's doctrine of \emph{Constructive Empiricism} \cite{Fraassen:1980}. Whether or not the Intuitionistic Type theory or some other system of constructive logic can be indeed applied in physics and other sciences remains an open question, which I shall not tackle here. Instead I present below an alternative informal analysis of Frege's   \emph{Venus} example, which shows that the non-standard notion of identity used in the Intuitionistic Type theory makes good sense in the real life too.   

Suppose that at certain point of history astronomers observe what they call the Morning Star ($MS$) and the Evening Star ($ES$). Behind these two names there is, of course, a strong epistemic assumption according to which every morning they observe one and the same object ($MS$) and every evening they observe another object ($ES$). In the following story I shall not try to make explicit possible grounds of this assumption but simply take it granted. 

So far our astronomers take $MS$ and $ES$ to be two different objects. Now suppose  that they get a new evidence, which suggests that $MS$ and $ES$ are different appearances of one and the same planet. It may be, for example, a data obtained from a new telescope. With these new data in hand the astronomers develop a new theory, which accounts for the new and the old data about $MS$ and $ES$ and on this basis provides a theoretical ground for the statement that $MS$ is identical to $ES$. On this basis they make a new linguistic convention by replacing the older names $MS$ and $ES$ by the new single name ``Venus''.

In terms of Intuitionistic Type theory the above setting can be formalized rather straightforwardly. First, one needs a type $A$ (for ``astronomy'') of observable brighting spots (``stars'') $s_{i}$ on the sky and some notion of identity of such things, which allow one to identify a given star $s$ observed today with the \emph{same} star observed yesterday. In spite of the fact that this identity is epistemically non-trivial  in the given context we take it to be \emph{definitional}. Then for every pair of stars $s_{1}, s_{2}$ we form another (dependent) type $Id_{A}(s_{1}, s_{2}$) elements (terms) of which (if any) are evidences (proofs) that $s_{1}$ and $s_{2}$ are identical. This latter notion of identity is \emph{propositional} and it should not be confused with the former (definitional) identity. When $s_{1} and s_{2}$ are definitionally identical they are also propositionally identical: in this case the definitional identity $s_{1} = s_{2}$ plays the role of evidence (proof) of their propositional identity. However $Id_{A}(s_{1}, s_{2}$ may be also inhabited, i.e., $s_{1}, s_{2}$ can be propositionally identical, when they are not definitionally identical. The \emph{Venus} example is a case in point: although $MS$ and $ES$ are definitionally different they are propositionally the same. Thus type  $Id_{A}(MS, ES)$ must be inhabited by an independent evidence like one obtained by our astronomers with the new telescope. For a further reference I denote this evidence $E_{1}$. 

We see how the distinction between the definitional identity and the propositional identity can apply to a physically meaningful context. This is however is not yet the end of our story: higher identity types of the Intuitionistic Type theory suggest the following development.  

Suppose that another group of astronomers makes independent observations trying either to confirm or to refute the claim of the first group that $MS$ and $ES$ is the same planet Venus. This other research group obtains evidence $E_{2}$, which in fact supports the claim of the first group. Now the question is whether the second group has simply repeated the observations made by the first group or obtained a genuinely new evidence supporting the claim. This latter question has an epistemic impact on the claim of identity of $MS$ and $ES$ and hence must be treated in the same context: if the two evidences are independent this provides a stronger support to the claim. The type  $Id_{A}(MS, ES)$ of evidences (of identity of $MS$ and $ES$) is equipped with a definitional identity of such evidences. However it may once again turn out that two evidences, which are definitionally different, turn to be propositionally the same. For that reason we form higher type  $Id_{Id_{A}(MS, ES)}(E_{1}, E_{2})$ and see whether or not it is inhabited. If it is inhabited by several (definitionally) different terms one may need to consider identity types of the third order, and so on. In \textbf{6.9} I show that the complex structure, which arises in this way, may be described as a \emph{groupoid of identities}. Here I want to suggest that it better reflects the complexity of empirical inquiry than Frege's universal notion identity, which wholly ignores this complexity.  

A Fregean may argue that what I discuss here is an epistemological issue (How we come to know an identity statement?) which is not Frege's problem. I don't think that this claim is historically correct. The following quote shows that epistemic concerns make part of Frege's inquiry into the notion of identity (equality):  

\begin{quote}
Equality gives rise to challenging questions which are not altogether easy to answer. Is it a relation? A relation between objects, or between names or signs of objects? In my \emph{Begriffsschrift} I assumed the latter. The reasons which seem to favor this are the following: $a = a$ and $a = b$ are obviously statements of differing cognitive value; $a = a$ holds a priori and, according to Kant, is to be labeled analytic, while statements of the form $a = b$ often contain very valuable extensions of our knowledge and cannot always be established a priori. The discovery that the rising sun is not new every morning, but always the same, was one of the most fertile astronomical discoveries. Even to-day the identification of a small planet or a comet is not always a matter of course. (\cite{Frege:1952}, 56)
\end{quote}

I cannot see that Frege provides a satisfactory answer to this epistemic concern. His theory of meaning and reference applies equally in the context where the meaning and the reference are known and in the context where  these things are merely assumed; it doesn't reflect the epistemic difference between the two situation. The multi-level type-theoretic propositional identity considered above does this job by requiring an explicit proof of (or reason) why certain things are equal and by making such proofs (reasons) into proper elements of the identity of the given entity. The topological interpretation of this construction presented in \textbf{6.9} will allow for an intuitive grasp of this construction.     

The above discussion provides an insight onto this crucial question: In which sense if any the concept of identity is fundamental?  First of all we should admit that unless  \emph{some} identities are fixed one cannot communicate (and arguably cannot even produce) any coherent thought. This is clear because any language - no matter how ``informal'' and how metaphoric - requires a \emph{local stability} of certain phonetic, written, syntactic and (last but not least) semantic patterns. In order to use a language one needs to recognize and reproduce patterns of all these sorts.  In any natural language such a local stability is combined with a non-trivial global dynamics: languages evolve in time and interact in space.  Fixing semantic patterns requires fixing identities of entities, which are not purely linguistic (like Sun, Morning Star and what not). The scientific discourse demonstrates a similar contrast between the local stability and the global dynamics but in this case the contrast is sharper: locally (say, in a given scientific publication) the meaning of terms is fixed more rigidly than the meaning of words in the common speech while globally the scientific language evolves more rapidly than the everyday language (since science itself evolves more rapidly than more traditional human institutions like religion, family, etc.). Thus any reasonable discourse and, in particular, any scientific discourse requires fixing some basic identities some of which are linguistic and some of which are not.      

These fixed identities must be locally stable in the following sense: once one assumes an identity $a = b$ one is forced to preserve it until the accomplishing of the given piece of reasoning. What exactly counts as an accomplished piece of reasoning is a tricky question, which I shall not treat here systematically. In the process of inquiry it often happens that one rejects certain earlier made assumptions (for example, in the light of new evidences). Such revisions may be or be not justified but in any event they are not regulated by rules of logical inference, which are supposed to draw some consequences from given assumptions without changing these assumptions. In this sense any logical reasoning is conservative: it preserves assumptions and reaches some new conclusions. (I understand now ``logic'' in the usual sense, which does not include the dialectical logic and similar non-conservative schemes.) Thus  there is no way to assume $a = b$ and then argue that $a$ and $b$ are different: this would be a sheer contradiction. 

This conservativity does not mean, however, that all identities are (locally) fixed ``once and for all''. Some identities may be not found among assumptions of the given reasoning but be obtained as conclusions of this reasoning. In this respect identities, or more precisely \emph{propositions} expressing identities, behave just like all other propositions. Now we should take into account that assumptions may play different roles in reasoning. Often one wants to chose as assumptions some fundamental propositions, which may serve as grounds for some further conclusions but cannot be themselves obtained as conclusions made on some ``deeper'' grounds (at least at the given stage of knowledge). In this case the assumptions qualify as ``first principles''. However in other cases assumptions serve rather for bootstrapping one's reasoning; in such cases they are no more ``fundamental'' than conclusions (even if their logical role is not the same as the role of conclusions). Definitional identities in the Intuitionistic Type theory seem to belong to this second category of assumptions. The fact that the Morning Star seen yesterday and the Morning Star seen today is one and the same ``star'' is an assumption, which on the one hand, fixes a linguistic convention, and on the other hand, identifies a continuing series of regularly observed phenomena as the appearance of one and the same object (the Morning Star). \emph{Mutatis mutandis} this applies to the Evening Star too. The fact that we take the identities of the Morning Star and of the Evening Star to be \emph{definitional} reflects a stage of knowledge: this is where the given reasoning begins; the reasons why one series of phenomena is made into the Morning Star and the other series of phenomena is made into the Evening Star are left behind. Then one obtains an evidence that the Morning Star and the Evening Star is one and the same planet Venus. This latter identity in the given context qualifies as propositional. It reflects a further stage of knowledge. Arguably it is more fundamental than the former definitional identities. (The true reason why the Morning Star observed yesterday is the same thing as the Morning Star observed today is the fact the Morning Star is Venus!). Notice that the progress of knowledge achieved with the new propositional identity ($MS = ES$) is strictly conservative in the sense that it does not require to revise the assumed definitional identities. (The discovery that $MS = ES$ does not contradict the fact that $MS$ observed yesterday and $MS$ observed today is one and the same object.) This allows one to incorporate the change of knowledge into an inferential  logical scheme.

\chapter{Identity Through Change, Category Theory and Homotopy Theory}
\section{Relations versus Transformations}

The replacement of the equivalence $xEy$ by the identity $x = y$ discussed by Frege (\textbf{5.5}) allows for an interpretation, which differs from Frege's. Namely, equivalence $E$ can be understood as an invertible  \emph{transformation} (rather than relation), which turns $x$ into $y$ and vice versa; then the identity = becomes the identity \emph{through} this transformation. If $E$ is the relation of Euclidean congruence then the corresponding transformation is the (Euclidean) motion; thinking  about $E$ as motion (rather than congruence) one says that  $y$ is the same object $x$ but subject to translation and/or rotation in the Euclidean space. Here $x$ and $y$ are said to be the same in the same sense of ``same'' in which, for example, an adult yesterday and today is the same person. So we think here geometrical figures in much the way we think of a substantial continuant - as an entity capable changing its states and/or positions. Such a ``substantialist'' interpretation works also for Frege's example of parallel lines
\footnote{For a more up-to-date account of the notion of substance and of identity through change see \cite{Wiggins:1980}.}.

The substantialist reinterpretation of mathematical relations may look like an exercise in oldfashioned
metaphysics but it appears surprisingly fruitful from the mathematical point of view.
For in mathematics the language of transformations is not formally equivalent to that of relations
as one might expect but is actually far richer. Given equivalence $xEy$ there are, generally speaking,
many distinguishable transformations turning $x$ into $y$ while $xEy$ only says that one such
transformation exists. So here the underlying naive metaphysics matters mathematically.
The difference becomes particularly evident in the case of (global) invertible  transformations of a
given geometrical space. In the language of relations the existence of such transformations
amounts only to the claim that a given space is equivalent to itself. But in fact such
transformations contain the most basic information about the corresponding space. This was first
recognized by Klein \cite{Klein:1872} when he formulated a new research program in geometry (known today as ``Erlangen Program'') as follows:
\begin{quote}
Given a manifold and a group of transformations on it one should investigate the structures on the
manifold with respect to those properties that respect the transformations of the group. (\cite{Klein:1872}, p.7, my translation from German)
\end{quote}

It is not the notion of a substantial form surviving through transformations that is the major
issue in the new framework for the study of geometrical structure proposed by Klein. Rather
there is something of a different sort, which also remain unchanged through the transformations.
That something is the structure(s) or forms of the transformations themselves. I refer to the fact
that invertible  geometrical transformations like Euclidean motions form algebraic groups under
composition. This fact remains completely hidden from view when one uses the language of
relations. Thus the traditional metaphysics of substance and form fulfills a mathematical need
which the new Frege-Russell metaphysics does not - whatever might be said in favor of the latter
against the former for philosophical reasons.

Let me next specify some terminology, which will be useful for what follows. We have considered
three different ways of thinking of what is involved in operating with an (arbitrary) equivalence
relation $xEy$:

\begin{enumerate}
\item Extension:  Consider equivalence classes formed of those things equivalent under the relation $E$
\item Abstraction: Replace the relation $xEy$ by identity $x = y$, and read $x, y$ anew as standing for a
(relational) property common to all and only members of the same equivalence class under $E$
\item Substantivation: Think of the given relation as an invertible  transformation of relata into each
other, and read $E$ as identity through this transformation
\end{enumerate}

In the case of Humean relation $H$ one may proceed from 1) to 3) through the following steps.
Given certain class of classes $x, y, ..$ equivalent by $H$

\begin{itemize}
\item think of the one-one correspondences between elements of given classes $x, y$ as invertible 
transformations (isomorphisms) $f, g...$ turning elements of $x$ into elements of $y$ and conversely
(the invertibility implies that different elements of $x$ turn into different elements of $y$ and vice
versa);
\item think of $x, y$ as different states of the same underlying substratum $X$, and think of
(auto)morphisms $f, g,...$ as \emph{changes} of $X$;
\item similarly identify all classes equivalent to $x$ and $y$ with $X$.
\end{itemize}

A non-trivial fact, which makes mathematical sense of this metaphysical exercise, is that the
automorphisms of $X$ form a \emph{group} called its \emph{permutation} group or \emph{symmetric} group. To see better what we gain and what we might lose in switching from relations to transformations consider the following table:

\begin {center}
\begin{tabular}{|p{8cm}|p{8cm}|}
  \hline
  \textbf{Extensional reading} & \textbf{Substantional reading} \\
  \hline
  Write $x \sim y$ for ``class $x$ is equivalent (isomorphic) to class $y$'' & Write $f: X\rightarrow X$ or simply $f$ for an isomorphism from a class $X$ to itself (automorphism) \\
   \hline
  $\sim$ is an equivalence relation. & Automorphisms of $X$ form a group. \\
   \hline
  $\sim$ is transitive, i.e., $x \sim y$ and $y \sim z$ implies $x \sim z$. & Given automorphisms $f, g$ there exists a  unique automorphism $fg$ resulting from the application of $g$ after $f$.\\
   \hline
  $\sim$ is reflexive, i.e., every class $x$ is isomorphic to itself: $x \sim x$. & There exist an identity automorphism $1$ such
that for any automorphism $f$ we have $1f = f1 = f$. \\
  \hline
  $\sim$ is symmetric, i.e.,  if $x \sim y$ then $y \sim x$. & every atomorphism $f$ has an inverse $f^{-1}$ such that $ff^{-1} = f^{-1}f = 1$.\\
 \hline
 \end{tabular}
\end {center}

Let me now comment on each raw of this table separately.  

\underline{Raw 1} \\
Classes $x, y$ from the left column are identified in the right column through  Frege's abstraction and denoted by the same symbol ``$X$''.
Notice that $x \sim y$ is a proposition but $f: X \rightarrow X$ is a (mathematical) object, namely a particular morphism (function). Proposition $x \sim y$ \emph{says} that there exists an isomorphism between $x$ and $y$, while $f$ \emph{is} such an isomorphism. It is helpful to forget for the moment about the abstraction and think of $f$ as an isomorphism of the form $x \rightarrow y$. Then the translation from the left to the right cell of this raw can be described as an \emph{instantiation}: while the left cell tells us that an isomorphism of certain form exists the right cell points to such an isomorphism. When $x$ and $y$ are identified through abstraction $f$ turns into an automorphism of the form $X \rightarrow X$. The instantiation of a given concept provides a concrete instance (concrete object) that falls under this given concept. This shows that translation from the left to the right column involves a double conceptual transformation, namely, it involves concretization (of the notion of isomorphism) along with abstraction (over given isomorphic classes). 

Given an object of certain type one may always claim that an object of this type \emph{exists} (in some appropriate sense of ``exists'' - we are now talking about the existence of mathematical objects and the present argument does not depend on any particular theory of mathematical existence). But one may be also in a position to describe further properties of the given object, including those properties, which this given object does not share with all other objects of the same type. In other words a given object may also have some detectable \emph{specific} properties.  As we shall briefly see this is a case in point: the language of transformations allows one not only to claim that certain isomorphisms exist but also to describe specific properties of such isomorphisms.  This shows that the right cell contains some information, which is not found the left cell. However if $x \sim y$ does \emph{not} hold we still have a proposition, which tells us something useful. Such information cannot be provided by means used in the right column: given no authomorphism of the appropriate type one has nothing to talk about here (unless one brings into consideration some further relevant objects like morphisms of more general sorts).  This shows that translation between the language of relations and the language of morphisms is not wholly transparent  in either direction.

\underline{Raw 2} \\
As the following comments make clear there is an interesting conceptual link between the notion of equivalence relation, on the one hand, and the notion of (algebraic) group, on the other hand. In order to show this I take the standard definition of equivalence relation (as transitive, reflexive and symmetric relation) and compare this definition with a category-theoretic (rather than standard set-theoretic) definition of group as a \emph{group of transformations}. The following three raws of the table establish a piecewise correspondence between the two definitions.   

\underline{Raw 3} \\
In order to see more precisely how works the translation from the left to the right cell it is once again helpful to begin with the instantiation. In the left cell we have $x \sim y$ and $y \sim z$; by instantiation we get $f: x\rightarrow y$ and $g: y\rightarrow z$ correspondingly. By the transitivity property we also have $x \sim z$ (in the left cell), which by instantiation gives us $fg: x\rightarrow z$ (in the right cell). Now identifying $x, y, z$ through abstraction we get the situation presented in the right cell. Observe however that  the transitivity property of $\sim$ does not reflect the fact that composition $fg$ is uniquely defined  by $f$ and $g$. 

\underline{Raw 4} \\
The reflexivity property of relation = amounts to the fact that for any given class $x$ admits at least one automorphism. Generally, classes admit many authomorphisms. While every element of any given abstract class $x$ is indiscernible from any other element of the same class, if any, (in the sense that the concept of class doesn't assume that every element of class has some particular properties, which may allow one to distinguish this given element of the given class from another element of the same class) automorphisms $A_{i}$ of a given class $x$ are not all alike (except the trivial case when family $\{A_{i}\}$ consists of a single element; notice that by reflexivity of relation $\sim$ family $\{A_{i}\}$ cannot be empty). For $\{A_{i}\}$ contains a distinguished automrphism called \emph{identity} automorphism, which sends every \emph{element} of $x$ into itself (but not only $x$ to itself). Other automorphisms of $x$ permute its elements and so don't have this property. The same property of the identity authomorphism can be described in terms of composition as this is shown in the right cell; this definition implies that the identity automorphism is indeed unique.  

Proposition $x \sim x$ in the left cell tells us that the family $\{A_{i}\}$ of automorphisms of class $x$ is not empty. The right cell provides a concrete instance of such automorphism, namely $1$. This concrete instance has a specific property (namely, the property of being \emph{identity} automorphism), which automorphisms of a given class do not have in general. This specific property is not chosen arbitrary because every  class has the identity automorphism but not not every class has automorphisms of other sorts (the empty class and the class consisting of a single element do not). Thus, once again, one can observe that the right cell contains more information than the left: while the left cell tells us only that every class admits an automorphism the right cell specifies that every class admits an automorphism with a specific property.
 
 Let's now see how abstraction works in this case. Given a family $\{x_{i}\}$ of isomorphic classes it is easy too check that families of automorphisms $\{A_{i}\}$ corresponding to these classes are also isomorphic in the sense that there exists a one-to-one correspondence between members of $\{A_{i}\}$ and $\{A_{j}\}$. However this notion of isomorphism turns to be too weak to make good mathematical sense in this latter case. When isomorphic classes $x_{i}, x_{j}$ are identified through abstraction it doesn't matter which element of $x_{i}$ is identified with which elements of $x_{j}$. In other words, such an identification can be made through any one-to-one correspondence between the elements of the two classes. This is because elements of a given class  are indiscernible in the sense that the concept of class doesn't assume that an element of class has some specific properties, which may allow one to distinguish this given element from any other element of the given class. However each family $A_{i}$ contains a distinguished automorphism $id_{i}$ called \emph{identity automorphism}, which sends every element of $x_{i}$ into itself. Otherwise this distinguished automorphism can be defined in terms of composition (without referring to elements) as this is shown in the right cell and called the \emph{unit} (of the given group of automorphisms). If families $A_{i}$ and $A_{j}$ are identified merely as classes through some arbitrary one-to-one correspondence between their elements the difference between identity automorphisms and other automorphisms is ignored, so what we get is again an abstract class rather than an abstract family of \emph{automorphisms}. To prevent this it is necessary to identify the identity $id_{i}$ of the first family with the identity $id_{j}$ of the second family rather than with any other member of this latter family. Unless an one-to-one correspondence between members of $\{A_{i}\}$ and $\{A_{j}\}$ has the specific property just mentioned it hardly deserves to be called isomorphism and cannot be used for an appropriate abstraction. This property (avoiding confusion of identity automorphisms with other automorphisms) is necessary but not sufficient for formulating the new notion of isomorphism. What one needs here is, of course, the notion of \emph{group isomorphism} defined in \textbf{8.5} below. Automorphisms of a given class is not just another class or a family but a \emph{group} in the sense explained in the right column of the table.

\underline{Raw 5} \\ 
The symmetry property of relation $\sim$ described in the left cell consists of the following: given an isomorphism of the form $x \rightarrow y$ there exists an isomorphism of the form $y \rightarrow x$. The notion of \emph{inverse} automorphism not only provides an instance of isomorphism of the latter form (modulo the identification of $x, y$ through abstraction) but also describes a specific property of this isomorphism.  

Does the approach outlined above provide any viable alternative to Frege's project of settling the
question of identity in mathematics by external logical means? Prima facie it seems that the notion
of identity through change (transformation) invoked here remains completely informal and not
likely to be helpful in avoiding paradoxes mentioned in \textbf{5.1}. However I claim we have here a
new formal concept of identity as the unity of a group of transformations. This group-theoretic
notion of identity meshes well with the metaphysical intuition that any changing entity contains a
core invariant through changes. Merging equivalence classes $x, y, ..$ into one ``class-substance'' $X$ indiscriminately, we recover a notion of identity as a particular transformation (and one unique
for a given group) which we may speculate is connected with the notion of \emph{repetition}. Here the notion of ``'repetition of the same'' is thought of as giving meaning to the notion of "the same'' rather than the other way round (cf. \cite{Deleuze:1968}). This group-theoretic identity is obviously relative in Geach's sense: objects can be identical up to a transformation of one type but different up to a transformation of a different type (\textbf{5.6}). It is not immediately clear whether this group-theoretic identity has anything to do with the logical notion of identity, which was Frege's concern. But at least we get a well-defined identity concept here, and one which makes the metaphysical intuitions behind it precise.

There are at least three objections, which can be brought against the suggestion that we should take this group-theoretic notion of identity as a serious candidate for the philosophical explication of the notion of identity either inside or outside mathematics.

\begin{enumerate}
\item
The logical (and metaphysical) notion of identity should apply to the widest possible domain
of entities, so one can say which things in a given domain are the same and which differ. But
group-theoretic identity is relevant to a single object, namely its group.
\item
Group-theoretic identity $1_{G}$ does not allow us to form propositions like $A = B$ ``$A$ is identical to $B$''. Generally, the group-theoretic identity $1_{G}$ like any other element of a given group $G$ is a particular mathematical object while identity is a basic logical concept. In the group theory like elsewhere in mathematics the role of logical identity is played by mathematical equality =. For the sake of the argument we can now ignore subtle differences between the logical identity and mathematical equality discussed in the last Chapter. Anyway $1_{G}$ and = have little if anything in common except the common name and some vague metaphysical intuitions behind it. For theoretical reason we need to distinguish the two things sharply and reflect the distinction in the terminology rather then allow ourselves to be led by confusing terms and the vague metaphysics.
\item
In particular the group-theoretic identity $1_{G}$ like any mathematical object needs certain identity
conditions. These identity conditions matter essentially, for example, when one proves the uniqueness of the identity of a given group. Proposition ``there exist $1_{G}$ such that for any $f\in G$ we have $1_{G}f = f1_{G} = f'$''  takes the logical identity (equality) = for granted. Hence one needs a prior logical notion of identity to cope with the group-theoretic notion of identity, so no way the latter can be a candidate for the normalization or mathematical explication of the more general notion.
\end{enumerate}

In what follows we will see that these problems can be partly fixed through generalizing the
concept of group up to that of category.

\section{How To Think Circle}

Arguably the best way to explain what is circle is to show one:

\begin{center}
\includegraphics[scale=0.4]{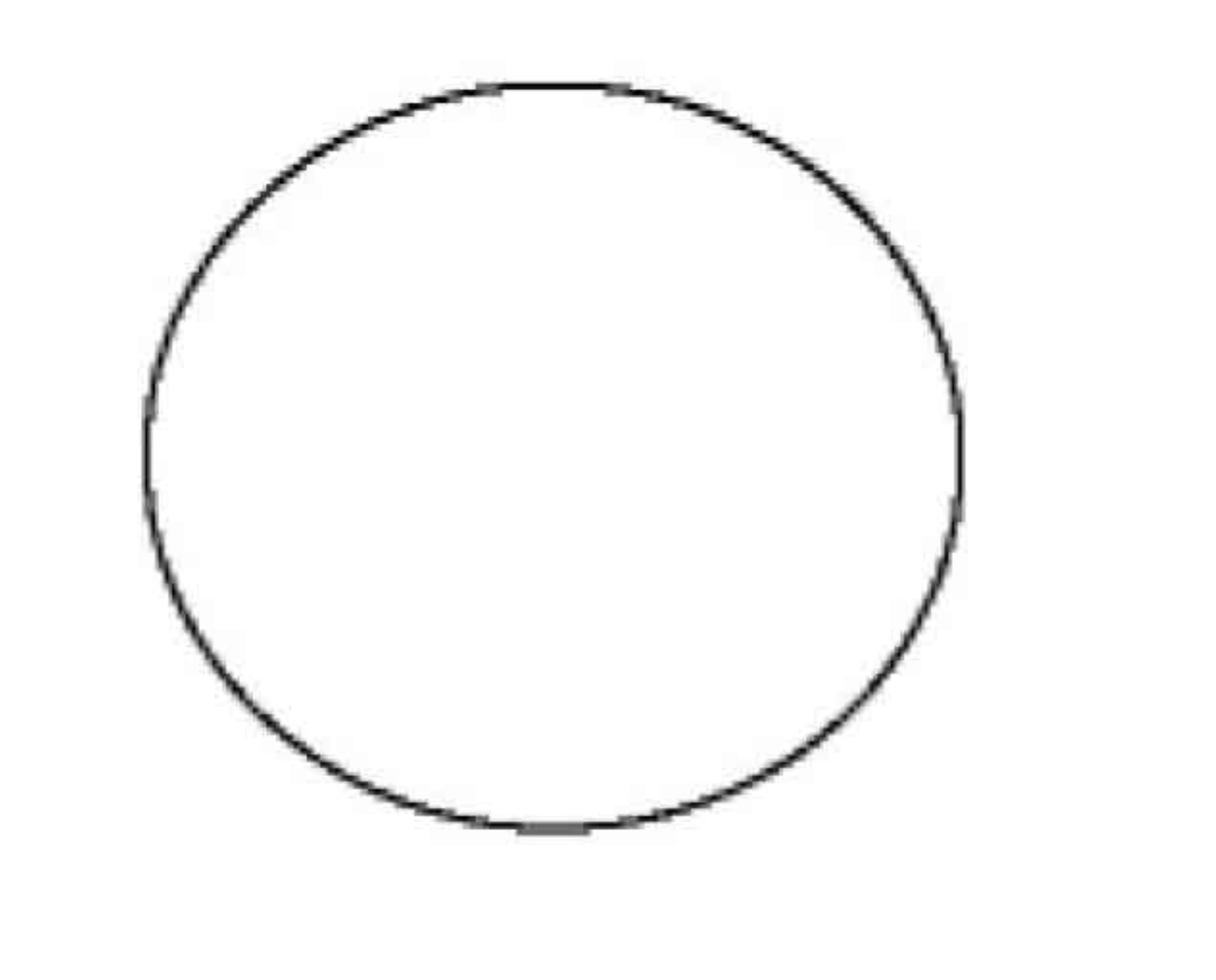}
\end{center}

\begin{center}
Fig. 6.1
\end{center}

However this wouldn't work for one who tries to grasp the notion of circle from scratch because the above circle has some specific features (like its size and its position on this book's page), which other circles may not have; so a beginner cannot possibly learn from this or any other single example which features of the given picture are generic and must be taken into consideration when one thinks and talks of circles, and which in the given context are superficial and thus must be ignored. The above picture alone doesn't allow one to see it \emph{as circle} unless one is already familiar with the \emph{concept} of circle, i.e., already knows how circles look like. 
Only when a learner is shown multiple figures and told which of them do qualify as circles (in spite of differences in size, color, position, etc.) and which do not qualify as such (in spite of some resemblance with circles like in the case of ovals) he may learn how correctly identify circles among figures of different shapes. 

The above story about learning the concept of circle reveals the following important feature of this and other similar concepts. Every particular circle represents the concept of circle just as well as any other particular circle; there is no sense in which a given circle $c_{1}$ is a better or worse circle than another given circle $c_{2}$: \emph{as circles}  $c_{1}$ and $c_{2}$ are strictly equivalent. Yet it is an essential feature of circles that such things are \emph{many} and so this feature is an essential element of the circle concept. It is equally essential that circles belong to a broader genus of things like a \emph{geometrical figure} and that this broader genus contains things of other types like ovals, etc\footnote{In Plato's view the distinction between better and worse ``copies'' of the generic ``ideal circle'' makes sense. However I assume here that at least among mathematical circles no similar distinction can be made: every mathematical circle is a circle, period.}.

What are general concepts and what are their particular instances and how concepts relate to their instances is a question that has been discussed in philosophy throughout it history at least since Plato's times. The modern distinction between \emph{types} and \emph{tokens} as well as the modern Fregean notion of abstraction (as the identification of members of equivalence classes modulo some equivalence relation) discussed in the last Chapter are rather traditional in this respect even if they come with some modern symbolic techniques. My purpose of mentioning here this very general issue is to point to a relatively new way of tackling it, which originates from the 19th century geometry.

Informally the idea is the following. Given two arbitrary items $c_{1}$ and $c_{2}$, which are supposed to instantiate the same concept $C$ (think about circles) consider transformations of the form $f: c_{1} \rightarrow c_{2}$ (which may be or be not invertible). Then by specifying a \emph{type} $T$ of such transformations and by choosing a single  \emph{generic} instance $c^{*}$ one may describe the whole extension $E_{C}$ of concept $C$ by saying that it consists of those and only those items, which are obtainable from $c^{*}$ through some transformation of type $T$. (Notice that if these transformations are all invertible then $E_{C}$ can be similarly obtained from any its element, i.e., every element of $E_{C}$ turns to be generic. Otherwise this, generally, is not the case.)  For example, it may be specified that by moving a given circle and scaling it (changing its size without changing the shape) one always obtains a circle while transformations of other sorts turn a circle into something else. Such a description (appropriately improved) is sufficient for understanding what is circle in general and thus it fairly presents the general concept of circle\footnote{
Since motions and scalings  are invertible any circle is generic. In order to see how this condition may brake assume a convention according to which a point counts as circle (of zero radius) and scaling (in the more liberal sense of the term than above) allows for shrinking a circle (of a non-zero radius) into a point. Such conventions are abundant in mathematics and there is nothing wrong with them. However it is clear that a zero-radius circle (point) is not generic because scaling (liberally understood) and moving of a point always brings back a point, and so no circle with a positive radius can be obtained from a point through these transformations.}.

One may object that the suggested way of thinking about the concept of circle through transformations is not independent because it requires the specification of an appropriate \emph{type} of transformations: thinking about circles as generated by a generic circle through appropriate transformations one no longer needs the circle-type but still needs a transformation-type. This is a faire point but there are in fact ways to tackle the problem. First, transformations can be specified in terms of composition with other transformations, see \textbf{6.1} above. So the notion of \emph{composition} of transformation provides means for distinguishing some transformations from some others. Second, one may think about types of transformations again in terms of transformation (of the second level). Consider transformation  $f: c_{1} \rightarrow c_{2}$ taking one circle into another and another transformation  $g: c_{1} \rightarrow c_{2}$ of the same type. Then consider a second-level transformation $\alpha : f \rightarrow g$ which takes $f$ into $g$. Geometrically $f$ and $g$ can be thought as paths in a geometrical space (for simplicity consider the case when $f$ and $g$ are motions); then $\alpha$ is a transformation of one path into another. Then one may specify an appropriate type for second-level transformations, point to a generic first-level transformation $f^{*}$ (or a number of such generic transformations) and finally describe first-level transformations as those transformations, which are obtained from $f^{*}$ through second-level transformations of the appropriate type. A similar procedure can be applied to second-level transformations. The fact that we get here an infinite regress from a mathematical viewpoint doesn't make the whole procedure pointless: one may think of an infinite structure involving transformations of $n$-th without assuming that $n$ is bounded. In \textbf{6.9} I describe a mathematical structure, which realizes this idea in a precise form. 

Let me now come back to the problem of identity. When an entity undergoes a transformation it always remains unclear  (at least at the linguistic level) whether this transformation produces a new different entity or only a different state of the same entity. This gives rise to paradoxes about identity similar to those mentioned in \textbf{5.1}. As I have already stressed in \textbf{5.1} in mathematics such ambiguities about identity are also ubiquitous. Instead of trying to fix this problem by imposing some external logical regimentation I shall rather explore the possibility of using transformations themselves for it. An \emph{identity transformation}, which leaves a transformed object as it is without producing any change in it may be thought of as an auxiliary formal notion like zero or empty set. However one may also tentatively think of it as a means, which \emph{determines} the identity of its object. More generally one may try to use transformations (in some mathematically refined sense of the term) for controlling identity. Think again about circles. As far as we are talking about Euclidean circles living on Euclidean plane there is a natural criterion of identity (leaving now aside the problem of \emph{coinciding} figures stressed in \textbf{5.1}), which in terms of transformations can be described as follows: rotations of a given circle about its center and the identity transformation count as \emph{self}-transformations while all other available transformations (always of type $T$) transform the given circle into another circle. If one now manages to distinguish identity transformations and appropriate rotations from other transformations of type $T$ without using identity conditions for circles established otherwise then such a distinction between different kinds of transformations may serve for introducing such identity conditions independently. Let me now bring into discussion some more mathematics and see whether this project is viable.   

\section{Categorification}
Circles together with their mutual transformations give us a toy example of a \emph{category} in the sense of Chapter \textbf{4}, i.e., in the general \emph{mathematical} sense of the term. 

Generally, a category comprises:
\begin{itemize}
\item
Class of its \emph{objects} $A, B, C, . . .$;
\item
For each ordered pair of objects $A, B$ class of \emph{morphisms}
$f: A \rightarrow B, g: A  \rightarrow B, . .$; given $f: A \rightarrow B$, $A$ is
called \emph{domain} of $f$ and $B$ is called \emph{codomain} of $f$;
\item
Composition $fg$ of morphisms $f, g$ such that the codomain of $f$
equals the domain of $g$ (see the diagram below); the composition
is associative : $h(gf) = (hg)f = hgf$;
\item
Identity morphisms $1_{A}$ associated with each object $A$ and defined
by the following condition : for all morphisms $f, g$, $1_{A} = f$ and $g1_{A} = g$
(provided the compositions $1_{A}f$ and $g1_{A}$ exist).
\end{itemize}
 
When in a categorical diagram  any arrow $A \rightarrow C$ equals to any other arrow between objects $A$ and $C$ obtained through composition of other arrows shown at this diagram  the
diagram is said to be \emph{commutative}. For example, saying this triangle
$$\xymatrix{&B\ar[dr]^g\\ A\ar[ur]^f\ar[rr]_h && C}$$
is commutative is simply tantamount to saying that $fg = h$. Morphisms resulting from composition of shown morphisms can be omitted at a commutative diagram when this doesn't lead to an ambiguity. For example, saying this square
$$\xymatrix{A \ar[r]^g & B \\ C \ar[r]_h \ar[u]^f & D\ar[u]_i}$$
is commutative is tantamount to saying that $fg = hi$.

Now the above construction with circles can be described as category $C$ where objects are circles and morphisms are mutual transformations of circles (motions and scalings) some of which transform a given circle into itself while some other transform a given circle into another circle. Composition of such transformations understood in the usual way obviously satisfies the conditions mentioned in the above definition of the notion of category. Notice that our category of circles $C$ has the following additional property not assumed in the general definition of category : all its morphisms
(transformations) are invertible . The invertibility (aka \emph{reversibility}) is a basic property
of all usual geometrical transformations (like motion or scale transformation) in virtue of which such transformations form groups. 

In the  category-theoretic terms just introduced the invertibility of transformation
(morphism) $f: A \rightarrow B$ amounts to existence of transformation
(morphism) $g: B \rightarrow A$ (called the \emph{inverse} of $f$) such that $fg = 1_{A}$ and
$gf = 1_{B}$. In category theory this property is taken as the definition
of isomorphism, so isomorphisms are invertible morphisms by definition. A category
like $C$ such that all its morphisms are isomorphisms is called \emph{groupoid}.
Thinking of objects of a groupoid ``up to isomorphism'' one gets a \emph{group}.
(So group is a category with only one object such that all its morphisms
are isomorphisms.) However such identification causes a lost of information,
namely the lost of distinction between morphisms of objects to
themselves (automorphisms) and morphisms of objects to other objects.
Thus groupoids provide an important counter-example against the widespread
belief according to which in categories all isomorphic objects can
be always viewed as identical.

The full strength of the notion of category is revealed through the case when morphisms between
objects are not all invertible, that is, are not all isomorphisms. A basic example is the category of
sets having sets as objects and functions between sets as morphisms (see \textbf{4.1} above). Further examples are
obtained through equipping sets with various structures like group structure or topological
structure. Then morphisms are required to ``preserve'' or ``respect'' the corresponding structure:
so in the category of groups morphisms are homomorphisms of groups, and in the category of
topological spaces morphisms are continuous transformations. (The precise definitions are given below in \textbf{8.5} where the idea of ``preservation of structure'' is critically reconsidered.) Using these common examples one should not forget that categories of structured sets don't cover all categories of interest as shows the example of Grothendieck topos from \textbf{4.9}. Our circle category $C$ also belongs to this latter sort (unless a circle is construed as a structured set of its points).

Thus the upgrade of the notion of group up to that of category involves two independent steps:
\begin{enumerate}
\item
introduction of multiple identities (multiple objects) instead of unique identity (unique object);
\item
allowing for non-invertible  morphisms. 
\end{enumerate}

This upgrade can be shown with the following diagram
\footnote{
In the standard set-theoretical setting monoid is defined as set $M$ provided with a binary
operation $\otimes$ and unit (identity element) $1$. The
existence of inverse elements is not required.
}:

\begin{center}
\includegraphics[scale=0.5]{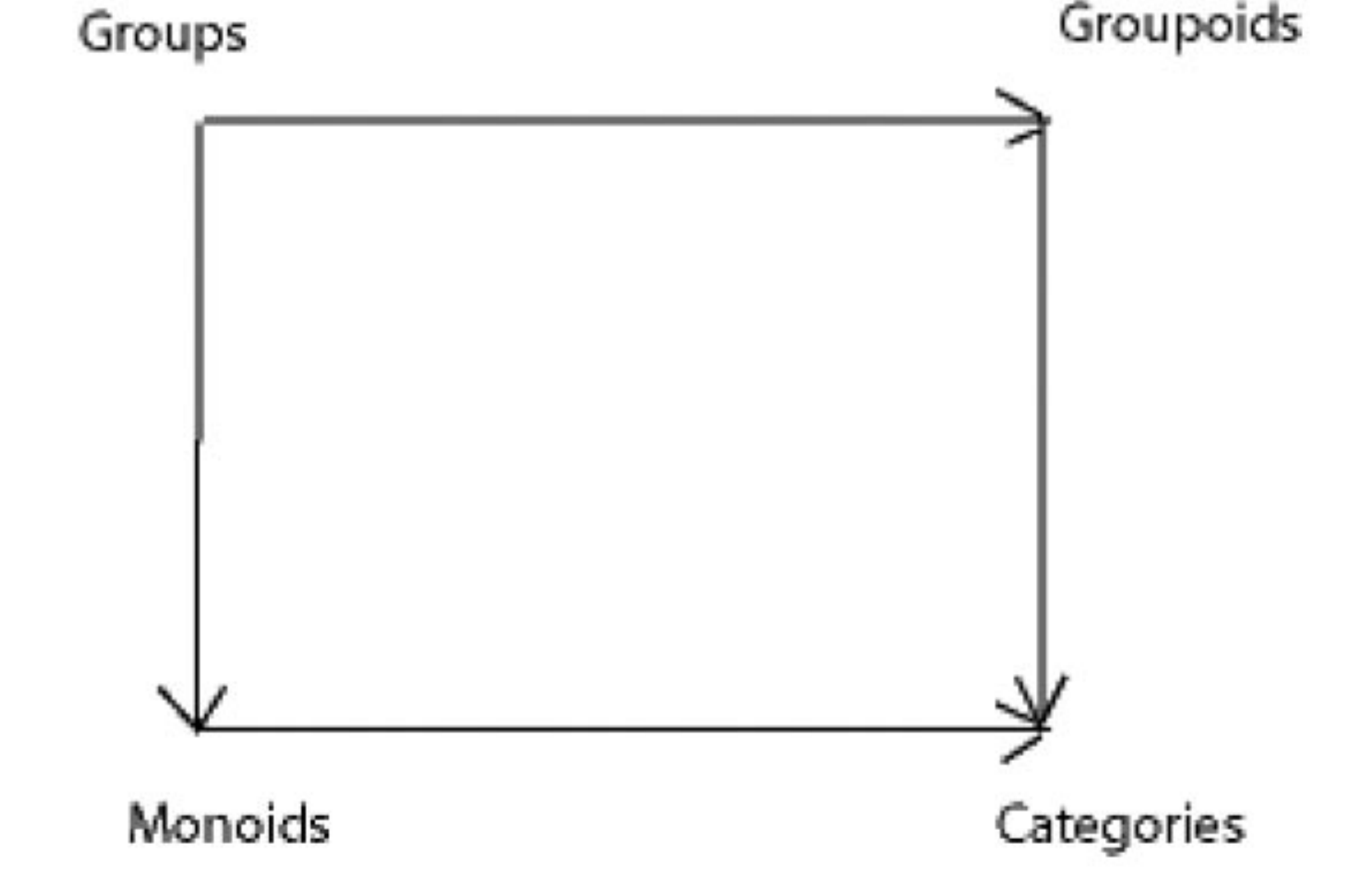}
\end{center}

Examples of categories given so far are concrete categories. This means that objects of such
categories are specified in advance (usually this means that they are construed \`a la Bourbaki as
structured sets), so a category could be seen as a structure over and above given class of its
specific objects. However category theory allows for a different approach: starting with the
general notion of category one specifies its algebraic properties to the effect that the structure of
morphisms between objects and their compositions determines properties of these objects. The
specification of given abstract category amounts to the requirement that certain morphisms exist
and certain diagrams commute. As it has been already explained in \textbf{4.1} a properly specified abstract category ``turns into'' the category of sets \cite{Lawvere:1964} in the sense analogous to that, in which logical variables in axiomatic systems like ZF turn into sets under its intended interpretation. 

At the early stage of category theory people often opposed categorical foundations of set theory to the standard foundations as ``external'' approach to ``internal''. The idea is that while in the case of standard foundations sets are reconstructed though their elements, that is, ``from inside'', in categorical foundations sets are taken as black boxes interacting through morphisms (functions), so what sets are is ultimately determined in ``sociological'' terms of their mutual behavior. This is a right point as far as it concerns basic intuitions about sets but from the formal point of view such intuitions are not essential. The relation of membership $x\in y$ taken as basic in ZF and its likes can be read in both senses -  from the left to the right and from the right to the left - and this makes no formal difference (although the intuition behind the extensionality axiom makes the former ``internal'' reading preferable). At the same time it is not correct that the categorical approach doesn't allow one to ``look inside an object'': in particular the relation of membership can be perfectly reconstructed by categorical means. In both cases objects of given theory (in particular sets) are first taken as abstract
individuals and then ``interactions'' between the objects tell us ``what these objects are''. A real
difference between the two approaches (and this is my second remark) concerns how exactly
these ``interactions'' are accounted mathematically. In the standard axiomatic approach they are
interpreted as relations, and relations in their turn are formalized as predicates (like the two-place
predicate $\in$). In the categorical approach the ``interactions'' are accounted for as morphisms
(transformations)
\footnote{
As Lawvere  puts it
\begin{quote}
The crystallized philosophical discoveries which still propel our subject include the idea that a
category of objects of thought is not specified until one has specified the category of maps which
transform these objects into one another and by means of which they can be compared and
distinguished. (\cite{Lawvere:1991}, p.1)
\end{quote}
Words ``category of objects of thought'' apparently paraphrase Cantor's famous definition of
set as ``a collection into a whole of [..] objects of our intuition or our thought'', see \cite{Cantor:1885}.
}.

Since categories represent concepts in a specific way the mathematical use of the term ``category'' is after all not in odd with how this term has been used throughout the history of philosophy.

\section{Are Identity Morphisms Logical?}

In the last Chapter we considered the complicated interplay between the ``usual'' mathematical equality and the ``usual'' (Fregean) logical identity. Then in \textbf{6.1} we introduced into the play the notion of identity transformation (as the unit of a group of transformations) and finally in \textbf{6.3} generalized this latter notion up to that of identity morphism in a category. In \textbf{6.1} we also suggested that the notion of identity \emph{transformation}, on the one hand, and the notion of identity \emph{relation}, on the other hand, are two different representations of the same pre-theoretical identity concept. We compared these two notions  on the basis of common metaphysical intuitions about identity including intuitions about the identity through change. Admittedly such a base for comparison is very shaky. One may argue then that the very idea to treat the identity relation and identity transformations (identity morphisms) on equal footing is obviously wrong. The argument goes as follows. The identity relation is a fundamental logical concept while identity morphisms are specific mathematical objects. In order to define the category-theoretic notion of morphism in general, and the notion of identity morphism in particular, one applies a logical notion of identity anyway. This allegedly shows that the identity relation makes part of foundations while identity morphisms do not
\footnote{For a general logicist critique of foundational pretenses of category theory see \cite{Feferman:1977}.}.

One way of responding to this critique is to emphasize the ``geometrical point of view'' \cite{Marquis:2009} on category theory against the ``logical point of view'' \cite{Quine:1953} without assuming that the two points of view are incompatible. I believe, however, that there is indeed a strong and important philosophical opposition between the two viewpoints, which needs to be resolved through a dialectical \emph{Aufhebung} rather than downplayed through some peaceful reconciliation of these viewpoints. And I also believe that we are now already well prepared for taking the logicist challenge. So I can reply as follows. 

True, the common intuitions about the identity through change cannot by themselves replace a rigorous logical concept of identity. However when these intuitions are put in a rigorous geometrical form they become strong candidates for representing objective features of our world and may help us to build a working identity concept appropriate for natural sciences. For reasons discussed in the Introduction and in the proceedings Chapters  I strongly reject the notion according to which some ``primitive'' logical notion of identity formed independently from our best science and our best mathematics can possibly qualify as a part of \emph{foundations} of these disciplines. I believe after Cassirer that the proper function of logic and mathematics ``is only within the empirical science itself'' (\cite{Cassirer:1907}, p. 43-44). This approach squares well with Hegel's notion of \emph{objective logic}, which inspires Lawvere's work in categorical logic (\textbf{4.8}). So what I propose to do with the naive intuitions about identity through time is not to ``clarify'' them with some independently designed logical machinery but rather use them for developing a new logical machinery that could better serve the needs of today's science. I claim that the relevance of these pre-theoretical intuitions is measured  \emph{only} by their role in the modern scientific practice but \emph{not} by their role in the everyday thinking and in the everyday language. Since the common intuition about the identity through change apparently continues to play a role in today's sciences (whether we are talking about the identity of living organisms or the identity of quantum particles), I take this notion seriously and believe that in an appropriately mathematized form it should make part of the scientific logic. Since category theory supports this intuition I consider this theory as an appropriate mathematical vehicle for it. However I also take into consideration the fact that the relevance of traditional intuitions in the modern science is very limited, and that in many important contexts such intuitions are manifestly inadequate and become obstacles for a further scientific progress (as in quantum physics). So my aim is not to save the old intuitions by putting them in a new mathematical form (as often do Analytic metaphysicians) but rather to modify these intuitions with the new mathematics in such a way, which would make them more adequate to today's science (see Chapter \textbf{7} below). Rather than take the \emph{logical} notion of identity as primitive and study its formal properties I want to construe this notion in terms of some geometrical categories, which qualify as candidates for the objective representation of space, time and quantity.    

Since the general notion of category can hardly qualify as such a candidate the categorical notion of identity morphism cannot  replace the standard logical identity relation. Even if this replacement makes some sense in the special context described in \textbf{6.1} it obviously makes no sense in other contexts including the context of foundations of category theory. So in order to push the project forward we need some neater ideas. In fact I have already described one such idea, namely, Lawvere's idea of internalization of the equality relation in the categorical logic (\textbf{4.5}). In the next Section I present a development of this idea in B\'enabou's work. In (\textbf{6.9}) we shall see a reappearance of this idea in a geometrical (more precisely, \emph{homotopic}) interpretation of Martin-L\"of's type theory called \emph{homotopy type theory}. In homotopy type theory morphisms of an appropriate category represent the logical notion of identity as it appears in Martin-L\"of's type theory, which is much richer than the standard one (see \textbf{5.11} above).      

\section{Fibred Categories}

The following discussion is based on \cite{Benabou:1985} by Jean B\'enabou. The idea is
the following. Recall that categories have been introduced in \textbf{6.3}
 as classes of a certain kind (classes of objects and morphisms). Which properties of classes are used in such a ``naive'' category theory? Let category $C$ be our ``object of study''
and category $B$ be our ``optical instrument'' for studying $C$. $B$ can be
thought of as category $S$ of sets; however we can also consider different
possibilities, in particular abstractly defined toposes. Following B\'enabou
I shall call objects of $B$ \emph{sets} (remembering that they could be somewhat
different than usual sets) and call classes of morphisms or objects of $C$
\emph{families}. (In what follows families will reappear as multiplicities of a
different sort than classes.) Now given a set $I$ (an object of $B$) we may
master category $C(I)$ called \emph{fiber over} $I$ whose objects and morphisms
are families of objects and morphisms of $C$ indexed by elements of $I$, that
is, families of the form $X = \{X_{i}\}$ and $f = \{f_{i}: X_{i} \rightarrow Y_{i}\}$ where $i \in I$.
B\'enabou remarks that speaking about categories naively we assume more
than this, so we cannot just fix some sufficiently large set $I$ and use it
for indexing every time when this is needed. Namely, we also assume the
possibility of re-indexing: given families $X = \{X_{i}\}, Y = \{Y_{j}\}$ in $C$ where
$i \in I, j \in J$ and morphism $u: J \rightarrow I$ in $B$ we assume that family of
objects $X_{u(j)}$ and family of morphisms $f = \{f_{j}: Y_{j} \rightarrow X_{u(j)}\}$ is uniquely
defined and ``behaves properly''. This allows us to extend $C(I)$ through
introducing new category $Fam(C)$ of families of $C$ where objects are
families of objects of $C$ indexed by different sets and morphisms are pairs
of the form $(u, f)$ where $u$ and $f$ are as just described. Morphisms of the
form $f = \{f_{i}: X_{i} \rightarrow Y_{i}\}$ we identify with $(id_{I}, f)$ where $id_{I}$ is identity
morphism of $I$ in $B$. The composition of morphisms in $Fam(C)$ is
defined in the obvious way. We equipe the construction with projection
functor $p_{C}$ which sends every family of objects of $C$ to the set by which
this family is indexed and every morphism $(u, f)$ between families to
morphism $u$ between sets: $p_{C}: \{X_{i}\} \rightarrow I, (u, f) \rightarrow u$.

Now suppose that we know what equality is in $Fam(C)$ and in $B$
but not in $C$. This implies that we cannot think of families (of objects
and morphisms of $C$) extensionally as usual. In particular given a family
$X = \{X_{i}\}$ where $i \in I$ and morphism $u: J \rightarrow I$ in $B$ we cannot
define another family $Y = \{Y_{j}\}$ by saying that $Y_{j} = X_{u(j)}$ because the
latter equality doesn't make sense for us! Nevertheless we can achieve
the same effect through requiring certain properties of $Fam(C)$ and $p_{C}$.
What we need for it is to characterize morphism $\phi_{(u,X)} = (u, \{id_{Y_{j}}\})$ in
$Fam(C)$ without using equality in $C$; $\{id_{Y_{j}}\}$ to be the family of identity
morphisms of objects $Y_{j} = X_{u(j)}$ in $C$.

Given $u: J \rightarrow I$ and $X = \{X_{i}\}$
$\phi$ is characterized up to unique isomorphism by the following property:\\
(i): for any morphism $\psi  = (v, g)$ with codomain $X$ in $Fam(C)$ and
any $v'$ such that $v = v'u$ in $B$ there exists in $Fam(C)$ a unique  $\psi '$
such that $\psi = \psi ' \phi$ and $p_{C}\psi ' = v'$.\\
In addition morphisms of the form $\phi = (u, \{id_{Y_{j}}\})$ satisfy the \emph{functoriality conditions}\\
(ii): $\phi_{(u,X)} = id_{X}$
(where $id_{X}$ is the identity morphism of family $X$), and $\phi_{(uv,X)} = \phi_{(u,X)}\phi_{(v, Y)}$ for each
$v: K \rightarrow J$.

Now we use these properties as definition of abstract functor $p: F \rightarrow B$ called \emph{Grothendieck fibration} over $B$ (or \emph{fibred category} over $B$) in the case when only the property (i) is taken into account, and called \emph{Grothendieck split fibration} over $B$ in the case when in addition for  each pair $(X, u: J\rightarrow p(X))$ one makes a particular choice of $\phi_{(u, X)}$ called \emph{splitting} satisfying the functoriality conditions (ii). Thus equality in a category can be defined as splitting of fibration over an appropriate base.  Noticeably given a fibration its splitting might not exist or be not unique. I refer the reader for further details to\cite{Benabou:1985}. Talking about fibrations in this Section I always mean Grothendieck fibrations. The relationships between this categorical notion of fibration and the geometrical notion of Hurewicz fibration mentioned in \textbf{4.5} above will be discussed in \textbf{6.8}.   

B\'enabou's theory of equality in categories allows for regarding objects and morphisms of a given category as families rather than bold individuals; these families can be occasionally split into elements through a (split) fibration in different ways dependently of the choice of base. Such splitting is the inverse operation with respect to the informal identification of isomorphic objects and morphisms, and unlike the latter it is performed more rigorously and ``more categorically''. This inversion is remarkable : it shows that given the definition of  equality through split fibration families are no longer thought of as extensional multiplicities, that is, as classes. Recall however that given a fibration $p: F\rightarrow B$ categories $F, B$ are construed in the usual way and, in particular, assume the usual equality of morphisms and objects, so the internalized equality relates only to hypothetical category $C$ such that $F = Fam(C)$. As B\'enabou stresses in the end of his paper such ``meta-equality'' is indispensable ``unless you do something different from Category theory''. In the end of this Chapter we shall see that this ``something different'' is a real option.

\section{Higher Categories}
Given abstract category $C$ consider class $Hom(A,B)$ of morphisms $f, g,...$ of the form $A \rightarrow B$.
Then turn $Hom(A,B)$ into a new category formally introducing morphisms of the form $\alpha: f \rightarrow g$, that is, morphisms between morphisms of $C$:
\begin{center}
\UseTwocells 
\[
$$\xymatrix{A \rtwocell^f_g & B}$$
\]
 \end{center}
Do this for all pairs of objects of $C$. Observe that
morphisms $\alpha$ can be composed in two different compatible ways shown at the below diagram

\UseTwocells 
\[
$$\xymatrix{A \rtwocell^f_g & B \rtwocell^h_i & C} $\longrightarrow$ \xymatrix{A \rtwocell^k_l & C}$$
\]
where $k = fh$ and $l = gi$ (horizontal composition)

\UseTwocells 
\[
$$\xymatrix{A \rtwocell^f_g &B} $$
$$\xymatrix{A \rtwocell^g_h &B}$$
\]
$\downarrow$

\UseTwocells 
\[
$$\xymatrix{A \rtwocell^f_h &B} $$
\]

(vertical composition)

Requiring now appropriate equational conditions (which say that certain diagrams commute),  we obtain a \emph{2-category}. It comprises objects $A,B,...$, morphisms $f, g,..$ between objects (the same as in $C$) called in this
context \emph{1-morphisms}, and morphisms between morphisms $\alpha, \beta, .. $ called \emph{2-morphisms}. See \cite{Leinster:2002} for precise definitions. (The reason why these definitions are many will be clear in a moment.)  
An example of 2-category which has been around from the very beginning of category theory (even it was not called by this name) is 2-category $2-Cat$ having (some, for example, all small) categories as objects, functors between these categories as 1-morphisms and natural transformations between the functors as 2-morphisms \cite{MacLane:1998}.

Let me now explain what 2-categories have to do with the internalization of identity (equality).
Remark that in a 2-category we have not only the usual composition of 1-morphisms
but also functor $F: Hom(A,B) \times Hom(B,C) \rightarrow Hom(A,C)$ (provided
that in category $Hom_{C}$ having $Hom$-categories of $C$ as objects Cartesian product $\times$ is available). On 2-morphisms this functor acts as their horizontal composition while in $Hom$-categories 2-morphisms are composed vertically. If functors of the form $F$ preserve identities in $Hom$-categories (2-identities) then equalities in $C$ may be omitted without any lost. This means that we don't even need to define $C$ as a category but may think of it as a class of objects and
morphisms between these objects, and then define composition of these morphisms ``from
above'' through functors like $F$. In this case one may speak indeed about ``replacement of relations
by morphisms'': 2-identites from $Hom$-categories make in $C$ the job of equalities. The situation
here is quite analogous to one we have seen in fibred categories: at the top ``meta-'' level of
construction (namely in $Hom$-categories and in the category $Hom_{C}$ of the $Hom$-categories) one uses the ``god-given'' equality but at the bottom level equalities are got rid of.

An apparent difference between the two approaches is this: in higher categories the notion of
class is used at all levels including the lowest one while in fibred categories this notion is used
only at the meta-level while at the lower level classes are replaced by non-extensional
families. But is the assumption that objects of 2-category form a class necessary?
Prima facie it is the case. For in order to compose 2-morphisms $\alpha, \beta$ in a $Hom$-category (that is, compose them vertically) we need to check that domain of $\alpha$ equals codomain of $\beta$. So 1-morphisms (objects of $Hom$-categories) should form a class of well-distinguishable elements provided with a notion of equality allowing for distinguishing them). If we take this view then the internalization of equality in $C$ just described will be only partial: it will apply to equalities of the form $fg = h$ but not to equalities of the form $f = f$. However it is easy to get around this point through identification of 1-morphisms with 2-identities, so all needed equational conditions could be written in terms of 2-morphisms. Then one may claim that in $C$ ``there is no equality'', and hence its elements don't need to form a class.

It should be noted that the usual interest of people working in higher category theory is not
internalization of equality as such but \emph{weakening} of equality, that is, finding a rigorous way of
replacing equalities with certain isomorphisms (see \textbf{8.2} below). This approach is quite natural in the given context since the requirement that functors of the form $F$ preserve all 2-identities is very strong. Since we are no longer obliged to think of $C$ as a category in the usual sense we get a room for playing. Instead of imposing on $Hom$-categories equational conditions implied by the assumption that $C$ is a category we can use weakened conditions which don't imply that 2-isomorphisms replacing equalities in $C$ are identities. Such weak 2-categories have been first introduced by B\'enabou in \cite{Benabou:1967} under the name of \emph{bicategories}. In bicategories 

\begin{itemize}
\item (i) the usual associativity law $(fg)h=f(gh)$  is replaced by the requirement of existence of
associativity (2-)isomorphism $a: (fg)h\rightarrow f(gh)$ (eventually called \emph{associator} by other authors);
\item (ii) the usual axioms of identity $1_{A}f = f$ and $f1_{B} = f$ (for $f: A\rightarrow B$ is replaced by the requirement of existence of unit 2-isomorphisms $l: 1_{A}f \rightarrow f$ and $r: f1_{B}\rightarrow f$. 
\end{itemize}

These isomorphisms are subjects of equational conditions called \emph{coherence laws}, which I shall not list here, see \cite{Leinster:2004}) where this idea is developed in several different variants. 

The notion of 2-category allows for a straightforward geometrical
analogy : think of objects as points, of 1-morphisms as oriented lines,
and of 2-morphisms as oriented surfaces bounded by the lines. Whether
we use this analogy (which is quite profound as we shall briefly see) or not the notion of 2-category calls for the inductive generalization to the notion of $n$-category for arbitrary $n$ and further
to $\omega$-category (leaving bigger ordinals apart). The strict (meaning
non-weakened) versions of the notions of $n$- and $\omega$-category look unproblematic
: the enrichment of a given category $C$ bringing the notion of
2-morphism explained in the beginning of this Section can be easily reformulated
as an inductive step bringing the notion of $k$-morphism provided
with $k$ different kinds of composition. However the
notions of weak $n$- and $\omega$-categories are much less transparent. There are many alternative definitions of weak $n$- and $\omega$-categories around ; ten
of them are presented in  \cite{Leinster:2002}. A specific obstacle for putting these
things into order, which is not irrelevant to the issue discussed in this
paper, is this : it is not immediately clear which notion of equivalence one should
apply to answer the question whether two given definitions of weak $n$-
category are equivalent or not.  In the rest of this Chapter I present a recent geometrical approach to weak higher categories, which sheds more light on the identity issue.    

\section{Homotopies}
The geometrical analogy mentioned in the end of the previous Section can be developed into a genuine example of $n$-category and then used for the introduction of this later concept \cite{Lurie:2009}. As we shall now see in such a context the idea of weak identity emerges naturally in the sense that one gets a weak category first and then needs some further efforts in order to construct a strict category out of it. This suggests a reverse of the usual conceptual order, which turns the so-called ``weak identity'' into  the ``true'' identity and leaves the whole strategy of ``weakening identity'' without its usual appeal. We shall also see that this geometrical setting is closely related (in a sense that will be made precise) to a logical calculus, which allows for discussing the non-standard notion of identity in precise formal terms. Finally we shall see that the categorical notion of fibration discussed in \textbf{6.5} also plays a crucial role in this setting.     

Consider a topological space $T$, points $A, B, ...$ of this space and \emph{paths} $f, g, ....$ between these point. A path is not to be confused with a curve. While a curve in a space is a subspace of this space a path in $T$ is a \emph{map} from a distinguished \emph{unit space} $I$ into $T$. Often $I$ is identified with the real interval [0,1], which can be thought of as the time interval corresponding to motion beginning at some point $A$ and ending at point $B$ and thus producing path $f$ from $A$ to $B$. Given two points $A, B$ there is, generally, many different paths from $A$ to $B$. Formally these different paths are represented by different continuous maps of the form   $f: [0,1] \rightarrow T$ such that $f(0) = A$ and $f(1) = B$ (see fig ???). Beware that different curves may represent the same curve: this happens every time when the corresponding functions $f$ have the same image. Such paths can be thought as being traced by a moving particle which moves between the same endpoints $A, B$ along the ``same route'' during the same time interval [0,1] and differ only in the character of their motion: for example particle $x$ may move with a constant speed and particle $y$ move faster than $x$ during a half of the time and slower than $x$ during the other half of the time. So what in the homotopy theory is called a ``path'' is indeed a motion, which is equipped with some notion of ``absolute time'', which allows one to localize a given moving particle at a given moment of time $t \in [0,1]$ at certain point $P$ of the given space $T$.   

\begin{center}
\includegraphics[scale=0.3]{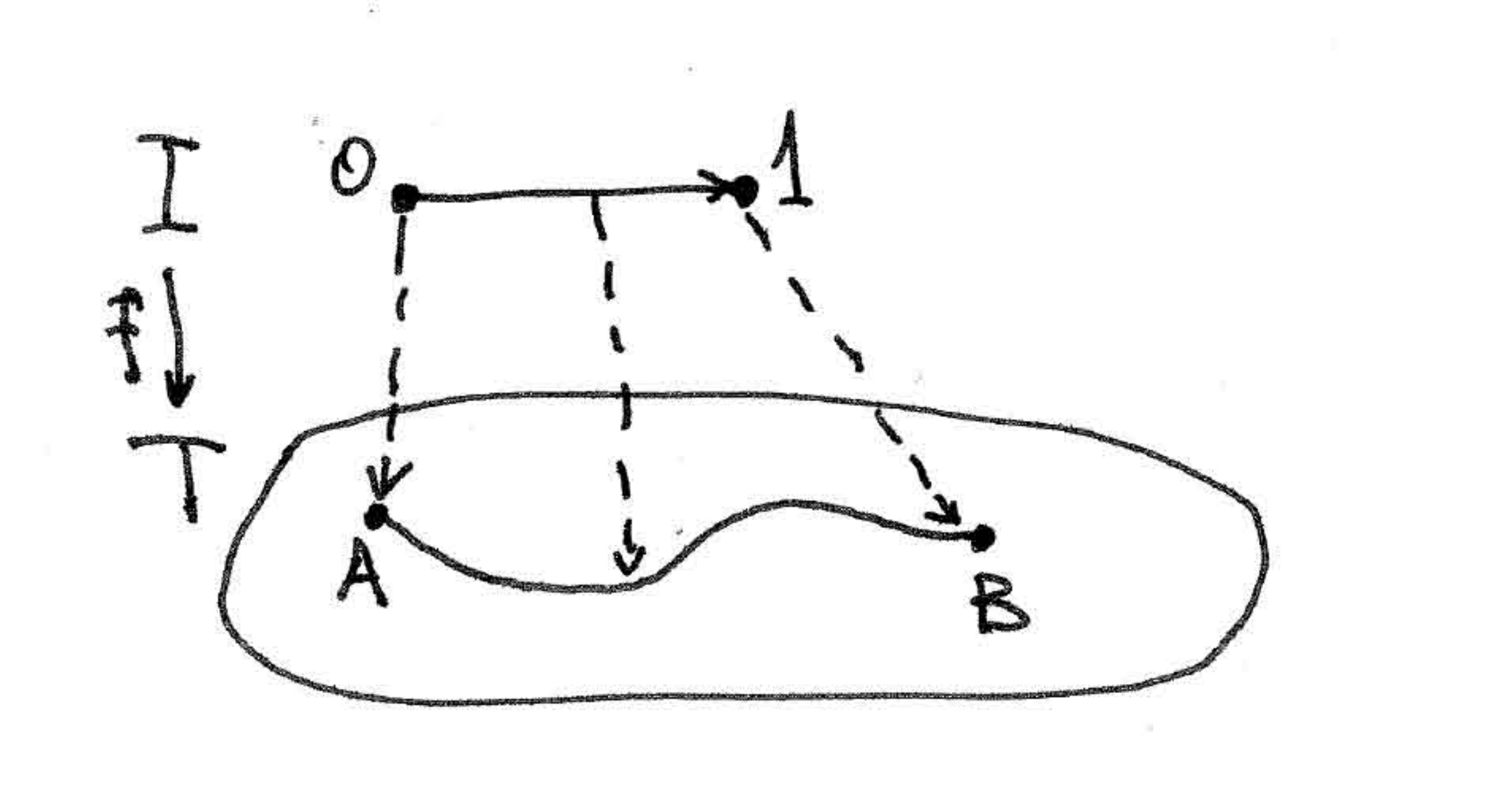}
\end{center}

\begin{center}
Fig. 6.2
\end{center}

Now consider a tentative category $P$ with objects $A, B, ...$ points of $T$ and morphisms paths $f, g, ...$ between these points. We assume that for every path $f: [0,1] \rightarrow T$, $f(0) = A$ and $f(1) = B$, there is the inverse path $g: [0,1] \rightarrow T$, $g(0) = B$ and $g(1) = A$, so $P$ is a \emph{groupoid}. It remains to define the composition of morphisms and check the axioms. Intuitively the composition of paths appears unproblematic: given path $f$ from $A$ to $B$ and another path $g$ from $B$ to $C$ one may easily imagine a composite path $h$ from $A$ to $C$, where $B$ is a midpoint. However an attempt of a more rigorous mathematical description reveals a hidden complexity of the situation. By composing $f: [0,1] \rightarrow T$ and $g: [0,1] \rightarrow T$ straightforwardly we get a map $h^{*}$ of the form $[0,2] \rightarrow T$, which is not a path in the sense of our definition. So in order to get $h = fg$ of the same form $h: [0,1] \rightarrow T$ we should do something different. Here the intuitive interpretation of paths as motions becomes helpful. Remind that [0,1] is the time taken by a continuous motion from the beginning of a given path to its end. Now we have two such motions - one from $A$ to $B$ and the other from $B$ to $C$ - each of which takes the same time [0,1]. We want to produce now a motion from $A$ to $C$, which takes the same time. Clearly, in order to do this one needs to move faster!  

\begin{center}
\includegraphics[scale=0.3]{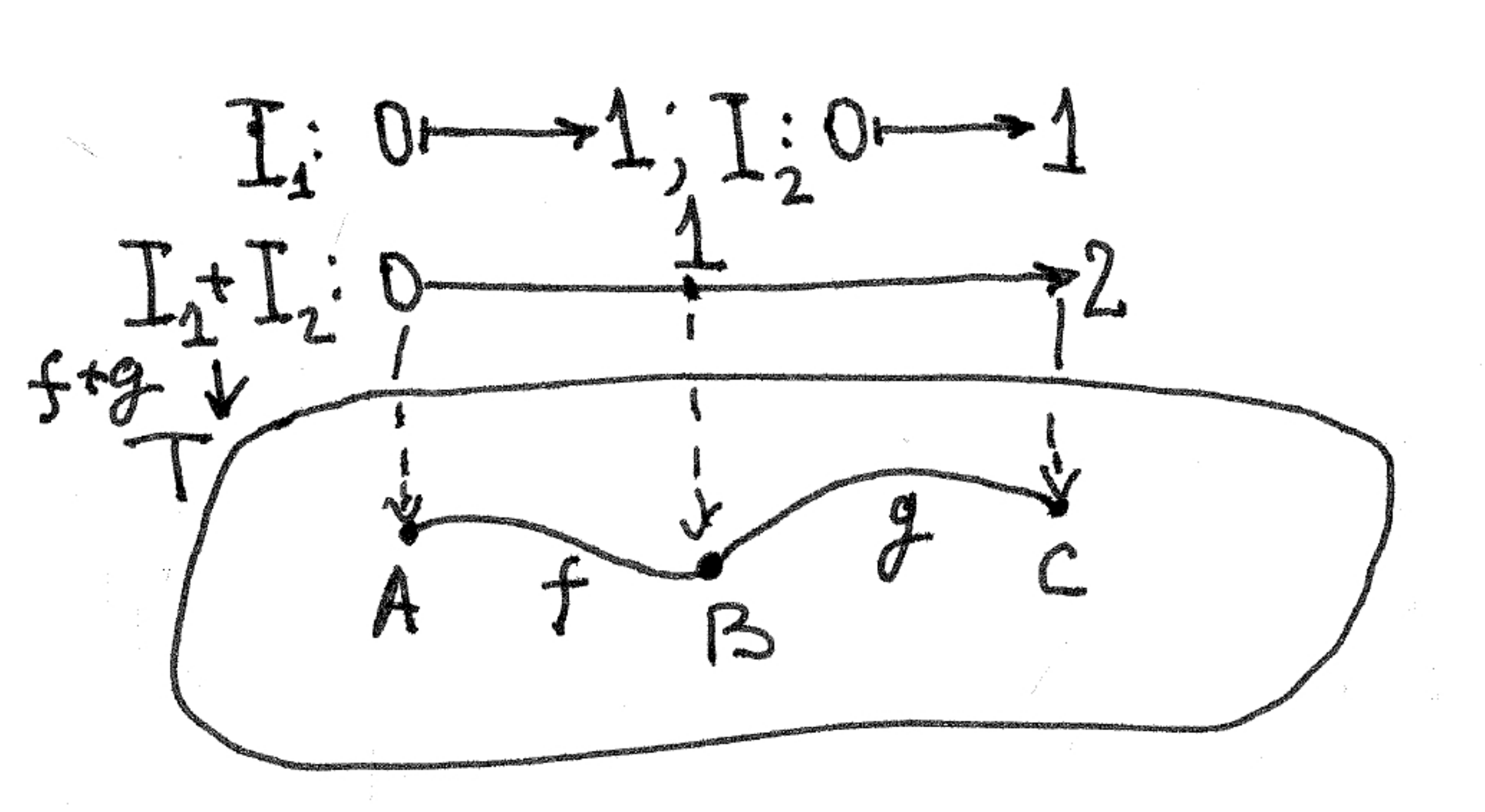}
\end{center}

\begin{center}
Fig. 6.3
\end{center}

Mathematically this intuitive idea can be realized as follows. First, we introduce a real variable $t \in [0,1]$ called in this context a \emph{parameter}, so that for every value of $t$ $f(t)$ represents a precise position (i.e. a point) on path $f$. Path $g$ is treated similarly. One doesn't need to introduce a different parameter for $g$ because motions along $f$ and $g$ can be seen as simultaneous, so $f(t)$ and $g(t)$ represent two different points. However in order to compose $f$ and $g$ one needs to put these two motions in the right order, and as I have just noticed, the composite motion must have a different speed. For this purpose we introduce a new parameter $t' = \frac{t}{2}$ (so the time is halved and hence the speed gets doubled) and using this new parameter define  the composition $h = fg$ as follows:
$$
h(t') = \left\{\begin{array} {l}f(2t'),  0 \leq t' \leq \frac{1}{2} \\ g(2t' - 1), \frac{1}{2} < t' \leq 1\end{array}\right.
$$

so $t'$ ranges from 0 to 1 as required. This trick is called a \emph{reparameterization}.

The problem is that there are many different reparameterizations, which can be used for this purpose, and which produce different results. So no particular reparameterization gives us a \emph{general} definition of path composition. Obviously the unit interval in   $h: [0,1] \rightarrow T$ can be cut not only into the two equal halves but also in any other proportion (so the motion along $f$ may become slower and the motion along $g$ become faster). Moreover the reparameterization needs not to be linear: the speed of motion along either of the two composed paths $f, g$ needs not to be constant. Generally, the reparameterization amounts to choosing a particular continuous ``scaling map''  $s: [0,1] \rightarrow [0,2]$ from the  space $[0,2]^{[0,1]}$ of such maps. Each particular choice of $s$ brings about a new definition of path composition. There is no way to make a ``right'' (or as mathematicians say \emph{canonical}) choice here. In the general situation the composition of paths can be pictured as follows (Fig. 6.4): 

\begin{center}
\includegraphics[scale=0.3]{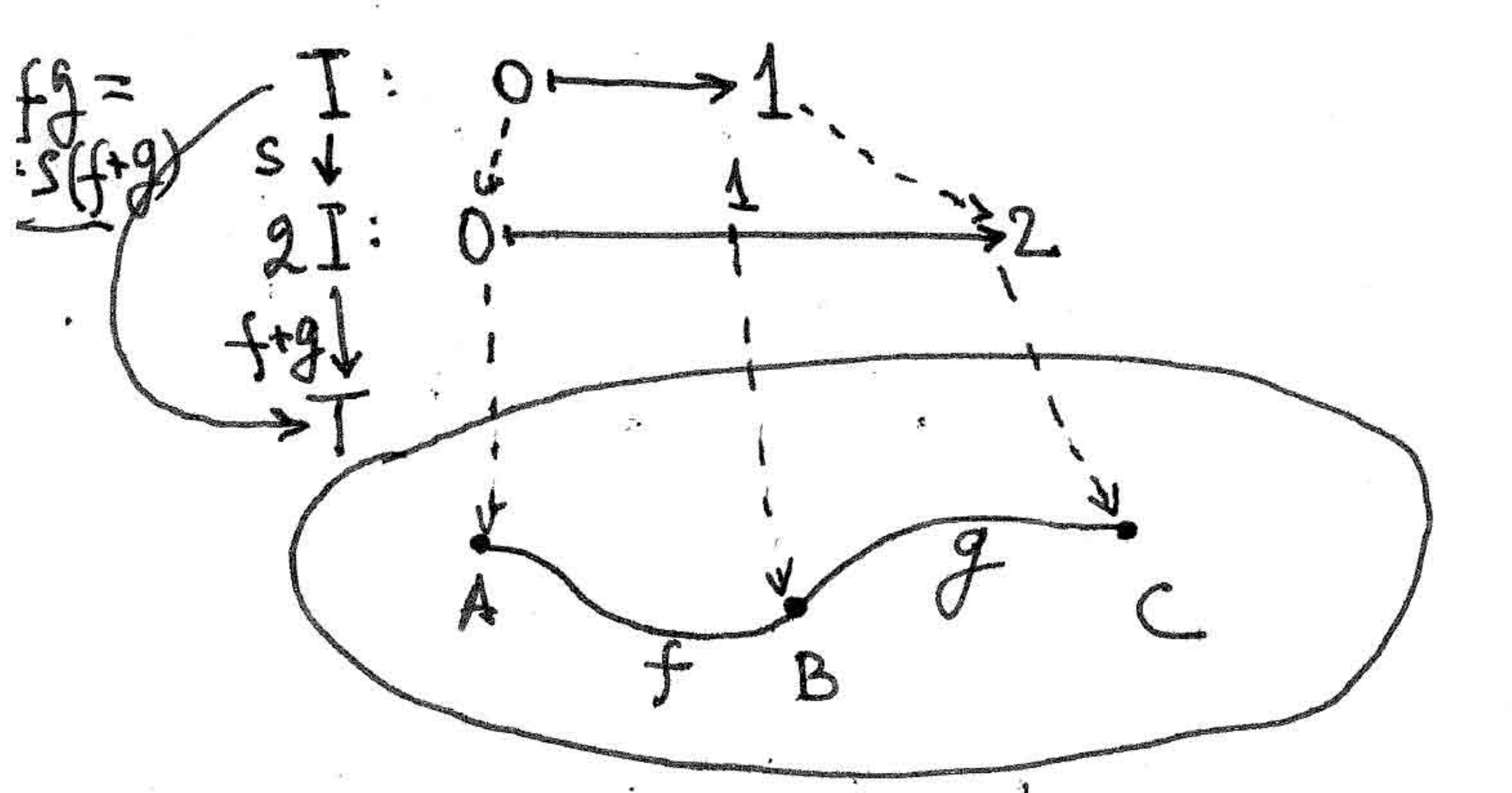}
\end{center}

\begin{center}
Fig. 6.4
\end{center}

Suppose now that some particular map $s: [0,1] \rightarrow [0,2]$ is chosen and the composition of paths $h = fg$ is defined as above. Let $s(t) = 2t$ for simplicity
\footnote{This is, of course, the same linear reparameterization ``by halves'' that we used earlier. However then we thought about it as a function $[0,2] \rightarrow [0,1]$. The above diagram shows that the right sense of the scale function $s$ is $[0,1] \rightarrow [0,2]$ because it must be composable with $(f+g): [0,2] \rightarrow T$. Generally, $s$ needs not to be invertible, so there can be no function transforming the old parameter into the new one.   
}. It remains to check whether points and paths composed according this rule indeed form a category. Identity morphism of a given point $A$ is defined as map sending each point of [0,1] to $A$ and it is straightforward to check that it behaves as expected. (Intuitively such identity paths can be thought of as ``staying in rest'' at $A$ during the whole time interval [0,1].) However checking the associativity property reveals a further problem: the composition of paths turns to be \underline{not} associative. This can be seen immediately at the below diagram: $(fg)h$ sends the middle point $\frac{1}{2}$ to $C$ while $f(gh)$ sends $\frac{1}{2}$ to $B$, which means that the two composite paths from $A$ to $D$ are different in spite of the fact that both are represented by the same curve (see Fig. 6.5). So $P$ defined as above is not a category!

\begin{center}
\includegraphics[scale=0.3]{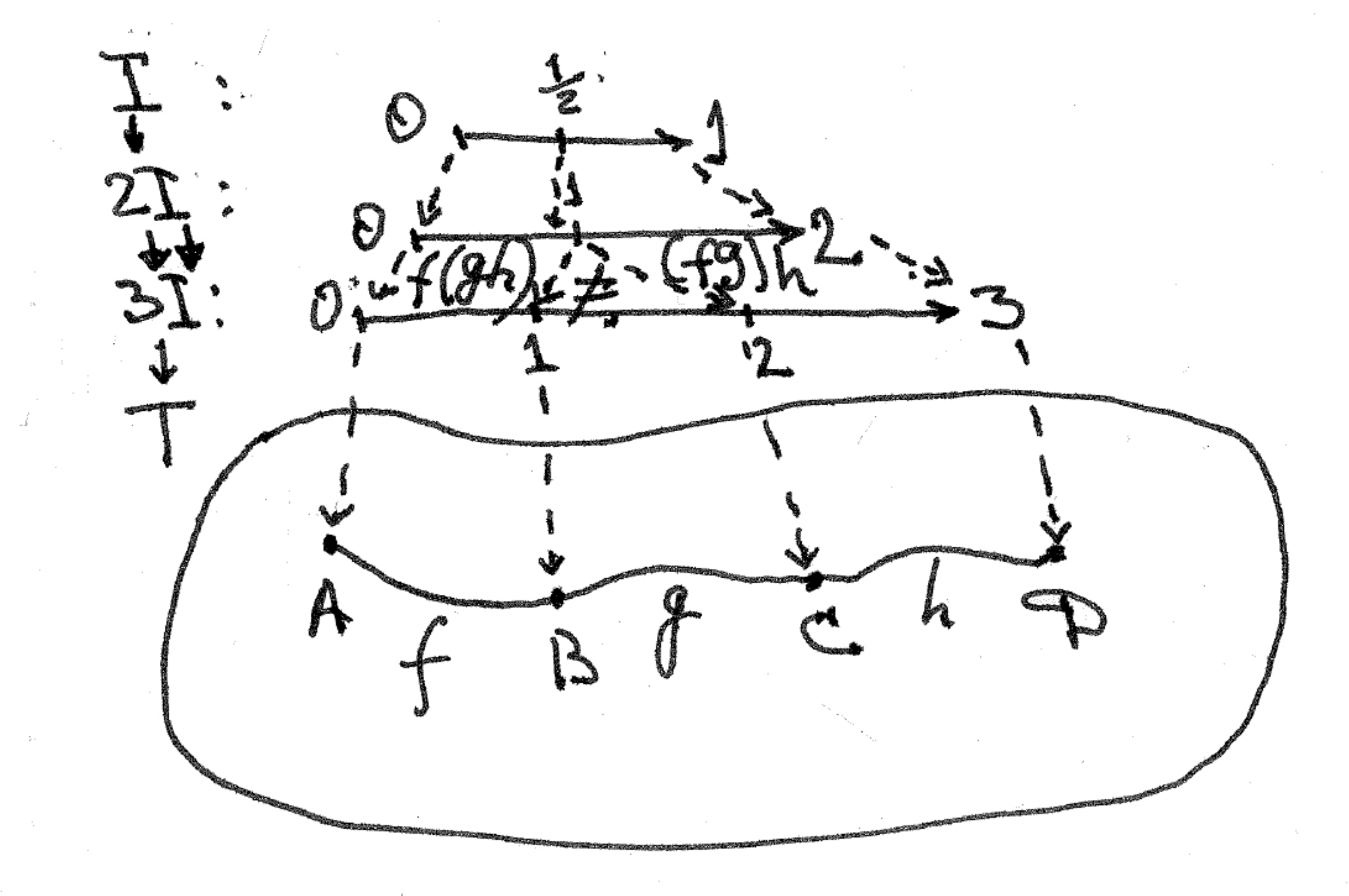}
\end{center}

\begin{center}
Fig. 6.5
\end{center}

Here we get an important conceptual choice. We may either fix the associativity with an additional effort or explore the obtained non-associative structure (which is not a category in the usual sense of the term). I shall now briefly consider both options; we shall see that the second option provides a more general view onto this situation and makes explicit some further structures, which otherwise ``remain hidden''.  

\underline{(a) Fixing the Associativity}  

In order to fix the associativity of path composition it is sufficient to redefine our category $P$ by taking its morphisms to be certain equivalence classes of paths rather than paths themselves. The appropriate equivalence relation is that of  \emph{homotopy}; two paths a called \emph{homotopic} when there is a \emph{homotopy} between them. Once again we encounter here a systematic terminological confusion when one and the same term is commonly used for denoting a relation and a map (see \textbf{6.1}). The notion of homotopy as a map is as follows: it is a continuous map  $u: [0,1]^{2} \rightarrow T$  (where $[0,1]^{2}$ is the real square) such that $u(t, 0) = f(t)$, $u(t, 1) = g(t)$ and  $u(0, r)= f(0) = g(0)$, $u(1, r) = f(1) = g(1)$ for all values of parameters $t, r$ from [0,1]. As we can see a homotopy between paths can be thought of as a ``path between paths'' or as a ``path of the second order''. Correspondingly, a path can be described as a ``zero-order homotopy''.  A homotopy can be represented with a two-dimensional \emph{surface} (cell) delimited by a pair of curves representing two homotopic paths in a way similar to which paths are represented by curves\footnote{
In the introductory literature the path homotopy is often also described as a continuous transformation of one path into another. This is somewhat misleading because such a description doesn't take into account the role of the interval [0,1] and suggests thinking about paths as if they were continuous curves. The difference between a path homotopy and a continuous transformation is similar to that between a path and a curve: a homotopy is not just a transformation but a parameterized transformation or ``transformation through time'', where the time has a particular shape of the unit interval (but can also be of some other specified shape). The representation of paths by curves is imperfect because different paths can be represented by the same curve. A similar remark concerns the representation of homotopies by 2-dimensional cells. When different paths are represented by the same curve the homotopies of these paths are not representable by 2-cells in the usual way.} (Fig. 6.6):

\begin{center}
\includegraphics[scale=0.3]{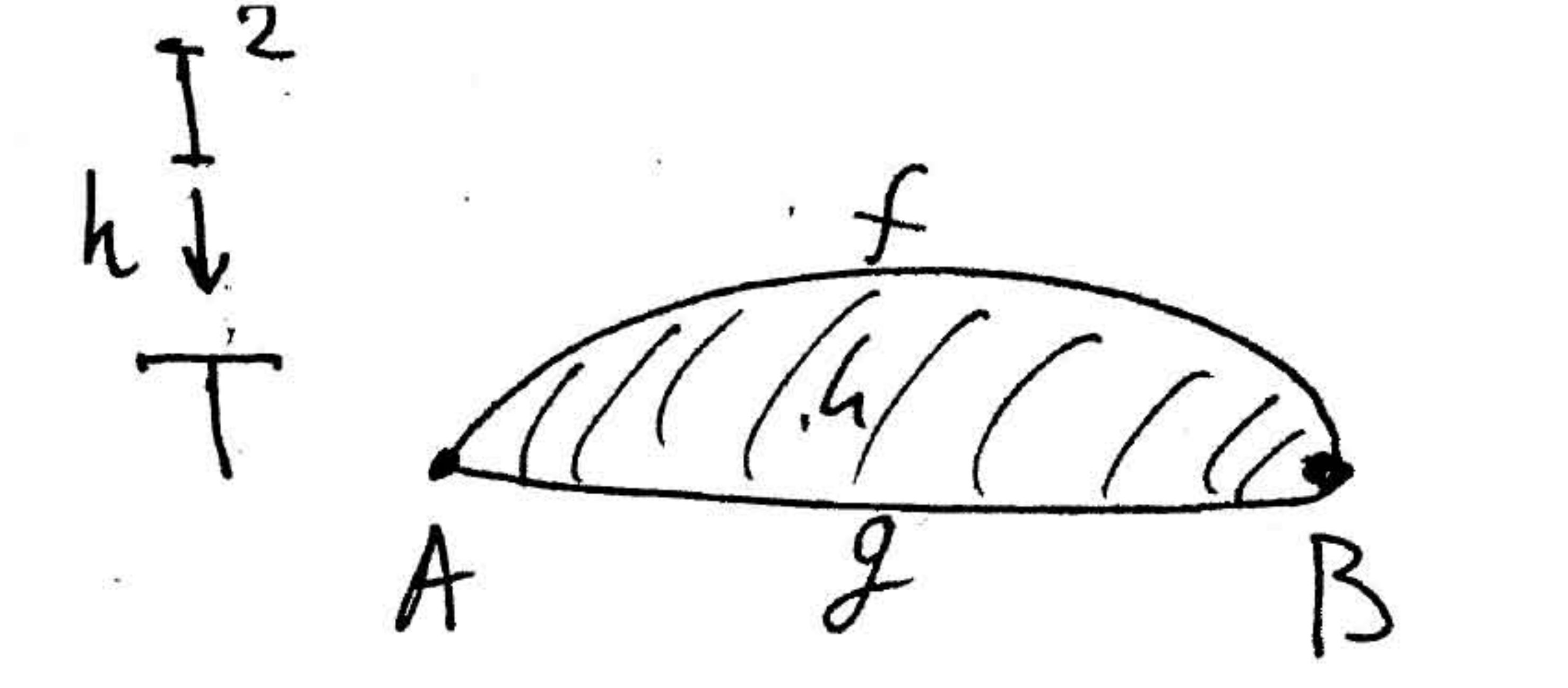}
\end{center}

\begin{center}
Fig. 6.6
\end{center}

Now observe that given two composable paths $f, g$ different choices of scale function $s$ always correspond to homotopic (albite, generally, not equal) compositions $h=fg$. The homotopy works in this case as a ``reparameterization of reparameterization'', which allows to transform one given reparameterization into another continuously along a new parameter $r$ ranging over the interval [0,1]. Since the morphisms of our category are no longer paths but homotopy classes of paths this observation allows us to define their composition \emph{up to homotopy}: just as we no longer distinguish between homotopic paths we shall not distinguish between homotopic compositions of paths. So the problem is partly fixed: the composition of morphisms is now defined in our category uniquely. Since for any choice of scale function $(fg)h$ is homotopic to $f(gh)$ the associativity is fixed too: instead of saying that  compositions $(fg)h$ and $f(gh)$ are homotopic we shall say that they are equal. Thus we get a well-defined category $P_{1}$ with objects points of $T$ and morphisms equivalence classes of homotopic paths. This category is called the \emph{fundamental groupoid}  of given topological space $T$. By collapsing groupoid $P_{1}$ into a group (through identification of all its objects) one gets the more familiar notion of \emph{fundamental group} of a topological space\footnote{There is a very interesting attempt to generalize the usual homotopy theory by allowing for non-invertible  paths; in this case the fundamental category is not necessary a groupoid, see \cite{Grandis:2001}.}. 

\underline{(a) Exploring the Non-Associativity} \\
 The idea is now the following: instead of forcing the uniqueness and the associativity of path composition through the identification of homotopic paths consider path homotopies as 2-morphisms of a higher category. Since the composition of paths is not unique and not associative the obtained 2-groupoid $P_{2}$ is  \emph{weak} in the sense of \textbf{6.6}. However if we want to complete the construction at this point we still need to fix the composition of 2-morphisms (i.e., the composition of path homotopies) precisely. This can be done by the trick already used for  $P_{1}$: the path homotopies proper are replaced by equivalence classes of such homotopies; more precisely, path homotopies are identified up to homotopy of the next order. Homotopies of the next order are ``homotopies between path homotopies'', i.e., maps of the form $[0,1]^{3} \rightarrow T$. 
 
Alternatively, we may take yet higher homotopies into consideration and proceed to build the ladder up until any $n$ and further ad infinitum. (What exactly happens when $n$ tends to infinity is a tricky question, which I informally discuss in \textbf{6.10}). The infinite-dimensional (weak) groupoid ($P_{\omega}$ so obtained is called  the \emph{fundamental} $\omega$-groupoid of $T$. This groupoid was first described by Grothendieck in a private letter \cite{Grothendieck:1983} dated by 19.02.1983 and then used by other authors as a basic motivating example of the very notion of (weak) $n$-category \cite{Leinster:2004}, \cite{Lurie:2009}. The idea according to which $P_{\omega}$ captures all the basic structure of the corresponding topological space $T$  is called sometimes the \emph{Grothendieck conjecture}.  

Considering the choice between the two above options we can again notice a controversy between using relations (equivalences), on the one hand, and using transformations, on the other hand, see \textbf{6.1} above. However since we are working now in a general categorical context (rather than the limited group-theoretic context) the ``substantialist'' interpretation of transformations developed in \textbf{6.1} becomes, generally, irrelevant. In what follows we shall see how this extended context allows for a more sophisticated identity concept than the traditional notion of substance preserving its identity through time and change.   

Similar constructions can be made with the category $Top$ of topological spaces (with morphisms continuous transformations) rather than a single topological space $T$ as above.  Two continuous  maps $\xymatrix{A\ar@<2pt>[r]^{f} \ar@<-2pt>[r]_{g}&B}$  between two spaces $A, B$ from $Top$ are called homotopical when there exist a homotopy between them. A  homotophy $h$ between two continuous maps $f, g: A \rightarrow B$ in $Top$ is a continuous map $h: A \times [0,1] \rightarrow B$ such that $h(0) = f$ and $h(1) = g$ (this later notion of homotopy generalizes upon the notion of path homotopy explained above). By identifying all homotopic maps in $Top$ one gets the \emph{homotopy category} $hTop$. However one may also consider homotopies of continuous maps in $Top$ as 2-morphisms. So one gets a two-dimensional version $Top_{2}$ of $Top$. Taking into consideration higher homotopies (i.e., continuous maps of the form  $h: A \times [0,1]^{k} \rightarrow B$ and proceeding inductively one gets an infinite-dimensional topological category $Top_{\omega}$.

In this latter setting there emerges another important homotopical concept, namely that of \emph{homotopy equivalence}, which is more general than the simple notion of being homotopic. Consider again a pair of topological spaces $A, B$ with two continuous maps, which this time go into the opposite directions:  $\xymatrix{A\ar@<2pt>[r]^{f} &B\ar@<2pt>[l]^{g}}$, and consider compositions $fg$ and $gf$. Remind that  if $fg = 1_{A}$ and $gf = 1_{B}$ then $f$ and $g$ are called mutually inverse. Now consider a weaker property, namely the situation when $fg$ is homotopic to $1_{A}$ and $gf$ is homotopic to $1_{B}$ (which is tantamount to saying that there exist homotopies $fg \times I \rightarrow 1_{A}$ and $gf \times I \rightarrow 1_{B}$). When this latter condition holds spaces $A$ and $B$ are called \emph{homotopy equivalent}, or interchangeably, belonging to the same \emph{homotopy type}. This definition straightforwardly generalizes (at the price of losing its concrete geometrical meaning) to general categories by replacing continuous maps by functors and homotopies by \emph{natural equivalences} of these functors, i.e., by isomorphisms of these functors, which are invertible 2-arrows. The obtained equivalence relation between general categories is called \emph{categorical equivalence} or \emph{equivalence of categories}; it has been already mentioned at several occasions in this book. 

\section{Model Categories} 
So far we considered the basic homotopy theory as a geometrical motivation of and/or a geometrical underpinning for the higher category theory. Now I want to outline another important idea that has emerged at the joint of category theory and homotopy theory. So far talking about homotopies we took the notion of topological space for granted and used it as a foundation of homotopy theory. We can however ask: What is a general category, which allows for doing homotopy theory in the same (or similar) way, in which we do it in $Top$? 

What we are now looking for is a categorical axiomatization of the homotopy theory, i.e., a list of axioms that makes an abstract category to behave like $Top$ with respect to its homotopic properties (although some reasonable generalization of these properties is also expected). Such axioms were suggested by Quillen \cite{Quillen:1967} in 1967; the resulting concept is called a \emph{model category} (the term is short for ``model of homotopy theory'', see  \cite{Quillen:1967}, p. 3). A model category is an abstract category with three distinguished classes of morphisms, each closed under composition: \emph{weak equivalences}, \emph{fibrations} and \emph{cofibrations}.  The weak equivalences play the role of generalized homotopy equivalences, fibrations are generalized fibrations and the cofibrations are maps, which are in an appropriate sense dual to fibrations (while fibrations are ``homotopic surjections'' cofibrations are ``homotopic injections''); morphisms of these three distinguished classes are the subject of several axioms. The notion of model category proves sufficient for developing an abstract homotopy theory, which comprises not only classical but also many non-classical examples and thus turns to be a fruitful generalization of the classical homotopy theory. While in the last Section we have shown how the (classical) homotopy theory helps to develop the higher category theory in concrete geometrical terms, now we are discussing how the (ordinary) category theory helps to develop the homotopy theory abstractly (which in the given context means \emph{axiomatically}). So we have here another case of fruitful dialectical interplay between the ``logical'' (axiomatic) and the ``geometrical'' points of view, which can be compared with the interplay between logic and geometry in Lawvere's topos theory (\textbf{4.9}). Quillen's book of 1967 \cite{Quillen:1967} presents in an axiomatic form geometrical results obtained in and around Grothendieck's school in early 1960ies.  The notion of fundamental $\omega$-groupoid is first suggested by Grothendieck in his letter to Quillen \cite{Grothendieck:1983} written in 1983. (So the order of my exposition in this case is not historical.) In the following two Sections we shall see how this fruitful dialectical interplay between logic and geometry  continues nowadays in the \emph{homotopy type theory} and how a new theory of identity emerges in it.

For a modern introduction into the theory of modal categories see \emph{Hovey:1999}. Here I consider only Quillen's key axiom, which expresses in an abstract form the so-called \emph{homotopy lifting property}. This property is used for defining (Hurewicz) fibrations in the classical homotopy theory (see \textbf{4.5} above) as follows. Consider in $Top$ a commutative square:

$$\xymatrix{X\ar[r]^{f}\ar[d]_{1_{X}\times 0} &E \ar[d]^{p} \\ X \times [0,1] \ar[r]_{h}\ar@{.>}[ur]|-{h*} & B}$$
       
where $h$ is a homotopy. Map $p$ is called \emph{fibration} when for any space $X$ there exist homotopy $h*$ (the dotted arrow) making both triangles commute. The idea is that fibrations ``lift'' homotopies from the base level ``upstairs''; recall from textbf{4.5} the example of M\"obius strip and observe that a given fibration may also happen to be a homotopy  equivalence.  

In Quillen's axiomatic setting the above property is reformulated and stated as an axiom in the following abstract form: given a commutative square as below

$$\xymatrix{X\ar[r]\ar[d]_{i} &A \ar[d]^{p} \\ Y \ar[r]\ar@{.>}[ur] & B}$$
 
where $i$ is a cofibration, $p$ is a fibration and at least one of these is also a weak equivalence, there exists the lifting homotopy shown by the dotted arrow. 

Beware that although the notion of \emph{Grothendieck fibration} (see \textbf{6.5} above) equally involves the intuitive idea of lifting and equally derives from the canonical example of fibered bundle, its relation to the notion of Hurewicz fibration used in the classical homotopy theory and in Quillen's axiomatic theory of model categories is not straightforward. No known model structure (i.e., the structure of model category) on $CAT$ (i.e., on ``the'' category of categories) turns Hurewicz fibrations into Grothendieck fibrations.

Nevertheless this works in the special case of the \emph{canonical} model structure on the category $Grpd$ of groupoids \cite{Anderson:1978}, which is the following: 

\begin{itemize}
\item
\underline{Weak equivalences} are equivalences of categories (see the definition in the last Section);
\item
\underline{Fibrations} are \emph{isofibrations} of groupoids, i.e., functors $p: E \rightarrow B$ such that for each object $e$ from $E$ and any isomorphism $i: p(e) \leftrightarrow b$ from $B$ there exists an isomorphism $j: e \leftrightarrow e'$ such that $p(j) = i$.  Any Grothendieck fibration is obviously an isofibration; the converse does not hold generally but it does hold in the case of groupoids where all morphisms are isomorphisms;
\item
\underline{Cofibrations} are \emph{faithful} functors between groupoids, i.e., functors injective on objects. 
 \end{itemize}   

Thus in this canonical model category fibrations (i.e., Hurewicz fibrations construed axiomatically) are also Grothendieck fibrations.      
 
\section{Homotopy Type Theory}
We have considered independently three non-standard identity concepts:  one (i) in Martin-L\"of's intuitionistic type theory (\textbf{5.11}), another (ii) in B\'enabou's theory of fibred categories (\textbf{6.5}), and finally  one (iii) in the higher category theory and homotopy theory (\textbf{6.6 - 6.7}). All of these approaches can be seen as different developments of Lawvere's 1970 insights about the internalization of the identity relation in the categorical logic and its geometrical realization (\cite{Lawvere:1970b}, see \textbf{4.5}) although the question about the precise impact of Lawvere's ideas obviously requires an accurate historical study in every particular case. In this Section I show how (i) and (iii) nicely combine together in a recently emerged theory called \emph{homotopy type} theory\footnote{This new theory emerged already after the author published an earlier version of this Chapter as a part of separate paper \cite{Rodin:2007}, so this present Section is an update.}. (ii) also turns to be relevant to the new theory both technically (as discussed in \cite{Warren:2008}, Section 2.4) and conceptually (see below); however it is, in my view, yet premature to discuss this latter connection in general philosophical terms.  
  
Remind from \textbf{4.4} that Martin-L\"of's Intuitionistic Type theory has a category-theoretic model called Locally Cartesian Closed category (LCCC). Any such model verifies the type-theoretic Extensionality Axiom (aka \emph{Reflection Principle}) according to which any propositional equality $Id_{A}(x, y)$ implies the definitional equality $x = y$. Since the \emph{intensional} version of the theory discussed in \textbf{5.11} is described in a merely negative way (as one where the Extensionality Axiom is not assured) a LCCC, formally speaking, is also a model of the intensional version of Martin-L\"of's theory. One may wonder, however, whether there exist intensional models of this theory, i.e., models, which do not verify the Extensionality Axiom.  

Using some earlier results, which I shall not treat here (see \cite{Streicher:1993} for further references), Hofmann and Streicher \cite{Hofmann&Streicher:1994}, \cite{Hofmann&Streicher:1998} suggested in 1994 an intensional  \emph{groupoid model} of Martin-L\"of's theory, where a given basic type $A$ (more precisely: a \emph{judgements} of the form $\vdash A: type$) is groupoid $A$, term $x$ of type $A$ ( judgement $\vdash x : A$) is object $x$ of groupoid $A$, and dependent type $B(x)$ (judgement $x : A \vdash B(x): type$) is fibration of the form  $B \rightarrow A$. (I abuse here the notation in the common way by using the same names for type-theoretic syntactic objects and for their corresponding semantic values.) Identity type $Id_{A}(x,y)$ in this model is the \emph{arrow groupoid} of groupoid $A$, which is a functor category of the form $[I, A]$ where $I$ is the connected groupoid having exactly two non-identical objects and a single non-identity isomorphism between these objects. Clearly this model allows groupoid  $Id_{A}(x,y)$ to be non-empty  when terms (objects) $x$ and $y$ are not definitionally equal. This shows that the Extensionality Axiom is not verified by the given model. 

For a more precise historical account I would like to reveal the fact that the above summary of Hofmann and Streicher's paper \cite{Hofmann&Streicher:1998} is somewhat anachronistic. Although the authors notice that their ``use of the term groupoid is in accordance with homotopy theory'' (\cite{Hofmann&Streicher:1998}, p. 91) they don't consider groupoids as topological fundamental groupoids but treat them as abstract categories with all morphisms invertible. Instead of talking about fibrations (as I did in the above summary without specifying the precise sense of the term) they use the notion of \emph{family} of groupoids indexed by groupoid developed in \cite{Dybjer:1996}. The idea of using families of objects and morphisms of one category as objects and morphisms of another category is similar to B\'nabou's (\textbf{6.5}) but in this case it involves further technicalities related to the purpose of modeling the type-dependence. The crucial idea of Hofman and Streicher was to replace families of sets indexed by sets by families of groupoids indexed by groupoids. This also makes an echo of B\'nabou's ideas about using some non-standard ``sets'' for indexing (even if B\'nabou means here some other toposes rather than $Grpd$). It is also interesting to observe that since in $Grpd$ Hurewicz fibrations are Grothendieck fibrations (\textbf{6.8}) - and since we also know how to describe Hofman\&Streicher's model in terms of Hurewicz fibrations \cite{Awodey&Warren:2009} - thinking about indexed families and fibrations in the sense of B\'nabou's \cite{Benabou:1985} turns to be wholly relevant in this context.  

In 2006 Voevodsky lectured and circulated an unpublished note   \emph{On the homotopy $\lambda$-calculus} \cite{Voevodsky:2006} where he described a large research program of studying type systems by homotopical methods. Independently Awodey\&Warren suggested in 2007 (the date of submission of \cite{Awodey&Warren:2009}) an intensional model of Martin-L\"of's type theory similar to Hofmann\&Streicher's model but more abstract and more explicitly homotopical in the sense that it used the setting of the general model category. Thus Awodey\&Warren revealed the homotopical aspect of Hofmann\&Streicher's groupoid model, and thus contributed by their research to Voevodsky's program.  By the same pattern Awodey\&Warren showed  that every model category admits an internal language, which is a form of Martin-L\"of's type theory. This latter result is of great interest for us on its own rights, because it shows that toposes are not exceptional in their role of joint between geometry and logic. Figuratively speaking, it shows that the geometrical operation of fibration is just  as significant \emph{from a logical point of view} as the geometrical operation of sheafification.

 A simple way of thinking about Hofmann\&Streicher's groupoid model geometrically is this: think of given groupoid $A$ (interpreting the corresponding type) as the fundamental groupoid of some topological space $T$; then terms $x, y$ present themselves as points of $T$ and the arrow groupoid $Id_{A}(x,y)$ is the  \emph{groupoid of paths} (in $T$) between points $x, y$ with morphisms path homotopies. This  geometrical viewpoint is not only intuitively appealing but also allows for a fruitful conceptual reconstruction of the initial groupoid model. As long as we are working with the usual ``flat'' groupoids (1-groupoids) $A$ and  $Id_{A}(x,y)$ are seen as two different objects. However as we already know from \textbf{6.7} the two things can be naturally made into elements of the same 2-categorical construction of weak 2-groupoid\footnote{
In the end of their paper Hofmann\&Streicher consider this possibility:
\begin{quote}
[I]t might be interesting to view equivalent categories as propositionally equal. This, however, would require ``2-level groupoids'' in which we have morphisms between morphisms and accordingly the identity sets are not necessarily discrete. We do not know whether such structures (or even infinite-level generalizations thereof) can be sensibly organized into a model of type theory. (\cite{Hofmann&Streicher:1998}, p. 108)
\end{quote}
}.  Moreover, we have seen that \emph{from a geometrical point view} there is no reason to stop at level 2: given a topological space $T$ the construction of its fundamental $n$-groupoids $P_{n}$ proceeds inductively ``up to'' the infinite-dimensional groupoid $P_{\omega}$.

Observe that the flat groupoid model of Hofman\&Streicher (or its 2-dimensional presentation) does not have resources for presenting the structure of higher identity types, i.e., types of the form $Id_{Id_{A}(x, y)}(x', y')$  and higher because $\alpha, \beta: Id_{Id_{A}(x, y)}(x', y')$ implies $\alpha = \beta$ (definitionally). Awodey \cite{Awodey:2010} calls this latter property the ``extensionality one dimension up'' ; a model of intensional type theory having this property interprets the non-trivial structure of first-order identity types like $Id_{A}(x,y)$ without reducing terms of this type to definitional identities (so it qualifies as intensional model) but it does reduce (or as mathematicians say \emph{truncate}) all the hierarchy of higher identity types. The aforementioned Awodey\&Warren's model\cite{Awodey&Warren:2009} has the same unwanted property. In his dissertation \cite{Warren:2008} Warren solves this difficulty by extending Hofman\&Streicher's model to higher dimensions with \emph{strict} $\omega$-groupoids and expresses the ``hope that these constructions can be generalized to yield models using other kinds of higher-dimensional structures such as, e.g., the weak $\omega$-groupoids of Kapranov and Voevodsky'' (p. 5) referring to \cite{Kapranov&Voevodsky:1991} where the authors introduce the notion of homotopic $\omega$-groupoid described in \textbf{6.7}. I shall not try to provide more details about the current state of the art in this active research area but stress the fact that homotopy type theory makes precise the idea of construing the notion of identity in terms of morphisms (paths) rather than relations (see \textbf{6.1} and \textbf{6.4}above). To conclude this Chapter I turn now to the foundational aspect of Voevodsky's program. 

\section {Univalent Foundations}
The ultimate ambition of Voevodsky's research program is building new foundations of mathematics, which he calls \emph{Univalent Foundations}. Here is how Voevodsky describes the purpose of this new foundations:

\begin{quote}
The broad motivation behind univalent foundations is a desire to have a system
in which mathematics can be formalized in a manner which is as natural as possible.
Whilst it is possible to encode all of mathematics into Zermelo-Fraenkel set theory,
the manner in which this is done is frequently ugly; worse, when one does so,
there remain many statements of ZF which are mathematically meaningless. This
problem becomes particularly pressing in attempting a computer formalization of
mathematics; in the standard foundations, to write down in full even the most
basic definitions - of isomorphism between sets, or of group structure on a set -
requires many pages of symbols. Univalent foundations seeks to improve on this
situation by providing a system, based on Martin-L\"of's dependent type theory,
whose syntax is tightly wedded to the intended semantical interpretation in the
world of everyday mathematics. In particular, it allows the direct formalization
of the world of homotopy types; indeed, these are the basic entities dealt with by
the system. (\cite{Voevodsky:2011}, p. 7)
\end{quote}

In addition to its theoretical aspect this foundational program pursue a practical purpose: to allow for a systematic use of computer proof-assistants like \emph{Coq} in the everyday mathematical practice. The success of this project would make 
 the \emph{formalization} of mathematical reasoning useful not only for the purposes of ``speculative foundations'' (Lawvere) but also for purely mathematical purposes like proof-checking.  The reason why the univalent foundations can be helpful in this respect is twofold. First, effective proof-assistants including \emph{Coq} are designed on the basis of Martin-L\"of's type theory or similar formal calculi; in fact the interaction between the constructive type theory and computer science has been active and very productive from the early days of both disciplines. Second, the univalent foundations purports to provide the contemporary everyday mathematics (which is quite unlike the everyday mathematics of early 1900-ies)  with a universal  language which has a formal aspect (of Martin-L\"of's type theory) and an intuitive geometrical (or more precisely, homotopical)  aspect. Proof-assistants  like \emph{Coq} become a part of the same combination. 

Before we discuss some epistemological aspects of this new tentative foundations let me describe the idea more precisely. To begin with I assume the $\omega$-groupoid homotopic interpretation of Martin-L\"of's type theory as our intended interpretation of this theory, and in Grothendieck's spirit identify the $\omega$-groupoids with the corresponding topological spaces (of which these groupoids are fundamental groupoids). Accordingly, I call types ``spaces'', call terms ``points'', write $paths_{A}(x, y)$ instead of $Id_{A}(x, y)$ (meaning the space of all paths from $x$ to $y$ in the space $A$) and so on. 

Consider now this inductive definition (see \cite{Voevodsky:2010}): 

\begin{itemize}
\item (i) Given space is called $A$ \emph{contractible} (aka space of $h$-level 0) when there is point $x:A$ connected by a path with each point $y:A$ in such a way that all these paths are homotopic. 
\item (ii) We say that $A$ is a space of $h$-level $n+1$ if for all its points $x, y$ path spaces $paths_{A}(x, y)$ are of $h$-level $n$. 
\end{itemize}
 This completes the definition.

Let us look at members of this hierarchy at low levels:
\begin{itemize}
\item \underline{Level 0}: up to homotopy equivalence there is just one contractible space that we call ``point'' and denote $pt$;
\item {Level 1}: up to homotopy equivalence there are two spaces at this level: the empty space $\emptyset$ and the point $pt$. (For $\emptyset$ condition (ii) is satisfied vacuously; for $pt$ (ii) is satisfied because in $pt$ there exists only one path, which consists of this very point.)  We call $\emptyset, pt$ \emph{truth values}; we also refer to types of this level as \emph{properties} and \emph{propositions}.  Notice that $h$-level $n$ corresponds to the logical level $n-1$: the propositional logic (i.e., the propositional segment of our type theory) lives at $h$-level 1.  
\item  {Level 2}: Types of this level are characterized by the following property: their path spaces are either empty or contractible. So such types are disjoint unions of contractible components (points), or in other words \emph{sets}  of points. This will be our working notion of set available in this framework.
\item  {Level 3}: Types of this level are characterized by the following property: their path spaces are sets (up to homotopy equivalence). These are obviously (ordinary flat) \emph{groupoids} (with path spaces hom-sets). 
\item  {Level 4}: Here we get 2-groupoids
\item
\item
\item {Level n+2}: $n$-groupoids
\item
\end{itemize}
  
The above informal definitions are straightforwardly written in the language of type theory; then one may formally prove some expected simple theorems. Let $iscontr(A)$ and $isaprop(A)$ be formally constructed types `` $A$ is contractible'' and ``$A$ is a proposition'' (for formal definitions see \cite{Voevodsky:2011}, p. 8). Then one formally deduces (= further constructs according to the same general rules) types $isaprop(iscontr(A))$ and $isaprop(isaprop(A))$, which are non-empty and thus ``hold true'' for each type $A$;  informally these latter types tell us that for all $A$ ``A is contractible'' is a proposition and ``A is a proposition'' is again a proposition. With the same technique one defines in this setting type $weq (A, B)$ of  \emph{weak equivalences} (i.e., homotopy equivalences) of given types $A, B$ (as a type of maps $e: A \rightarrow B$ of appropriate sort) and formally proves its expected properties. These formal proves involve a \emph{different} type $isweq (A, B)$ of $h$-level 2, which is a proposition saying that $A, B$ are homotopy equivalent, i.e., that type $weq (A, B)$ is inhabited.)

Remark that the non-trivial homotopic identity (i.e., in Martin-L\"of's terminology the \emph{propositional} identity) works so far only for terms of each given type but not for types themselves. This means that we still treat the identity of types only as ``god-given'' (B\'enabou), i.e., definitional. However in fact we have already smuggled an \emph{intrinsic} identity condition for types when we reasoned above about these types (under the name of ``spaces'') up to homotopy equivalence: this condition is, of course, the homotopy equivalence itself. If we simply replace types by equivalence classes of types in the usual way then we immediately fall back into the circle of classical foundational problems about classes and relations, which we have discussed in Chapter \textbf{5}. And we also abandon in this case at this basic level of our construction the geometrical view according to which sets and classes are homotopy types of special sort (of $h$-level 2). This is why Voevodsky suggests here a different solution. Consider ``big'' type $U$ (for ``universe'') of all ``small'' types $A, B,..$. As terms of $U$ the small types can be compared by path spaces $paths_{U} (A, B)$. Since the small types include equivalence types $weq (A, B)$ we can now produce map $\theta: paths_{U} (A, B) \rightarrow weq (A, B)$, which tells us that propositionally identical types are homotopy equivalent. We want however to have also the inverse map, which would allow us to see homotopy equivalent types as identical. For this reason (and some other reasons, which I leave aside) Voevodsky introduces at this point the \emph{axiom of univalence}, according to which for all types $A, B$ the above map $\theta$ is itself a (homotopy) equivalence. Informally Voevodsky states this axiom (which gives the Univalent Foundations its brand name) as follows:

\begin{quote}
[T]he homotopy theory of the types in the universe should be fully and faithfully reflected by the equality on the universe. \cite{Voevodsky:2011}, p. 8)
\end{quote} 

The axiom of univalence is the only geometrically motivated principle, with which Voevodsky extends Martin-L\"of's type theory. Given this axioms the interchangeable terms ``type'' and ``space''  can be also interchanged with the term ``homotopy type''. Thus in Voevodsky's univalent foundations homotopy types turn to be the elementary bricks for constructing the whole of the mathematical universe \underline{including its logic}.  (The axiom of univalence also has more specific non-trivial consequences, which I shall not consider here.) 

Let us now compare Voevodsky's style of formalization with the standard Hilbert-style formalization. In both cases we have a symbolic syntax and an interpretation of this syntax but the relationships between the syntax and its interpretation is not the same. Remind from \textbf{2.4} that in Hilbert's case the interpretation works separately for logical and non-logical symbols: while the semantic of \emph{logical} symbols is fixed the semantic of non-logical symbols is variable.  In Voevodsky's case we also make difference between the symbolic syntax of type theory and the geometrical interpretation of this syntax (even if Voevodsky's strategy is to merge the two things as close as possible). However Voevodsky's approach does not involve the same difference between logical and non-logical symbols: on his account the same symbols and the same constructions are interpreted both in logical and geometrical terms \footnote{In his lecture \cite{Voevodsky:2010} Voevodsky says that ``logic concerns the $h$-level 1 because this is where the truth values live'' (min. 12 of the recording). This, in my sense, must be understood in the sense that logic necessarily involves objects of the $h$-level 1 (since it uses truth values) but not in the sense that it concerns \emph{only} types of this level. 
}; in particular, in the univalent setting \emph{propositions} are objects (of $h$-level 1) among other mathematical rather than some sort of Fregean metaphysical eternal entities, which live in a separate domain. In other words logic in this setting is wholly \emph{internalized}.

This sheds a new interesting light on Hilbert's distinction between the real and the ideal, and actually on the whole of Formal Axiomatic Method. In his lecture Voevodsky \cite{Voevodsky:2010} criticizes the standard set-theoretic foundations by saying that such foundations require a reconstruction of types of all $h$-levels from types of one particular $h$-level (namely, from types of $h$-level 2, i.e., sets) while his novel approach allows one to reach types of all levels directly. I think that a similar critical argument can and should be made with respect to the ``logical'' $h$-level 1; in this case the critique touches upon the Formal Axiomatic Method itself rather than only upon the idea of using sets in foundations. Observe that Hilbert's mathematical world in fact comprises three kinds of entities: (1) ``real'' finite strings of symbols; (2) propositions and (3) ``ideal'' mathematical objects (that is, all mathematical objects except finite strings of symbols). It doesn't matter whether propositions count as real or ideal on this account; what does matter is that fact that propositions serve here as the only joint between (1) and (3), i.e., between the real and the ideal. Remind Hilbert's idea: any well-formed string of symbols, which represents a system of consistent propositions, determines some ideal mathematical object (or, typically, a system of such objects), which makes all these propositions true. Think of Hilbert's axioms for Euclidean geometry and the Euclidean space (however represented). 
Thus with the exception of the limited domain of finitary combinatorial objects HIlbert's Axiomatic Method is a method of putting every mathematical theory into a propositional form. This propositional reductionism is a tenet of what I have called above the \emph{mathematical logicism} in the broad sense, and it is technically built into Hilbert's Formal Axiomatic Method.  
Now in Voevodsky's vein I remark that propositions form just one particular level of the infinite hierarchy of homotopy types (namely, the $h$-level 1), and that the propositional reductionism is both unreasonable and unnecessary. In any event this propositional reductionism does not make part of Voevodsky's approach. This feature makes Voevodsky's Axiomatic Method quite unlike Hilbert's Axiomatic Method. In Chapter \textbf{9} I generalize upon Lawvere's and Voevodsky's approaches to axiomatization, suggest on this basis a general description of New Axiomatic Method and contrast this new method to Hilbert's method.

Lawvere's tripartition \emph{Formal - Conceptual - Semantic}  \cite{Joyal:2011} suggests thinking about the homotopic interpretation of Martin-L\"of's type theory as the ``conceptual layer'' of this theory; in that sense Voevodsky's univalent foundations can be described as a \emph{conceptual} foundations. To extend this Lawverian perspective on homotopy type theory it seems me also appropriate to consider this theory as a form of \emph{objective logic} in Hegel\&Lawvere's sense. The fact that the homotopy type theory can arguably serve as a foundation of today's everyday mathematics is a strong argument in favor of its objectivity. However as I have already argued in \textbf{6.4} a further test of objectivity involves testing the relevance of this system of logic in natural sciences. I conclude this Chapter with some preliminary remarks, which suggest that this crucial test may be positive. 

Remind from \textbf{5.11} the discussion about the relevance of identity types of Martin-L\"of's type theory in empirical contexts, and in particular in Frege's \emph{Venus} example. I argued then that type $Id_{P}(MS, ES)$ provides an adequate logical form of a proposition stating that the Morning Star and the Evening Star are the same planet Venus; term $e$ of this type is an evidence (aka proof) that the proposition is true. Let us now reconsider this example using the homotopical interpretation of identity types as path groupoids. This interpretation tells us how evidence $e$ should look like: it is a \emph{continuous path} from $MS$ to $ES$. Surprisingly this is exactly what the naive intuition about objects enduring through time suggests: given such a path between $MS$ and $ES$ we qualify these two things as two different appearances of the same thing! (Let me for simplicity limit now the time of observations by 24 hours.) Even if one cannot observe the whole of this path directly one should find other theoretical and empirical reasons of its existence. For the existence of such a path constitutes the very \emph{meaning} of the claim that $MS$ and $ES$ are the same object (Fig. 6.7). As we see the homotopy interpretation of identity types performs really well against this basic notion of identity through time, which is responsible for identity of particles in the classical mechanics and by an extension also in the classical astronomy, biology and many other areas of science. 

\begin{center}
\includegraphics[scale=0.3]{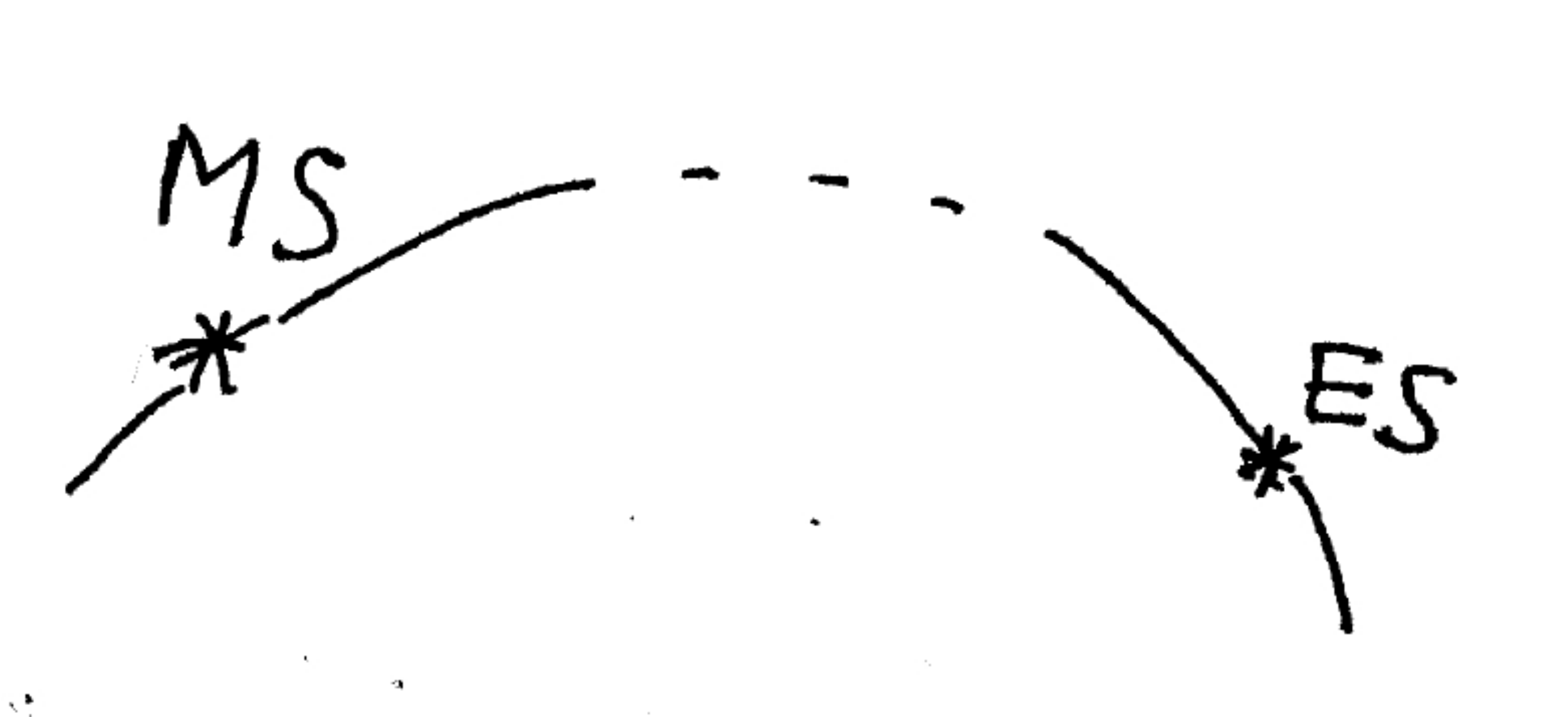}
\end{center}

\begin{center}
Fig. 6.7
\end{center}

The classical mechanical setting, which we have just considered, implicitly assumes the ``extensionality one dimension up'' and leaves no place for path homotopies. However in the spacetime of General Relativity path homotopies become relevant. Consider the phenomenon of \emph{gravitational lensing}: a remote light source $S$ is placed behind massive galaxy $G$ with respect to the observer $O$;  the galaxy bends the light rays in such a way that the observer sees two ``false images'' $S', S''$ of the source (Fig. 6.8).

\begin{center}
\includegraphics[scale=0.3]{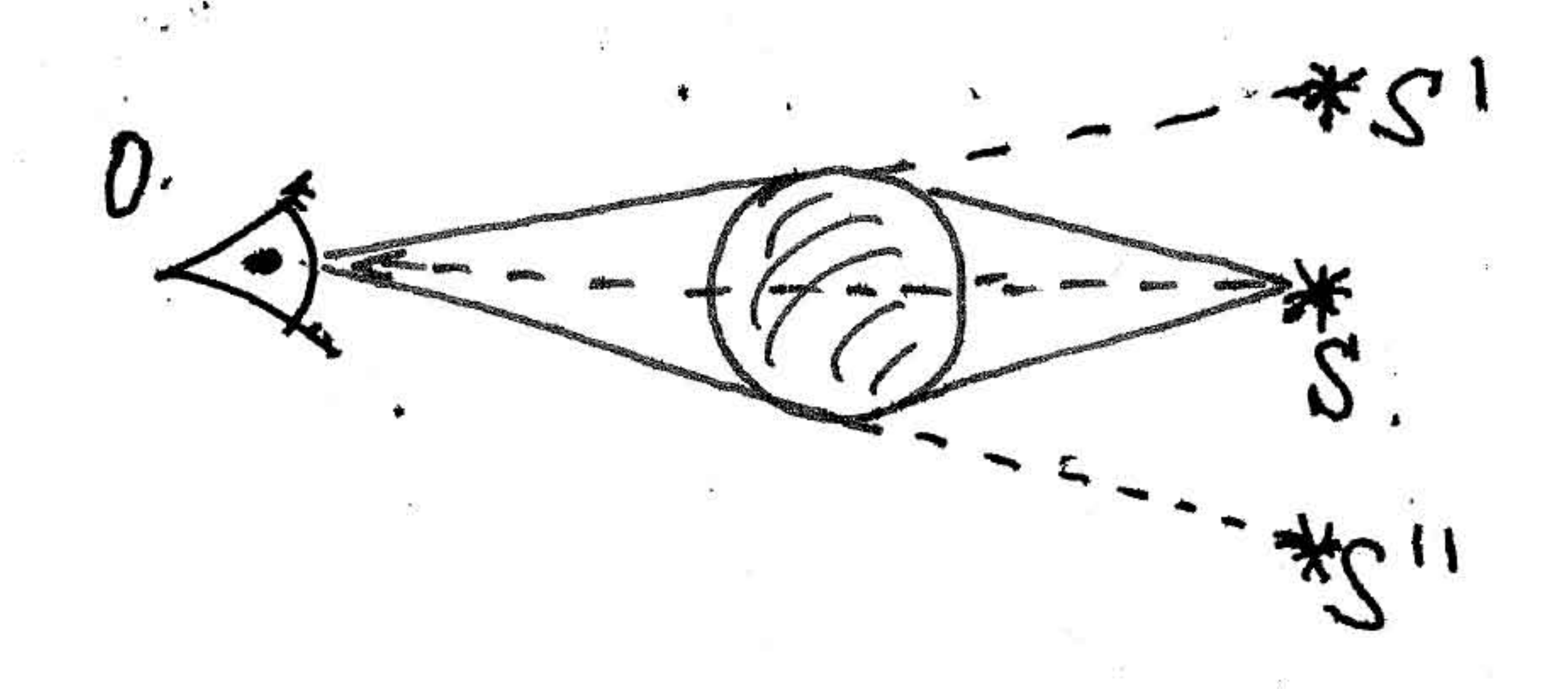}
\end{center}

\begin{center}
Fig. 6.8
\end{center}

In order to identify the source $S$, i.e., to prove that $S'$ and $S''$ is one and the same object, it is not sufficient to trace paths (i.e., light rays) $p_{1}, p_{2}$ between the observer and the source but it is sufficient to prove that $p_{1}, p_{2}$ are homotopic (by performing some actual homotopy or otherwise). However this latter condition is not necessary. What happens when the galaxy is replaced by a spacetime wormhole $W$? Under some realistic conditions the lensing effect persists. But paths $p_{1}, p_{2}$ are no longer homotopic, and so we get in this case a different path groupoid (Fig.6.9).  

\begin{center}
\includegraphics[scale=0.3]{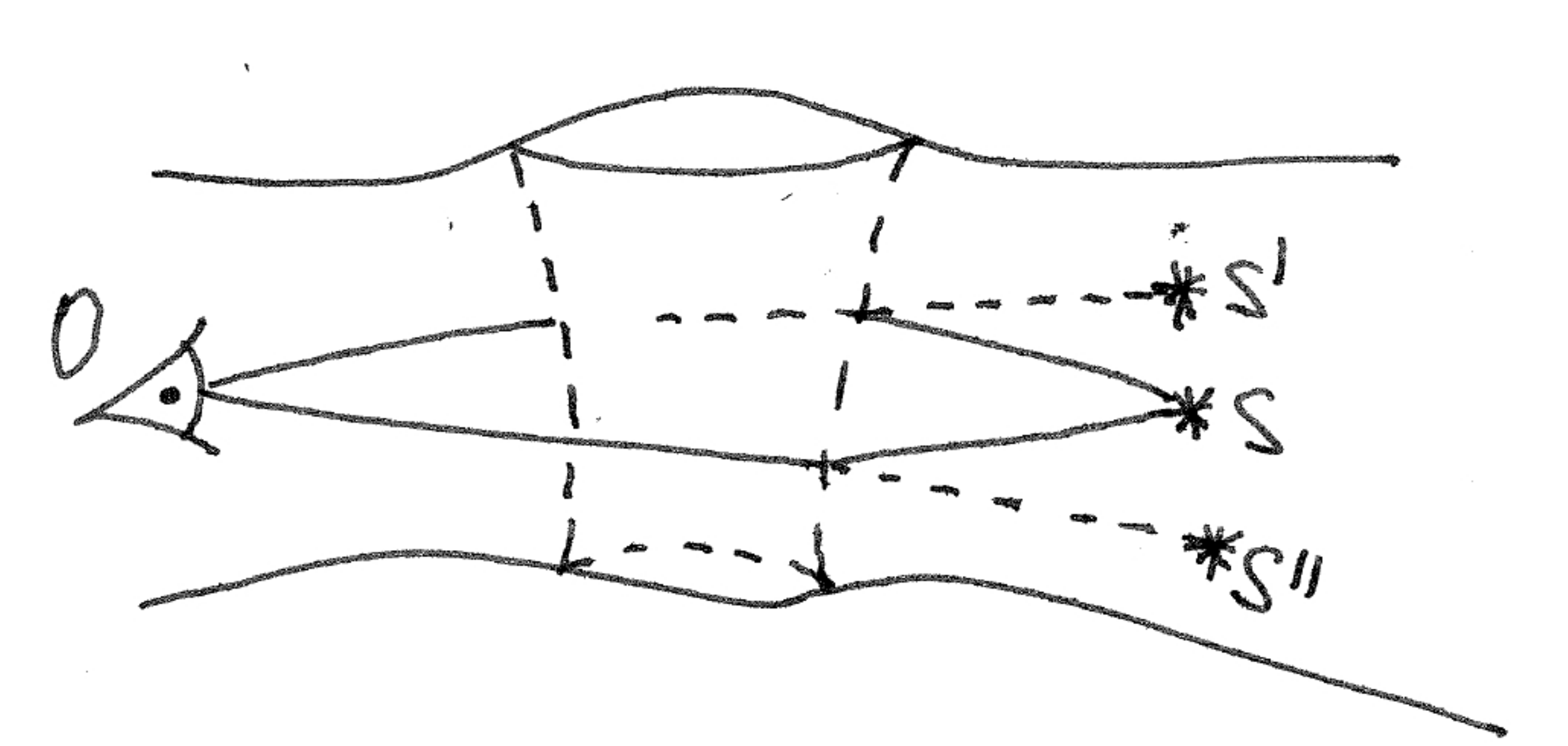}
\end{center}

\begin{center}
Fig. 6.9
\end{center}

It is known that in the quantum field theory homotopical effects have a bearing on the statistics of particles \cite{Suzuki:2011}, and that different statistics are naturally interpreted in terms of different identity conditions of particles \cite{French&Krause:2006}. This suggests that homotopies in the spacetime of general relativity we may have similar effects, so that the replacement of $G$ by $W$ may have a bearing on the identity types of the light source. (In our toy model the corresponding identity types are different but it remains unknown whether this difference accounts for any real physical effect.) These, of course, are bold speculations, which cannot prove anything. However they suggest, in my view, a sound research strategy. Considering the fact that in the past important breakthroughs in science always involved (albeit never reduced to) a renewal of mathematics at the foundational level I assume that such a renewal will necessarily make part of a substantial scientific progress in the future.  Since genuinely original foundational proposal in mathematics are rare, I believe that each one deserves an empirical test. Even if Lobachevsky's attempt to apply his new geometry to the real physical space did not bring immediately any valuable physical result, it established an epistemic strategy, which proved successful later.          

\addcontentsline{toc}{chapter}{Conclusion of Part 2}
\chapter*{Conclusion of Part 2}
We have seen that category theory provides a new and genuinely original mathematical treatment of the identity concept. If one needs to describe this new approach to identity with a single word then, I think, the right word is ``internalization''. The strategy of internalization of identity can be described as follows: instead of assuming after Frege that the identity concept must be fixed before one begins doing mathematics, try to reconstruct this concept \emph{within} mathematics using all appropriate means; after performing such a bootstrap think how to make the newly designed identity concept self-standing by trying to get rid of any pre-established ``external'' notion of identity in the foundations of your science. 

In the category theory such an approach to identity was initiated by Lawvere in 1970 \cite{Lawvere:1970b}. Independently, and wholly outside category theory, Martin-L\"of developed in early 1980ies his version of type theory \cite{Martin-Lof:1984} where he also treated the identity concept internally. The complexity of the identity concept in Martin-L\"of's theory was not designed in purpose but was discovered later by Hofman and Streicher \cite{Hofmann&Streicher:1998} who refuted the so-called Uniqueness of Identity Proof conjecture (aka UIP conjecture) according to which axioms of Martin-L\"of's type theory allegedly implied that all identity types are simple  in the following sense: given two proofs $p_{1}, p_{2}: Id_{A}(x, y)$ there is always a proof $t: Id_{Id_{A}(x, y)}(p_{1}, p_{2})$ showing that $p_{1}, p_{2}$ are identical. So it is fair to say that the intrinsic complexity of Martin-L\"of's identity types was discovered rather than put there by hand. From the retrospect it is clear that since type-theoretic considerations were also behind Lawvere's ideas about identity the two lines of research had to meet at some point. But when this really happened several years ago it was quite surprising anyway.  

In Chapter \textbf{5} I tried to sketch the history of thinking about identity and equality in mathematics (without following the historical order of events) paying a special attention to the view, which has emerged in the end of the 19th and the beginning of the 20th century. In the following Chapter \textbf{6} I showed how the identity concept is treated in category theory. Two things seem me clear: first, that the new categorical approaches to identity are philosophically relevant (and sometimes also strongly philosophically motivated) and, second, that the standard Frege-Russell's view on identity is not adequate to these new approaches. We have seen that the contemporary mathematics played a crucial role in Frege's and Russell's then-new thinking about identity: these philosophers improved on the traditional logic with new mathematical concepts (like that of function) and, reciprocally, tried to improve on their contemporary mathematics with a new logical regimentation. It appears to me that the first part of this project was by far more fruitful. Having this in mind I tried in Chapter \textbf{6} to learn some philosophical lessons from recent mathematical theories about identity rather than impose any particular philosophical agenda upon mathematics. Nevertheless I believe that in a longer run such philosophical lessons may have a feedback effect. The main philosophical lesson about the identity concept that I have learnt from the categorical logic in general, and the homotopy type theory, in particular, is this: the identity concept (as well as other basic logical notions) must be thought of on equal footing with all other scientific concepts but not as a prerequisite for doing science.      
        
\part{Subjective Intuitions and Objective Structures}
\chapter{How Mathematical Concepts Get Their Bodies}
\section{Changing Intuition}
By Kant's popular word  ``thoughts without content are empty, intuitions without concepts are blind'' \emph{Critique of Pure Reason} (B75). In \textbf{1.3} we have seen how this general claim applies to Euclid's geometry, that is, how exactly concepts and intuitions work here closely together. In this Chapter I consider a question that Kant himself never systematically studied, namely the question of how mathematics develops. It is obvious that mathematical concepts change through time: some of them change their content preserving the name (and something like ``general idea''), some new mathematical concepts regularly emerge and some old concepts stop growing and die off.  Today's philosophers of mathematics like biologists of the 18th century still often replace a study of living populations by a study of dead specimen preserved with some special techniques.  But even if the evolutionary view in the philosophy of mathematics is not common the phenomenon of 
conceptual change in mathematics is so evident that nobody can deny its very existence.

What about mathematical intuitions? Do mathematical intuitions change along with mathematical concepts or not? This question does not have a similar obvious answer. For one thing it is not clear how one may notice an intuitive change if any. (Hereafter I use the expression ``intuitive change'' formed by analogy with the established expression ``conceptual change'' in the sense of \underline{changing intuition}, not in the sense of change, which has an intuitive meaning or has some other special relation to intuition.) Indeed, the emergence of a new concept or a change of content of some earlier existing concept can be precisely described on the basis of textual and other relevant symbolic evidences, in particular, by comparing the corresponding mathematical definitions. But how one can possibly access and study raw intuitions? One may even argue that the notion raw intuition is ill-formed because if such a thing would exist it could not be rationally known. 

In addition to this general epistemological difficulty there are two more specific reasons why the issue of intuitive change in mathematics has been never systematically studied, so that the very existence of such a phenomenon remains questionable. One reason has to do with some interpretations of Kant's philosophy, according to which mathematical intuitions are somehow wired into our brain and cannot possibly change without making a change in our biological constitution as \emph{homo sapiens}. The other reason is related to the former. Since throughout the 19th century mathematical concepts manifestly changed, and since the new changed concepts were no longer compatible with older intuitions, some philosophers including Russell convinced themselves and tried to convince other people that Kant was altogether wrong about the fundamental role of intuition in mathematics; they developed a new philosophy of mathematics, which didn't reserve for intuition any significant epistemic role. The influence of this philosophy, which in the previous Chapters I called mathematical logicism, is also a reason why the phenomenon of intuitive change in mathematics remains largely unnoticed by philosophers. 

In this Chapter I shall not try to develop a systematics philosophical critique of Kant's and Russell's views on mathematics but instead try to demonstrate at some historical examples that the intuitive change in mathematics is a real historical phenomenon. My general picture is this. In the beginning of the 20th century many people overdramatized new tendencies in mathematics and convinced themselves that mathematics during several decades wholly changed its nature. I don't believe that such a radical change really happened. A broader historical outlook shows that mathematics develops both through the invention of new concepts and the invention of new intuitions. However it often happens that one of these parallel developments goes ahead the other. When the intuitive development goes ahead mathematics is full of poorly conceptualized intuitions and so mathematicians work on their conceptualization. When the conceptual development goes ahead mathematics is full of poorly intuited concepts and mathematicians work on building the appropriate intuitions. Both kind of situations may co-exist in the historical time. However in a given historical period (and in a given mathematical community) situations of one of these two sorts may dominate. In other words, in the real history of mathematics we find a dialectical interplay between conceptual and intuitive developments rather than the stable Kantian harmony\footnote{
In fact Kant himself did not describe the relationships between concepts and intuitions as a perfect harmony. His notion of \emph{regulative idea} accounts for the situation when ``the concepts go ahead''. However Kant did not treat the question of historical development of mathematics systematically.}.
 The first half of the 20-th century was the time when the concept-building clearly dominated over the building of new intuitions. The development of new intuitions was not wholly suppressed but it became a prerogative of a narrow circle of creative mathematicians who invented new mathematical concepts and solved with them real mathematical problems; the rest of the community received these concepts mostly in the sterilized form of formal and quasi-formal (like in Bourbaki's case) axiomatic theories. Hence the popular thesis about the abstract character of the new mathematics. However already in the second half of the 20th century this general tendency began to change; the homotopic approach in the category theory described in \textbf{6.9} is just one evidence of this change among many others. Thus I claim that the tendency towards a ``higher abstraction'' of the 20th century mathematics is nothing but a local effect comparable with similar tendencies taking place in other historical periods (see below); it does not represent a global tendency in the historical development of mathematics.

Before I substantiate this controversial claim with some historical evidences let me make two methodological points. First, I should answer the question ``How one can study raw intuitions?''. My answer is that one can \underline{not} study raw intuitions just one cannot study bare concepts whatever this might mean. Here I am wholly with Kant. What one, however, \emph{can} do is to reconstruct the historical dialectics of concepts and intuitions and specify some situations when ``intuitions go ahead'' and some other situations when ``concepts go ahead''. I'm not going to represent here in this way the overall history of mathematics but I shall demonstrate what I have in mind with several examples. Second, I should say more precisely what I mean by intuition. I borrow my basic notion of mathematical intuition from Kant - both in the sense of Kant's Transcendental Aesthetics and in the sense of Kant's Doctrine of Method (see \textbf{1.3} above). One may argue that Kant's notion of intuition does not make sense outside the mathematical context, which Kant had in mind when he developed this notion. I cannot see why this can be true. As far as I can see in an appropriately generalized form Kant's notion of mathematical intuition applies to mathematics throughout its history. The following historical examples demonstrate this secondary claim too.

\section{Form and Motion}
It is a truism that in the geometrical Books of his  \emph{Elements} Euclid develops a systematic theory. It is perhaps less obvious that this systematic approach concerns not only geometrical concepts but also geometrical intuitions. The intuitive geometrical content of Euclid's \emph{Elements} is organized as follows. It comprises 

\begin{itemize}
\item primitive intuitions corresponding to concepts of \emph{point}, \emph{straight line} and \emph{circle};
\item complex intuitions built with the primitive intuitions. Such combination of primitive intuitions is colloquially called a ``construction by ruler and compass''. In Euclid's  \emph{Elements} such constructions are regulated by Postulates 1-3 (see \textbf{1.4}).
\end{itemize}

Constructions by ruler and compass allow for an intuitive grasp on complex geometrical configurations, which cannot be so grasped immediately. It is not uncommon in the literature to identify the notion of intuition with a primitive \emph{immediate} intuition but this restriction is not justified. Mathematics in general and geometry in particular requires complex intuitions as well as complex concepts. Let me demonstrate the notion of complex intuition with an example. Think of a regular 1024-gon (polygon with 1024 sides) and proceed as follows. Inscribe a square into a circle. Then by halving angles between diagonals of the square construct an octagon inscribed into the same circle as shown at the below diagram (fig 7.1).

\begin{center}
\includegraphics[scale=0.6]{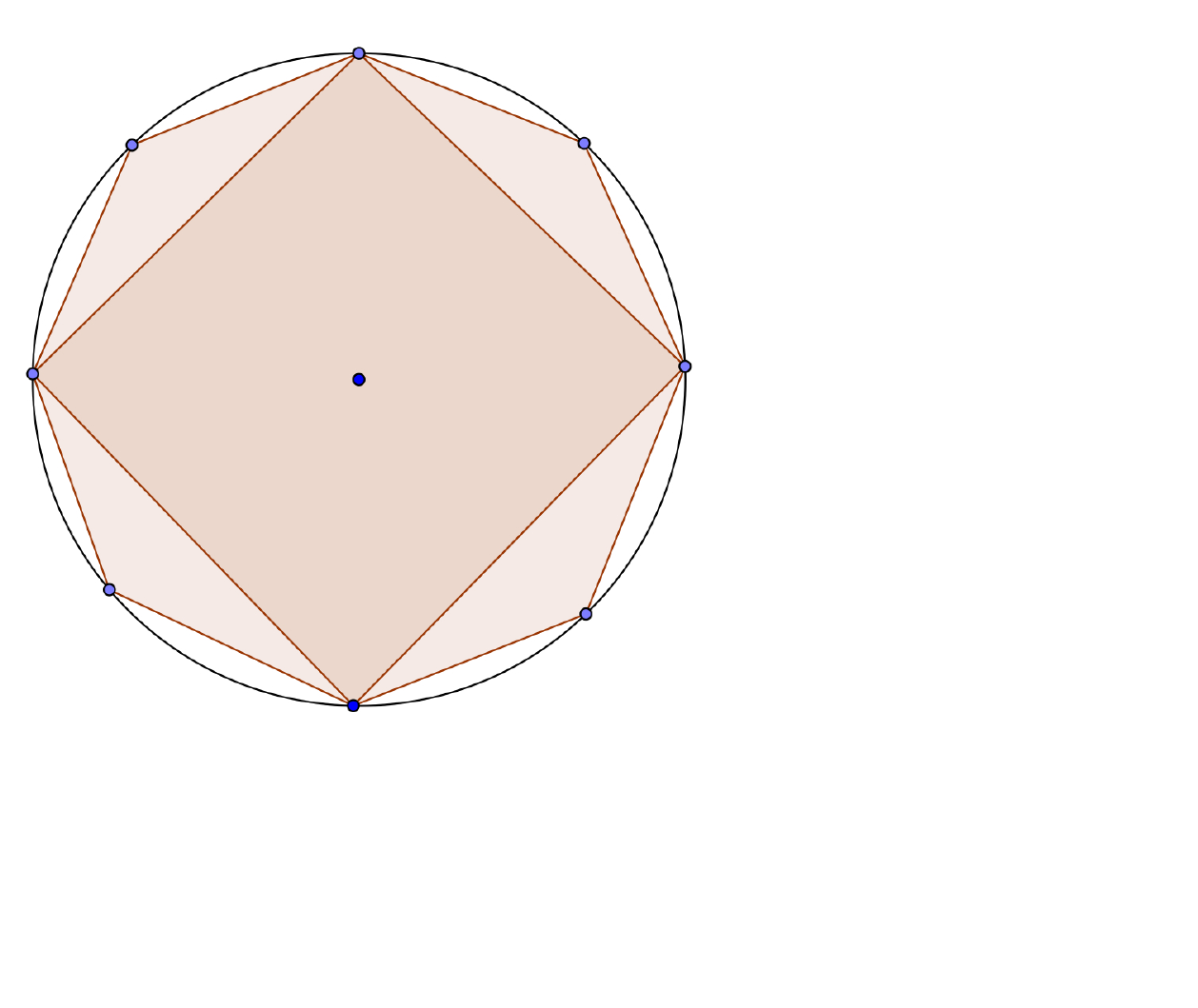}
\end{center}

\begin{center}
Fig. 7.1
\end{center}

Then halve the angles between all the neighboring diagonals of the octagon, get 16-gon and repeat the same procedure 5 more times. (I leave obvious details of this construction to the reader.)  Although the resulting figure cannot be grasped immediately it is now intuitively grasped through the above procedure, which involves only a few primitive geometrical intuitions. In addition to these primitive geometrical intuitions it involves an intuition lying behind the very notion of \emph{procedure} and, more specifically, behind the idea of \emph{iteration}. This latter intuition is temporal rather than spatial. The construction of 1024-gon also involves some basic arithmetical intuitions, which allow one to calculate powers of 2.

A chiligon (1000-gon) is often referred to (after Descartes) as an example of geometrical object, which can be thought of but cannot be imagined. The above example shows that the big number of sides of a polygon does not always make an intuitive grasp impossible. However the method of intuitive grasping used in this example does not work for the chiligon because it cannot be constructed by ruler and compass. A similar difficulty concerns a regular heptagon, doubled cube and many other relatively simple geometrical concepts. This difficulty shows that the intuitive and the conceptual legs of Euclid's geometry don't always match each other as expected\footnote{As it has been first conjectured by Gauss and proved in 1837 by Wantzel \cite{Wantzel:1837} a regular polygon is constructible by ruler and compass if and only if its number of sides is a product of a power of 2 and a Fermat prime. Fermat primes are primes of the form $2^{2^n} + 1$ where $n$ is a positive integer.}.

The problem was actually twofold. On the one hand, Greek geometers came across a number of geometrical concepts, which they couldn't provide with an intuitive support by standard means, i.e. by performing a construction by ruler and compass.  On the other hand, they produced a number of intuitively appealing constructions by other means than ruler and compass but didn't know how to treat them theoretically; I am talking about a cycloid, a spiral and other so-called \emph{mechanical} curves. (Unlike the circle such curves could not be described without referring to a mechanical setting.) Thus they got a number of problematic cases some of which could be described as poorly intuited concepts while some other could be described as poorly conceptualized intuitions. They applied mechanical curves for solving problems unsolvable by ruler and compass; such solutions pointed to new ways of matching concepts with intuitions. However they didn't develop any alternative systematic theory which could compete with the standard theory based on ruler and compass.

This situation changed significantly only in the Early Modern times. At this point of history the ancient notion of primitive geometrical form lost its earlier appeal, and the idea of construction by ruler and compass was already commonly seen as a pure convention like today. In his \emph{Geometry} of 1637 \cite{Descartes:2008} Descartes describes an alternative instrument (sometimes called \emph{Cartesian compass}) able to produce a large variety of curves other than circles;  he argues that constructions made with this new mechanism and the traditional ruler should be treated on equal footing with traditional constructions by ruler and (standard) compass\footnote{Descartes does not give up the ancient distinction between mechanical and (properly) geometrical curves altogether but changes its sense: on Descartes' account only non-algebraic curves qualify as mechanical.  Following the ancient pattern Descartes claims that such mechanical curves do not belong to the subject-matter of pure geometry. However this ban of transcendental curves was hardly ever respected in Modern mathematics.} 
. The Cartesian compass is, of course, only a by-product of the Early Modern reform of mathematics but it well illustrates the fact that mechanical intuitions, which remained in Greek geometry marginal, in the new geometry became central.

In the beginning of his \emph{Geometrical Lectures} of 1735 \cite{Barrow:2006} Barrow claims that all geometrical objects are ultimately generated by motion. Although today we can describe Euclid's and Barrow's geometrical universes by the same term ``Euclidean space'' Barrow's universe is a great extension of Euclid's universe. Curves produced by arbitrary continuous motions were somehow ``around'' in Euclid's universe but they did not belong to it properly. Obviously in order to treat such curves mathematically Barrow, Newton and other people had to invent a wholly new conceptual apparatus. But they also had to clarify the geometrical intuitions that fitted the new conceptual apparatus. This new geometrical intuition no longer supports the view on geometrical objects as combinations of primitive forms but supports the view on geometrical objects as traces of continuous motions. 

\section{Non-Euclidean Intuition}

The invention of Non-Euclidean geometries in 19th century is often referred as a crucial moment of history when geometry lost its traditional link with empirical spatial intuitions. But this view is plainly wrong. This invention involved the development of new spatial intuitions, not abandoning of spatial intuitions altogether. The idea to cut geometry from spatial intuitions came about later, in the very end of 19th century. It soon developed into a powerful trend, which we have discussed in Chapter \textbf{2}. The biased view on the 19th century geometry  mentioned in the beginning of this paragraph is a by-product of this later trend
\footnote{
After explaining the basics of Hilbert's Axiomatic Method in a college-level textbook its author writes \cite{Greenberg:1974}
\begin{quote} 
The formalist viewpoint just stated is a radical departure from the older notion that
mathematics asserts ``absolute truths'', a notion that was destroyed once and for all
by the discovery of Non-Euclidean geometry. This discovery has had a liberating effect
on mathematics, who now feel free to invent any set of axioms they wish and deduce
conclusions from them. In fact this freedom may account for the great increase in the
scope and generality of modern mathematics. (p. ??)
\end{quote}
Leaving aside the issue of mathematical truth I would like only to stress that Greenberg's account of events is wrong historically.  Although it is not unreasonable to consider the ``discovery of Non-Euclidean geometry'' as a historical cause, which gave rise to the Formal Axiomatic Method, there is a historical distance of at least 40 years between these events (I count from Lobachevsky's German publication of 1840 \cite{Lobachevsky:1840} to Pasch's \emph{Foundations of New Geometry} of 1882 \cite{Pasch:1882}, which was the most important predecessor of influential Hilbert's \emph{Foundations} of 1899 \cite{Hilbert:1899} 
One may suggest that during this period of time the new geometry remained in an obscure state, so that before Pash and Hilbert treated it axiomatically nobody could make good sense of it. This is again historically false: there were a lot of important work in the new geometry during this period, which include Riemann's pathbreaking generalization \cite{Riemann:1854} and Klein's synthesis \cite{Klein:1871}, \cite{Klein:1873} that gave to Lobachevsky's geometry its modern name ``hyperbolic''. If these argument are not convincing for one who believes that the  Hilbert-style axiomatization is the only way to make a mathematical theory reasonable, I may remark that my opponent should apply this principle also to the Euclidean geometry. This view implies that the Euclidean geometry was in an obscure state before Hilbert too.   
}. 

It is well-known that the discovery of non-Euclidean geometries was a result of long attempts to prove Euclid's 5th Postulate: failing to get a  \emph{reductio ad absurdum} of the negation of this Postulate some people guessed that they were exploring a new vast territory rather than approach the expected dead end. But why did the 5th Postulate attract so much attention to begin with? An obvious answer is that 5th Postulate does not have the same intuitive appeal as other Euclid's first principles; this is why since Antiquity this Postulate was widely seen as problematic. This remark shows that the often repeated claim that the human ``natural'' spatial intuition is intrinsically Euclidean does not stand against obvious historical evidence. Were this true geometers would always take the 5th Postulate as a self-evident principle along with other Euclid's Postulates and Axioms and wouldn't spent significant efforts trying to prove this particular Postulate on the basis of other principles. It was the mathematical education rather than human nature that later made the majority of people to take the Euclidean intuition for granted, so only the tiny minority of experts could still see and address the problem. 

Today Lobachevsky's hyperbolic geometry is usually taught in a Hilbert-style axiomatic form provided with some \emph{models} including Beltrami pseudosphere model, Klein disc model and Poincar\'e semisphere model \cite{Greenberg:1974}. I recognize that this axiomatic presentation of geometry has certain merits; in any event it provides a pattern of axiomatic thinking, which is crucially important for the 20th century mathematics. However such a presentation of hyperbolic geometry  suggests a wrong idea that there are only these two conceptions of this geometrical theory: one purely \emph{formal} and the other one \emph{interpreted} with one of the aforementioned models.  I claim that there is yet another way of thinking about the hyperbolic geometry, which is Lobachevsky's own. Lobachevsky associates geometrical intuitions to geometrical concepts directly in the same way in which anyone does this when one studies Euclidean geometry. Since Lobachevsky thinks about his new geometry as a generalization of the traditional geometry rather than as a separate theory
\footnote{
Here is an explicit statement:
\begin{quote}
The principle conclusion, to which I arrived .... was the possibility of Geometry in a
broader sense than it has been [earlier] presented by Euclid the Founder. This extended
notion of this science [=of Geometry] I call Imaginary Geometry; Usual [=Euclidean]
Geometry is included in it as a particular case. (\cite{Lobachevsky:1949}, Introduction, my translation from Russian)
\end{quote}
}
 his geometrical intuition (which I call \emph{hyperbolic}) is an extension of the Euclidean intuition rather than a wholly new type of geometrical representation. I stress this point in order to confront the widespread philosophical myth according to which the invention of Non-Euclidean geometry required an abrupt departure from the ``usual'' geometrical intuition. In what follows I demonstrate Lobachevsky's hyperbolic intuition at work (after \cite{Lobachevsky:1840}); I treat Lobachevsky's case more attentively than other historical cases considered in this Chapter because of its importance for the history of the Axiomatic Method. 
 
Instead of Euclid's Fifth Postulate (P5) Lobachevsky uses the Axiom of Parallels (AP, aka Playfair Axiom) known to be equivalent to P5 since the Antiquity:

\begin{quote}(AP) Given a line and a point outside this line there is unique other line which is parallel to the given line and passes through the given point.\end{quote}

Here the term ``parallels'' stands as usual for straight lines having no common points. We shall shortly see how Lobachevsky changes this Euclidean terminology. For a
terminological convenience I shall call a given straight line \emph{secant} of another given
straight line when the two lines intersect (in a single point). Let's now make the
required construction and see what the intuition shows about it. Although the intuition does \emph{not} show whether AP is true or not it shows several other important things (all of which can be proved from Euclid's first principles without using P5 unless it is explicitly stated otherwise):

\begin{itemize}
\item (i) Parallel lines exist (unlike round squares); moreover through each given point $P$ outside any given straight line $l$ passes at least one parallel line $m$. 

\item (ii) Given a straight line and point $P$ outside this line there exist secants of the given
line passing through the given point. To construct a secant take any point of the given
line and connect it to the given point outside this line.

\item (iii) Let $PS$ be perpendicular to $l$ and $A$ be a point of $l$. Consider a straight line $PR$ such that angle $SPR$ is a proper part of angle $SPA$ (and hence is less than angle $SPA$). I shall call line $PR$ \emph{lower} than line $PA$ and call $PA$ \emph{upper} than $PR$). Beware that this definition involves the perpendicular $PS$, and so depends on the choice of $P$. Then $PR$ intersects $l$ in some point $B$, i.e. it is a secant. Thus a line, which is lower than a given secant, also is a secant (Fig.7.2). (In order to prove (iii) rigorously one needs Pasch's axiom which Lobachevsky never mentions but often uses tacitly. This axiom first introduced in \cite{Pasch:1882} reads: Given a triangle and a straight line intersecting one of the triangle's sides but passing through none of the triangle's apexes the given line intersects one of the two other sides of the given triangle. To apply this axiom to the given case one needs a simple additional construction that I leave to the reader. 

\begin{center}
\includegraphics[scale=0.4]{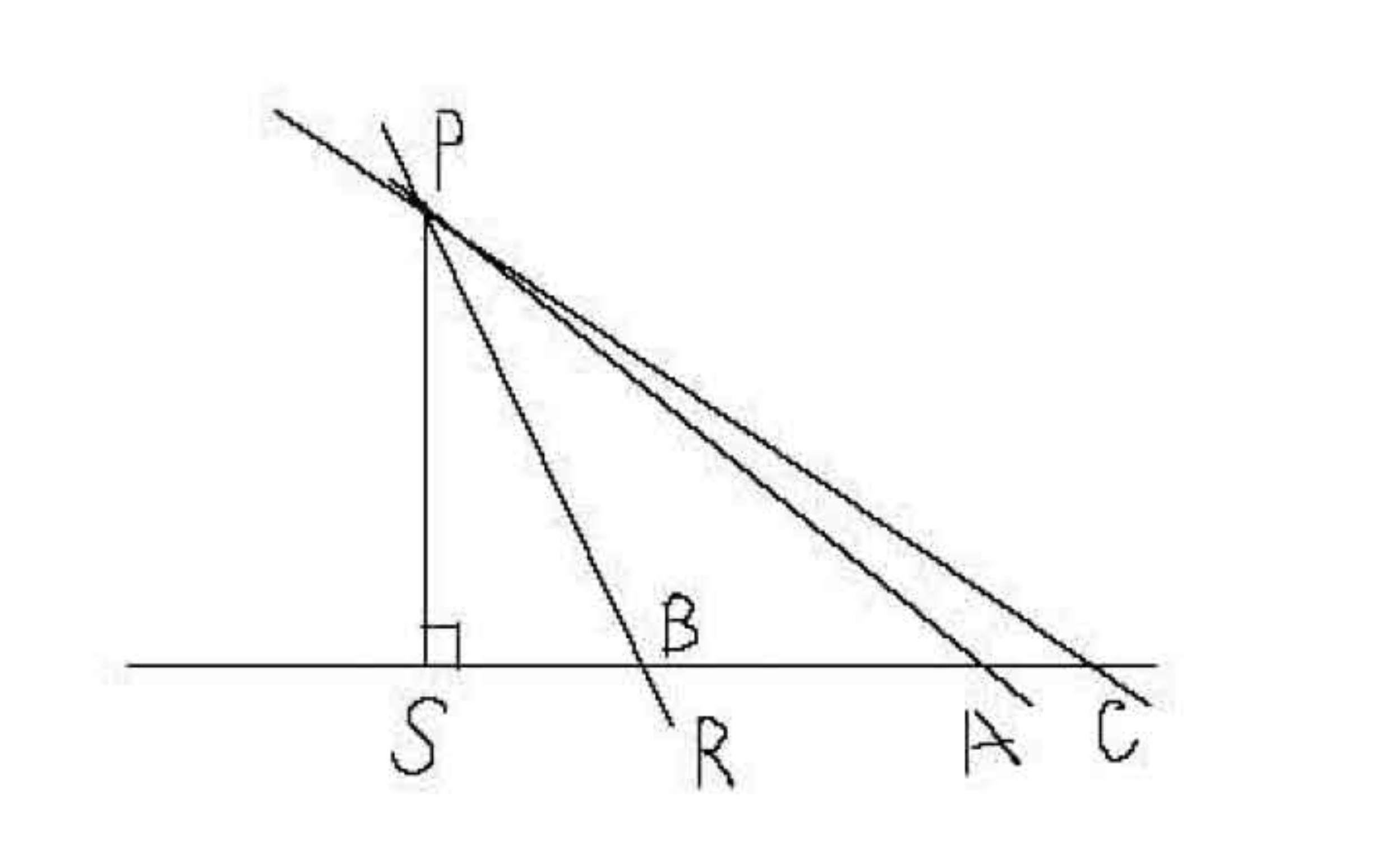}
\end{center}
\begin{center}
Fig.7.2
\end{center}

\item (iv) There exist no upper bound for secants of a given line passing through a given
point outside this given line. For given some secant $PA$ one can always take a further
point $C$ such that $A$ will lay between $S$ and $C$ and so secant $PC$ be upper than the given secant $PA$ (Fig. 7.2).

\item(v) Let $m$ be parallel to $l$ , which is constructed as in (i). Let $n$ be another parallel to $l$ passing through the same point $P$. Suppose that $n$ is lower than $m$ (obviously this condition doesn't restrict the generality). Then any straight line which is upper than $n$ and lower than $m$ is also parallel to $l$ (Fig. 7.3)

\begin{center}
\includegraphics[scale=0.4]{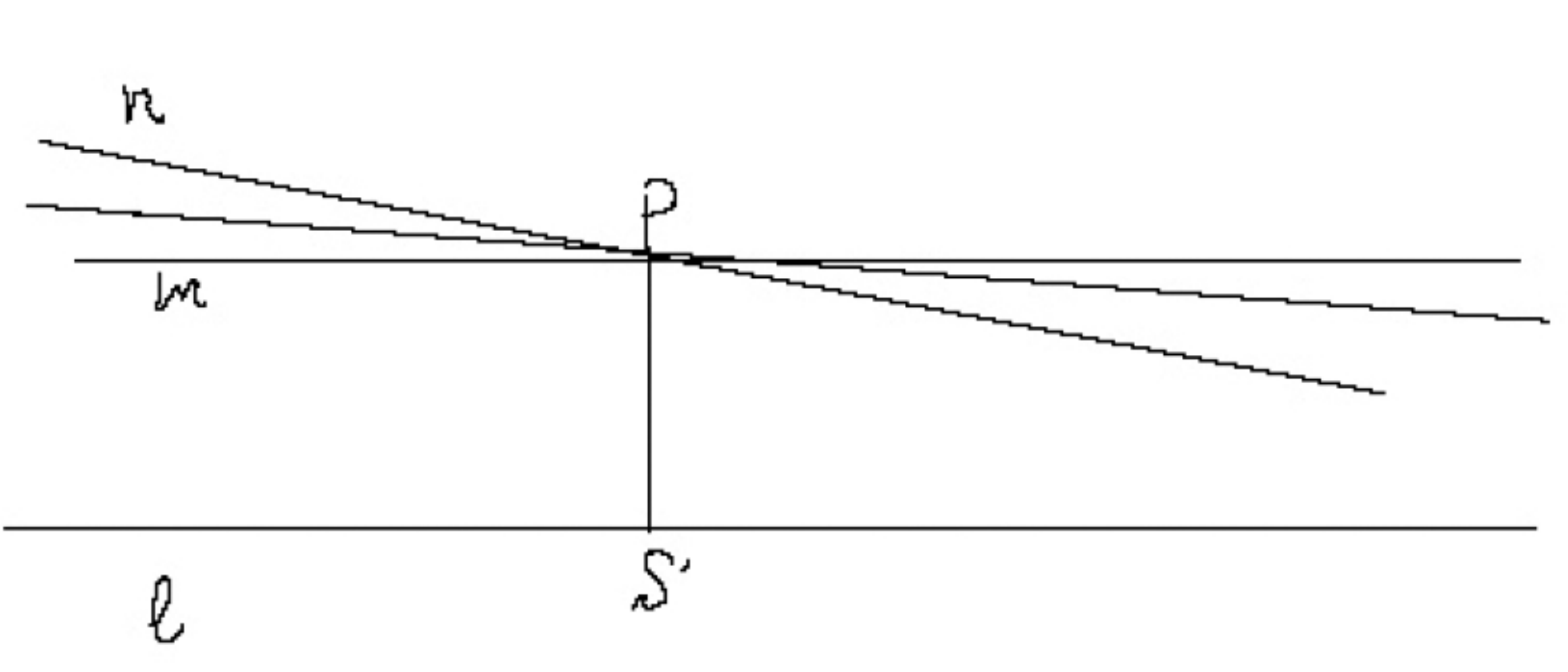} 
\end{center} 
\begin{center}
Fig. 7.3
\end{center}

\item (vi) Parallels to a given straight line passing through a given point have a lower bound.
(In order to prove this rigorously one needs some continuity principle like one asserting the existence of Dedekind cuts. Lobachevsky doesn't state such a principle explicitly. 

\item (vii) Any straight line $PA $- a secant or a parallel - passing through point $P$ as shown at Fig. 2 is wholly characterized by its characteristic angle $SPA$. In particular this
concerns the lowest parallel mentioned in (vi). Let the measure of $SPA$ corresponding
to the case of the lowest parallel be $\alpha$. Now it is clear that by an appropriate choice of
$l$ and $P$ one can make $\alpha$ as close to $\pi/2$ as one wishes. For given any angle $SPA < \pi/2$ it is always possible to drop perpendicular $AT$ on $PS$ (Fig. 7.4). Then $PA$ is a secant of $AT$ and so by (iii) all parallels to $AT$ including its lowest parallels are upper than $PA$. Hence the value of $\alpha$ corresponding to straight line $AT$ and point $P$ outside this line is between $SPA$ and $\pi/2$. Since the only variable parameter of the configuration is the distance $d$ between the given straight line and the given point outside this line the angle $\alpha$ is wholly determined by this distance.

\begin{center}
\includegraphics[scale=0.4]{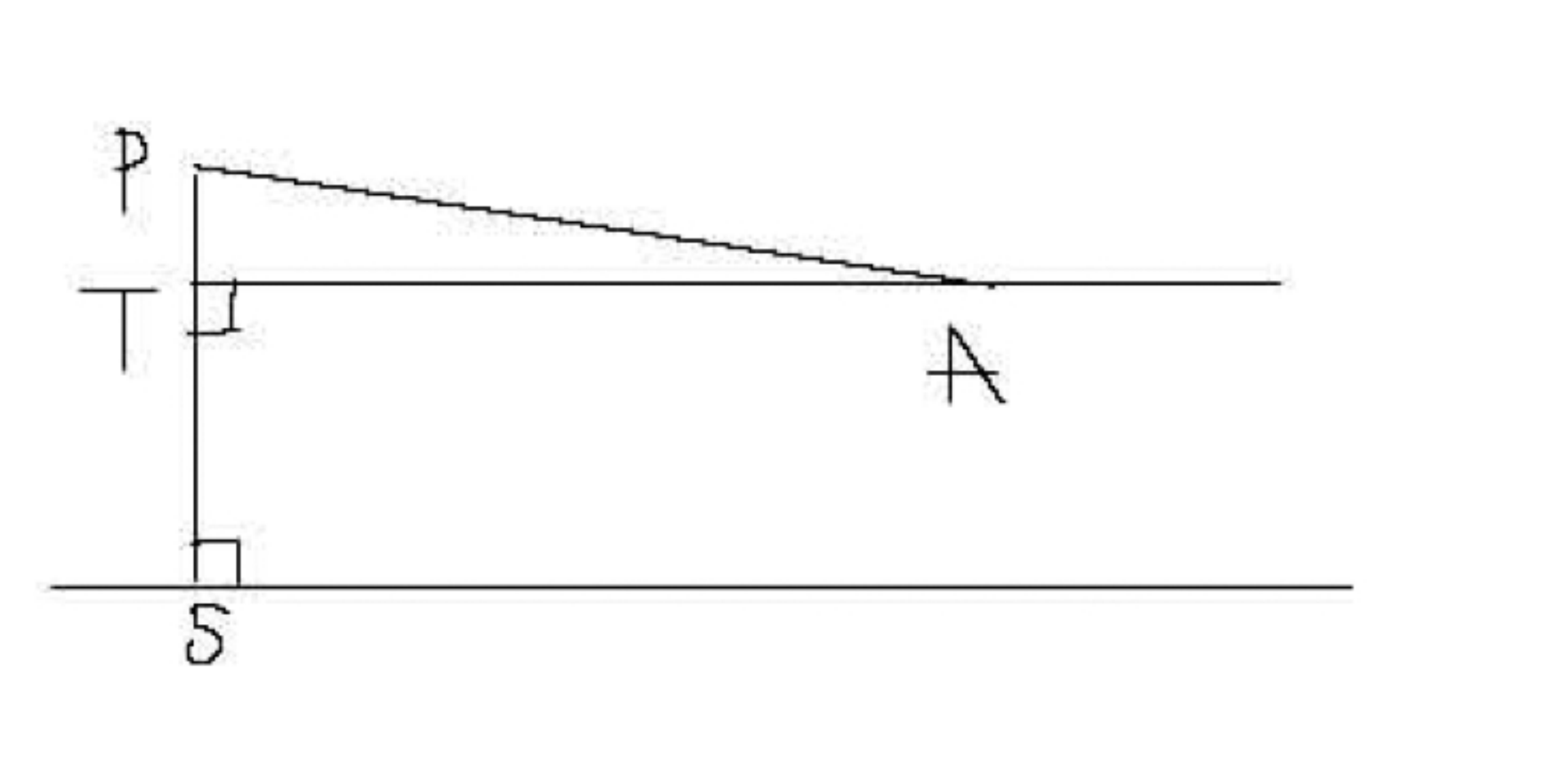} 
\end{center} 
\begin{center}
Fig. 7.4
\end{center}

\end{itemize}

Propositions (i - vii) provide an intuitive basis for Lobachevsky's  geometry (see \cite{Lobachevsky:2007}, propositions 7, 16, 21). Then he proceeds as follows. First, he makes a terminological change: he calls ``parallels'' (not just non-intersecting straight lines but) the two boundary lines which separate secants from non-secants (i.e. parallels in the usual terminology) passing through a given point. So in Lobachevsky's terms there exist exactly \emph{two} parallels to a given straight lines passing through a given point, which may eventually coincide if AP holds (i.e. in the Euclidean case). For further references I
shall call these two parallels \emph{right} and \emph{left} (remembering that this assignment of parity is arbitrary). Since Lobachevsky's definition of parallels involves the choice of $P$ it is not immediately clear that the parallelism so defined is an equivalence. So Lobachevsky must show that the property of being parallel (in his new sense) to a given straight line is independent of this choice (\cite{Lobachevsky:2007}, proposition 17), and
that the relation of being parallel is symmetric and transitive (while reflexivity is taken for
granted by the usual convention) (\cite{Lobachevsky:2007}, propositions 18, 25). For the obvious reason the transitivity applies here only for parallels of the same orientation, i.e., separately for right and for left parallels. Lobachevsky provides the required proofs making them in the traditional synthetic Euclid-style manner. Then Lobachevsky proves some further properties of parallels, in particular the fact that the angle $\alpha$ characterizing a parallel (see (vii) above) can be made however close to $\pi/2$ and however close to zero (\cite{Lobachevsky:2007}, proposition 23). This immediately implies that if AP does not hold then given an angle $ABC$, however small, there always exist a straight line $l$ laying wholly inside this angle and intersecting none of its two sides (Fig. 7.5). This is already by far more counterintuitive than (i-vii) but still not counterintuitive enough to rule out this construction as absurd and on this ground to claim a proof of AP.

\begin{center}
\includegraphics[scale=0.4]{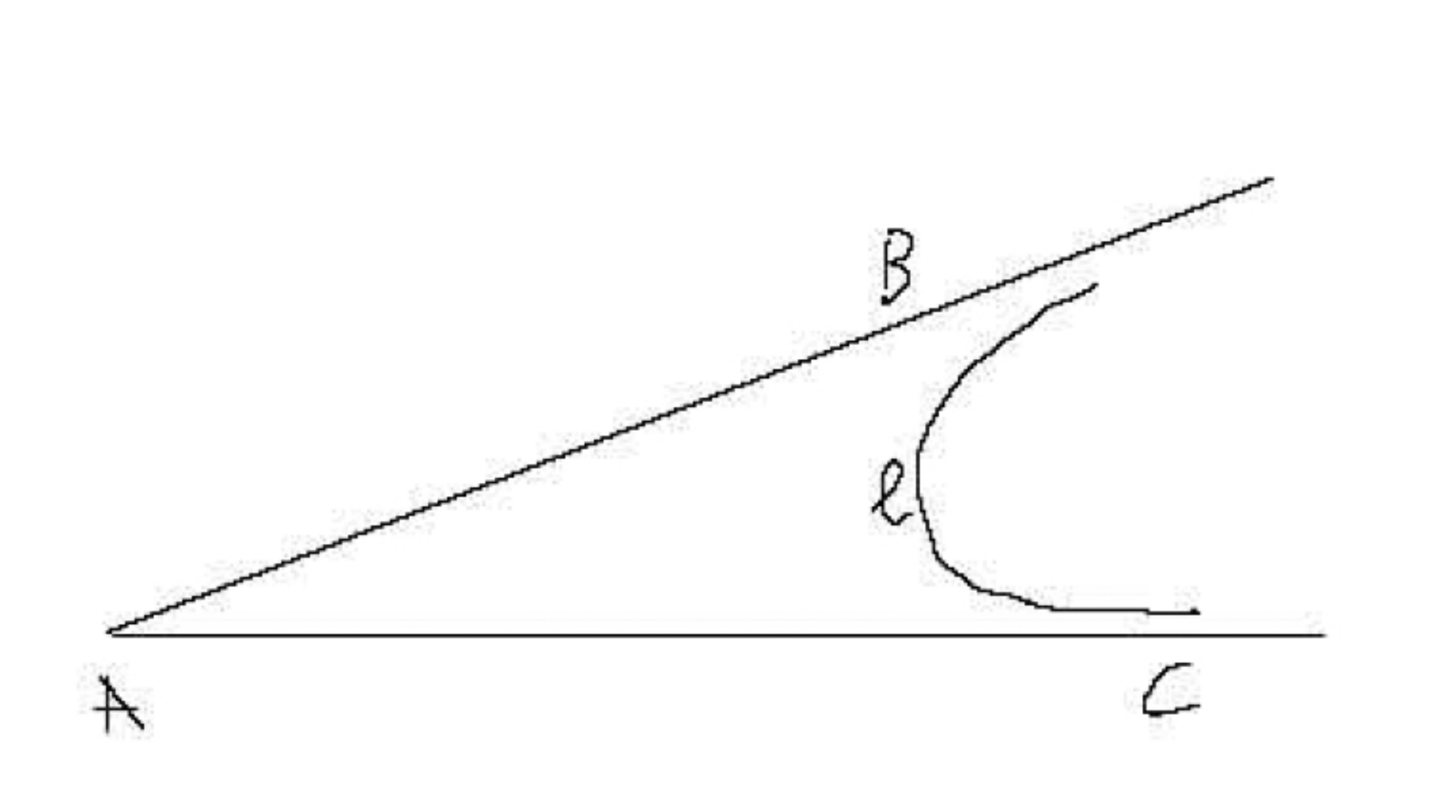} 
\end{center} 
\begin{center}
Fig. 7.5
\end{center}

It is appropriate to notice that the above elementary intuitive considerations were known to people working on the Problem of Parallels long before Lobachevsky (see  \cite{Bonola:1955}). So the hyperbolic intuition has been actually developed (in the minds of very few intelligent people) long before the official birth of Non-Euclidean geometry. Lobachevsky first managed to conceptualize these intuitions (without leaving them out!) in a way that allowed him to develop on this basis a sound mathematical theory. I provide some further details in \textbf{8.8} below. 

Let me now briefly mention another important intuitive change occurred in the 19th century geometry, which is related to  the new concept of geometrical space put forward in 1854 by Riemann \cite{Riemann:1867} on the basis of earlier geometrical works by Gauss \cite{Gauss:1828};  see \cite{Bonola:1955} for historical details. Today this concept of space is known as \emph{Riemanian manifold}; besides its importance in the pure mathematics it plays an important role in the contemporary physics by providing the mathematical foundation for the theory of General Relativity, which remains to the date the best available theory of physical spacetime. The concept of Riemanean manifold needs a special discussion that I cannot include here; although without this concept the corresponding intuition cannot be properly communicated I nevertheless provide below its very rough description.        
 
One may get a sense of the Riemanian intuition through the experience which once led Gauss to his geometrical discoveries: move around in a hilly environment and try to forget about the third dimension (which is not particularly difficult when the explored territory is sufficiently large). If you need to reach a destination, which is only few feet away, and ask yourself which way is the shortest, the answer is obvious: this is the old good Euclidean straight line. But if the destination is few miles away you need to think hard to find the shortest path and it may turn out that the optimal solution is not unique. Thus one needs here a new kind of spatial intuition, which does not cancel the Euclidean intuition but qualifies it as \emph{local} and extends it with new \emph{global} features.

Even if the notion of Riemann manifold is dispensable when one needs to map the hilly surroundings of Hanover (this was Gauss' practical charge that helped him to make his great theoretical discoveries) it is certainly indispensable for mapping the surroundings of our planet at the cosmological scales. So there is no point in saying that Riemann's generalized notion of geometrical space is ``more abstract'' than the traditional notion of Euclidean space. People knew something about curve spaces long before Riemann and Einstein! What they didn't know was how to treat such spaces mathematically: their spatial intuition remained blind in this respect. Thanks to Gauss, Lobachevsky, Riemann and other great mathematicians of the past this intuition now became much sharper and helps us to orientate ourselves in the physical spacetime at large scales where the Euclidean intuition becomes manifestly insufficient (albeit not wholly irrelevant).

\section{Lost Ideals}
The above examples of intuitive change are examples of sharpening pre-mathematical intuitions and transformation of earlier developed mathematical intuitions with newly invented geometrical concepts. Now I am going to provide a historical example of different sort when a mathematical concept was put in place before an appropriate intuition became available. As we shall see this particular example has a general significance because it demonstrates a modern way of intuiting mathematical objects, which applies very broadly. 

Given a commutative ring 
\footnote{
A commutative ring is an algebraic structure which has two operations, which are subject to certain axioms like in the case of algebraic group. A canonical example of commutative ring is the ring of integers. 
}
$(K, \oplus, \otimes)$ an \emph{ideal} (in the sense of modern algebra) is a subset $I$ of $K$ such that $I$ is (i) closed under $\oplus$ and (ii) closed under $\otimes$ in a stronger sense than usual: for all $a$ from $K$ and all $b$ from $I$  $a\otimes b$ is in $I$.  Learning this standard definition a student may wonder where the funny name ``ideal'' comes from. The definition itself gives no clue to it but historical works provide a clear answer (see \cite{Edwards:1980}). The term stems from the notion of \emph{ideal number} (\emph{Idealzahl}) introduced by Kummer in 1847. Kummer observed that in some rings of cyclotomic integers 
\footnote{A ring of cyclotomic integers is the ring of elements $m + \sqrt[p]{1}$ where $m$ is a usual integer and $\sqrt[p]{1}$ is a complex $p$-th root of unity.}  
usual integers extended with  $p$th root of 1 where $p$ is a prime number. 
the unique factorisation into primes fails. Kummer stipulated, i.e., ``put by hand'', numbers of new type that provided the unique prime factorisation in such problematic cases. He called these new numbers ``ideal'' because he did not have any justification of the \emph{existence} of such things. This was not something quite unusual in mathematics: ideal points in geometry, imaginary roots in algebra and even the familiar negative numbers were introduced into mathematics in a similar way. 
It goes without saying that Kummer's work did not reduce to this bold stipulation: assuming that the ideal numbers provide for the unique factorization he proved a number of non-trivial theorems about these stipulated things.   

It is important to keep in mind that the intended meaning of ``ideal'' depends on what in a given context counts as ``real''. Although our present discussion concerns only mathematical contexts and does not touch upon general philosophical notions of real and ideal, there are important differences between various mathematical contexts too.  Remind from \textbf{2.4} that for Hilbert  only strings of symbols are real, while all other mathematical objects are ideal. According to another popular view, real numbers are real but complex numbers are not. For Kummer all complex numbers are real but his newly invented ``ideal'' numbers are not. (In the last paragraph I deliberately confuse the technical and the non-technical senses of term ``real number'' but I hope that the reader can easily distinguish them properly.)  

The relativity of the real/ideal distinction may suggest that it is nothing but a matter of subjective attitude. However we have just seen that what is really at stake here is the question about the existence of mathematical objects of certain sort. And this latter question has a bearing on the question of legitimacy of mathematical reasoning that involves such objects. In the given case the question concerns the legitimacy of Kummer's theory of ideal numbers. Is this theory sound or not? By saying that this is a matter of subjective attitude one loses any epistemic criterion in mathematics whatsoever. 

Talking about the existence of mathematical objects we must also restrict the context appropriately. Mathematicians may disagree about  metaphysical views concerning the nature of mathematical objects (some may believe that such objects sit at the Platonic heaven, some that they are fictions, etc.) and still agree about which mathematical theory is sound and which is not. I leave now the metaphysics wholly aside and focus on the epistemology. In order to make this clear I change the terminology: instead of asking ``do ideal numbers really exist or not?'' I shall ask ``are ideal numbers well-formed or not?''. I assume that if a given mathematical theory is sound then all its objects are well-formed. (I use now the word ``sound'' not in the special logical but in the general epistemological sense - instead of saying ``sound'' I could say ``good'' or ``acceptable''.)  

We already know Hilbert's way of answering such questions: formalize the given theory and check its formal consistency. If Kummer's theory of ideal numbers is consistent then it is sound and Kummer's numbers exist  or are well-formed if one prefers.  We also know that the check of consistency is generally more problematic than Hilbert believed when he designed this method. But in any event my present purpose is not to apply Hilbert's approach (or some of its modernized versions) to Kummer's theory but rather look at this theory and its further development from a Kantian perspective paying attention both to its conceptual content and to its intuitive content. My talk of ``forming objects'' is, of course, of Kantian origin. 

In Kantian terms mathematical ideal objects  (aka ``ideal elements'') can be described as ``thoughts'', which are not wholly ``empty'' (for otherwise they wouldn't qualify as mathematical objects at all) but which are only poorly supported by intuitions.  
Like all modern mathematical objects the ideal mathematical objects are supported by some \emph{symbolic} intuition (\textbf{2.4}), which comes into the play as soon as one uses symbols for denoting such objects. But unlike mathematical objects of other sorts the ideal objects lack any other intuitive support. For example the concept of number 3 is supported both (i) by the intuition associated with symbol ``3'' and (ii) by a complex of intuitions associated with counting. In the case of ideal objects such a further intuitive support is missing. 

Now I am going to argue that such a situation can be improved: an ``ideal'' object first introduced purely formally may acquire an additional intuitive support and thus, to put it metaphorically,  \emph{get a body}. A classical example of this sort is that of the so-called ``imaginary'' number $i = sqrt{-1}$, i.e., the number satisfying $i^{2} = -1$) This number has been first introduced in the 16th century as a formal device for solving cubic equations (by ``solving equation'' I mean here finding its real roots). In the 19th century Gauss and other people developed a geometrical interpretation of $i$ and other complex numbers 
representing these numbers by points on the Euclidean plane \cite{Nahin:1998}. Ever since no research mathematician (leaving philosophers aside) doubts that the imaginary numbers are just as real as the real numbers (Kummer's example is representative in this sense);  today the very term ``imaginary number'' is no longer much in use and  in any event it wholly lost its original appeal. 
\footnote{
Since complex numbers ``officially'' don't count as geometrical objects one may argue that complex number $c$ and its geometrical representation with point $p$ are two different things. Indeed in this case the concept and the intuition are more detached than, say, in the case of Euclidean triangle. We shall shortly see that in the case of algebraic ideals the appropriate intuition and the concepts fit each other more closely.}. 

Let us now see how a similar improvement has been made with Kummer's ideal numbers. In order to understand how ideals lost their ideal character and became full-fledged mathematical object it is useful once again to look into the history of this notion. This notion was given its modern shape by Dedekind who also introduced the term (the term ``ideal'' used as a noun), see \cite{Edwards:1980} and the historical references therein. Dedekind definition of ideal differs from ours in only one respect: instead of an abstract ring $K$ Dedekind considers cyclotomic rings and other similar rings based on complex numbers. Then he shows that Kummer's ideal numbers can be identified with certain infinite classes of complex numbers, which fall under Dedekind's general definition.  
In comparison with our modern notion of ideal Dedekind's strikes as too special and too ``concrete''. One should not however forget that considering infinite collections of mathematical objects (like complex numbers) as mathematical objects in their own rights is a relatively recent idea. In Dedekind's times it was very new and even revolutionary. Anyway there is a sense in which Dedekind ``constructs'' algebraic ideals from complex numbers in a way analogous to which Euclid constructs triangles from straight segments. Admittedly this sense of ``construction'' is overtly non-constructive (in any of usual technical senses of being constructive). Yet the term ``construction'' is, in my view, wholly appropriate for describing a mental act by which one collects some given objects into an infinite collection and then operates with this infinite collection as a new object (compare \textbf{3.2} above). If elements of the infinite collection are full-fledged mathematical object provided with a non-symbolic intuitive support (as in the case of complex numbers) such a (liberally understood) construction brings about a new full-fledged object, which is equally supported by the intuition. This intuition about infinite collections allows us to think of Dedekind's ideals as full-fledged mathematical objects rather than purely formal entities like Kummer's ideal numbers.

In order to get from Dedekind's notion of algebraic ideal to the modern notion defined in the beginning of this Section one more step is needed: to replace infinite collections of  complex numbers by infinite collections of abstract elements  (or ``pure units'' by Cantor's word) and define abstract algebraic operations correspondingly. This step can be justly described as a generalizing \empty{abstraction} because the modern notion of ideal generalizes on various ``concrete'' examples of ideals including Dedekind's examples. For this reason the modern algebra is sometimes referred to as \emph{abstract} algebra \footnote{This abstract style of doing algebra follows the pattern of Waerden's course \cite{Waerden:1930-31} first published in 1930 (vol. 1)  - 1931 (vol.2). For the history of modern algebra see \cite{Corry:2004}.}. 
I suggest, however, that this move towards a higher abstraction in mathematics is accompanied by an opposite development, which involves a new modification of mathematical intuition; this new intuition allows for thinking about abstract algebraic structures \emph{in concreto} without restricting the generality and without using familiar examples. Indeed an algebraic ideal in the sense of modern definition can be called abstract \emph{only} against the background of a variety of concrete examples like Dedekind's. Although the knowledge of such examples helps one to appreciate the power of the modern concept of ideal, this knowledge is not absolutely necessary for its understanding. The modern structural mathematics (see Chapter \textbf{8} below) uses abstract sets in a way similar to which Euclid uses points, straight lines and circles, namely as a basic material for building all further constructions (in the liberal sense of ``construction'' explained above). This develops a specific set-theoretic intuition which represents elements of abstract sets as primitive ``dots'' and interprets abstract operations in terms of manipulations with dots.  Since the relevant sets are, generally, infinite the everyday experience with small finite collections often becomes misleading in this context. However this common experience serves as a basis for a further refinement and further modification of the set-theoretic intuition; such a refinement makes part of the modern mathematical education. The fact that the resulting intuition turns to be problematic in many respects is not a reason for denying its very existence. It is rather a reason for building new mathematical concepts and developing new mathematical intuitions, which are more satisfactory.

The set-theoretic intuition developed in the 20th century mathematics stems from Cantor's ``naive'' set theory. In order to analyze the intuitive aspect of this theory it is useful to apply Russell's distinction between the ``extensional genesis'' and the  ``intensional genesis'' of (infinite) classes (see \textbf{5.8} above). Remind that the former is a form of infinite enumeration and the latter amounts to using the \emph{comprehension principle} according to which for each property $P(x)$ there exist (= it is possible to generate) class $X$ of all $x$ such that $P(x)$. I shall consider some intuitive aspects of both these methods in turn.  

In his \emph{Foundations of a general theory of manifolds} of 1883 \cite{Cantor:1883a} Cantor purports to extend the ordinary natural numbers series with new ``transfinite'' numbers. Cantor gives the name of \emph{first principle of generation} to the usual operation (+1) which allows one to build every finite number $n$ starting with zero. According to the\emph{second principle of generation} the infinite iteration of the first principle has a \emph{limit}. When one applies the first principle iteratively starting with zero the corresponding limit $\omega$ can be described as the smallest number, which is bigger than every natural $n$. Applying again the first principle to $\omega$ one gets $\omega+1$, $\omega+2$, and so on; then applying again the second principle one gets $2\omega$, then again by the first principle $2\omega+1$, $2\omega+2$, .., then again combining the two principles one gets $3\omega$, .. , $\omega ^2$,  $\omega^{\omega}$, etc.. 

As we see only the last mentioned step brings one beyond the domain of countable collections. It is often argued that this latter step is crucial: while the countable infinity can be accepted as a reasonable idealization of finitary enumeration the domain of uncountable cannot be grasped with the extensional genesis (i.e., with an idealized enumeration) in principle.  Cantor suggests a different analysis. The crucial step is where the \emph{second principle of generation} first comes into the play and brings about the limit value $\omega$. The usual intuition about counting does not support the idea of such a limit. However a relevant \emph{geometrical} setting makes the idea more clear. Consider the series of fractions $\frac{1}{n}$ (where $n$ is a natural number) against the background of the \emph{real line} shown at the below diagram Fig. 7.6

\begin{center}
\includegraphics[scale=0.3]{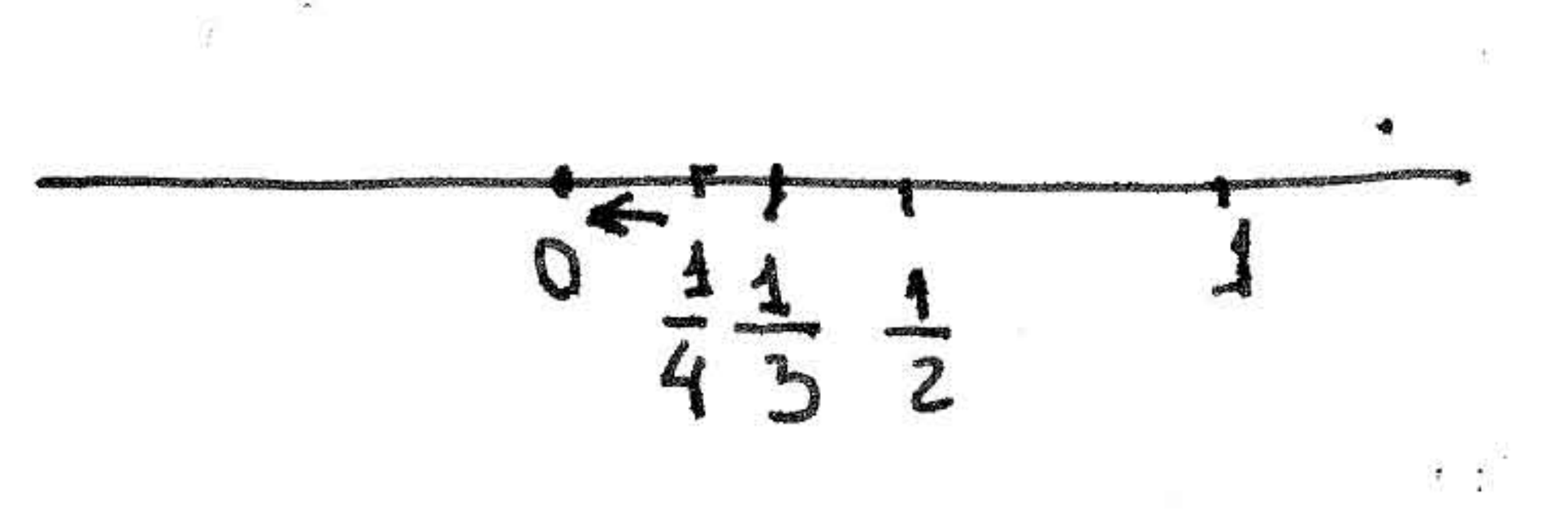}
\end{center}

\begin{center}
Fig. 7.6
\end{center}

It provides a clear picture of how this series tends to the limit 0.  The proposition saying that $\frac{1}{n}$ tends to zero can be strictly formulated in terms of the linear order relation between fractions extended with 0: 0 is the biggest number, which is smaller than each number $\frac{1}{n}$. In the case of the natural series $n$ we have a similar linear order: the order of natural numbers $n$ is obtained from the order of fractions $\frac{1}{n}$ by the simple reversal. As long as the series of fractions is taken together with its limit 0 the reversal of order translates 0 into Cantor's $\omega$. Having in mind this reversal one may use the above diagram or its like for representing $\omega$. Although this geometrical representation does not directly apply beyond $\omega$ (beware that Cantor's theory does not allow for infinitesimals!) it provides an intuitive support for the \emph{second principle}, which can be applied iteratively. 

Turning now to the intensional genesis I propose the following question: Is there anything similar to the idea of building mathematical objects using the comprehension principle in the pre-Cantorian mathematics? After all it is very improbable that this way of building mathematical objects, which today is used widely (with appropriate restrictions of the comprehension principle) would not have any predecessor in the older mathematics. At least one such example immediately suggests itself: it is building geometrical objects as {loci}. A simple example of geometrical \emph{locus} is given by Euclid's definition of circle (Definition 15 of Book 1 of the \emph{Elements}) according to which

\begin{quote}     
A circle is a plane figure contained by a single line [which is called a circumference], (such that) all of the straight-lines radiating towards [the circumference] from a single point lying inside the figure are equal to one another. \cite{Euclid:2011}
\end{quote} 

Reading this definition as a definition of geometrical locus amounts to the following paraphrase \footnote{The concept of geometrical locus was certainly already known to Euclid.  The term ``locus'' is the Latin translation of Greek word ``topos'', which was systematically used in Greek geometry as a technical term with the same meaning. For the history of the notion and the term see \cite{Allman:1889}.}: 

\begin{quote} 
A circle is a line $C$ such that for all straight-lines $OA$ and $OB$ where $O$ is a given fixed point and $A,B$ are two points of $C$, we have $OA = OB$
\end{quote} 
 
Today we can further paraphrase Euclid's definition as follows:
\begin{quote} 
A circle is a \emph{set of points} equidistant from a given point
\end{quote} 

This latter definition would hardly make sense to Euclid since for him a circle is first of all a \emph{line} (with such-and-such properties) and he never suggests the idea of equating lines with their points. Anyway this modern rendering of the ancient notion of geometrical locus suggests a convenient way of thinking about infinite sets formed with the comprehension principle in terms of familiar geometrical objects. So in the case of intentional genesis a properly modified traditional geometrical intuition is also relevant. A further modification of this intuition may play a role in ``abstract'' contexts
\footnote{
\begin{quote}
Some years ago I began an introductory course on Set Theory by attempting to explain the invariant content of the category of sets [..]. I was concerned to present an ideological vision of the significance of the objects of this category, which I called \emph{abstract sets}. I emphasized that an abstract set may be conceived of as a bag of dots which are devoid of properties apart from mutual distinctness. (\cite{Lawvere:1994a}, p. 5)
\end{quote} 
}
, for example, when one defines an algebraic ideal in the modern structural way explained above. 

It is interesting to notice that Cantor's general theory of sets emerged as a generalization of a theory of point-sets, on which Cantor worked earlier.  There are also other important historical connections between Cantor's set theory and the 19th century geometry \cite{Ferreiros:1999}, \cite{Lawvere:1994a},  which I cannot discuss here. These historical facts suggest that Cantor smuggled into his general notion of set,  which he famously described as ``a collection into whole of definite, distinct objects of our intuition or our thought'' \cite{Cantor:1885}, quite a lot of geometrical content. It remains an open question whether Cantor's set theory can be reasonably reconstructed on some appropriate geometrical basis. Voevodsky's suggestion to treat sets as a particular homotopy type (\textbf{6.10}) is a possible way of solving this problem. Without going into details I would like only stress here a foundational significance of this endeavor. In the standard set-theoretic foundations of mathematics sets are primitive and geometrical spaces are built as structured sets through the Formal Axiomatic Method (\textbf{3.2}). Making sets into a special \emph{geometrical} concept reverses this conceptual order and makes geometrical spaces more basic. In the next Section I explain why this conceptual reversal is, in my view, desirable.

\section{Are Intuitions Fundamental?}
I hope having convinced the reader by the above examples that mathematical intuitions change at least as fast as mathematical concepts. Another recent intuitive change of the same sort is connected to the category theory: I have already discussed the intuitive geometrical aspect of this theory earlier in this book (see particularly \textbf{6.2 - 6.3}) and I shall discuss it again in the next Chapter (\textbf{8.8}). Let me now try to answer a more general question: Is the mathematical intuition fundamental (as Kant believed) or it is only an auxiliary instrument, which helps mathematicians to do their science? I claim that it \emph{is} fundamental. The core argument goes as follows. I assume after Cassirer (i) that  ``[t]he principle according to which our concepts should be sourced in intuitions means that they should be sourced in the mathematical physics'' and (ii) that ``[l]ogical and mathematical concepts must no longer produce instruments for building a metaphysical 'world of thought' '' (\cite{Cassirer:1907}, p. 43-44). I also assume that this sourcing of mathematics from empirical sciences is a fundamental issue. The conclusion follows.

It is appropriate to ask what the examples of mathematical intuition treated in this Chapter have to do with the empirical science. As long as Euclid and Barrow are concerned the answer is obvious: both geometrical theories were closely related to the contemporary empirical research. In Euclid's times, when the mathematical physics in the modern sense of the word did not yet exist, the relevant empirical research was mostly limited to astronomy, geography and optics; the new geometry of curves pioneered by Barrow (among other people) provided a geometrical basis for Newton's physics. The case of non-Euclidean geometries may look less convincing in this respect - and the fact that a non-Euclidean geometry provided a basis for Einstein's theory of relativity may look just as a stroke of luck. However in fact Lobachevsky from the very beginning of his geometrical studies tentatively thought about his newly invented non-Euclidean spaces as physical spaces. He does not take for granted the idea according to which the geometry of physical space is unique, and does not think about the spectrum of new geometries as a spectrum of merely mathematical possibilities.  See, for example, the following passage that reads today as a precipitation of Einstein's General Relativity and modern field theories:

\begin{quote}
[T]he assumption according to which some natural forces follow one Geometry while some other forces follow some other specific Geometry, which is their proper Geometry, cannot bring any contradiction into our mind. (\cite{Lobachevsky:1949}, p. ??, my translation)
\end{quote}

This shows that Lobachevsky's ``hyperbolic intuition'' qualifies as intuition in Kant's and Cassirer's sense: it is not a merely subjective play of imagination but a suggestive way of thinking about ``natural forces''. When this aspect of the history of the 19th century geometry is taken into consideration the success of the new geometry in the 20th century physics no longer looks as a miracle \cite{Gray:1999}.   

In the case of the \emph{set-theoretic} intuition the situation is more complicated. It is often argued that since there is no infinite sets in the nature the set theory (except its finitely fragment) is irrelevant to natural sciences. It follows that set-theoretic intuitions (if any) cannot be relevant in natural sciences either. In my view this argument points to a real problem, but in order to better illuminate the problem this argument must be formulated differently. What does it mean that there is no infinite sets in nature, how we can know this? And how this fact makes the set theory irrelevant in natural sciences? After all there is no perfect Euclidean points in the nature either: at best a point can only \emph{represent} a physical particle in appropriate empirical contexts. Why an infinite set cannot do a similar job? 

The real problem, as I see it, is that in the 20th century the set theory developed wholly independently from empirical sciences. Interestingly, Cantor himself thought of possible empirical applications of his theory \cite{Grattan-Guinness:1970}. However since Zermelo in 1908 \cite{Zermelo:1908} turned the research in set theory into a formal axiomatic research (\textbf{3.1})  Cantor's approach was disqualified as ``naive'' and no longer developed (until in 1994 Lawvere \cite{Lawvere:1994a} took Cantor's original approach seriously and tried to develop it in a novel way). Many mathematicians have been convinced by the argument according to which the pure logic and mathematics somehow extends beyond any possible experience into the putative universe of infinite sets; the trick of Hilbert-style formalization offers itself as a means allowing to reach this ``metaphysical world of thought'' (Cassirer) aka ``Cantorian Paradise'' (Hilbert) by finitary means. In this context the symbolic intuition, which allows for manipulation with finite strings of symbols according certain rules, is supposed to be sufficient for treating infinities mathematically (as ideal objects in Hilbert's sense). As I have already argued in Chapter \textbf{3} this letter proposal has  never been popular in the mainstream mathematics. Outside the formal axiomatic studies of sets mathematicians still think about sets ``naively'' and support this thinking by non-symbolic intuitions, which I have described in the last Section. Although these intuitions hardly support an empirical application of set-based mathematics directly, they help to connect the set-based mathematics to the older mathematics, which is still widely used in physics and other natural sciences (think of differential equations and the theory of dynamical systems, for example).  

The question about the role of mathematical intuitions in empirical sciences needs a separate discussion that I cannot provide here. Evidently the standard Kantian understanding of this role does not immediately apply to Quantum Mechanics and other parts of the contemporary physics where empirical observations and measurements are very indirect. However taking it for granted that the empirical science provides the only epistemic access to reality, I claim after Kant and Cassirer that all sound mathematics is sourced from this type of science. One may further ask whether the pure mathematics is ultimately itself an empirical science or it is rather a non-empirical (a priori) ingredient of empirical sciences. I leave now this question aside because it has no bearing on my principle claim according to which the pure mathematics cannot be founded independently for the mathematized empirical science.

Cassirer calls the link between the pure mathematics and the empirical science \emph{intuition} without specifying what it is and how it works in the modern context. One may argue that since the role of direct observations by the naked eye no longer plays a fundamental role in science, the symbolic intuition is after all the best possible interface between the pure mathematics and the empirical science.  In order to evaluate this argument we need first of all to distinguish between two different claims: one according to which (i) the symbolic intuition is the basic sort of intuition in the contemporary mathematics, and the other one according to which (ii) the symbolic intuition \emph{as a part the Formal Axiomatic Method} plays (or should play) a fundamental role in the contemporary mathematics\footnote{(ii) is Hilbert's mature view, see \textbf{2.4, 2.6} above.}. (ii) implies (i) but not conversely. As I have already stressed in \textbf{3.3} the way in which symbols are used in theories built by the Formal Axiomatic Method is quite unlike the way in which symbols are used elsewhere in mathematics; I have also argued in \textbf{3.3} that this particular way of using symbols makes the formal mathematics not apt for being applied in the modern empirical science. Since after Cassirer I want to understand the mathematical intuition as a link between the pure mathematics and the empirical science I rule (ii) out. But this argument has no bearing on (i), which remains plausible. 

As for (i) I refer the reader to the earlier discussion on the notion of symbolic intuition in \textbf{2.4 - 2.5}. As long as by the symbolic syntax one understands the standard strings of symbols, one may remark that this syntax is very good for some purposes and not so good for some other purposes. In any event we cannot grant to this type of syntax - and hence to the symbolic intuition, which allows one to manipulate with symbols using this syntax - any exceptional epistemic role. So we should treat the symbolic intuition in this restricted sense of the term on equal footing with mathematical intuitions of other types including, for example, the geometrical intuition behind the homotopy theory (\textbf{6.7}).

If, however, one understands the notion of symbol in a more general sense, namely, as any representation, which does not involve the idea of resemblance between a given symbol and what this symbol stands for, then (i) is indeed justified. However in this case the notion of symbolic intuition becomes unspecific. For Kant, saying that the traditional geometry provides a representation of the physical space that \emph{resembles} this physical space, would not make sense because \emph{things in themselves} and their representations do not admit such a relation; one may rather talk here about resemblances between different representations. Saying that the traditional geometrical diagrams resemble the geometrical objects, which those diagrams represent, does not make sense either, because for Kant the geometrical objects \emph{are} representations.  Thus given the above general notion of symbol, all mathematical intuition turns to be symbolic (which is Cassirer's view) and (i) turns into a general claim about the symbolic character of mathematical intuition.

The relevant difference between older and newer types of empirical representation is this. A traditional diagram representing the Solar system resembles the Solar system. A mechanical model of Solar system (I mean a real mechanism used for educational or other purposes) also resembles the Solar system.  Now I am talking about the resemblance between two material objects (one of which is natural and the other artificial), not about the Kantian issue of mental representation. Clearly, no reasonable model of atom resembles the actual atom in a similar way. Remark, however, that this relation of resemblance plays no role in Kant's theory of scientific representation. Hence the Kantian approach to scientific representation (albeit not all specific details of Kant's theory) applies to the Quantum physics after the model of Newtonian  physics. Thus Kant already showed us how to do the modern science without using the naive notion of resemblance between physical objects and their images. He did this for the Newtonian physics but his general approach applies to today's physics as well. Moreover, in today's physics this approach becomes even more pertinent because the recent physics makes the older notion of resemblance between physical objects and their images manifestly absurd.

To repeat, I do not want to say that Kant's original approach is fully adequate to today's science. It is obviously not. But I do want to suggest that any adequate approach to scientific representation should build on Kant's in one way or another.  
True, many basic intuitions playing the central role in the Newtonian physics (like the intuition of mechanical motion stressed by Barrow) turn to be irrelevant in the new context, or at least they need to be conceptualized very differently. However these differences are differences between different intuitions and their corresponding concepts rather than differences between different types of relationships between concepts and intuition. This is why I think that Kant's general notion of intuition is suitable for today's physics as well as to today's mathematics. I recognize that Hegel was right when he criticized Kant's ``subjective bent'' even if I don't quite satisfied by his proposed solution of this problem. In my view this problem can and should be addressed rather in a naturalistic framework providing an objective scientific representation of human subjectivity as a part of an objective representation of human natural environment.  It goes without saying that this topic needs a more specific discussion, which I cannot provide in this book. 

To conclude this Chapter I would like to suggest that the argument according to which all intuitions about infinite collections are doomed to be irrelevant in the empirical sciences, is wrong. So I do not defend here some form of \emph{finitisme} according to which infinities must be banned in mathematics because they cannot make any empirical sense. I believe after Cantor \cite{Grattan-Guinness:1970} that they really can. But this requires a different set theory; Lawvere's ETCS and Vovodsky's homotopical theory of sets give an idea how such a new theory of set may look like. The traditional Scholastic distinction between the potential and the actual infinity used by Cantor in \cite{Cantor:1883a} must not be taken for granted and used then as the philosophical basis for mathematical theories of the infinite. Relevant empirical consideration \cite{Rashevskii:1973} are both more illuminating and more suggestive in this respect. The metaphysical idea of mathematics penetrating into the non-empirical domain of the infinite must be definitely abandoned as non-scientific. 

\chapter{Categories versus Structures}

In this Chapter I shall discuss a notion, which have been already used throughout this book, namely, the notion of \emph{mathematical structure}. I shall also discuss a philosophical view known as \emph{structuralism} and analyze its relationships with the category theory. According to a popular opinion the category theory wholly justifies the structural approach in mathematics and provides a framework for developing the structural mathematics. I shall argue that the relationships between the category theory and the mathematical structuralism are more involved: although the category theory helps the mathematical structuralism to win over some of its traditional rivals (including all varieties of mathematical \emph{substantialism}) it also transforms the traditional structuralism into something very different. Although many people in this context prefer to talk about a new form of  structuralism \cite{Makkai:1998}, \cite{Hellman:2001}, \cite{Hellman:2003}, \cite{Awodey:2004}, \cite{McLarty:2004}
\footnote{A structuralist view on category theory also underlies the recent monograph \cite{Marquis:2009} even if the author does not discuss structuralism explicitly in this book.}
 I find this terminological choice rather unfortunate because it only points to the continuity of the conceptual change but does not help one to describe this change itself. I wholly recognize the fact that the category theory grew up on the basis of the structuralist mathematics, and that structuralist motivations played a major role in the historical development of this theory.  Nevertheless I claim that the category theory supports a new vision of mathematics and science, which is incompatible with the traditional structuralism. Of course, this claim is defensible only if the structuralist view is sharpened accordingly. Some people may argue that my interpretation of structuralism is too restrictive. Anyway in this Chapter I shall try to persuade the reader  that there is a sharp conceptual difference between the structuralist view and the new view on mathematics emerging with the category theory.

\section{Structuralism, Mathematical}
In its most popular usage the name of structuralism refers to a broad intellectual movement in France during the time period from 1945 to late 1960-ies or early seventies, which was very influential not only in philosophy but also in mathematics, linguistics, sociology, economics, political science and other humanities \cite{Dosse:1997}. Bourbaki's structuralist manifesto \cite{Bourbaki:1950} was first published in French during the same year 1948 when Levi-Strauss defended his thesis \emph{The Elementary Structures of Kinship} \cite{Levi-Strauss:1969}, which soon afterwards became a canon of the structural approach in humanities. Today one can reasonably ask whether the notion of mathematical structure as formulated by Bourbaki and the notion of structure used by Levi-Strauss was indeed the same notion \cite{Rabouin:2011}. Although the answer to this question is far from being obvious it is a well-documented historical fact that the structural approach in mathematics, on the hand, and the structural approach in humanities, on the other hand, closely and fruitfully interacted (in particular, through the application of apply new mathematical methods in linguistics and other humanities). This, in my view,  is a sufficient reason to count the mathematical structuralism of Bourbaki as a proper part of the French post-war structuralist movement.

I insist on this point for the following reason. Just about the time when according to the official history  \cite{Dosse:1997} the French structuralism declines (the late 1960-ies) and gets succeeded by a variety of new intellectual trends, which North-American academic experts classify under the title of \emph{post-structuralism}, the structuralism has a new birth in North America. This new structuralism started by Sneed \cite{Sneed:1971} and Suppe \cite{Suppe:1977}  limits its scope by the philosophy of science and forcefully denies any relationships with French structuralism with the noticeable exception of Bourbaki's mathematical structuralism. The program of this new structuralism can be roughly described as application of Bourbaki's method for a philosophical reconstruction of scientific theories; the purpose of such reconstruction is the \emph{structural analysis} of theories, which is a special sort of formal analysis, see  \cite{Balzer&Moulines&Sneed:1987} for further details. Independently,  yet another version of North-American structuralism establishes itself in 1990-ies as an influential philosophy of mathematics; this latter structuralism is known under the name of \emph{mathematical structuralism}, see Hellman's \cite{Hellman:2006a} for a definition and the references. Mathematical structuralism equally recognizes Bourbaki's structuralism but not the rest of French post-war structuralism. In this context the historical connection between Bourbaki's mathematical structuralism and the rest of French structuralism needs to be mentioned specially. 

To complete my short list of structuralisms I would like to mention the \emph{structural realism} emerged in UK in early 2000ies  \cite{Ladyman:2009}. This kind of ontological structuralism in the philosophy of science does not have any direct connection to Sneed's structuralism but does have many direct connections to the mathematical structuralism in Hellman's sense.

Although I don't discuss in this book French structuralism (except Bourbaki's mathematical structuralism) I borrow something important from its spirit and from its general intellectual attitude. The story of decline of structuralism in one part of the world followed by the reappearance of (a different) structuralism in another part of this world tells us something important about different ways of doing philosophy. French post-war structuralism functioned as a suggestive idea, which prompted significant developments in mathematics and humanities (I don't know about sciences); it never construed itself as a philosophical doctrine supposed to be learnt and transmitted through generations in its original form. It's early decline in the late 1960-ies could be expected because French structuralism never construed itself as long-living. It had something of an intellectual fashion - and for this a reason its  esteem in certain academic circles was very poor. The fact that in the late 1960-ies structuralism in France was succeeded by new influential trends (the so-called post-structuralism) is an evidence that the decline of structuralism was not a symptom of intellectual stagnation. 

The North-American structuralism in any of its multiple varieties, on the contrary, construed itself in the form of philosophical doctrine to begin with. This is why it is much easier to describe this later structuralism as a coherent system of beliefs. According to  Hellman the contemporary mathematical structuralism is    
 
\begin{quote}[..] a view about the subject matter of mathematics according to which what matters are structural relationships in abstraction from the intrinsic nature of the related objects. Mathematics is seen as the free exploration of structural possibilities, primarily through creative concept formation, postulation, and deduction. The items making up any particular system exemplifying the structure in question are of no importance; all that matters is that they satisfy certain general conditions - typically spelled out in axioms defining the structure or structures of interest - characteristic of the branch of mathematics in question. (\cite{Hellman:2006a})
\end{quote} 

(Beware that the above definition describes only the core mathematical structuralism, which branches then into multiple varieties, which I leave aside.) 

Mathematical structuralism in the above sense generalizes upon well-established patterns of the 20th century mathematics including Bourbaki's \emph{Elements}; in what follows I shall explain in more detail its relevance to the Formal Axiomatic Method and other distinctive features of mathematics practiced during this century.  However Hellman's mathematical structuralism presents itself as a view on mathematics \emph{simpliciter} without any explicit historical reflection; it simply tells us what mathematics \emph{is},  not what it used to be in the past,  nor what it will be in the future. Does this mean that Hellman and other adherents of North-American mathematical structuralism really believe that mathematics has some sort of eternal nature, which they try to describe? Perhaps some of them do. However there is also a more sophisticated way of tackling this issue, which perhaps better fits this people's intellectual attitude:

\begin{itemize}  
\item (1) Since (1.1) today's mathematics is the best known mathematics it follows that (1.2) in order to answer the question What is mathematics? it is sufficient to take into consideration only today's mathematics. The history of mathematics is an important subject but it does not have a bearing on this question. 
\item (2) Talking about today's mathematics it is sufficient to take into consideration only well-established patterns of modern mathematical reasoning and leave aside all recent ideas and new approaches, which so far remain controversial in eyes of mathematicians themselves. In the unlikely event that in the near future mathematics once again will change its shape (like it happened in early 1900-ies) a philosopher should change her mind accordingly. But as long as a philosophy of mathematics fits the content of standard mathematical textbooks there is no reason to promote such changes. 
\end{itemize}  

Let me now express some criticisms and explain my own position with respect to French and North-American structuralisms. 

I agree with premiss (1.1) but I disagree that the conclusion (1.2) follows. The fact that mathematics rapidly develops is its essential feature, i.e., an essential part of the answer to the question What is mathematics? (\emph{Mutatis mutandis} this is true about any science.) One cannot understand what mathematics is and how it develops without studying its history and without anticipating its future. 

(2) points to a crucial difference between French and North-American structuralisms and explains why the latter structuralism emerges when the former structuralism dies off. To repeat, the French structuralism worked as a suggestive idea, which died when the relevant structuralist suggestions were already successfully used in mathematics and elsewhere, so that the idea was no longer suggestive but became somewhat commonsensical. Exactly at this point North-American philosophers could take the opportunity of putting this idea in the form of philosophical doctrine. My own project is mixed. Following the pattern of French structuralism I shall consider the mathematical structuralism as a suggestive idea, not as a doctrine. Thus I reject (2) and not only discuss some well-established patterns of modern mathematical thinking but also discuss some history and some work in progress. However I also take the advantage of using Hellman's precise definition of structuralism for making explicit a non-structuralist aspect of category theory (see \textbf{8.6}) below). Thus my following critique of structuralism does not aim at showing that the mathematical structuralism is wrong as a doctrine; it rather aims at showing that the structuralism is outdated and no longer suggestive. This critique also helps me to formulate a different idea, which in my view is both more modern and more suggestive. Although this alternative idea is motivated more by recent mathematics than by recent philosophy I don't mind if someone describes it as a form of mathematical \emph{post-structuralism}. 
  
\section{What Replaces What?}

Here is how a working mathematician describes an example of mathematical structure for a philosophical reader:

\begin{quote}All infinite cyclic groups are isomorphic, but this infinite group appears over and over again - in number theory, in ornaments, in crystallography, and in physics. Thus, the ``existence" of this group is really a many-splendored matter. An ontological analysis of things simply called ``mathematical objects" is likely to miss the real point of mathematical existence. \cite{MacLane:1996}\end{quote}

 This example well demonstrates the generalizing power of the structural approach: the concept of cyclic group
\footnote{Remind the definition of general group from \textbf{6.1}. An infinite cyclic group is a group with an infinite number of elements and such that any of its elements is generated  by some distinguished element $g$ and its inverse $g^{-1}$. A group is said to be \emph{generated} by a set of its distinguished elements (called \emph{generators}) when every element of this group is a product (in the sense of the group operation) of the generators. A canonical example of an infinite cyclic group is the additive group of whole numbers, which is generated by numbers 1 and -1.} can be instantiated by apparently very different mathematical constructions belonging to different areas of mathematics; in this way this concept establishes a conceptual link between these different areas. Mac Lane's example also illustrates the controversy between the mathematical structuralism and the rival view of set-theoretic substantivism. According to this latter view a group or any other mathematical object construed as a structured set \`a la Bourbaki (see \textbf{3.2}) is a \emph{set} at the first place. The above quote suggests that this base set does not really matter; what really matters is only the \emph{structure} supported by this (or any other) set. Before I shall try to make it more clear what is meant here by structure let me remind a passage from Burbaki's structuralist manifesto of 1948, which has been already quoted in \textbf{3.2}. After discussing logical difficulties of set theory Dieudonn\'e (under the name of Bourbaki) remarks:  

\begin{quote}The difficulties did not disappear until the notion of set itself disappears ... in the light of the recent work on the logical formalism. From this new point of view mathematical structures become, properly speaking, the only ``objects" of mathematics. (quoted by English translation  \cite{Bourbaki:1950}, p. 225, footnote)
\end{quote}

It is not clear what Dieudonn\'e exactly means by the ``disappearance" of sets, and I doubt that in 1948 he could really justify his claim. However the above quote clearly demonstrates the intention. 

Let us now see more precisely what this controversy is really about. For the sake  of my further argument I modify Mac Lane's example as follows: I replace the words ``infinite cyclic group" by the words ``number three" and the word ``isomorphic" by the word ``equal". So we get this statement

\emph{All threes are equal but this number appears over and over again - in number theory, in ornaments .... Thus the ``existence" of this number is really a many-splendored matter.}

which describes the traditional ``wobble about identity'' in mathematics, which I treated in Part \textbf{II} of this book. Thus the above modification reveals a traditional aspect of structuralism, which often remains unnoticed when people stress the novelty of this approach. Indeed the familiar number three is just as promiscuous as the infinite cyclic group or perhaps even more promiscuous. The number three equally ``appears" (to use MacLane's word) both inside and outside mathematics: in a trio of apples, a trio of points, a trio of groups, a trio of numbers or a trio of anything else. As in Mac Lane's original example, there is a systematic ambiguity between the plural and the singular forms of nouns in our talk about numbers. (Notice Mac Lane's talk about ``all infinite cyclic groups" and ``this infinite group" in the same sentence; in my paraphrase I talk similarly about a number.) This shows that the notion of ``many-splendored existence'' (i.e., of multiple instantiation) is not at all specific for the 20th century structural mathematical thinking. Thus in order to understand what is indeed specific for the structural thinking one should look elsewhere.  

Comparing Mac Lane's example with its modified version one can see that in Mac Lane's example the notion of \emph{isomorphism} plays the same role as the notion of equality (as distinguished from identity) plays in the traditional mathematics. The idea that isomorphic objects can be treated as equal is crucial for structuralism - at least if we are talking about structuralism as a trend in mathematics rather than a philosophical doctrine about mathematics. 

However we need to be again more specific here. Remind from \textbf{6.1} that the idea of ``taking isomorphism for equality'' branches into two very different ideas. The first idea can be dubbed the \emph{replacement of isomorphisms by equalities}, the second the \emph{replacement of equalities by isomrophisms}. To see clearly the difference between the two approaches one should first of all to disambiguate the term ``isomorphism'', which may denote either a relation or a map (transformation). The two meanings of the term are related as follows: the \emph{relation} of isomorphism $x \simeq y$ holds if and only if there is a \emph{map} $f: x\rightarrow y$, which is isomorphism (or in other words, which is invertible). Notice that the notion of isomorphism as a map stands first in the conceptual order: one needs this notion in order to define the further notion of isomorphism as a relation. 

The \emph{replacement of isomorphisms by equalities} works as follows. One takes the relation of equality $=$ for granted or introduces it independently from the notion of isomorphism. Then one introduces the relation of isomorphism $\simeq$ (using $=$ if needed). Finally one ``replaces'' $\simeq$ by $=$ through the Fregean abstraction (\textbf{5.5}).  Thus the replacement of $x \simeq y$ by $x = y$ brings an abstract object that one may denote $x$ or $y$ or by a new symbol $X$. The  \emph{replacement of equalities by isomrophisms} is more involved because it doesn't rely upon usual ways of thinking about equality but purports to construe this notion anew using isomorphisms. Since the notion of isomorphism as relation depends on the notion of isomorphism as a map this idea amounts to a construal of equalities (identities) in terms of maps. We have explored this idea above in Chapter \textbf{6}; remind that in the most developed form existing to the date this idea has been realized in the homotopy type theory where the relevant isomorphisms are paths and their homotopies (including all higher homotopies) in some topological space (see \textbf{6.9}). 

Which of the two approaches is more helpful for explaining what is a mathematical structure? At the first glance the  \emph{replacement of isomorphisms by equalities} is exactly what one needs for this purpose. Suppose one notices a similarity between (i) an ornament, (ii) an arithmetical construction and (iii) a piece of mathematical physics, which one wants to describe as a \emph{structural} similarity. In order to make this observation more precise one constructs mutual isomorphisms between (i), (ii), (iii); the existence of such isomorphisms is the intended precise sense of the expression ``structural similarity''.  Then one wants to
study the structural aspect of (i), (ii), (iii) independently from their other specific properties (in a way similar to which one may want to study the arithmetical aspect of trio of apples, trio of points and of any other trio disregarding all specific properties of apples, points, etc.).  For this purpose one forms an appropriate notion of \emph{structure} thinking about this structure as a thing, which (i), (ii), (iii) share in common and which makes (i), (ii), (iii) structurally similar. The \emph{replacement of isomorphisms by equalities}, which is a special case of Fregean abstraction, brings one the notion of an abstract group $G$ of the appropriate type (let it be the infinite cyclic group of Mac Lane's example).  Now one may forget about irrelevant properties of  (i), (ii), (iii) and study this group independently as an abstract algebraic structure. So one no longer cares about patterns of ornament, about numbers or about physical bodies but studies the ``structural relationships in abstraction from the intrinsic nature of the related objects. [..] The items making up any particular system exemplifying the structure in question are of no importance.'' (Hellman) So it looks like the \emph{replacement of isomorphisms by equalities} is a right way of making mathematical structures.   

Let me now explain why this way of making mathematical structures is not quite satisfactory. Suppose we have two abstract groups $G_{1}, G_{2}$ formed as above; these groups are construed as sets $S_{1}, S_{2}$ of abstract elements provided with abstract binary operations. The notion of group isomorphisms, which has been used for treating (i), (ii), (iii), makes perfect sense in the new abstract setting too. Suppose that the given groups are isomorphic, in symbols $G_{1} \simeq G_{2}$. We want now to say that groups $G_{1}, G_{2}$ are the same (or equal) but it turns out that this statement has some unwanted consequences.  Namely,  $G_{1} = G_{2}$ implies $S_{1} =  S_{2}$. The latter equality is unwanted because the idea is that different sets $S_{1}, S_{2}$ may support the same structure. So we want here a further abstraction, which would allow us to consider the group structure independently from any underlying set. This explains why Dieudonn\'e's wishes that sets disappear.  However the \emph{replacement of isomorphisms by equalities} is no longer helpful. It allows one to abstract from irrelevant features of (i), (ii), (iii) but it does not allow one to get rid of the fact that each of (i), (ii), (iii) is a collection of well-distinguishable elements. Unless one takes this fact into consideration one cannot build isomorphisms between  (i), (ii), (iii), so the \emph{replacement of isomorphisms by equalities} becomes irrelevant.  Thus the notion of abstract set, which one wants to use as a precise mathematical expression of the idea of collection of ``abstract elements having no intrinsic nature'' (Hellman), appears to be an obstacle for the structuralist way of thinking. One says that those abstract ``items [..] are of no importance'' (Hellman) but one does not know how to get rid of them altogether. Set-theoretic substantivism is a view according to which this difficulty is fundamental, so one should accept some notion of set and some fundamental identity relation before talking of abstract mathematical structures. Structuralism suggests to do something else, typically without making this suggestion more specific (like in Bourbaki's case).    
   
What about the \emph{replacement of equalities by isomrophisms}? Since this approach did not exist in a precise mathematical form until very recently, it had no significant influence on the mathematical structuralism as a philosophical doctrine. But as a general idea (not always well distinguished from the idea of \emph{replacement of isomorphisms by equalities}) it pushed these  developments in mathematics, which finally made a precise mathematical expression of this idea possible. In \textbf{8.5} I show how the \emph{replacement of equalities by isomrophisms} and, more generally, the \emph{categorification} (in the broad sense of \textbf{6.3})  indeed allow for a further development of the structuralist line of thinking. However in \textbf{8.6} I show how the same development transforms the structuralist way of thinking about mathematical (and other) matters into something different.

\section{Erlangen Program and Axiomatic Method}

Hilbert's Formal Axiomatic Method (in the sense of \emph{Foundations}Êof 1899 rather than  \emph{Foundations}Êof 1934-39) is one of the main tenets of mathematical structuralism.  Remind Hilbert's words according to which

\begin{quote}
[I]t is self-evident that every theory is merely a framework or schema of concepts together with their necessary relations to one another, and that basic elements can be construed as one pleases. [..] [E]ach and every theory can always be applied to infinitely many systems of basic elements. For one merely has to apply a univocal and invertible one-to-one transformation and stipulate that the axioms for the transformed things be correspondingly similar. (quoted by \cite{Frege:1971}, see the longer quote in \textbf{2.1})
\end{quote}
 
This Hilbert's description of a mathematical theory perfectly squares with Hellman's description of the subject-matter of mathematics. When Hellman describes mathematics as ``the free exploration of structural possibilities, primarily through creative concept formation, postulation, and deduction'' he most certainly has Hilbert's Axiomatic Method in mind. Noticeably,   Hilbert explicitly mentions in the above quote the notion of isomorphism (under the name of univocal and invertible one-to-one transformation), which is absent from Helmann's description but crucial for the structuralism. Since Hilbert talks here about isomorphisms as transformations, a (formal) theory in Hilbert's intended sense can be described as an invariant of such transformations, namely, as an invariant of all invertible transformations of each model of this theory into each other model. A geometrical origin of this idea becomes evident when one compares it with Klein's \emph{Erlangen Program} \cite{Klein:1872}, which, remind, is the idea of studying geometrical spaces through the identification of groups of automorphismes of these spaces and their invariants (\textbf{6.1})\footnote{On the history of Klein - Hilbert partnership see \cite{Rowe:1989}.} Now I would like to stress a foundational and unificational aspects of Klein's program, which I did not discuss earlier. 

The \emph{Erlangen Program} was formulated by Klein as a part of his research aiming at common foundation and unification of Euclidean and recently discovered Non-Euclidean geometries. Two major publications relevant to this research are Klein's papers \cite{Klein:1871} of 1871 and  \cite{Klein:1872} of 1873; the \emph{Program} was first published just in between, in 1872. Klein's idea (which he partly borrows from Cayley, see  \cite{Shenitzer:1991}) is to use the group-theoretic \emph{Erlangen} approach for presenting Lobachebsky's hyperbolic geometry and Euclidiean geometry (which Klein calls \emph{parabolic}) as special cases of \emph{projective} geometry. Klein, of course, did not build a projective  space as an abstract structure like we do this today.  Instead, he first describes  projective transformations and their invariants in the traditional way using Euclidean geometry (in the analytic form),
and then specified some transformations of Euclidean lengths, which did not affect the group of projective transformations but brought what we would call today \emph{models} of Non-Euclidean spaces of two sorts: of constant negative curvature (aka hyperbolic or Lobachevskian) and  of constant positive curvature (elliptic). Using this construction Klein presents the projective geometry as a foundation of a whole family of metric geometries including the Euclidean geometry  (as the case of null curvature). From the modern viewpoint Klein's approach to foundations of geometry looks peculiar or even wholly misguided because it limits itself to the case of differentiable manifolds of constant curvature and does not involve any explicit topological consideration. Leaving now aside the general question concerning the role of Klein's approach in the history of geometry, I shall talk only about a special relationship between Klein's and Hilbert's approaches to foundations of geometry. 

Klein describes the projective structure as the core geometrical structure, which unifies geometry (or at least its fragment) as described above. By the same pattern Hilbert points to the \emph{logical structure} as the core structure that unifies all theories including geometrical theories. In order to make this analogy more precise I shall distinguish between two different senses of the expression ``logical structure'', both of which are relevant in this context. First, one may talk about the logical structure of certain theory $T$. Logical structure of $T$ is made explicit by the formal version $FT$ of $T$. So a logical structure in this first sense is a formal theory (or the axiomatic structure of a given formal theory if one prefers). In Klein's spirit such formal theory $FT$ can be described as an invariant of transformation of any given model $T$ of this theory into any other model of the same theory: remind Hilbert's ``univocal and invertible one-to-one transformation'' in the above quote. Second, by the logical structure one may mean the universal logical structure of the world, which ``is invariant under all possible one-one transformations of the world onto itself'' ( \cite{Tarski:1986}, p. 149). The latter quote is Tarski's definition of logicality, see \textbf{2.2} above. Taski's ``world'' is, technically, a typed universe, which comprises individuals, classes of individuals, classes of classes of individuals, etc.; ``transformations of the world onto itself''  are automorphisms of this universe, which respect typing. I shall comment also on the philosophical aspect of Tarski's  notion of world shortly. For the sake of my argument I assume hereafter that Tarski's notion of logicality correctly explicates Hilbert's implicit notion of logicality.  

Saying that Hilbert considers the logical structure as the core structure of all theories I am talking about the logical structure in the second sense. So Hilbert's axiomatic construction of formal theories like Klein's construction has two levels: it comprises logic at the base level and (formal) theories at the next level. This Hilbert's two-level construction reflects Klein's two-level construction as follows: for Klein each (metric) geometry is the projective geometry provided with an additional metric structure; for Hilbert each (axiomatic) theory is logic provided with an additional axiomatic structure (say, with the axiomatic structure of Euclidean geometry). Geometrical structures and logical structures are construed here similarly in terms of invariants of corresponding transformations.  
This analogy between Klein's and Hilbert's foundations is visualized at the below diagram where $PG$ stands for projective geometry, $MG$ for metric geometry (in the restricted sense of the Riemanian geometry of manifolds with constant curvature), $s_{1}, s_{2}$ are functors forgetting some structure (the metric structure for $s_{1}$ and the non-logical axiomatic structure for $s_{2}$) and, finally, $a_{1}, a_{2}$ for my suggested two-level analogy between Klein and Hilbert:

$$\xymatrix{MG\ar[r]^{a_{1}} \ar[d]_{s_{1}}&Theories \ar[d]^{s_{2}} \\ PG \ar[r]_{a_{2}} & Logic}$$

Freudeltal \cite{Freudental:1960} is greatly impressed by Hilbert's \emph{Foundations} of 1899 and describes the preceding history anachronistically from Hilbert's axiomatic viewpoint. In his view  Klein's approach in foundations of geometry is fundamentally confused.  According to Freudental Klein ``did not understand the logical function of the model'' and as a result  ``[t]he logic of geometry was obscured rather than clarified by [Klein's] discovery of a model for non-Euclidean geometry''; thus ``this discovery did not contribute to clarifying the foundations of geometry'' (\cite{Freudental:1960}, p. 614). Freudental stresses the fact that Hilbert's axiomatic approach helps one to avoid a  vicious circle in Klein's reasoning who first uses Euclidean geometry for introducing projective geometry and then suggests the latter as a foundation of the former (Klein himself denies that the circle is vicious) (\cite{Freudental:1960}, p. 614). A similar point can be made about many other logical difficulties of the 19th century geometry easily found, in particular, in Lobachevsky's works. However these advantages of Hilbert's Axiomatic Method come with a price about which Freudental is well aware and which he readily accepts:
\begin{quote} 
[According to the view popular in the end of the 19th century,  i]n any case geometry deals with real space - Pasch, Enriques, Veronese, Pieri, Klein stressed this, and at the eleventh hour (1987) Russell wrote his philosophy [of] \emph{Foundations of Geometry}, which reveals the faint footstep of Kant rather then the paw of the lion. \\
\emph{Wir denken uns drei vershiedene Systeme von Dingen.} [Let us consider  three distinct systems of things.] - the bond with reality is cut. Geometry has become pure mathematics.  (\cite{Freudental:1960}, p. 618)
\end{quote}

``Wir denken uns drei vercshiedene Systeme von Dingen..'' is the beginning  of the first chapter of Hilbert's \emph{Foundations} of 1899 \cite{Hilbert:1899}, see \textbf{2.1}.  Unlike Greenberg \cite{Greenberg:1974} Freudental locates the ``cut of the bond with reality'' precisely: the cut was made not by the discovery of Non-Euclidean geometries but by Hilbert's axiomatic treatment of these geometries. Freudental describes Pasch's philosophical work aiming at the empirical justification of geometrical axioms as a ballast, which didn't allow Pasch's axiomatization of geometry  of 1882 \cite{Pasch:1882} to become as influential as Hilbert's axiomatization of 1899.  Freudental quotes Einstein \cite{Einstein:1921} saying that

\begin{quote} 
The progress entailed by axiomatics consists in the sharp separation of the logical form' and the realistic and intuitive contents. [..] The axioms are voluntary creations of the human mind. [..] To this interpretation of geometry I attach great importance because if I  had not been acquainted with it, I would never have been able to develop the theory of relativity. (quoted after  \cite{Freudental:1960}, p. 619)
\end{quote}

Stressing the ``liberating effect'' (Greenberg) and the success of Hilbert's Axiomatic Method in ``the lobbies of science and philosophy'' ( \cite{Freudental:1960}, p. 619) Freudental does not provide any philosophical argument justifying the idea of pure mathematics cut from any ``realistic content''. He rather suggests that the public success of this idea talks for itself. However today this success appears today to be more problematic than it appeared a century ago. Obviously the liberating effect of Hilbert's axiomatic approach triggered important developments in mathematics and physics, which wholly changed the shape of these sciences. However neither in the pure mathematics nor in physics Hilbert's Axiomatic Method is systematically applied today in its original form. The liberation from outdated patterns of thinking is necessary for scientific progress but it is not sufficient, so we need today to think about new ways of linking mathematics with experience and with reality. 

A drawback of Hilbert's axiomatic revolution, which I have already stressed in \textbf{2.2} and judged unacceptable, is the revival of the traditional metaphysical thinking (that Kant calls \emph{dogamtic}) in the new cloths of formal logical methods. Tarski's analysis of logicality makes this return of metaphysics very explicit. Tarski's \emph{world} (which, remind, serves Tarski for defining logicality) is a ``metaphysical world of thought'' (Cassirer), which is not simply prior to any experience but which (unlike Kant's world of pure intuitions) is also wholly independent of the world of experience. Tarski's world of individuals and classes would be still there even if the empirical science would not exist and would not be possible. Whether this metaphysical world is believed to exist independently on some Platonic heaven, or believed to exist in the human or some super-human mind, is a secondary metaphysical issue, which I shall not go into.

A more specific critical argument against Hilbert's axiomatic approach and in favor of Klein's approach consists of saying that using logic is \emph{the} universal invariant structure of all theories for the unification of geometry is an overkill, which leaves unnoticed some important details the hierarchy of relevant geometrical structures. This explains why we still call the hyperbolic geometry by this name invented by Klein. The axiomatic classification of geometries into Euclidean and Non-Euclidean on the basis of Euclid's Fifth Postulate reflects the history of the question rather than the structure of the subject-matter itself. As Weyl put this in 1923

\begin{quote}
The question of the validity of the ``fifth postulate'', on which historical development started its attack on Euclid, seems to us nowadays to be a somewhat accidental point of departure. \cite{Weyl:1923} p.??
\end{quote}    
 
This Weyl's remark squares with Lawvere's general critique of Formal Axiomatic Method according to which Hilbert-style axiomatic presentations of theories are ``subjective'' (\textbf{4.3}). In Chapter \textbf{4.3} we shall see how this problem is fixed with the New Axiomatic Method. 

Let me summarize this Section. We have seen that Hilbert's Axiomatic Method combines two ideas: one is the idea of mathematical structure as an invariant of some group of transformations, and the other is the idea of logic as a ultimate foundations of scientific thought on a par with metaphysics. The former idea stems from the 19th century geometry, and more specifically from Klein's  \emph{Erlangen Program}  \cite{Klein:1872}, the latter - from Aristotle and the following Scholastic tradition. The two ideas are combined by Tarksi's method: logic is thought of as a structure invariant under all possible transformations of the metaphysical world unto itself. None of these two ideas entails the other. However their mixture, which one may call \emph{structural logicism}, became in the 20th century popular and influential.  Mathematical structuralism as described by Hellman \cite{Hellman:2006a} shares Hilbert's logicist assumption and so qualifies as a version of structural logicism. This logicist assumption also makes part of the set-theoretic substantialism but in this latter case it is not combined with the structural view on logic.

\section{Objective Structures}

Kant's theory of mathematical representation provides a link between the pure mathematics and the mathematically-laden empirical science. Euclid has his own theory of representation, which he develops on the basis of his geometry and presents in his \emph{Optics} (\cite{Euclides:1883-1886}, vol. 7; English translation \cite{Euclid:1945})\footnote{
Euclid's \emph{Optics} is a geometrical theory of vision rather an optical theory in the modern sense of the term. Latin word \emph{perspectiva} is the translation of Greek \emph{optika} first used by Boethius (\cite{Andersen:2007}, p. xx)
}. Euclid's \emph{Optics}, Ptolemy's \emph{Optics} \emph{Ptolemy:1996} and later geometrical theories of spatial representation developed under the name of  theories of \emph{perspective} \cite{Andersen:2007} are interesting for the history of physics because they provide geometrical models of observation and explicitly introduce the observer into a mathematical picture. Although Kant's theory of mathematical representation and mathematical theories of perspective are very different in their character it seems me appropriate to consider them from a common viewpoint. 

According to Kant the origin of mathematical objectivity lies in the fact that all mathematical objects are constructed according to certain rules, which can be either explicit as in the case of Euclid's Postulates (\textbf{1.4}) or implicit. This is why for Kant all sound mathematics is objective, and there is a single idealized subject (that Kant calls \emph{transcendental subject}) who is doing the sound mathematics. This Kantian view squares with Newton's notion of \emph{absolute} space and time, which allegedly lies behind any subjective \emph{relational} representation of space and time (say, in same frame of reference). So Kant's notion of objectivity accounts for the objectivity of Newton's absolute space and time rather than for a relational notion of objectivity that may arise from comparing various frames of reference\footnote{
In fact Kant's attitude to Newton's theory of absolute space and time was critical and complicated; it changed during Kant's career. What I say about it here is a rough approximation, for details see \cite{Friedman:1992}.}.   
Today one cannot, of course, use Newton's outdated views on space and time for defending the Kantian view on objectivity. On the contrary, one can refer to the success of \emph{relational} theories of space and time (most importantly in Einstein's theory of relativity) as a reason for taking seriously mathematical theories of objectivity, which explicitly introduce the observer into the picture. Leaving the history of perspective for another study I shall now describe a recent mathematical theory of objectivity, which is closely related to the structuralist trend in mathematics and physics. 

Consider a geometrical object $S$, for example a stature thought of as a solid in the Euclidean 3-space, and a series of images $P_{i}$ of $S$ taken from different perspectives. Considering changes of perspective as  transformations of one such image into another ($P_{i} \rightarrow P_{j}$)  one may describe $S$ as an invariant of the given group of transformations. Under appropriate conditions such a description determines $S$ uniquely. This way of thinking about geometrical objects is pertinent when one uses \emph{coordinate systems}. For example circle $C$ with radius $R$ and centre $O$ is represented with the orthogonal coordinate system with the origin in $O$ by the equation $x^{2}+y^{2} = R^{2}$. But in a different orthogonal coordinate system this very circle $C$ is represented by a different equation. So if one wants to study the geometry of this circle using orthogonal coordinates one needs to learn how to switch from one such coordinate system to another and extract the information, which is invariant under the change of coordinate system (which in the given example is easy but in other cases is not). As long as coordinate systems represent physical frames of reference this problem is also relevant to physics. A classical example is given by Einstein's Special Relativity (SR). In the Newtonian physics it is assumed that spatial lengths and time intervals are \emph{objective} in the sense that they don't depend on the particular frame of reference, in which they are measured. SR tells us that the Newtonian assumption is only approximately correct for frames of reference which move with respect to each other with velocities significantly smaller than the speed of light. So in SR spatial lengths and time intervals are no longer objective in that sense: here measurements of lengths and times by different observers brings results, which do not agree with each other in the usual Euclidean way. However the spacetime of SR is also equipped with the notion of \emph{spatiotemporal interval}, which measures the spatiotemporal distance between \emph{events};  in SR spatiotemporal intervals are invariant under the change of reference frame, and in this sense are objective. Since the spatiotemporal interval interval between some observed events can be calculated on the basis of measurements of lengths and time intervals made in any given frame of reference, different observers can again agree with each other about what they observe. 

Generalizing upon the last example one can attach the following physical sense to the \emph{Erlangen Program}: given a space think about the automorphisms of this space as transformations of one point of view on this space (= reference frame) into another. Invariants of these transformations represent objective features of this space (and the objective features of stuff found ``in'' this space) as distinguished from subjective representations of these features in each reference frame. Such transformations, which transform a given reference frame into another reference frame, are called in physics \emph{passive} transformations. So the invariance under passive transformations is the most obvious meaning of being objective, which is formulated in mathematical and physical rather than general philosophical terms. However the invariance under some  \emph{active} transformations, which transform the given space itself rather than views on this space, also contribute to the colloquial notion of physical objectivity. (Beware that the distinction between passive and active transformations comes from physics and has no purely mathematical sense: we are now talking about the same geometrical transformations but interpret them differently.) It turns out that for purely mathematical reasons the assumption about certain \emph{symmetries} of space and time ) implies the existence of certain quantities globally conserved in all physical processes (Noether theorem). In particular, the assumption according to which a world living some time ahead or some time behind our actual world in the absolute Newtonian time would not be distinguishable from our world by its physical properties, implies the existence of the conserved quantity that we call \emph{energy}. The principle of conservation of energy (which survives in the relativistic model without the absolute time) is ``more objective'' than any contingent spatiotemporal state of affairs in the sense that unlike the spatiotemporal facts this principle can be  independently tested in very different (closed) physical systems at different places in different times by different people \cite{Nozick:2001}. 

So far we were talking about invariants of transformations as \emph{quantities}, which in physics are usually represented by real or complex numbers. When Klein formulated his  \emph{Erlangen Program} he meant a more general notion of invariant but, as I have already mentioned, he did not use our modern mathematical notion of invariant \emph{structure};  in particular, he did not know how to describe metric and projective spaces as abstract structures. The notion of objectivity as invariance under transformations that does not involve the notion of structure can be called \emph{pre-structural}\footnote{On the basis of this pre-structural notion of objectivity coming from physics Nozick  \cite{Nozick:2001} develops a general philosophical theory of objectivity, which covers the objectivity of truth, the objectivity of knowledge and the ethical objectivity.}. One can however go further and introduce the modern notion of mathematical structure into physical contexts and provide relevant structures with some sense of physical objectivity.  I cannot say which  role if any this idea plays in today's physics but it does play a role in some recent philosophy of physics: it is the core idea of the recent \emph{structural realism} and the recent structuralist trend in the philosophy of physics \cite{Debs&Redhead:2007}\cite{Ladyman:2009}. 

In Hellman's terms the structuralist idea of objectivity can be formulated as follows: invariant structures are objective, exchangeable particular items instantiating these structures are not. When this idea is interpreted ontologically (so that the objectivity implies existence) the resulting ontology resembles Plato's ontology where only enduring invariant Forms really exist while material items, which ``partake'' these forms don't  really exist but rather \emph{become} (\textbf{5.4}). As Nozick puts this (discussing conserved quantities in physics)
\begin{quote}
[S]omething whose amount in this universe cannot be altered, diminished or augmented  should count as (at least tied for being) the most objective thing there actually is. (\cite{Nozick:2001}, p. 81) 
\end{quote}

In \textbf{8.8} I describe a more dynamical non-structuralist way of thinking about objectivity suggested by category theory. 

\section{Types and Categories of Structures}
Remind from \textbf{2.1} that Hilbert's idea of Formal Axiomatic Method did not work out quite as expected: soon after publishing in 1899 his \emph{Foundations of Geometry} Hilbert realized that his system of axioms for geometry has non-isomorphic models, i.e., that this theory is not \emph{categorical}. As long as Hilbert's Axiomatic Method is used as intended the lack of categoricity is seen as a problem, and one wants to enforce the categoricity by some additional means (like Hilbert's second-order axiom of completeness used in the second edition of \emph{Foundations of Geometry} published in 1903 \cite{Hilbert:1903}. However as I have already stressed in Chapter \textbf{3} this Axiomatic Method is rarely (if ever) used as intended by Hilbert; at least, it is not used as intended in the 20th century \emph{structural} mathematics. 

A canonical example of structural axiomatic theory is \emph{group theory} (see \textbf{3.2} above); for the reader's convenience I repeat here the relevant axioms: 

\textbf{G1}: $x\circ (y\circ z) = (x\circ y)\circ z$  (associativity of $\circ$)

\textbf{G2}: there exists an item \emph{1} (called \emph{unit}) such that for all $x$  $x\circ 1 = 1\circ x = x$ 

\textbf{G3}: for all $x$ there exists $x^{-1}$ (called \emph{inverse} of  $x$) such that $x\circ x^{-1} = x^{-1}\circ x =  1$.

The fact that there are non-isomorphic groups means that the formal theory \textbf{GT} determined by axioms \textbf{G1-3} is not categorical. However \textbf{GT} is not the group theory in the usual sense of the word.  The group theory in the usual sense of the word is a theory of \emph{models} of \textbf{G1-3} called groups. When axioms \textbf{G1-3} are called \emph{axioms of group theory} then either the relevant notion of theory is not Hilbert's or the expression ``group theory'' does not have its usual sense. The categoricity of \textbf{GT} can be enforced with additional axioms, for example, with axioms, which specify that the intended model of our new formal theory \textbf{GC} is the infinite cyclic group. So \textbf{GC} pins down this particular group. But the group theory studies all groups but not only this or that particular group. (Talking about a ``particular group'' I mean, of course, ``particular up to isomorphism''.) Thus the lack of categoricity of theory \textbf{GT} is not a problem but rather an advantage: axioms \textbf{G1-3} determine a \emph{type} of structures of  interest and by adding new axioms one may specify particular structures (like the infinite cyclic group) and study these structures in the context of other structures of the same type. This makes the content of the ``pure'' group theory. However groups are also studied in  the context of structures of some other types like in the case of homotopy theory. In such cases one leaves the boarders of the ``pure'' group theory and applies this theory in other parts of mathematics. A \emph{homotopy group} (i.e., a group of homotopies) is a mixed structure, which combines a group structure with an appropriately decorated topological structure.   
The possibility to combine structures of different types is crucial for Bourbaki's project of unifying mathematics  and designing its global ``architecture'' \cite{Bourbaki:1950}.  The general idea of Bourbaki's architecture is this: to specify a small number of  ``great types of structures'' and then reconstruct the rest of mathematics by combining structures of these ``great types''. Thus Bourbaki's structural mathematics can be described as a joint model theory of a whole bunch of formal theories like \textbf{G}, which combine with each other in various ways. The artificial character of this description reflects the artificial character of Hilbert's notion of formal theory. In the context of modern structural mathematics it is more natural to call \textbf{GT} and its likes \emph{definitions} rather than theories. This remark reveals a traditional aspect of the modern structural mathematics: just like in Euclid's mathematics one introduces here some general concepts through definitions, specifies these general concepts through some further definitions, instantiates these concepts with certain objects (which in the structural mathematics are called structures) according to certain rules and, finally, applies some further constructions for proving  non-trivial theorems (which, generally, do not follow directly from the definitions).

Let me now use the example of  group theory for being more specific. Saying that groups $G_1$ and $G_2$ are isomorphic is tantamount to the following:

\textbf{I1}: elements of $G_1$  are in one-to-one correspondence with elements of $G_2$; 

\textbf{I2}: for all elements $a_1$, $b_1$, $c_1$ from $G_1$ such that $a_1\oplus b_1 = c_1$ the corresponding elements $a_2$, $b_2$, $c_2$ from $G_2$  satisfy  $a_2 \otimes b_2 = c_2$  where  $\oplus$ is the group operation in $G_1$ and $\otimes$ is the group operation in $G_2$.

A one-to-one correspondence between elements of two given groups that satisfies \textbf{I2} is called (group) isomorphism. Groups are isomorphic if and only if there exists isomorphism between them. One may think about a one-to-one correspondence between elements of groups $G_1$ and $G_2$, which satisfies condition \textbf{I2}, as a map or transformation  \emph{i}:  $G_1\rightarrow G_2$ of one group into another group. Since a one-to-one correspondence is a symmetric construction the choice of $G_1$ as the source and $G_2$ as the target of this transformation is arbitrary. In other words, one and the same isomorphism-\emph{qua}-correspondence gives rise to two isomorphisms- \emph{qua}-transformations  \emph{i}: $G_1\rightarrow G_2$  and \emph{j}: $G_1\rightarrow G_2$ , which run in opposite directions and cancel each other on both sides
\footnote{``[C]ancel each other'' means exactly this:  the composition transformation $i\circ j$ resulting from the application of transformation \emph{j} after transformation \emph{i} sends every element of $G_1$ into itself and composition transformation $j\circ i$ sends every element of $G_2$ into itself (beware that none of the two conjuncts implies the other). Given these conditions each of transformations  \emph{i} and  \emph{j} is called the  \emph{inverse} of the other. The inverse map can be defined without referring to elements of groups if these groups are thought of as objects of a category: maps $i, j$ are called mutually inverse when $i \circ j = 1_{G_1}$ and $j \circ i = 1_{G_2}$ where  $1_{G_1}$ is the identity morphism of $G_1$ and $1_{G_2}$  is the identity morphism of $G_2$.  
}. As I have already remarked, this terminology is slightly confusing but it is too common for trying to change it.

Invertible maps aka isomorphisms are not the only sort of maps between groups suggested by definition \textbf{G1-3}. Instead of one-to-one correspondence between elements of  $G_1$, $G_2$, one considers a more general kind of correspondence that is allowed to be many-to-one (but not one-to-many). In other words, one considers \emph{functions} \emph{f}: $S_1 \rightarrow S_2$ from the set $S_1$ of elements of $G_1$ to the set $S_2$ of elements of $G_2$.   Condition \textbf{I2} remains the same; notice that it can be satisfied when elements $a_1$, $b_1$ are different but elements $a_2$, $b_2$ are the same. A maps satisfying these relaxed conditions is called a (group) \emph{homomorphism}. 

Similar general notion of map can be defined for any mathematical structure. However there is no canonical way for doing this as can be seen at the example of topological spaces. Remind that a topological space is the set $T = P(S)$ of subsets of some given set $S$ some of which are called \emph{open} while some other are called  \emph{closed}, which satisfy certain axioms. An isomorphism of topological spaces  (aka \emph{homeomorphism} aka invertible continuous transformation) \emph{i}: $T_1\rightarrow T_2$ takes each open of $T_1$ into an open of $T_2$ and each closed of $T_1$ into a closed of $T_2$. Now one wants to define a \emph{continuous} (but not necessarily invertible) transformation of topological spaces as a general transformation that transforms one topological space into another (similarly to homomorphisms for groups). One may think of two possibilities: either to take opens to opens without assuming the corresponding condition for closed, or the other way round, to take closed to closed without assuming the corresponding condition for opens. The axioms of topological space treat opens and closed differently but they don't allow one to rule out one of these two options in favor of the other on some formal reasons. So the definition of continuous transformation is a matter of further specification of the notion of topological space. The standard definition of continuous transformation corresponds to the second option: the \emph{inverse} image of an open set is always open while the direct image of an open set can be also closed.

Group homomorphisms and similar maps between structures of other types are colloquially called ``structure preserving" (or ``structure reflecting'' as in the above example of continuous transformation). This is somewhat misleading because if such maps preserve anything at all it is a \emph{type} of structure but not a particular structure. Think about this trivial example: for all groups $G_1$, $G_2$ there exist a homomorphism \emph{h}: $G_1 \rightarrow G_2$  which sends every element of $G_1$ to the unit of $G_2$.  This homomorphism ``destroys all information" about $G_1$ reducing its image to a single element;  it doesn't provide any information about $G_2$ either (except the fact that $G_2$ is a group). Thus homomorphisms, generally, don't allow for invariants in anything like the same sense in which isomorphisms do so. The colloquial talk of ``preservation'' of structure apparently stems from the early days of group theory when group homomorphisms were not well distinguished from group isomorphisms and were seen as a sort of ``imperfect'' isomorphisms \cite{Wussing:2007}.

In order to see that thinking about homomorphisms as imperfect isomorphisms is misleading try to replace isomorphisms by homomorphisms in the process of abstraction described in \textbf{5.5}. One might expect to get in this way a generalized notion of structure but this doesn't work. Recall the first step: given class \textbf{\emph{G}} of groups we have divided it into equivalence subclasses of isomorphic groups. Two groups are isomorphic if and only if there exists isomorphism (i.e., an invertible  transformation) between them; clearly this is an equivalence relation. Let me (for the sake of argument) call two groups \emph{homomorphic} if and only if there is a homomorphism between them. Although this latter relation is also an equivalence, one can see the difference: since \emph{all} groups are homomorphic (see the above example of group homomorphism) one cannot use this equivalence for dividing \textbf{\emph{G}} into equivalence subclasses. Saying that two given groups are homomorphic is tantamount to saying that the given groups are groups. So the relation of homomorphism just introduced (not to be confused with the standard notion of homomorphism as transformation) doesn't make sense.

We see that homomorphisms cannot do the same job as isomorphisms: the invertibility condition stressed by Hilbert in the above quote turns out to be crucial for structural abstraction. One cannot reason ``up to homomorphism" in anything like the same way in which people reason up to isomorphism doing structural mathematics. Since ``invariant" in the given context is just another word for structure it is clear that homomorphisms, generally, don't have invariants in anything like the same sense in which isomorphisms and groups of isomorphisms do so.

Let me now explain how the category theory enters into the picture. The emergence of category theory in he 1940s and its further development in the context of structural mathematics was related to a growing awareness of the role of general maps (not only isomorphisms). This is a simple theorem \cite{MacLane:1996} that a class of structures of any fixed type provided with an appropriate notion of general map form a category. A non-mathematical reader is advised to refresh the definition of category given in \textbf{6.3} above. In this present Section I present some basic mathematical context in which this notion has been first invented and proved useful. This is the context of structural mathematics. For a detailed historical account of early days of category theory see \cite{Kromer:2007}.

Remind from \textbf{8.2} the Bourbaki's intention to get rid of the set-theoretic background of structural mathematics and work with ``pure structures'' whatever this might mean. Apparently the notion of category makes this structuralist dream true. Consider the category $G$ of groups for example.  In order to build this category in the usual way one needs, first, to introduce groups (through \textbf{G1-3} or with an equivalent system of axioms), second, introduce group homomorphisms and, finally put these elements together and, in particular, explain how the homomorphisms compose with each other. In this construction the set-theoretic background is indispensable. However one may try to reverse the conceptual order and introduce $G$ from the stretch, first, as an abstract generally category, and then (through some additional axioms) as a specific category distinct from any other category by some distinctive categorical properties.  As soon as this is done one may claim that the notion of group is fully accounted for in terms of maps (including identity maps aka objects), and so the set-theoretic background is effectively dispensed. 
  
We have seen in \textbf{4.1} how Lawvere realized such a project for the category of sets \cite{Lawvere:1964}. This and other similar achievements  show that the idea of categorification (in the broad sense of \textbf{6.3}) works out. We shall now see, however, that the better it works the less relevant becomes its structuralist philosophical underpinning. Since this philosophical underpinning concerns the Axiomatic Method, it is not just a matter of personal philosophical taste. In the next Chapter I shall argue that Hilbert's structural Axiomatic Method is not appropriate for the categorical mathematics and suggest a replacement.

\section{Invariance versus Functoriality}
In their seminal paper that marks the official birth of category theory Eilenberg and Mac Lane describe this theory as a continuation of the \emph{Erlangen Program}:

\begin{quote}
This may be regarded as a continuation of the Klein Erlanger Program, in the sense that
a geometrical space with its group of transformations is generalized to a category with its
algebra of mappings. ( \cite{Eilenberg&MacLane:1945}, p. 237) 
\end{quote}

Commenting on this passage Marquis says:
\begin{quote}
Klein's fundamental idea was that to study a geometry, one had to look at its
group of transformations and, furthermore, the geometric properties of that geometry
are those which are invariant under the group of transformations. This seems to
be the core of the generalization that Eilenberg and Mac Lane had in mind. (\cite{Marquis:2009}, p. 10) 
\end{quote}
  
I disagree with Mariquis' understanding of the core of this generalization. In my view, the core of the generalization is another feature of the \emph{Erlangen Program}, which Marquis also mentions, namely the very idea of considering mathematical objects together with their transformations (mappings). The categorical generalization of this idea amounts to including non-invertible mappings into the picture and treating them on equal footing with invertible ones (\textbf{6.3}). The notion of invariant, which is crucial for Klein's original approach, does not play the same major role in a categorical setting: there are no invariants of categories of geometrical transformations similar to invariants of groups of such transformations. Algebras of mappings don't play the same role: they are rather abstract algebraic objects similar to abstract groups (as distinguished from groups of concrete geometrical transformations). The categorical counterpart of invariance is rather \emph{functoriality} that splits into \emph{covariance} and \emph{contravariance}. Saying this I don't try to follow Eilenberg\&MacLane verbatim. Apparently these authors indeed have the idea of ``preservation of structure''  in their minds when they say:

\begin{quote}
The invariant character of a mathematical discipline can be formulated in these terms. Thus, in group theory all the basic constructions can be regarded as definitions of co- and contra-variant functors, so we may formulate the dictum: The subject of group theory is essentially the study of those constructions of groups which behave in a covariant or contravariant manner under induced homomorphisms. More precisely, group theory studies functors defined on well specified categories of groups, with values in another such category.
\end{quote}

I cannot see that the subject-matter of group theory construed as ``those constructions of groups which behave in a covariant or contravariant manner under induced homomorphisms'' indeed expresses the ``invariant character'' of this discipline. In this passage Eilenberg\&MacLane describe the co- and contra-variance (i.e., the functoriality) as a generalized invariance. But this is misleading.  What has been said above about non-invertible group homomorphisms \emph{mutatis mutandis} applies to non-invertible functors: they don't leave anything  invariant, not even the type of structure because the domain and the codomain categories of a given functor $e: A \rightarrow B$ can be categories of structures of different types. The \emph{covariance} of functor $e$ amounts to the commutativity of all diagrams of the form

  $$\xymatrix{x\ar[r]^{e_{x}} \ar[d]_{f}&e(x) \ar[d]^{e(f)} \\ y \ar[r]_{e_{y}} & e(y)}$$

were $f: x \rightarrow y$ is morphism in $A$, $e(f): e(x) \rightarrow e(y)$ is the image of morphism $f$ under functor $e$ in $B$, and $e_{x}, e_{y}$ are called \emph{components} of functor $e$. 

A \emph{contravariant} functor $e'$ reverses arrows as follows:  

$$\xymatrix{x\ar[r]^{e'_{x}} \ar[d]_{f}&e'(x)  \\ y \ar[r]_{e'_{y}} & e'(y)\ar[u]_{e'(f)}}$$

Assuming after Klein, Eilenberg, MacLane and Mariquis that morphisms and functors are transformations one can see that the above square diagrams are wholly dynamic: everything changes in them and nothing remains invariant. These squares  represent not some invariant structures surviving through changes but certain \emph{coherences} between different changes, which make these diagrams to commute. This coherence of transformation in mathematics is called functoriality. Unless the relevant functors are invertible functoriality does not imply invariance. The tendency of thinking of functoriality as generalized invariance is the same tendency by which people think of homomorphisms as imperfect isomorphisms. This tendency can be described as a case of conceptual inertia, which prevents one from making the full justice of a new concept.   

According to Eilenberg\&MacLane's dictum the genuine group-theoretic features are those features of groups, which vary functorially. This dictum provides a clear epistemic criterion for distinguishing epistemically significant features of groups from all their other features. But contra  Eilenberg\&MacLane I claim that this new epistemic criterion does not reduce to the older   platonic or structuralist criterion according to which only invariant features are epistemically significant, while all variable features are not.   

The switch from the structuralist thinking in terms of invariance to the new categorical thinking in terms of covariance and contravariance (i.e., functoriality) signifies a decisive brake with the structuralist viewpoint, which is not yet commonly realized. 
This conceptual change remains, however, an open philosophical and mathematical problem rather than a mere psychological problem. Since the concept of category is born in the context of structural mathematics it comes equipped with a structural method of theory-building (to wit Hilbert's Formal Axiomatic Method) and with the structural definition of category given by Eilenberg and Mac Lane in \cite{Eilenberg&MacLane:1945}. In order to establish the new point of view one needs to remake this basic framework on new grounds. Before the job is done the new idea survives through various compromises with the established view; one of such compromises I am now going to describe.

\section{Are Categories Structures?}
Axioms for category theory suggested by Eilenberg and Mac Lane \cite{Eilenberg&MacLane:1945} provide a structural definition of category in the same way in which axioms \textbf{G1-3} provide a structural definition of group. This suggests considering the notion of category as just another type of structure on a par with groups, rings, topological spaces, etc. . However the intended application of this notion described in \textbf{8.5} does not quite fit this idea as we shall now see. The category $G$ of groups or the category of structures of any other type \emph{embodies} the corresponding type of structures in a way similar to which the set $N$ of natural numbers embodies the (extension of the) concept of natural number. In the pre-Cantorian mathematics one could embody the concept of natural number only by instantiating this concept with some particular numbers; after Cantor one can also represent the full extension of this concept by set $N$. Similarly, in the pre-categorical mathematics one could embody the general concept of group only by instantiating it with some particular structure like that of infinite cyclic group
\footnote{    
Notice that I talk here about the instantiation of the group concept not in the same sense in which Hellman talks about the instantiation of structures by particular ``systems''. Each particular group structure can be instantiated by such a system of concrete ``items'' and concrete operations with these items. This is Hellman's instantiation. I am talking here about the instantiation of the general group concept by a particular group structure. One may object that what I call instantiation is rather a specification because a particular group structure is not an individual object but another concept. I do count particular structures as particular objects determined up to isomorphism.  
};
after Eilenberg and Mac Lane one can also represent the full extension of group concept by category $G$. One may remark that the concept of proper class is sufficient for it. However as I have already stressed in \textbf{5.8} the proper classes mathematically sterile in the sense that one cannot do any further construction with them. Categories, on the contrary, are very productive in this sense: given $G$ as above one may square it ($G^{2}$ is the category of functors from the two-point set to the category of groups), power with it some other category (consider, for example category $Vect^{G}$ of linear presentations of groups), etc. Thus unlike proper classes categories are effective vehicles for types of structures just like structures are effective vehicles for concepts like that of the infinite cyclic group. 

Using Hilbert's distinction between mathematics and metamathematics one may describe a category of structures  as a \emph{metastructure} meaning that all structures of a given type are models of the same formal theory (like theory \textbf{GT}), and so the given category belongs to the metatheory (namely, the model theory) of this theory. Does this metastructure qualify as a structure in the usual sense of the word? The answer is in negative because of the \emph{size problem} (\textbf{5.8}): since the collection of all groups or all structures of any other fixed type is not a set but a proper class it is problematic to think about the categories of such structures on a par with these structures themselves. 

A category in which all morphisms (including identity morphisms) form a set (as distinguished from a proper class) is called \emph{small}. Small categories can be thought of as structures on their own. Some well known mathematical structures are small categories. In particular, a (single) group is a category with only one object and all morphisms invertible;  a partial order is a category having at most one morphism going from one given object to another given object. Eilenberg and MacLane  realize that small categories are the only categories, which are fully legitimate (as long as one uses their axiomatic definition of category) - and they also realize that the reduction of categories to small categories makes the concept of category uninteresting. So they suggest an ``intuitive'' point of view:

\begin{quote}
We remarked [..] that such examples as the``category of \emph{all} sets'', ``category of \emph{all} groups'' are illegitimate. The difficulties and antinomies here involved are exactly those of ordinary intuitive \emph{Mengenlehre} [i.e. the naive set theory]; no essentially new paradoxes are apparently involved. Any rigorous foundation capable of supporting the ordinary theory of classes would equally well support our theory. Hence we have chosen to adopt the intuitive standpoint, leaving the reader to insert whatever  type of logical foundation (or absence thereof) he may prefer. ( \cite{Eilenberg&MacLane:1945}, p. 246)
\end{quote}

In \textbf{4.4} we have seen how Lawvere  \cite{Lawvere:1966a} tries to solve this problem with his idea of category of categories as foundation (CCAF). He first uses Eilenberg\&MacLane's axioms for defining an abstract category $CAT$, about which he thinks intuitively without considering model-theoretic issues (the \emph{elementary theory}); then he strengthens this theory specifying by \emph{internal} means of $CAT$ each particular object $C$ of $CAT$ as a category (the \emph{basic theory}). Lawvere does not refer in CCAF to classes but one may argue after Mayberry \emph{Mayberry:2000 } that some primitive notion of class or collection is  involved in this construction anyway because otherwise the Formal Axiomatic Method used at the first step wouldn't work.  In order to use Eilenberg\&MacLane's axioms, so the argument goes, one needs to think of collection of individuals and relations between these individuals, which satisfies these axioms under an appropriate interpretation; since the given system of axioms is obviously non-categorical one should also distinguish between different collections satisfying these conditions. Since the Formal Axiomatic Method applies everywhere in mathematics and not only in category theory one cannot identify such a primitive notion of collection of individuals with a category from the outset but is obliged to consider a category as a special type of mathematical structure on a par with other types of structure.  

I agree with Mayberry that a primitive notion of collection is an indispensable ingredient of the Formal Axiomatic Method. But unlike Mayberry I don't believe that this method is itself indispensable in mathematics. Instead of thinking about categories as structures of special sort I suggest to think about structures as categories of special sort\footnote{Lawvere suggests (in person) to identify structures with functors from small categories to a large background category like that of sets. In order to make this categorical notion of structure self-sustained one needs to treat the size issue categorically without referring to set theory; I leave now this important subject aside.}. Since Hilbert's Formal Axiomatic Method brings about nothing but structures (\textbf{8.3}) in order to justify the view on mathematical structure as a special type of category one needs to build category theory and the rest of mathematics with a different method.

\section{Objects Are Maps}
Remind Gauss' idea of \emph{intrinsic geometry}, which I mentioned talking about the intuitive background of Riemann's geometry in (\textbf{7.3}). In his \cite{Gauss:1828} Gauss identified some geometrical properties of surfaces, including the intrinsic notion of curvature, that do not depend on the way in which a given surface is embedded into the outer 3-space. This idea later allowed Einstein (after Riemann) to conceive of a curve space independently of any embedding of this space in some outer flat space. The notion of intrinsic geometry has a clear \emph{invariant} character: the Gaussian curvature and other intrinsic characteristics (including topological characteristics) of a given surface do not change when this surface is folded or unfolded. However I would like to stress a different aspect of Gauss' idea, which allowed for a new way of thinking about objects and spaces in geometry. 

Given a surface one can think of it (i) in the usual way as a two-dimensional \emph{object} living in the Euclidean 3-space and (ii) as a 2-space on its own rights (characterized by the intrinsic properties of the given surface), which is a home for its points, lines, triangles, etc.. Generally, a geometrical \emph{object} can be described in this context as a map of the form $s: B \rightarrow C$ where $B$ is a \emph{type} of the given object and $C$ is a \emph{space} where the given object lives and \emph{instantiates} (or \emph{represents}, which in the given context is the same) its type. This way of thinking about spaces and objects in spaces can be represented by this diagram:

$$\xymatrix{TYPE\ar[r]^{object} & SPACE}$$

It is suggestive also to think about a general categorical morphism in this way. As we can see it agrees with the way in which morphisms are thought of in logical categories (\textbf{4.3}) and also illustrates the dialectical unity of geometry and logic stressed by Lawvere at the example of topos (\textbf{4.9}). Since we interpret all domains as spaces and all codomains as types these notions are relational in the given context (each type serves as a space for incoming morphisms and each space serves as a type for outgoing morphisms). I shall illustrate this way of thinking about objects, spaces and types at some elementary geometrical examples. 

Let us distinguish between two different notions of Euclidean plane: (i) the domain of Euclidean Planimetry and (ii) an object living in the Euclidean 3-space ($ESPACE$). I shall write $EPLANE$ for Euclidean plane in the first sense, and write $eplane$ for Euclidean plane in the second sense. Then an $eplane$ can be presented as a map:

$$\xymatrix{EPLANE\ar[r]^{eplane} & ESPACE}$$ 

Such maps are many (there are many planes in the space) but they all ``are of'' the same type;  this type in its turn is inhabited (as a space) by objects of different types: 

 $$\xymatrix{CIRCLE\ar[r]^{circle} & EPLANE}$$ 
 
 A more interesting example I borrow from Lobachevsky \cite{Lobachevsky:1840}. Although Lobachevsky reasoned about the hyperbolic space intuitively without using Euclidean models he actually used what in modern term can be described as a non-standard hyperbolic model of Euclidean plane. Namely, he found in the hyperbolic space a special surface that he called the \emph{horisphere} and showed that intrinsically the geometry of this surface is the plane Euclidean geometry. (This helped Lobachevsky to develop the hyperbolic trigonometry and on this basis build an analytic model for his geometry.) Thus we have got an object of type $EPLANE$ that does not look like $eplane$:

$$\xymatrix{EPLANE\ar[r]^{horisphere} & HSPACE}$$ 
 
($HSPACE$ stands for hyperbolic space). We can see that the idea tho classify geometrical objects into types by their shapes and forms is misleading because it works only when the background space is fixed. However one can learn about any geometrical type by studying it intrinsically as a space, i.e., by studying objects of all types living in it. My suggested approach unlike the Riemanian approach does not privilege the intrinsic description against the extrinsic one: the fact that a horisphere is intrinsically an Euclidean plane (in the sense of being of type $EPLANE$) is just as significant as the fact that this horisphere is an object in the hyperbolic 3-space ($HSPACE$): when one studies geometrical objects there is, generally, no epistemic reason for privileging their types over their spaces or privileging their spaces over their types. 

The geometrical objects so construed are composeable in the obvious way. Here is an example of composite object:

 $$\xymatrix{CIRCLE\ar[r]^{circle_{1}} \ar[d]_{circle_{2}} & EPLANE \ar[dl]^{eplane} \\ESPACE}$$
 
In the given situation we tend to identify $circle_{1}$ living on $EPLANE$ with  $circle_{2}$ living in $ESPACE$. However if $ESPACE$ is projected back onto $EPLANE$ and this projection turns $circle_{2}$ into an oval the difference becomes obvious.

I leave it to the reader to check that the composition of objects is associative.
 
Developing this toy categorical geometry it is suggestive to think about spaces as places where objects meet: 

 $$\xymatrix{\ar[dr] &\\& {\bullet} \\ \ar[ur] &}$$
 
think about types as places where objects split: 

 $$\xymatrix{&\\{\bullet}\ar[ur]\ar[dr] &\\ &}$$
 
 and think about the composition of objects as an operation that glues spaces and types together: 
 
$$\xymatrix{\ar[dr]&&\\&{\bullet}\ar[ur]\ar[dr] &\\ \ar[ur]&&}$$

and thus form new objects capable for ``self-representation'': 

$$\xymatrix{\ar[dr]&&\\&{\circlearrowright}\ar[ur]\ar[dr] &\\ \ar[ur]&&}$$

Among such self-representing objects there is one that we call \emph{identity object} and denote 1; the identity object  is distinguished by the usual conditions ( $f¡1 = f$ for each object $f$ represented in the same space, and $1¡g = g$ for each object  $g$ of the same type), which in the given context are read as the conditions of being \emph{neutral} with respect to the composition of objects. Think about EPLANE for example. We know how EPLANE represents objects of various types (circles, triangles and the like) and we also know how EPLANE is represented in its turn in various other spaces. Now in order to make sense of saying that all these representations are representations in and of the same thing one should think of this thing itself as an object that represents itself in a way, which stabilizes the dynamics of all inner (incoming) and outer (outgoing) representations. 

The associative composition of objects and the above assumptions about the identity objects makes these objects into a category, which I denote $Geo$ for further references.  Clearly $Geo$ is not a well-defined specific category but rather a way of thinking about categories geometrically and of thinking about geometry categorically\footnote{For an introduction into the modern categorical geometry I refer the reader to \cite{Gelfand&Manin:2003}. }.

As the reader have seen what in category theory is usually called morphism I call object and what in category theory is usually called object I call identity object. This suggested terminological change is not without a reason. The distinction between objects and morphisms is useful in the structural mathematics because it helps to construct categories from structures of certain types (like groups) and appropriate morphisms of these structures (like group homomorphisms). $Geo$ can be also construed in this way as the category of differentiable manifolds and differentiable maps (assuming that differentiable manifolds are construed as structured sets) or as some other similar category. However I suggest a different way of thinking about $Geo$ and about categories in general. Before I shall try to clarify this different way of thinking let me remind that the usual distinction between objects and morphisms of categories is formally dispensable: since with each object $A$ of given category $C$ is associated a unique identity morphism $1_{A}$, one may formally identify objects of $C$ with their corresponding identity morphisms and thus consider objects as morphisms of special sort. Thus, formally, a general category can be described as a class of things called morphisms provided with a (partial) binary associative operation called composition. I claim that the name of \emph{objects} is more appropriate for these things than the name of morphisms. Saying this I do \underline{not} mean that any such thing can be called object in the most general sense of the term (as suggested among others by Parsons \cite{Parsons:2008}). Instead I have in mind a particular notion of object, which implies that objects form categories. This way of thinking about objects can be expressed by the slogan ``objects are maps''. 

One may wonder why I am not happy with the established mathematical terminology and don't want to get rid of the term ``object'' altogether and talk about maps or morphisms. The reason is that the term ``object'' does not belong exclusively to mathematics but has also a philosophical meaning; as a philosophical notion the notion of object is closely related to the notion of objectivity. Although my notion of object is not standard for the 20th century philosophy and for the 20th century mathematics, it is rooted in an earlier philosophy and earlier mathematics, in particular, in the Kantian philosophy and in the 19th century geometry. Since these historical links play a significant role in my proposal, I am not ready to scarify the term ``object'' by replacing it by the term ``map''. In addition, since the term ``object'' is used in category theory anyway, I want to contrast in this context my proposed notion of object against the usual notion, which underpins the established terminology.

In the beginning of his \emph{Founadations} of 1899 \cite{Hilbert:1899} Hilbert proposes thinking about geometrical objects as abstract ``things'' (of several basic \emph{types}) and about geometrical spaces as ``systems of things'' (\textbf{2.1}). I claim that the Kantian idea according to which an object is not simply a thing but is a \emph{represented} thing, deserves to be taken seriously. The principal modification of Kant's original viewpoint, which I suggest, is the following: while Kant assumed that all objects are represented in the same space I allow for representations in different spaces. We have seen that the 19th century geometry provides us with relevant examples. The example of Lobachevsky's \emph{horisphere} is particularly useful in this respect because it shows that geometrical objects, generally, are determined not only by their types (i.e., by their intrinsic properties) but also by spaces in which they are represented. We have also seen that the modern mathematical notion of category provides a suitable framework for such objects if one identifies these objects with morphisms of some category. 

Further, I believe that the notion of object is epistemically prior to the notions of type and space. We encounter objects first and only then classify them into types and organize them into spaces. That types of geometrical objects are many has been well understood already in Euclid's times; the fact that there are also many representation spaces, and that each geometrical type may serve as a representation space, has been well understood by mathematicians only in the 19th century (although the development of this idea can be traced through the history of mathematical studies of perspective back to Euclid's \emph{Optics} \cite{Euclid:1945},\cite{Andersen:2007}). 

I suggest in the Kantian vein that this way of thinking about objects is relevant not only to objects of special sort that we call geometrical objects but to all mathematical objects and to all other objects. This claim goes on a par with the popular slogan of category-theorists according to which objects cannot be adequately studied unless their corresponding maps are also taken into consideration. I add that those maps are true objects, and that what the category-theorists call objects are objects of very special sort.

The notion of \emph{objecthood} described in this Section suggests a revision of the structuralist notion of \emph{objectivity} described in Section \textbf{8.4} above. In the 20th century philosophy this structuralist notion of objectivity is usually opposed to the logicist notion of objectivity as the universal logical validity, which stems from Frege. However as soon as structuralism is reconciled with logicism (\textbf{8.3}) there is a space for agreement: the proponents of Fregean logical objectivity and proponents of more liberal structuralist objectivity agree that the key to objectivity is \emph{invariance}. Then the controversy  reduces to questions about different kinds of invariance: some people defend the invariance of Form, some other defend the invariance of Substance, etc.. Nozick (\cite{Nozick:2001}, ch.2) provides a general solution for such debates by allowing the objectivity to have degrees, which correspond to different kinds of transformations and different kinds of invariants. 

The categorical notion of object as a map and the idea of functoriality as a generalized invariance ( \textbf{8.6}) suggests thinking about objectivity in terms of functoriality (i.e., co- and contravariance) rather than invariance. This approach sublates the debate about invariance and refuses the Platonic viewpoint (taken by Nozick as a matter of course) according to which the most objective things are those, which are the most stable and the lest capable for change. The functorial objectivity amounts to the \emph{coherence} of changes rather than to the lack of change. Think about $Geo$ as a category of observers observing each other. One does not need Leibniz' notion of pre-established harmony (which sounds like a structural notion) for achieving a harmony (i.e., coherence) between observations made by different observers. What one needs is a set of appropriate global properties of $Geo$ like the property of being Cartesian closed; such properties can be thought of as emerging properties like elsewhere in science. If $Geo$ is apt for supporting an internal logic this logic is ``shared'' by all observers without being construed as an invariant structure. If this logic supports truth-values it enables our observers to agree about objective truths. Beware that objects in this picture are, generally, not the observers themselves but mutual representations of the observers. The observers (which can be also called \emph{subjects}) are objects of special sort capable for the (objective) self-representation.  

Interestingly, this picture allows different observers to play different roles, and thus does not require any global symmetry of the universe. Remark that the example of Special Relativity used by Nozick for the justification of his structuralist notion of objectivity as invariance does not immediately generalize to the General Relativity.  Although the \emph{general covariance}, which is the key concept of General Relativity, can be locally expressed in terms of invariant quantities, the general covariance does not imply the existence of global invariants. At the same time the general covariance can be always expressed in terms of of functoriality \cite{Baez:2004}. I leave this and other relevant physical examples for a future study. 

Let me summarize. 

In \textbf{8.3} I have argued that Hilbert's Axiomatic Method not only provided the Euclidean and Non-Euclidean geometry with some particular foundations but also suggested a way of unifying various geometrical theories on the same logical basis. Ultimately Hilbert aimed at the unification of all sciences by the same method ( \textbf{2.1, 2.4}). The specific way in which Hilbert thinks about logic and about axiomatic theories allows one to interpret Hilbert's Axiomatic Method as a far-reaching generalization of Klein's project of unification of geometry on the basis of projective geometry. Then in  \textbf{8.6} I have shown that although the category theory as conceived by Eilenberg and MacLane in 1945  \cite{Eilenberg&MacLane:1945} continues Klein's \emph{Erlangen Program} in a way, it does not give the same place to the structural notion of invariance but replaces it by a more general and conceptually very different notion of functoriality. Finally, in the present Section I showed how geometrical types, spaces and objects can be organized into a category, and on this basis suggested a non-structural way of thinking about objecthood and objectivity. In the next concluding Chapter I show how this way of thinking affects the Axiomatic Method.

\chapter{New Axiomatic Method (instead of conclusion)}
In the following long promised presentation of the New Axiomatic Method I shall use as a guide Lawvere's description of Axiomatic Method as ``unification and concentration'' (\cite{Lawvere:2003}, p. 213) and generalize upon some examples of axiomatic thinking due to Lawvere and Voevodsky. I begin with the unification, then turn to the concentration and, finally, discuss the place and the special character of logic in the New Axiomatic Method. 

\section{Unification}

Categorical logic and categorical geometry suggest a unification strategy, which essentially differs from the structuralist unification strategy used by Klein and Hilbert. Instead of looking for a core invariant structure shared by all geometrical spaces one studies maps between these spaces (i.e., objects in the sense of \textbf{8.8}) and organizes the universe of these maps/objects into a category.  Given such category $Geo$ (as in \textbf{8.8}) there is, generally, no way and no sense to dispense with its objects in a way similar to which one may dispense with concrete examples of the infinite cyclic group and conceive of this group as an abstract structure. 

One may ask what to do if some multiplicity of objects that one wishes to account for by category-theoretic methods does not really form a category. Since I am trying now to describe only a very general epistemic strategy (alternative to the structuralist strategy) I don't stick to the standard definition category. I already mentioned one modification of this standard notion, namely the notion of (weak) $n$-category; some other modification of the original concept are found in the literature and some other may appear in the future. The key unificatory notion offered by the New Axiomatic Method is that of composition (which is, generally, non-commutative). The general unificatory strategy is this:

\textbf{(Ia)} \emph{Take relevant objects and try to compose them:}

 $$\xymatrix{\ar[r]&{\bullet}\ar[r]&}$$

 \textbf{(Ib)} \emph{Use joints} $\bullet$  \emph{between composable objects as types and spaces. This works as follows. Fix some object $x$ and try to compose it with other objects from the left; classify objects composable with $x$ in this way into the same type:}

$$\xymatrix{&&&\\ \ar[r]&{\bullet}\ar[ur]\ar[dr] &\\ &&&}$$

 \emph{Then fix some object $y$ and try to compose it with other objects from the right; put all objects composable with $y$ in this way into the same space:}

$$\xymatrix{\ar[dr]&&\\&{\bullet}\ar[r]&\\ \ar[ur]&&}$$

\textbf{(Ic)} \emph{Think of joints between objects as neutral self-represented objects, i.e., by the identity objects.}     
  
I add this latter item having in mind to organize the multiplicity of objects into a category. Since the standard notion of category may develop this latter item can be understood liberally.

Let me now compare the unificatory strategy \textbf{(Ia-c)} with the unificatory strategy offered by the Formal Axiomatic Method. As a basis for unification of mathematics and all science Hilbert tacitly assumes a notion of universal logic. So the Hilbertian unification works top-down. As soon as the idea of universal logic is given up, and one is converted into a logical pluralist, the unifying power of the Formal Axiomatic Method is lost. The categorical unification just described works in a bottom-up way through the double-sided composition of objects. One may ask whether this method of unification is appropriate for ``big unifications'', which involve objects of different categories. In fact such a unification can be achieved in the same bottom-up way through a reiteration of principles \textbf{(Ia-c)}.  Given two categories $A, B$ consider ``interdisciplinary'' or ``mixed'' objects of the form  $A \rightarrow B$ and then build a category of such mixed objects. For a mathematical example (for the mathematical reader) think about the homology groups understood as contravariant functors from the topological category of chain complexes to the algebraic category of abelian groups. Clearly such mixed objects play a major role in the modern structural mathematics and in its global architecture. Thus the New Axiomatic Method unlike the Formal Axiomatic Method bridges different fields directly through appropriate objects (and categories of such objects) rather than tries to put all such fields on some common invariant ground.  

\section{Concentration}
\textbf{(Ia-c)} are elements of a strategy rather than algorithmic steps because in the real world very few objects are immediately available, so one needs further efforts for constructing them. Let me first point to some relevant examples and then suggest a generalization. My first example is Lawvere's functorial semantics (\textbf{4.6}) Prima facie the functorial semantics simply reproduces the Hilbert-Tarski scheme (I mean the notion of formal theory and its interpretation) by categorical means. However the application of categorical technique is not an innocent step as we shall now see. In the language of category theory an interpretation of given formal theory $T$ is construed as functor $T \rightarrow B$ where $B$ is a base contentual theory. After Tarski Lawvere takes for $B$ the category $Set$ of sets (supported by ETCS). The first step that takes Lawvere away from Hilbert's structuralism is his observation that in most situations of interest the \emph{categoricity} of theory $T$ (in logicians' sense of the uniqueness of its model up to isomorphism) is not only unrealistic but also undesirable because what really matters is a set of properties of category $M$ of models, which cannot obtain unless $M$ is sufficiently reach. Lawvere's 2004 commentary \cite{Lawvere:2004}, according to which $T$ can be seen as a special \emph{generic} object of $M$, i.e. as a generic model (see Lawvere's quote in \textbf{4.6} above), is a further step in the same direction.  

Remarkably, this last step brings us back to the traditional way of thinking about first principles of geometry found in Euclid, who introduces several basic geometrical shapes (point, straight line, circle) and describes several procedures (in Postulates 1 - 3) that generate  from these basic shapes all further constructions. Independently from Lawvere's functorial semantics this traditional geometric idea reappears in category theory in 1960s in  Ehresmann's  \emph{sketch theory}. In the case of sketch theory the Euclidean geometrical analogy is even more precise because sketches are not logical categories that express propositional axioms but rather geometrical ``generic figures";  in some approaches sketches are not even categories but directed graphs with an additional structure, see  \cite{Wells94sketches:1994} for an overview and further references. Independently from Ehresmann Lawvere used a similar approach in the   \emph{basic theory} making part of his CCAF \cite{Lawvere:1966a}; remind from \textbf{4.2} that Lawvere's internal description of general  category involves a sketch-like construction with four basic generic shapes: point for objects, straight segment for morphisms, triangle for composition of morphisms and square for the associativity of composition.  

Yet more closely the Euclidean genetic approach is reproduced by Voevodsky \cite{Voevodsky:2010} in the first geometrical part of his Univalent Foundation: like Euclid Voevodsky begins constructing his hierarchical universe of homotopy types with a point and then applies a simple inductive procedure for generating from this point the rest of this universe (\textbf{6.10}).  

Remind also that the traditional genetic approach also plays a role in Hilbert's Formal Axiomatic Method, namely in its more advanced symbolic version (\textbf{2.3 - 2.4}). Every student of mathematical logic learns today this traditional genetic method with the notion of well-formed formula generated by a given alphabet of symbols through several  iteratable syntactic operations.  The New Axiomatic Method lifts Hilbert's restriction according to which the genetic method of theory-building and object-building must be reserved to strings of symbols. Accordingly the New Method allows for generic procedures other than manipulations with symbols (albeit it takes advantage of reducing such procedures to manipulations of symbols when such reduction is possible and useful like in the case of computer implementation of Martin-L\"of's type theory). Which generic procedures are allowed in mathematics and which are not, in my understanding, is not a general question about the Axiomatic Method and so it does not have a general answer.  The idea to restrict the spectrum of possible generic procedures once and for all on some epistemic, pragmatic or other grounds cannot be justified. 

Thus the second element of the New Axiomatic Method is this: 

\textbf{(II)} \emph{Specify basic objects and basic constructions that concentrate the current knowledge about the given field. Specify procedures that generate all other relevant objects and relevant constructions from the basic ones.} 

Admittedly this rule represents a traditional rather than original feature of the new method. Nevertheless its revival in today's mathematics is quite remarkable. It is a form of atomism which remains alive in the modern science. Looking at Hilbert's formal mathematics from a historical perspective one may ironize that Hilbert apparently took too straightforwardly the metaphor of the Ancient atomists according to which atoms were letters of Greek alphabet. The New Axiomatic Method uses the atomistic genetic principle in a more general sense without reducing its scope to syntactic issues. 

The reader may remark at this point that although we have already spoken about unification and concentration, which is supposed to cover the whole of Axiomatic Method, we did not yet say a word about logic. Indeed \textbf{(I-II)} cover only the \emph{geometrical} part of this method. These two steps prepare (so far solely by intuitive geometrical means) a field of objects, which constitutes an \emph{objective background} for the further logical organization of a given theory. Without such a background the logical organization of theories cannot be \emph{objective} in Hegel-Lawvere's sense (\textbf{4.8}). Once again let me stress that this intuitive non-logical element of Axiomatic Method is also present in Hilbert's Formal Axiomatic Method albeit in a rudimentary form that involves only syntactic constructions (\textbf{2.4}). Hilbert recognizes other forms of mathematical intuition \cite{Hilbert&Chon-Vossen:1932} but leaves them outside the ready-made \emph{axiomatic} mathematics; the New Axiomatic Method brings them back in and thus significantly extends the scope of axiomatic thinking. 

\section{Internal Logic as a Guide and as an Organizing Principle}
As soon as the objective geometrical background $G$ of a given theory is fixed and organized into a category with \textbf{(I-II)}  it makes sense to think about the ``subjective logic of inference between statements'' (\cite{Lawvere:1994}, p. 16, see \textbf{4.8}), i.e., about logic in the common sense of the word. The appropriate technical notion of logic is that of \emph{internal logic} $L$ of category $G$. As soon as $L$ is taken together with $G$ but not in an abstract form it qualifies not only as subjective but also as objective (this squares with the idea that every subject is an object of sort, see \textbf{8.8} above). $L$ allows one to reproduce the geometric steps \textbf{(I - II)} (unification and concentration) in a new \emph{logical} form; this logical reproduction can be understood as an axiomatic \emph{reflection} upon the preceding ``naive'' geometrical stage. However in our case unlike Hilbert's case this reflection does not involve any attempt to abolish geometrical intuitions or to reduce these intuitions to some special form. Like in the last Section I begin the following presentation with some examples and then suggest a generalization.           

We have already learned about two sorts of geometrical categories supporting an internal logic: toposes (\textbf{4.9}) and (higher) homotopy categories (\textbf{6.9}). I refer specifically to Lawvere's paper \cite{Lawvere:1970a} where the author provides his axiomatic treatment of topos theory, and the recent Voevodsky's lecture \cite{Voevodsky:2010} where the author suggests new Univalent Foundations of mathematics (\textbf{6.10}). These works are relatively rare examples of successful application of logical methods in the recent mainstream mathematics: in both cases logical considerations greatly simplify and clarify otherwise very difficult mathematical concepts and ideas. The role played by Lawvere's axiomatic topos theory with respect to the earlier geometrical research on toposes and related subjects in the Grothendieck school can be justly compared with the role of Hilbert's \emph{Foundations of Geometry} with respect to the 19th century geometry. Voevodsky's Univalent Foundations  play a similar role with respect to the modern algebraic geometry. 

Lawvere's and Voevodsky's axiomatic approaches share a common feature, which concerns the relationship between geometry and logic in their theories. Lawvere describes this relationships in Hegelian terms as a dialectical contradiction.  Voevodsky uses the title of ``direct formalization'' for describing the fact that the relevant fragment of homotopy theory models Martin-L\"of's type theory without any additional axiom. (The only geometrically motivated axiom used by Voevodsky in the Univalent Foundations is his Axioms of Univalence.) In both cases the relevant system logic is internal with respect to the corresponding geometrical theory. This makes a sharp difference with Hilbert's axiomatic approach, where logic is always external, i.e., fixed wholly independently from geometrical or any other specific mathematical theories, which this logic supports. Hilbert's approach assumes an asymmetric relationship between geometry and logic: geometry is supposed to be logical while logic is not supposed to be geometrical. Lawvere's and Voevodsky's approaches make this relationship symmetric: their geometrical theories are logical and their logic is geometrical. Lawvere's axiomatic topos theory and Voevodsky's homotopy type theory qualify both as theories of logic and as theories of geometry. Although logical and geometrical aspects of these theories can be distinguished, this distinction is not the usual distinction between the logical form and the geometrical content. The relevant distinction between logical and geometrical aspects is rather a matter of interpretation. For example, one and the same element of a given topos Lawvere interprets in turn (i) logically as a set of truth-values (more precisely, a truth-value object) and (ii) geometrically as a particular sheaf. Similarly in the homotopy type theory one and the same element is interpreted either (i) logically as identity type or (ii) geometrically as groupoid of homotopies. Talking in this context about interpretation one should beware that the relevant notion of interpretation is not standard either. The standard model-theoretic notion of interpretation (or more precisely, its most basic version) assumes the relationships between a theory and its models, which can be pictured in that way:

 $$\xymatrix{ &Theory \ar[dl] \ar[d] \ar[dr] &  \\ Model_{1} & Model_{2} & Model_{3}}$$

In topos theory and homotopy type theory there is no such intermediate element, which can be (and also can be not)  interpreted either logically or geometrically: one \emph{immediately} interprets geometrical notions in logical terms (Voevdsky's ``direct formalization'') and, reciprocally, interprets logical notions in geometrical terms. (While the geometrical interpretation of logical notions qualifies as interpretation in the sense of model theory, the reciprocal interpretation does not.) This can be visualized as follows: 

 $$\xymatrix{ Logic \ar@<2pt>[r] & Geometry \ar@<2pt>[l]}$$

An elementary duality between geometry and logic can be noticed already in the toy geometrical category $Geo$ from \textbf{8.8} where morphism $A \rightarrow B$ is a geometrical object of \emph{type} $A$ represented in \emph{space} $B$. Types are turned into spaces and reciprocally by the reversal of arrows. Since the notion of space belongs to geometry and the notion of type belongs to logic this duality of types and spaces represents in a rudimentary form the dialectical duality of geometry and logic.

Given geometrical category $G$ with internal logic $L$ one may ``do mathematics in $G$'' using the standard Formal Axiomatic Method, i.e., take $L$ as basic logic and built further theories in the usual way by considering various non-logical axioms and systems of such axioms. This gives rise to the notion of \emph{internal theory} (internal with respect to $G$). I would like to stress once again that Lawvere's axiomatic topos theory as presented in \cite{Lawvere:1970a} and Voevodsky's homotopy type theory as presented in \cite{Voevodsky:2010} are \underline{not} internal theories in the above special sense. This is why I claim that Lawvere and Voevodsky apply a different Axiomatic Method, which I call the New Axiomatic Method. In the New Axiomatic Method the internal logic $L$ serves for building the geometrical theory of $G$ but not for building some further theories ``in'' $G$. The theory of $G$ is \underline{not} a metatheory, which provides some additional geometrical support to $L$, but a full-fledged axiomatic geometrical theory. $G$ and $L$ can be also described as two different aspects of the same theory, one of which is geometrical and the other is logical. 

As I have already mentioned at several occasions in Chapter \textbf{4} a similar kind of relationships between geometry and logic is found in Euclid's \emph{Elements}, which by Friedman word is a general ``form of rational representation'' (\textbf{2.5}) rather than a specific mathematical theory comparable with the theory of Euclidean geometry presented in Hilbert's \emph{Foundations of Geometry}. Lawvere's topos theory and Voevodsky's homotopy type theory also qualify as such general forms of rational representation rather than specific theories treating some small fragments of the Platonic mathematical universe. This return of the older pattern of axiomatic thinking appears natural when one considers the history of Hilbert's Axiomatic Method during the passed century. Hilbert granted  the freedom to invent new mathematical axioms but not the freedom to invent new systems of logic. This second degree of axiomatic freedom has been developed afterwards and changed the sense of Axiomatic Method.  Hilbert's thinking about this method is this: first, one develops an intuitive mathematical theory, and then axiomatizes this theory using a ready-made system logic as a tool, which helps one to put ideas in the right order. Today one has a choice between different ready-made logical tools and a freedom to construct new logical tools appropriate to the given task. Thus one is no longer in a position allowing for relying on logic as something given; the epistemic requirement according to which one must ``reason logically'' in the new context means that one must pay attention to logical issues like truth and rules of inference but not that one must stick to some particular logical rules. This is why every modern axiomatic theory, or at least every axiomatic theory, which is supposed to be used as a foundation in mathematics and science, can be today only a self-sustained form of rational reasoning and rational representation. 

How internal logic $L$ of given geometrical category $G$  may help one to build an axiomatic theory of $G$? In order to see how it works remind of Hilbert's notion of \emph{axiomatization of logic} as distinguished from the axiomatization of particular mathematical theories like the theory of Euclidean geometry (\textbf{2.3}). Remind also of Russell's radical mathematical logicism according to which all mathematics is logic and each mathematical theorem is a logical tautology. The alleged tautological character of mathematics does not make mathematical theorems trivial and useless. Similarly, the fact that $L$ qualifies as a system of logic rather than a non-logical formal theory like Hilbert's formal theory of Euclidean geometry, does not imply that all tautologies of $L$ are trivial and non-interesting. Moreover since $L$ is fully interpreted in $G$ each tautology of $L$ translates into a contentual geometrical proposition. If $L$ is presented axiomatically, i.e., if its tautologies are generated from a  set of distinguished tautologies called axioms, the axiomatic presentation of $L$ directly translates into an axiomatic presentation of $G$. 
Since $L$ unlike $G$ makes explicit forms of judgement and forms of inference of judgements from other judgements $L$ allows one to complete the ``geometrical axiomatization'' \textbf{(I - II)} described earlier in this Chapter with a full-fledged ``logical axiomatization'', which involves the logical specification of forms of correct reasoning in $G$   

\footnote{Although Voevodsky calls the association of $L$ (Matrin-L\"of's type theory) to $G$ (homotopy theory) the ``direct formalization'' (of $G$ by $L$) it is clear that $L$ does not represent the \emph{logical form} of $G$ in the usual sense of the word. Let $x, y$ be individual variables and $R$ be a binary predicate variable. The expression $xRy$ is said to represent the logical form of linguistic expression \emph{Mary loves John} in the sense that this latter expression obtains from $xRy$ when variables $x, y, R$ are given the appropriate semantic values; by giving these variables some other semantic values one may obtain other contentual expressions like \emph{Peter hates Paul} or $1 > 0$, which have the same logical form but different contents. Hilbert's formal theory of Euclidean geometry represents the logical form of the traditional intuitive Euclidean geometry in the same traditional sense of ``logical form'' (\textbf{2.1}) However homotopy theory is not just one contentual intuitive interpretation of Matrin-L\"of's type theory among a bunch of other such interpretations. Matrin-L\"of's type theory and the appropriately adjusted homotopy theory are rather two different legs of the same theory called homotopy type theory. The basic relationship between these two legs is not that of form and content.}.  
  
So let me formulate the last element of the New Axiomatic Method, which concerns logic:

\textbf{(III)} \emph{Given a category of objects formed according to} \textbf{(I - II)}  \emph{use these objects for building a system of logic, which is sound in this field, and which expresses essential facts about this field in the form of tautologies. Then use the obtained logic for a further axiomatic organization of the given field through the logical unification and logical concentration.}

Lawvere's notion of quantifier as adjoint to substitution (\textbf{4.3, 4.5}) is a logical supplement of the basic geometrical  unificatory mechanism described in \textbf{9.1}: given object $f: A \rightarrow B$ that provides the basic geometrical unification of areas $A, B$, the quantification along $f$ is a logical operation that moves propositions from one of these areas to the other.  The logical mode of concentration is the familiar choice of axioms for $L$, which complements the geometrical concentration described in \ \textbf{9.2}. The Curry-Howard correspondence in Cartesian closed categories (\textbf{4.4}) is the correspondence between these two modes of concentration.

It is appropriate to ask what to do if obtained geometrical category $G$ does not support any internal logic. In this case the theory of $G$ cannot be an independent axiomatic theory and thus must be included into a larger theory with internal logic. Generally, the New Axiomatic Method does not allow one to make up axiomatic theories by a fiat because unlike the Formal Axiomatic Method it doesn't use ready-made external logical tools. I don't see this feature of the New Method as a disadvantage because the purpose of this method is to build axiomatic theories having some objective significance rather than support the pure speculation. Even if the pure mathematical intuition alone cannot provide a clear-cut boundary between science and speculation it does this job as a part of the empirical intuition. This is why in the New Axiomatic Method the intuition plays its traditional role of join between the pure mathematics and the empirical science.  The fact that Hilbert's Axiomatic Method in its original form contributed very little to the 20th century science must be taken seriously as a lesson. I suggest that the New Axiomatic Method outlined in this concluding Chapter will perform better in this respect. I leave the justification of this claim for another book.  

\bibliographystyle{plain} 
\addcontentsline{toc}{chapter}{Bibliography}
\bibliography{catax} 
\end{document}